%% file: paper.tex
\documentclass[showlinks,onefignum,onetabnum]{siamart250211}

\input{shared}

\ifpdf
\hypersetup{
  pdftitle={Trace estimation for parameter-dependent matrices},
  pdfauthor={F. Matti, H. He, D. Kressner, H. Lam}
}
\fi

\begin{document}
\maketitle

\begin{abstract}
Stochastic trace estimation is a well-established tool for approximating the trace of a large symmetric matrix $\mtx{B}$. Several applications involve a matrix that depends continuously on a parameter $t \in [a,b]$, and require trace estimates of $\mtx{B}(t)$ for many values of $t$. This is, for example, the case when approximating the spectral density of a matrix. Approximating the trace separately for each matrix
$\mtx{B}(t_1), \dots, \mtx{B}(t_m)$ clearly incurs redundancies and a cost that scales linearly with $m$. To address this issue, we propose and analyze modifications for three stochastic trace estimators, the Girard-Hutchinson, Nyström, and Nyström++ estimators. Our modification uses \emph{fixed} randomization across different values of $t$, that is,
every matrix $\mtx{B}(t_1), \dots, \mtx{B}(t_m)$ is multiplied with the \emph{same} set of random vectors.
When combined with Chebyshev approximation in $t$, the use of such constant random matrices allows one to reuse matrix-vector products across different values of $t$, leading to significant cost reduction.
Our analysis shows that the loss of stochastic independence across different $t$ does not lead to deterioration. In particular, we show that $\mathcal{O}(\varepsilon^{-1})$ random matrix-vector products suffice to ensure an error of $\varepsilon > 0$ for Nyström++, independent of low-rank properties of $\mtx{B}(t)$. We discuss in detail how the combination of Nyström++ with 
Chebyshev approximation applies to spectral density estimation and provide an analysis of the resulting method. This  improves various aspects of an existing stochastic estimator for spectral density estimation. Several numerical experiments from electronic structure interaction and neural network optimization validate our findings.
\end{abstract}



\input{intro.tex}

\input{analysis.tex}

\input{spectraldensity.tex}

\input{numresults.tex}
\input{conclusion.tex}

\input{appendix.tex}

\section*{Acknowledgments}
The authors gratefully acknowledge the help of Lin Lin in reproducing the example in \cref{subsec:hamiltonian}. They thank David Persson for many enlightening discussions relating the proof in \cref{subsec:nystrom-pp}. They would also like to thank the reviewers for their detailed comments, which helped to improve the manuscript.

\bibliographystyle{siamplain}
\bibliography{bibliography}

\clearpage


\end{document}

%% file: shared.tex
\usepackage{lipsum}
\usepackage{amsfonts}
\usepackage{graphicx}
\usepackage{epstopdf}
\usepackage{algorithmic}
\ifpdf
  \DeclareGraphicsExtensions{.eps,.pdf,.png,.jpg}
\else
  \DeclareGraphicsExtensions{.eps}
\fi

\newsiamremark{remark}{Remark}
\newsiamremark{hypothesis}{Hypothesis}
\crefname{hypothesis}{Hypothesis}{Hypotheses}
\newsiamthm{claim}{Claim}
\newsiamremark{fact}{Fact}
\crefname{fact}{Fact}{Facts}

\headers{Trace estimation for parameter-dependent matrices}{F. Matti, H. He, D. Kressner, H. Lam}

\title{Stochastic trace estimation for parameter-dependent matrices applied to spectral density approximation\thanks{Submitted to the editors Feb 27, 2025.}}

\author{Fabio Matti\thanks{Institute of Mathematics, EPFL, 1015 Lausanne, Switzerland
  (\email{fabio.matti@epfl.ch}, \email{haoze.he@epfl.ch}, \email{daniel.kressner@epfl.ch}, \email{hysan.lam@epfl.ch}).}
\and Haoze He\footnotemark[2]
\and Daniel Kressner\footnotemark[2]
\and Hei Yin Lam\footnotemark[2]}

\usepackage{amsopn}

\usepackage{graphicx}
\usepackage{tikz}
\usepackage{subcaption}
\usepackage{booktabs}

\usepackage{xcolor}

\definecolor{darkblue}{HTML}{2F455C}
\definecolor{mainblue}{HTML}{0e437c}
\definecolor{darkorange}{HTML}{F98125}
\colorlet{lightblue}{darkblue!10!white}
\colorlet{lightorange}{darkorange!10!white}
\colorlet{linkcolor}{black}

\newcommand{\mtx}[1]{\boldsymbol{#1}}
\newcommand{\vct}[1]{\boldsymbol{#1}}

\DeclareMathOperator{\Trace}{Tr}
\newcommand\Hutch[1]{\Trace_{#1}}
\newcommand\Nystr[2]{\widehat{#2}_{#1}}
\newcommand\Nystrpp[2]{\Trace_{#2, #1}}

\newcommand\sfrac[2]{{}^{#1}{\mskip -5mu/\mskip -3mu}_{#2}}

%% file: intro.tex
\section{Introduction}
\label{sec:introduction}

Spectral distributions of matrices and linear operators reveal important properties in numerous applications across physics, chemistry, engineering, and data science: For example, in electronic structure calculations, eigenvalues represent the energy levels that electrons occupy~\cite{drabold-1993-maximum-entropy, ducastelle-1970-moments-developments, haydock-1972-electronic-structure, lin-2017-randomized-estimation}, in machine learning they are indicative of the topology of the loss landscape~\cite{ghorbani-2019-investigation-neural, yao-2020-pyhessian-neural}, and in graph processing they can uncover hidden graph motifs~\cite{huang-2021-density-states}. Given a symmetric matrix $\mtx{A} \in \mathbb{R}^{n \times n}$, the distribution of its eigenvalues $\lambda_1, \dots, \lambda_n \in \mathbb{R}$ can be represented by the spectral density $\phi(t) = n^{-1} \sum_{i=1}^{n} \delta(t - \lambda_i)$, which places a rescaled Dirac delta distribution $\delta$ at each eigenvalue. Clearly, assembling this expression amounts to computing all eigenvalues of the matrix; an operation that is often prohibitively expensive. Several techniques have been developed for approximating $\phi$, including moment matching~\cite{cohen-steiner-2018-approximating-spectrum, braverman-2022-sublinear-timea} and polynomial approximation, either implicitly by a Lanczos procedure \cite{lin-2016-approximating-spectral, chen-2021-analysis-stochastic}, or explicitly by expansion \cite{weisse-2006-kernel-polynomial, lin-2016-approximating-spectral}. Some of these methods have also been combined with a deflation step to yield even more efficient algorithms \cite{lin-2017-randomized-estimation, bhattacharjee-2025-improved-spectral}.

Being composed of Dirac delta distributions, the function $\phi$ itself is hard to approximate. However, as many applications only require a rough estimate of $\phi$, it usually suffices to instead approximate a smoothed spectral density of the form
\begin{equation}
    \phi_{\sigma}(t) = \frac{1}{n} \sum_{i=1}^{n} g_{\sigma}(t - \lambda_i) = \frac{1}{n} \Trace(g_{\sigma}(t \mtx{I}_n - \mtx{A})),
    \label{equ:smoothed-spectral-density}
\end{equation}
for some fixed smoothing kernel $g_{\sigma}$, typically a Gaussian~\cite{lin-2016-approximating-spectral, lin-2017-randomized-estimation} or a Lorentzian~\cite{haydock-1972-electronic-structure, lin-2016-approximating-spectral}. We thus arrive at the problem of computing the trace of the parameter-dependent matrix function $\mtx{B}(t) \equiv g_{\sigma}(t \mtx{I}_n - \mtx{A}) \in \mathbb{R}^{n \times n}$. Evaluating this matrix function \emph{exactly} requires the diagonalization of $\mtx{A}$, which does not seem to yield any gains compared to the original problem of computing its spectral distribution. The crucial point is that we can now \emph{estimate} the trace by matrix-vector products with $\mtx{B}(t)$, which --- thanks to the smoothness of $g_\sigma$ --- can be well approximated by a few matrix-vector products with $\mtx{A}$.

For a \emph{constant} matrix $\mtx{B}$, one of the most popular trace estimators is the Girard-Hutchinson estimator \cite{girard-1989-fast-montecarlo, hutchinson-1990-stochastic-estimator} along with variance reduction techniques~\cite{gambhir-2017-deflation-method, saibaba-2017-randomized-matrixfree, lin-2017-randomized-estimation, meyer-2021-hutch-optimal, persson-2022-improved-variants, chen-2023-krylovaware-stochastic, epperly-2024-xtrace-making}. Suitable extensions to parameter-dependent matrices have been considered, e.g., in~\cite{lin-2017-randomized-estimation,chen-2023-krylovaware-stochastic}, but we are not aware of an analysis providing rigorous justification and insight of these extensions. In passing, we note that dynamic trace estimation~\cite{dharangutte-2024-dynamic-trace,woodruff-2024-optimal-query} is an efficient technique for subsequently estimating the traces of matrices $\mtx{B}(t_1), \dots, \mtx{B}(t_m)$ when the increments $\mtx{B}(t_{i+1}) - \mtx{B}(t_i)$ are relatively small in norm. The potential of dynamic trace estimation appears to be limited in our setting because $\mtx{B}(t)$ may change rapidly close to eigenvalues, with $g_{\sigma}$ approximating a Dirac delta function.

The methods considered in this work apply to any matrix $\mtx{B}(t)$ that continuously depends on a parameter $t$ and are based on the following simple idea: Apply an existing randomized trace estimator to $\Trace(\mtx{B}(t))$ with \emph{fixed} random vectors, that is, the same randomization is used for each value of the parameter $t$. For example, the Girard-Hutchinson estimator becomes $\mathsf{est}(t):=n_{\mtx{\Psi}}^{-1} \sum_{j=1}^{n_{\mtx{\Psi}}} \vct{\psi}_j^{\top} \mtx{B}(t) \vct{\psi}_j$ for $n_{\mtx{\Psi}}$ fixed Gaussian random vectors $\vct{\psi}_1, \dots, \vct{\psi}_{n_{\mtx{\Psi}}}$. We will use the $L^1$-norm $\int_{a}^{b} | \Trace(\mtx{B}(t)) - \mathsf{est}(t) |~\mathrm{d}t$ to measure the error of such an estimator for all $t$ in an interval by $[a,b]$; similar to the approach taken in~\cite{kressner-2024-randomized-lowrank} for low-rank approximation.

\paragraph{New contributions} We analyze three well-established randomized trace estimators when they are applied to arbitrary parameter-dependent matrices $\mtx{B}(t)$ with fixed randomness (with respect to $t$): The Girard-Hutchinson estimator \cite{girard-1989-fast-montecarlo, hutchinson-1990-stochastic-estimator}, the trace of the Nyström low-rank approximation \cite{gittens-2013-revisiting-nystrom}, and the Nyström++ estimator \cite{persson-2022-improved-variants}. Combined with Chebyshev approximation, this allows one to reuse the majority of computations across different values of the parameter $t$, making the estimators scale favorably with respect to the number of parameter evaluations. We propose several additional modifications to the methods from \cite{lin-2017-randomized-estimation} for the special case of approximating the smoothed spectral density $\phi_{\sigma}$ from~\cref{equ:smoothed-spectral-density}. In particular, we use a more rigorous approach that preserves non-negativity for approximating matrix functions in terms of Chebyshev polynomials. The resulting Chebyshev-Nyström++ estimator 
unifies all methods from~\cite{lin-2017-randomized-estimation} in a single algorithm.
Combining our analysis of the randomized trace estimators with standard error bounds for Chebyshev approximation leads to~\cref{thm:chebyshev-nystrom}, which shows that $\mathcal O(\varepsilon^{-1})$ matrix-vector products suffice to guarantee an error $\varepsilon$ in the $L^1$-norm, which parallels a similar result for the Wasserstein-1 metric \cite{braverman-2022-sublinear-timea}.

\paragraph{Structure} In \cref{sec:analysis} we will analyze the error of three standard trace estimators --- the Girard-Hutchinson, Nyström, and Nyström++ estimators --- when they are applied to parameter-dependent matrices with fixed random matrix. In \cref{sec:application}, we propose and analyze the Chebyshev-Nyström++ estimator, which uses these estimators to approximate the smoothed spectral density of real symmetric matrices. Finally, we will illustrate and validate our developments for examples from various applications in \cref{sec:results}.

\paragraph{Reproducibility} \input{re-pro-badge.tex}

%% file: re-pro-badge.tex
The code for our numerical experiments can be found on \url{https://github.com/FMatti/parameter-trace}. The figures and tables were generated from the commit \href{https://github.com/FMatti/parameter-trace/tree/fcee24d}{fcee24d} on 2026-02-19 at 17:05:44 UTC.

%% file: analysis.tex
\section{Randomized trace estimators for parameter-dependent matrices}
\label{sec:analysis}

Consider a general parameter-dependent matrix
$\mtx{B}(t) \in \mathbb{R}^{n \times n}$ --- not necessarily of the form $g_{\sigma}(t\mtx{I} - \mtx{A})$, as will be used later in spectral density estimation (see \cref{sec:application}). Each
entry $b_{ij}(t)$ is assumed to be a 
continuous function on a bounded interval $[a,b]$. In this section, we 
describe and analyze methods that estimate, for several values of $t$, the 
trace $\Trace(\mtx{B}(t)) = b_{11}(t)+ \cdots + b_{nn}(t)$
from matrix-vector products of $\mtx{B}(t)$ with \emph{constant} random vectors. When $\mtx{B}$ is constant, our methods reduce to known algorithms. When $\mtx{B}$ is variable, important differences arise from the fact that the random vectors do not depend on $t$. On the one hand, as we will see in~\cref{sec:application}, this allows one to reuse computations. On the other hand, this significantly complicates the analysis, which will be the main focus of this section.

\subsection{Methods}
\label{subsec:methods}

\paragraph{Girard-Hutchinson estimator} One of the simplest methods for trace estimation, the Girard-Hutchinson estimator~\cite{girard-1989-fast-montecarlo,hutchinson-1990-stochastic-estimator} proceeds by sampling $n_{\mtx{\Psi}}$ independent Gaussian random vectors $\vct{\psi}_1,\dots, \vct{\psi}_{n_{\mtx{\Psi}}} \in \mathbb{R}^{n}$
and computes the approximation 
\begin{equation}
    \Trace(\mtx{B}(t)) \approx \Hutch{\mtx{\Psi}}(\mtx{B}(t))
    := \frac{1}{n_{\mtx{\Psi}}} \sum_{j=1}^{n_{\mtx{\Psi}}} \vct{\psi}_j^{\top} \mtx{B}(t) \vct{\psi}_j
    = \frac{1}{n_{\mtx{\Psi}}} \Trace( \mtx{\Psi}^{\top} \mtx{B}(t) \mtx{\Psi}).
    \label{equ:hutchinson-trace-estimator}
\end{equation}
Here, $\mtx{\Psi} := [\vct{\psi}_1 ~ \cdots ~ \vct{\psi}_{n_{\mtx{\Psi}}}] \in \mathbb{R}^{n \times n_{\mtx{\Psi}}}$ is a Gaussian random matrix, that is, 
its entries are independent standard normal random variables.
There are other sensible choices for the probability distribution, such as the Rademacher distribution or the uniform distribution on the $(n-1)$-sphere of radius $\sqrt{n}$ \cite{epperly-2024-xtrace-making,chen-2025-randomized-matrixfree}. Our theoretical developments, however, are intended for the Gaussian case.

The Girard-Hutchinson estimator is a popular choice for spectral density estimation; it has been used in the Haydock method \cite{haydock-1972-electronic-structure}, as well as the Delta-Gauss-Legendre and Delta-Gauss-Chebyshev methods from~\cite{lin-2016-approximating-spectral, lin-2017-randomized-estimation}. However, it shares the disadvantage of Monte Carlo methods: Its error is proportional to $n^{-\sfrac{1}{2}}_{\mtx{\Psi}}$, making it expensive in applications that require decent accuracy.

\paragraph{Nyström estimator} Given a good low-rank approximation of a matrix, the trace can be cheaply estimated by simply taking the trace of that low-rank approximation. In the context of spectral density estimation, the use of a non-negative smoothing kernel implies that $\mtx{B}(t)$ is symmetric positive semi-definite (SPSD) for every $t$. This makes the Nyström method
from~\cite{gittens-2013-revisiting-nystrom} well suited for constructing a (randomized) low-rank approximation.
Given a Gaussian random matrix $\mtx{\Omega} \in \mathbb{R}^{n \times n_{\mtx{\Omega}}}$, Nyström computes the approximation
\begin{equation}
    \mtx{B}(t) \approx \Nystr{\mtx{\Omega}}{\mtx{B}}(t) := (\mtx{B}(t) \mtx{\Omega}) (\mtx{\Omega}^{\top} \mtx{B}(t) \mtx{\Omega})^{\dagger} (\mtx{B}(t) \mtx{\Omega})^{\top}.
    \label{equ:nystrom-approximation}
\end{equation}
Related randomized low-rank approximations for general parameter-dependent matrices have been analyzed in~\cite{kressner-2024-randomized-lowrank}.
Exploiting the invariance of the trace under cyclic permutations, taking the trace of~\cref{equ:nystrom-approximation} gives rise to the Nyström estimator
\begin{equation}
    \Trace(\mtx{B}(t)) \approx \Trace(\Nystr{\mtx{\Omega}}{\mtx{B}}(t)) = \Trace( (\mtx{\Omega}^{\top} \mtx{B}(t) \mtx{\Omega})^{\dagger} ( \mtx{\Omega}^{\top} \mtx{B}(t)^2 \mtx{\Omega})).
    \label{equ:nystrom-trace-estimator}
\end{equation}
This estimator forms the basis of the so-called spectrum sweeping method \cite{lin-2017-randomized-estimation}. 

\paragraph{Nyström++ estimator}
The Nyström estimator~\cref{equ:nystrom-trace-estimator} relies on the feasibility of low-rank approximation, that is, the singular values of $\mtx{B}(t)$ need to decay sufficiently quickly.
If this property cannot be ensured, it is preferable to use an estimator that combines the Girard-Hutchinson and Nyström estimators.
The combination suggested in~\cite{lin-2017-randomized-estimation} applies Girard-Hutchinson to correct the error of the Nyström approximation~\cref{equ:nystrom-approximation}:
\begin{equation}
    \Nystrpp{\mtx{\Omega}}{\mtx{\Psi}}(\mtx{B}(t)) = \Trace(\Nystr{\mtx{\Omega}}{\mtx{B}}(t)) + \Hutch{\mtx{\Psi}}(\mtx{B}(t) - \Nystr{\mtx{\Omega}}{\mtx{B}}(t)),
    \label{equ:nystrompp-trace-estimator}
\end{equation}
where $\mtx{\Psi} \in \mathbb{R}^{n \times n_{\mtx{\Psi}}}$, $\mtx{\Omega} \in \mathbb{R}^{n \times n_{\mtx{\Omega}}}$ are independent Gaussian random matrices. 

For constant $\mtx{B}$, the estimator~\cref{equ:nystrompp-trace-estimator} coincides with the Nyström++ estimator from~\cite{persson-2022-improved-variants}, which is based on the Hutch++ estimator~\cite{meyer-2021-hutch-optimal}. In this situation (constant $\mtx{B}$), these estimators were both shown to achieve a relative $\varepsilon$-error with $\mathcal{O}(\varepsilon^{-1})$ matrix-vector products only, independent of the singular value decay of $\mtx{B}$.

At this point, we acknowledge the existence of the XNysTrace estimator from \cite{epperly-2024-xtrace-making} which oftentimes outperforms the Nyström++ estimator. Moreover, a straightforward extension of the XNysTrace estimator to the parameter-dependent setting seems to be in reach. However, what distinguishes \cref{alg:nystrom-chebyshev-pp} for estimating spectral densities are two key observations --- using the cyclic invariance of the trace and the affine linear form of Chebyshev approximations (see \cref{subsubsec:chebyshev-nystrom-implementation} for details) --- which it exploits to significantly speed up the computation. Unfortunately, we do not see a way of marrying these observations with the efficient implementation of XNysTrace described in \cite[Section 2.2]{epperly-2024-xtrace-making}, which limits its suitability for spectral density estimation.

\subsection{Error bounds}

In this section, we derive bounds on the $L^1$-error (with respect to $t$) for each of the trace estimators from \cref{subsec:methods}. Our bounds parallel existing results in the constant case. As in~\cite{kressner-2024-randomized-lowrank}, our proofs incorporate parameter dependence by proceeding through moments $\mathbb{E}^p[X] = (\mathbb{E}[|X|^p])^{\sfrac{1}{p}}$ of certain random variables $X$.

\subsubsection{Girard-Hutchinson estimator}
\label{subsec:hutchinson}

We first derive an $L^1$-error bound for the Girard-Hutchinson estimator~\cref{equ:hutchinson-trace-estimator} that is relative to the $L^1$-norms of the Frobenius and spectral norms of $\mtx{B}(t)$, defined by
\begin{equation*}
 \|\mtx{B}(t)\|_{F,1} := \int_{a}^{b} \|\mtx{B}(t)\|_{F}~\mathrm{d}t, \quad  \|\mtx{B}(t)\|_{2,1} := \int_{a}^{b} \|\mtx{B}(t)\|_{2}~\mathrm{d}t.
\end{equation*}

\begin{theorem}[Girard-Hutchinson estimator for parameter-dependent matrices]\label{thm:hutchinson}
For a symmetric matrix $\mtx{B}(t) \in \mathbb{R}^{n \times n}$ that depends continuously on $t \in [a, b]$, consider the  Girard-Hutchinson estimator~\cref{equ:hutchinson-trace-estimator} with an $n\times n_{\mtx{\Psi}}$ Gaussian random matrix $\mtx{\Psi}$. Then for any $p \geq 1$, $p \in \mathbb{R}$, and $\gamma \geq 1$, the following bound holds with probability at least $1 - \gamma^{-p}$:
\begin{equation}
    \int_{a}^{b} |\Trace(\mtx{B}(t)) - \Hutch{\mtx{\Psi}}(\mtx{B}(t))|~\mathrm{d}t \leq 4 \gamma \Big( \sqrt{\frac{p}{n_{\mtx{\Psi}}}}  \|\mtx{B}(t)\|_{F,1} + \frac{2 p}{n_{\mtx{\Psi}}} \|\mtx{B}(t)\|_{2,1} \Big).
    \label{equ:girard-hutchinson-bound}
\end{equation}
In particular, given $\varepsilon \in (0, 1)$ and $\delta \in (0, \sfrac{1}{2}]$, the bound $\int_{a}^{b} | \Trace(\mtx{B}(t)) - \Hutch{\mtx{\Psi}}(\mtx{B}(t)) | ~\mathrm{d}t \allowbreak \leq \varepsilon \|\mtx{B}(t)\|_{F,1}$ holds with probability at least $1-\delta$ if $n_{\mtx{\Psi}} = \mathcal{O}(\varepsilon^{-2} \log(\delta^{-1}))$.
\end{theorem}
\begin{proof}
    Let us first consider $t$ fixed. From the proof of \cite[Theorem 1]{cortinovis-2022-randomized-trace}, we know that $\Trace(\mtx{B}(t)) - \Hutch{\mtx{\Psi}}(\mtx{B}(t))$ is a $(2 \lVert \mtx{B}(t) \rVert _F^2 / n_{\mtx{\Psi}}, 2 \lVert \mtx{B}(t) \rVert _2 / n_{\mtx{\Psi}})$-sub-gamma random variable. This observation allows us to apply
    \cref{lem:sub-gamma-moments}, which spells out moment bounds for sub-gamma random variables and yields
    \begin{equation}
        \mathbb{E}^{p}[\Trace(\mtx{B}(t)) - \Hutch{\mtx{\Psi}}(\mtx{B}(t))] \leq 4  \sqrt{\frac{p}{n_{\mtx{\Psi}}}} \lVert \mtx{B}(t) \rVert _F + 8 \frac{p}{n_{\mtx{\Psi}}} \lVert \mtx{B}(t) \rVert _2.
        \label{equ:moment-bound-hutchinson}
    \end{equation}

    To address the parameter-dependent case, we first note that the continuity assumption implies that $\Trace(\mtx{B}(t)) - \Hutch{\mtx{\Psi}}(\mtx{B}(t))$ is measurable. Therefore, we can apply Minkowski's integral inequality~\cite[Theorem 202]{hardy-1952-inequalities} to conclude from the moment bound~\cref{equ:moment-bound-hutchinson} that
    \begin{align*}
        \mathbb{E}^{p}\left[ \int_{a}^{b} |\Trace(\mtx{B}(t)) - \Hutch{\mtx{\Psi}}(\mtx{B}(t))|~\mathrm{d}t  \right]
        &\leq \int_{a}^{b} \mathbb{E}^{p}\left[ \Trace(\mtx{B}(t)) - \Hutch{\mtx{\Psi}}(\mtx{B}(t)) \right]~\mathrm{d}t \notag \\
        &\leq 4 \sqrt{\frac{p}{n_{\mtx{\Psi}}}} \lVert \mtx{B}(t) \rVert_{F,1} + 8 \frac{p}{n_{\mtx{\Psi}}} \lVert \mtx{B}(t) \rVert_{2,1}.
    \end{align*}
    This implies~\cref{equ:girard-hutchinson-bound} by Markov's inequality.

    We now use $\lVert \mtx{B}(t) \rVert _2 \leq \lVert \mtx{B}(t) \rVert _F$ and determine $n_{\mtx{\Psi}}$, such that for some $\gamma, p \geq 1$
    \begin{equation}
    4 \gamma \left( \sqrt{\frac{p}{n_{\mtx{\Psi}}}} + \frac{2p}{n_{\mtx{\Psi}}} \right) \leq \varepsilon.
    \label{equ:girard-hutchinson-bound-constants}
    \end{equation}
    The left-hand side of \cref{equ:girard-hutchinson-bound-constants} monotonically decreases with $n_{\mtx{\Psi}}$. Equality is achieved for
    \begin{equation*}
    n_{\mtx{\Psi}} = 4 \left( 1 + \sqrt{1 + 2 \varepsilon/\gamma} \right)^2 \varepsilon^{-2} p \gamma^2.
    \end{equation*}
    We define the failure probability $\delta = \gamma^{-p} \in (0, \sfrac{1}{2}]$. Since $\varepsilon / \gamma \leq 1$, we can bound $(1 + \sqrt{1 + 2 \varepsilon/\gamma})^2 \leq (1 + \sqrt{3})^2$. The factor $p \gamma^{2} = p \delta^{-\sfrac{2}{p}}$ takes its minimum when $p = 2 \log(\delta^{-1})$, at which point it assumes the value $2 \log(\delta^{-1})e$. With these observations, we see that inequality \cref{equ:girard-hutchinson-bound-constants} still holds for
    \begin{equation*}
    n_{\mtx{\Psi}} \geq 8 e \big( 1 + \sqrt{3} \big)^2 \varepsilon^{-2} \log(\delta^{-1}) = \mathcal{O}(\varepsilon^{-2} \log(\delta^{-1})).
    \end{equation*}
\end{proof}

When setting $\mtx{B}(t) \equiv \mtx{B}$ for a constant matrix $\mtx{B}$, \cref{thm:hutchinson} exactly matches an existing bound for the Girard-Hutchinson estimator; see \cite[Lemma 2.1]{meyer-2021-hutch-optimal}. 

\begin{remark}
We note the following:
\begin{enumerate}
\item A variant of \cref{thm:hutchinson} can be derived when assuming that $\mtx{B}(t)$ is SPSD and nonzero for every $t \in [a, b]$. After dividing both sides of \cref{equ:moment-bound-hutchinson} by $\Trace(\mtx{B}(t))$, one can apply the remaining steps of the proof to deduce that
\begin{equation*}
    \int_{a}^{b} \frac{| \Trace(\mtx{B}(t)) - \Hutch{\mtx{\Psi}}(\mtx{B}(t)) |}{\Trace(\mtx{B}(t))} ~\mathrm{d}t < \varepsilon
\end{equation*}
holds with probability at least $1 - \delta$ if $n_{\mtx{\Psi}} = \mathcal{O}(\log(\delta^{-1}) \varepsilon^{-2} \rho^2)$. Here, 
the quantity $
    \rho := \frac{1}{b-a} \int_{a}^{b} \frac{\lVert \mtx{B}(t) \rVert _F}{\Trace(\mtx{B}(t))} ~\mathrm{d}t \in [1,n]
$
is small when the singular values of $\mtx{B}(t)$ decay slowly for the majority of $t\in [a,b]$ and large otherwise. This parallels a similar result for constant matrices; see~\cite[Remark 2]{cortinovis-2022-randomized-trace}.
\item Using the observations of \cite{epperly-2024-note-self}, \cref{thm:hutchinson} can be restated in the same form for random vectors $\vct{\psi}_1, \dots, \vct{\psi}_{n_{\mtx{\Psi}}} \in \mathbb{R}^{n}$ drawn uniformly from the $(n-1)$-sphere of radius $\sqrt{n}$. In the bound, $\mtx{B}(t)$ is replaced by $\overline{\mtx{B}}(t) := \mtx{B}(t) - n^{-1} \Trace(\mtx{B}(t)) \mtx{I}$.
\end{enumerate}
\end{remark}

\subsubsection{Nyström estimator}
\label{subsec:nystrom}

As a consequence of~\cite{kressner-2024-randomized-lowrank}, the following result confirms the intuition that the Nystr\"om estimator is very well suited when the singular values of $\mtx{B}(t)$ decay quickly. More specifically, it shows that the $L^1$-error of the estimator stays on the level of the best rank-$r$ approximation error of $B(t)$, measured in the nuclear norm.

\begin{theorem}[Nyström estimator for parameter-dependent matrices]\label{thm:nystrom}
    Let $\mtx{B}(t) \allowbreak \in \mathbb{R}^{n \times n}$ be SPSD and continuous for $t \in [a, b]$. For integers $r \geq 2$ and $n_{\mtx{\Omega}} \geq r + 4$, consider the Nystr\"om approximation $\Nystr{\mtx{\Omega}}{\mtx{B}}(t)$ defined in~\cref{equ:nystrom-approximation} with an $n\times n_{\mtx{\Omega}}$ Gaussian random matrix $\mtx{\Omega}$. Then, for any $\gamma \geq 1$, the bound 
    \begin{equation*}
        \int_{a}^{b} \big| \Trace(\mtx{B}(t)) - \Trace(\Nystr{\mtx{\Omega}}{\mtx{B}}(t)) \big| ~\mathrm{d}t
        \leq \gamma^2 (1 + r) \int_{a}^{b} \sum_{i = r+1}^{n} \sigma_i(\mtx{B}(t)) ~\mathrm{d}t.
    \end{equation*}
    holds with probability at least $1 - \gamma^{-(n_{\mtx{\Omega}} - r)}$.
\end{theorem}
\begin{proof}
Using that $\mtx{B}(t) - \Nystr{\mtx{\Omega}}{\mtx{B}}(t)$ is SPSD (see, e.g.,~\cite[Lemma 2.1]{frangella-2023-randomized-nystrom}) we obtain
    \begin{align*}
        \big| \Trace(\mtx{B}(t)) - \Trace(\Nystr{\mtx{\Omega}}{\mtx{B}}(t)) \big|
        & = \big| \Trace(\mtx{B}(t) - \Nystr{\mtx{\Omega}}{\mtx{B}}(t)) \big|
        = \lVert \mtx{B}(t) - \Nystr{\mtx{\Omega}}{\mtx{B}}(t) \rVert _{\ast} \\
        & = \lVert (\mtx{I}_n - \mtx{\Pi}_{\mtx{B}(t)^{\sfrac{1}{2}} \mtx{\Omega}}) \mtx{B}(t)^{\sfrac{1}{2}} \rVert _F^2,
    \end{align*}
    where $\lVert \cdot \rVert _{\ast}$ denotes the nuclear norm and we used \cite[Theorem 1]{gittens-2011-spectral-norm} in the last equality. This allows us to apply \cite[Theorem 5]{kressner-2024-randomized-lowrank}, which immediately implies the statement claimed by the theorem.
\end{proof}

Compared to the constant case~\cite[Theorem 8.1]{tropp-2023-randomized-algorithms}, the bound of \cref{thm:nystrom} features an additional factor $1+r$, which cannot straightforwardly be compensated for with oversampling.

\subsection{Nyström++ estimator for parameter-dependent matrices}
\label{subsec:nystrom-pp}

\cref{thm:hutchinson} shows that the Girard-Hutchinson estimator achieves a relative $\varepsilon$-error when using $\mathcal{O}(\varepsilon^{-2})$ queries. The aim of this section is to show that the Nyström++ estimator~\cref{equ:nystrompp-trace-estimator} improves this to $\mathcal{O}(\varepsilon^{-1})$ queries (that is, matrix-vector products with random vectors), \emph{without} requiring any assumption on the singular values of $\mtx{B}(t)$. This parallels existing results for the constant case~\cite{meyer-2021-hutch-optimal,persson-2022-improved-variants}. In particular, \cite[Theorem 3.4]{persson-2022-improved-variants} shows that $|\Trace(\mtx{B}) - \Nystrpp{\mtx{\Omega}}{\mtx{\Psi}}(\mtx{B})| \leq \varepsilon \Trace(\mtx{B})$ holds with high probability for (constant) SPSD $\mtx{B}$ when $n_{\mtx{\Psi}} = n_{\mtx{\Omega}} = \mathcal{O}(\varepsilon^{-1})$ in the Nyström++ estimator $\Nystrpp{\mtx{\Omega}}{\mtx{\Psi}}(\mtx{B})$ from~\cref{equ:nystrompp-trace-estimator}. 

Our analysis mimics the strategy from~\cite{meyer-2021-hutch-optimal}, by first establishing a relative Frobenius norm error proportional to $n_{\mtx{\Omega}}^{-\sfrac{1}{2}}$ for the parameter-dependent Nyström approximation.

\begin{lemma}[Nyström approximation for parameter-dependent matrices]\label{lem:nystrom}
    Let $\mtx{B}(t) \in \mathbb{R}^{n \times n}$ be SPSD and continuous for $t \in [a, b]$. Consider the Nyström approximation $\Nystr{\mtx{\Omega}}{\mtx{B}}(t)$ with a Gaussian random matrix $\mtx{\Omega} \in \mathbb{R}^{n \times n_{\mtx{\Omega}}}$ and \emph{even} $n_{\mtx{\Omega}} \geq 4$.
    Then, for any $\gamma \geq 1$, the bound
    \begin{equation}
        \int_{a}^{b} \lVert \mtx{B}(t) - \Nystr{\mtx{\Omega}}{\mtx{B}}(t) \rVert _F~\mathrm{d}t \leq \gamma \frac{c}{\sqrt{n_{\mtx{\Omega}}}} \int_{a}^{b} \Trace(\mtx{B}(t))~\mathrm{d}t
        \label{equ:nystrompp-theorem-bound}
    \end{equation}
        holds with probability at least $1 - \gamma^{-\sfrac{n_{\mtx{\Omega}}}{4}}$ for 
        $c = 154$. In particular, given $\varepsilon > 0$ and $\delta \in (0, 1)$, the bound $\int_{a}^{b} \lVert \mtx{B}(t) - \Nystr{\mtx{\Omega}}{\mtx{B}}(t) \rVert _F~\mathrm{d}t \leq \varepsilon \int_{a}^{b} \Trace(\mtx{B}(t))~\mathrm{d}t$ holds with probability at least $1-\delta$ if $n_{\mtx{\Omega}} = \mathcal{O}(\varepsilon^{-2} + \log(\delta^{-1}))$.
\end{lemma}

\begin{proof}
For fixed $t$, consider the spectral decomposition $\mtx{B}(t) = \mtx{U} \mtx{\Lambda} \mtx{U}^{\top}$ with  
$\mtx{\Lambda} = \operatorname{diag}(\lambda_1, \dots, \lambda_n)$ and $\lambda_1 \ge \dots \ge \lambda_n \ge 0$. For $k < n_{\mtx{\Omega}}$ (chosen arbitrarily at the moment), partition
    \begin{equation*}
        \rule[\dimexpr-2ex-\ht\strutbox]{0pt}{\dimexpr2ex+4ex+\baselineskip}
        \mtx{U} = \begin{bmatrix}
            \smash{\underbrace{\mtx{U}_1}_{n \times k}} & \smash{\underbrace{\mtx{U}_2}_{n \times (n-k)}}
        \end{bmatrix}
        \quad \text{and} \quad
        \mtx{\Lambda} =
        \begin{bmatrix}
            \smash{\overbrace{\mtx{\Lambda}_1}^{k \times k}} & \\ & \smash{\underbrace{\mtx{\Lambda}_2}_{(n-k) \times (n-k)}}
        \end{bmatrix},
    \end{equation*}
    we set
    $\mtx{\Omega}_1 := \mtx{U}_1^{\top} \mtx{\Omega} \in \mathbb{R}^{k \times n_{\mtx{\Omega}}}$ and $\mtx{\Omega}_2 := \mtx{U}_2^{\top} \mtx{\Omega} \in \mathbb{R}^{(n - k) \times n_{\mtx{\Omega}}}$, which are independent Gaussian random matrices.
Applying \cite[Theorem B.1]{persson-2023-randomized-lowrank} for $f(x) = x$, see also proof of \cite[Corollary 8.2]{tropp-2023-randomized-algorithms}, yields the bound
    \begin{equation}
        \lVert \mtx{B}(t) - \Nystr{\mtx{\Omega}}{\mtx{B}}(t) \rVert _F 
        \leq  \lVert \mtx{\Lambda}_2 \rVert _F + \lVert \mtx{\Lambda}_2^{\sfrac{1}{2}} \mtx{\Omega}_2 \mtx{\Omega}_1^{\dagger} \rVert _{(4)}^2,
        \label{equ:nystrom-proof-persson-bonud}
    \end{equation}
    where $\lVert \cdot \rVert _{(4)}$ denotes the Schatten-4 norm. Following the proof of~\cite[Lemma 3]{meyer-2021-hutch-optimal}, we replace the the first term by a simpler bound:
    \begin{equation}
        \lVert \mtx{\Lambda}_2 \rVert _F
        \leq \sqrt{ \lambda_{k+1} (  \lambda_{k+1} + \cdots + \lambda_{n})}
        \leq \sqrt{ \Trace(\mtx{B}) / k \cdot \Trace(\mtx{B})}
        \leq \Trace(\mtx{B}) / \sqrt{k}.
        \label{equ:nystrom-proof-frobenius-trace}
    \end{equation}
    
To treat the second term in~\cref{equ:nystrom-proof-persson-bonud}, we bound 
its $q$th moment using a standard Schatten norm bound and submultiplicativity of the spectral norm along with the independence of $\mtx{\Omega}_1$ and $\mtx{\Omega}_2$:
    \begin{equation}
        \mathbb{E}^{q}\left[ \big\| \mtx{\Lambda}_2^{\sfrac{1}{2}} \mtx{\Omega}_2 \mtx{\Omega}_1^{\dagger} \big\|_{(4)}^2 \right]
        \leq \mathbb{E}^{q}\left[ \sqrt{k} \big\| \mtx{\Lambda}_2^{\sfrac{1}{2}} \mtx{\Omega}_2 \mtx{\Omega}_1^{\dagger} \big\| _2^2 \right]
        \leq \sqrt{k} \mathbb{E}^{q}\left[ \big\| \mtx{\Lambda}_2^{\sfrac{1}{2}} \mtx{\Omega}_2  \big\|_2^2\right] \mathbb{E}^{q}\left[ \big\| \mtx{\Omega}_1^{\dagger} \big\|_2^2  \right].
        \label{equ:nystrom-proof-processed-tail}
    \end{equation}
    To continue from here, we set $k = n_{\mtx{\Omega}}/2$ and $q = n_{\mtx{\Omega}}/4$. By the calculations after~\cite[Eqn~(B.7)]{tropp-2023-randomized-algorithms}, one has
    \begin{equation*}
     \mathbb{E}\left[ \big\| ( \mtx{\Omega}_1 \mtx{\Omega}_1^{\top} )^{-1} \big\|_2^{\sfrac{n_{\mtx{\Omega}}}{4}} \right]
        \leq
        \left(1 + \frac{n_{\mtx{\Omega}}}{2}\right)
        \left( \frac{3}{4} n_{\mtx{\Omega}}\right)^{\sfrac{n_{\mtx{\Omega}}}{4}}
        \big( ( n_{\mtx{\Omega}} / 2 + 1)!\big)^{-\frac{n_{\mtx{\Omega}}}{2+n_{\mtx{\Omega}}}}.
    \end{equation*}
    Therefore,
    \begin{equation}
      \mathbb{E}^{\sfrac{n_{\mtx{\Omega}}}{4}}\left[ \big\| \mtx{\Omega}_1^{\dagger} \big\|_2^2 \right]
        = \mathbb{E}\left[ \big\| ( \mtx{\Omega}_1 \mtx{\Omega}_1^{\top} )^{-1} \big\|_2^{\sfrac{n_{\mtx{\Omega}}}{4}} \right]^{\sfrac{4}{n_{\mtx{\Omega}}}}
        \le \frac{3}{4}  
        \frac{e^4 n_{\mtx{\Omega}}}{(n_{\mtx{\Omega}} / 2+ 1)^2} \le 
        \frac{3}{4}  \frac{e^4}{n_{\mtx{\Omega}}},
        \label{equ:pinv-spectral-norm-bound}
    \end{equation}
    where we used $(1 + m)^{\sfrac{1}{m}} \leq e$ and $(m!)^{-\sfrac{1}{m}} \leq e/m$. Note that, in contrast 
    to the result of~\cite[Lemma B.3]{tropp-2023-randomized-algorithms}, this inequality is valid for arbitrarily large $n_{\mtx{\Omega}}$, at the expense of a slightly larger constant.

    The other factor in~\cref{equ:nystrom-proof-processed-tail} is treated using moment bounds for the spectral norm of a non-standard Gaussian random matrix (\cref{lem:spectral-norm-moment} with $\mtx{A} = \mtx{\Lambda}_2^{\sfrac{1}{2}}$):
    \begin{equation*}
        \mathbb{E}^{\sfrac{n_{\mtx{\Omega}}}{4}}\left[ \big\| \mtx{\Lambda}_2^{\sfrac{1}{2}} \mtx{\Omega}_2 \big\|_2^2 \right]
        \leq \frac{5}{4} n_{\mtx{\Omega}} \Big( 2 \big\| \mtx{\Lambda}_2^{\sfrac{1}{2}} \big\|_2^2 + \frac{1}{n_{\mtx{\Omega}}} \big\| \mtx{\Lambda}_2^{\sfrac{1}{2}} \big\|_F^2 \Big).
    \end{equation*}
    Inserting this inequality and~\cref{equ:pinv-spectral-norm-bound} into~\cref{equ:nystrom-proof-processed-tail} gives
    \begin{equation*}
        \mathbb{E}^{\sfrac{n_{\mtx{\Omega}}}{4}}\left[ \lVert \mtx{\Lambda}_2^{\sfrac{1}{2}} \mtx{\Omega}_2 \mtx{\Omega}_1^{\dagger} \rVert _{(4)}^2 \right]
        \leq \frac{15 e^4}{16}  \sqrt{n_{\mtx{\Omega}}} \Big( 2 \lVert \mtx{\Lambda}_2^{\sfrac{1}{2}} \rVert _2^2 + \frac{1}{n_{\mtx{\Omega}}} \lVert \mtx{\Lambda}_2^{\sfrac{1}{2}} \rVert _F^2 \Big).
    \end{equation*}
    Bounding $\lVert \mtx{\Lambda}_2^{\sfrac{1}{2}} \rVert _2^2 = \lambda_{k+1}  \leq \Trace(\mtx{B})/k$ and $\lVert \mtx{\Lambda}_2^{\sfrac{1}{2}} \rVert _F^2 = \Trace(\mtx{\Lambda}_2) \le \Trace(\mtx{B})$ (recall that $k = n_{\mtx{\Omega}}/2$) yields
    \begin{equation}
        \mathbb{E}^{\sfrac{n_{\mtx{\Omega}}}{4}}\left[ \lVert \mtx{\Lambda}_2^{\sfrac{1}{2}} \mtx{\Omega}_2 \mtx{\Omega}_1^{\dagger} \rVert _{(4)}^2 \right]
        \leq \frac{15 e^4}{16}  \sqrt{n_{\mtx{\Omega}}} \Big( \frac{2}{n_{\mtx{\Omega}}} \Trace(\mtx{B}) + \frac{1}{n_{\mtx{\Omega}}} \Trace(\mtx{B}) \Big)
        \leq  \frac{154}{\sqrt{n_{\mtx{\Omega}}}} \Trace(\mtx{B}).
        \label{equ:nystrom-proof-tail-bound}
    \end{equation}
    
    Inserting~\cref{equ:nystrom-proof-tail-bound} along with \cref{equ:nystrom-proof-frobenius-trace} in \cref{equ:nystrom-proof-persson-bonud}, letting $c=154$, and using the triangle inequality for $\mathbb{E}^{\sfrac{n_{\mtx{\Omega}}}{4}}[\cdot]$, we obtain
    \begin{equation}
        \mathbb{E}^{\sfrac{n_{\mtx{\Omega}}}{4}} \left[\lVert \mtx{B} - \Nystr{\mtx{\Omega}}{\mtx{B}} \rVert _F \right]
        \leq \mathbb{E}^{\sfrac{n_{\mtx{\Omega}}}{4}} \left[ \lVert \mtx{\Lambda}_2 \rVert _F \right] + \mathbb{E}^{\sfrac{n_{\mtx{\Omega}}}{4}} \left[ \lVert \mtx{\Lambda}_2^{\sfrac{1}{2}} \mtx{\Omega}_2 \mtx{\Omega}_1^{\dagger} \rVert _{(4)}^2 \right]
        \leq \frac{c}{\sqrt{n_{\mtx{\Omega}}}} \Trace(\mtx{B}).
        \label{equ:frobenius-moment-bound}
    \end{equation}
    As in \cite{kressner-2024-randomized-lowrank}, one can show that the error $\lVert \mtx{B}(t) - \Nystr{\mtx{\Omega}}{\mtx{B}}(t) \rVert _F$ is measurable.
    As in the proof of~\cref{thm:hutchinson}, the claimed bound~\cref{equ:nystrompp-theorem-bound}
    follows from~\cref{equ:frobenius-moment-bound} using Minkowski's integral inequality and Markov's inequality.

    Fixing $\gamma = e$ and letting $n_{\mtx{\Omega}} = \lceil(c e \varepsilon)^{-2}  + 4 \log(\delta^{-1}) \rceil$ establishes the second part of the theorem.
\end{proof}

Combining the bound on the Girard-Hutchinson estimator from~\cref{thm:hutchinson} with the result of~\cref{lem:nystrom} finally gives the desired error bound for the Nyström++ trace estimator.
\begin{theorem}[Nyström++ trace estimator for parameter-dependent matrices]\label{thm:nystrom-pp}
    Let $\mtx{B}(t) \in \mathbb{R}^{n \times n}$ be SPSD and continuous for $t \in [a, b]$. Choose $\delta \in (0, 1]$ and $\varepsilon \in (0, 1)$.
    Then, for $n_{\mtx{\Psi}} = n_{\mtx{\Omega}} = \mathcal{O}(\varepsilon^{-1} \log(\delta^{-1}))$,
    the inequality 
    \begin{equation*}
        \int_{a}^{b} | \Trace(\mtx{B}(t)) - \Nystrpp{\mtx{\Omega}}{\mtx{\Psi}}(\mtx{B}(t)) |~\mathrm{d}t
        \leq \varepsilon \int_{a}^{b} \Trace(\mtx{B}(t))~\mathrm{d}t
    \end{equation*}
    holds with probability at least $1 - \delta$.
\end{theorem}
\begin{proof}
    Choosing $n_{\mtx{\Psi}} = n_{\mtx{\Omega}} = \mathcal{O}(\tilde{\varepsilon}^{-2} \log(\tilde{\delta}^{-1}))$, \cref{thm:hutchinson} and~\cref{lem:nystrom} imply that
    \begin{align*}
        \int_{a}^{b} | \Trace(\mtx{B}(t)) - \Nystrpp{\mtx{\Omega}}{\mtx{\Psi}}(\mtx{B}(t)) |~\mathrm{d}t 
        &= \int_{a}^{b} | \Trace(\mtx{B}(t) - \Nystr{\mtx{\Omega}}{\mtx{B}}(t)) - \Hutch{\mtx{\Psi}}(\mtx{B}(t) - \Nystr{\mtx{\Omega}}{\mtx{B}}(t)) |~\mathrm{d}t  \\
        \leq\ & \tilde{\varepsilon} \int_{a}^{b} \lVert \mtx{B}(t) - \Nystr{\mtx{\Omega}}{\mtx{B}}(t) \rVert _F ~\mathrm{d}t       \leq \tilde{\varepsilon}^2 \int_{a}^{b} \Trace(\mtx{B}(t)) ~\mathrm{d}t 
    \end{align*}
    holds with probability at least $1 - 2\tilde{\delta}$. Taking $\varepsilon = \tilde{\varepsilon}^2$ and $\delta = 2 \tilde{\delta}$ concludes the proof.
\end{proof}

Comparing \cref{thm:nystrom-pp} with the analogous result for constant matrices \cite[Theorem 3.4]{persson-2022-improved-variants}, we notice an additional factor of $\sqrt{\log(\delta^{-1})}$ in the required number of random vectors $n_{\mtx{\Omega}}$ and $n_{\mtx{\Psi}}$.

%% file: spectraldensity.tex
\section{Application to spectral density approximation}
\label{sec:application}

We now get back to the task of computing the spectral density $1/n \sum_{i=1}^{n} \delta(t - \lambda_i)$ of a real symmetric matrix $\mtx{A} \in \mathbb{R}^{n \times n}$ with eigenvalues $\lambda_1, \dots, \lambda_n$, where $\delta(\cdot)$ denotes the Dirac delta distribution. As explained in the introduction, finding accurate approximations to this distribution is challenging, and we will therefore work with
a smoothed version.

\subsection{Smoothed spectral density}
\label{subsec:spectral-density}

We recall that the smoothed spectral density $\phi_{\sigma}$ is defined as
\begin{equation}
    \phi_{\sigma}(t) = \frac{1}{n} \sum_{i=1}^{n} g_{\sigma}(t - \lambda_i) = \frac{1}{n} \Trace(g_{\sigma}(t\mtx{I}_n - \mtx{A}))
    \label{equ:smooth-spectral-density}
\end{equation}
for a smoothing kernel $g_{\sigma}$ parametrized by a smoothing parameter $\sigma > 0$. In this work, we choose a Gaussian smoothing kernel of width $\sigma > 0$:
\begin{equation}
    g_{\sigma}(s) = \frac{1}{\sigma \sqrt{2\pi}} e^{-\sfrac{s^2}{2\sigma^2}}.
    \label{equ:smoothing-kernel}
\end{equation}
The expression~\cref{equ:smooth-spectral-density} gives rise to the task of estimating the trace of the parameter-dependent SPSD matrix $g_{\sigma}(t\mtx{I}_n - \mtx{A})$, which has been analyzed in \cref{sec:analysis}.

Choosing the smoothing parameter $\sigma$ results in the following trade off: A larger $\sigma$ usually makes it  easier to approximate the matrix function $g_{\sigma}(t\mtx{I}_n - \mtx{A})$, whereas a smaller $\sigma$ allow one to stay closer to the original spectral density $\phi$. To quantify the latter, one can measure the error between $\phi$ and its smoothened version $\phi_{\sigma}$ by
\begin{equation}
    \sup_{f \in \mathcal{S}} \int_{a}^{b} f(t) (\phi(t) - \phi_{\sigma}(t))~\mathrm{d}t,
    \label{equ:error-metric}
\end{equation}
where $\mathcal{S}$ represents an appropriately chosen space of test functions. Among others, the reference~\cite{lin-2016-approximating-spectral} uses $\mathcal{S} = \{ f: f(t) \equiv g_{\sigma}(t - s), s \in [a, b]\}$, \cite{chen-2021-analysis-stochastic} uses $\mathcal{S} = \{f : f(t) = \Theta(s - t), s \in [a, b] \}$ with the Heaviside step function $\Theta$, and \cite{braverman-2022-sublinear-timea,bhattacharjee-2025-improved-spectral} use $\mathcal{S} = \{f : |f(t) - f(s)| \leq |t - s| \}$ for which the error \cref{equ:error-metric} is equivalent to the Wasserstein-1 distance between $\phi$ and $\phi_{\sigma}$.
When choosing 
the last metric,
assuming that the eigenvalues are somewhat uniformly distributed and $\sigma$ is rather small, a simple calculation allows one to conclude that the smoothing error \cref{equ:error-metric} is of the order $\sigma$. Consequently, to ensure that the approximation deviates by at most a factor of $\varepsilon > 0$, one needs to choose $\sigma \approx \varepsilon$.

Having discussed the error incurred by smoothing the spectral density $\phi$, we now focus our analysis on the $L^1$-norm approximation error for the smoothed spectral density $\phi_{\sigma}$ itself.

\subsection{Chebyshev approximation of smoothing kernel}
\label{subsec:chebyshev-approximation}

To estimate the trace in \cref{equ:smooth-spectral-density} efficiently, we approximate the function 
$s\mapsto g_{\sigma}(t- s)$ by a
Chebyshev approximation $g_{\sigma}^{(m)}(t-s) = 
\sum_{l=0}^{m} \mu_l(t) T_l(s)$
with the approximation coefficients $\mu_l(t) \in \mathbb R$ and the Chebyshev polynomials
$T_l(s) = \cos(l \cdot \arccos(s))$ for $l = 0,1,\ldots,m$. We choose the coefficients $\mu_l(t)$ in such a way that $g_{\sigma}^{(m)}(t-s)$ interpolates $g_{\sigma}(t-s)$ at the Chebyshev nodes $s_i = \cos(\pi i / m),~i=0,1,\dots,m$. In this case, $g_{\sigma}^{(m)}$ is often referred to as the Chebyshev interpolant~\cite[Chapter 2]{trefethen-2020-approximation-theory}.
Assuming that the spectrum of $\mtx{A}$ is contained in $[-1,1]$, this gives rise to
\begin{equation}
    g_{\sigma}(t\mtx{I}_n - \mtx{A}) \approx g_{\sigma}^{(m)}(t\mtx{I}_n - \mtx{A}) = \sum_{l=0}^{m} \mu_l(t) T_l(\mtx{A}).
    \label{equ:matrix-approximation}
\end{equation}
Inserting~\cref{equ:matrix-approximation} into~\cref{equ:smooth-spectral-density} yields the approximate smoothed spectral density
\begin{equation}
    \phi_{\sigma}(t) \approx \phi_{\sigma}^{(m)}(t) =  \frac{1}{n} \Trace(g_{\sigma}^{(m)}(t \mtx{I}_n - \mtx{A})).
    \label{equ:expanded-spectral-density}
\end{equation}
\begin{remark}
    We note the following:
    \begin{enumerate}
    \item For a general symmetric matrix $\mtx{A}$, one first estimates an interval $[a,b]$ containing the spectrum of $\mtx{A}$~\cite{zhou-2011-bounding-spectrum}. The usual affine transform
    $\tau(t) = \frac{2t - a - b}{b - a}$
    then lets us define a matrix $\bar{\mtx{A}} = \tau(\mtx{A})$ for which the spectrum is contained in $[-1, 1]$ and, hence, the approximation \cref{equ:matrix-approximation} is sensible.\footnote{Note that one has to adjust the smoothing parameter $\bar{\sigma} = 2 \sigma (b - a)^{-1}$ to obtain an undistorted approximation of the original spectrum when inverting the transformation.} This allows us to assume, from now on and without loss of generality, that the spectrum of $\mtx{A}$ is contained in $[-1, 1]$.
    \item The Chebyshev approximation approach more generally applies to estimating the trace of any parametrized matrix function $f_t(\mtx{A})$ for which $t \mapsto f_t(s)$ is continuous and $s \mapsto f_t(s)$ analytic. One such example is $\exp(-t \mtx{A})$, whose trace can be used to, e.g., describe the partition function of a thermodynamic system \cite{chen-2023-krylovaware-stochastic}. However, to obtain explicit bounds for spectral density estimation, we specialize our derivations to $f_t(\mtx{A}) = g_{\sigma}(t \mtx{I} - \mtx{A})$.
    \end{enumerate}
\end{remark}

\subsubsection{Approximation error}
\label{subsubsec:approximation-error}

We now analyze the impact of the Chebyshev approximation of the Gaussian smoothing kernel $g_{\sigma}$ in the approximation of the spectral density $\phi_{\sigma}$. We start with a standard error estimate concerning the approximation of $g_{\sigma}$.

\begin{lemma}[Chebyshev approximation error]\label{lem:chebyshev-error}
    The Chebyshev approximation $g_{\sigma}^{(m)}$ of the Gaussian smoothing kernel $g_{{\sigma}}$ from  \cref{equ:smoothing-kernel} satisfies
    \begin{equation}
        \sup_{s,t \in [-1, 1]} \big| g_{\sigma}(t - s) - g_{\sigma}^{(m)}(t - s) \big| \leq 2\sqrt{\frac{2e}{\pi}} \frac{1}{\sigma^2} (1 + \sigma)^{-m} =: E_{\sigma, m}
        \label{equ:chebyshev-interpolation-sup-error-kernel}
    \end{equation}
    for every degree $m \in \mathbb{N}$ and smoothing parameter $\sigma > 0$.
\end{lemma}
\begin{proof}
    The proof of this result is standard. Let $\mathcal{E}_{\chi} \subset \mathbb{C}$ denote the Bernstein ellipse  with foci $\{-1, +1\}$ and the sum of semi-axes equal to $\chi > 1$. By \cite[Theorem 8.2]{trefethen-2020-approximation-theory},
    \begin{equation*}
        \sup_{s \in [-1, 1]} \big| g_{\sigma}(t - s) - g_{\sigma}^{(m)}(t - s) \big| \leq \frac{4}{\chi^m (\chi - 1)} \max_{z \in \mathcal{E}_{\chi}} |g_{\sigma}(t - z)|,
        \label{equ:bernstein-bound}
    \end{equation*}
    where used the analyticity of $g_{\sigma}(t- \cdot)$.
    For $z = x + \mathrm{i} y \in \mathcal{E}_{\chi}$, we have
    \begin{equation*}
    | g_{\sigma}(t - z) | 
    = \frac{1}{\sigma \sqrt{2 \pi}} \big| e^{-\sfrac{(t - z)^2}{2\sigma^2}} \big|
    = \frac{1}{\sigma \sqrt{2 \pi}} e^{-\sfrac{(t - x)^2}{2\sigma^2}}e^{\sfrac{y^2}{2\sigma^2}}
    \le \frac{1}{\sigma \sqrt{2 \pi}} e^{\sfrac{(\chi - \chi^{-1})^2}{8 \sigma^2}},
    \end{equation*}
    where the last inequality uses 
    $e^{-\sfrac{(t - x)^2}{2\sigma^2}} \leq 1$ and $|y| \le (\chi - \chi^{-1}) / 2$.
    Choosing $\chi = 1 + \sigma$ and using $\chi - \chi^{-1} \leq 2\sigma$ completes the proof of~\cref{equ:chebyshev-interpolation-sup-error-kernel}.
\end{proof}

\cref{fig:chebyshev-heatmap} confirms that \cref{lem:chebyshev-error} not only predicts the qualitative behavior of the error with respect to $\sigma$ and $m$, but it also provides a relatively good quantitative estimate. A simple algebraic manipulation shows that choosing a degree $m = \mathcal{O}(\sigma^{-1} \log(\sigma^{-2} \varepsilon^{-1}))$ ensures an error $\mathcal{O}(\varepsilon)$. The lemma also allows us to bound the $L^1$-error of the spectral density approximation
$\phi_{\sigma}^{(m)}$ defined in~\cref{equ:expanded-spectral-density}:
\begin{align}
    \int_{-1}^{1} \left| \phi_{\sigma}(t) - \phi_{\sigma}^{(m)}(t) \right|~\mathrm{d}t 
    &= \int_{-1}^{1} \left| \frac{1}{n} \sum_{i=1}^n \big(g_{\sigma}(t - \lambda_i) - g_{\sigma}^{(m)}(t - \lambda_i)\big) \right|~\mathrm{d}t \notag \\
    &\leq \int_{-1}^{1} \max_{i = 1, \dots, n} \big| g_{\sigma}(t - \lambda_i) - g_{\sigma}^{(m)}(t - \lambda_i) \big|~\mathrm{d}t    \leq 2 E_{\sigma, m}.
    \label{equ:chebyshev-interpolation-spectral-density}
\end{align}

\begin{figure}[ht]
    \centering
    \scalebox{0.8}{\input{plots/chebyshev_heatmap.pgf}}
    \caption{Supremum norm error of the degree $m$ Chebyshev interpolation of $g_{\sigma}$. The error is approximated by evaluating the interpolant in $n_t = 1000$ uniformly spaced values on $[-1, 1]$. The dashed line is the contour corresponding to an error of $10^{-6}$. The hatched region shows the combinations of $m$ and $\sigma$ for which the error bound $E_{\sigma, m}$ of \cref{lem:chebyshev-error} remains below $10^{-6}$.}
    \label{fig:chebyshev-heatmap}
\end{figure}

\subsubsection{Preservation of non-negativity} \label{sec:preservnonneg}

In general, the Chebyshev approximation $g_{\sigma}^{(m)}$ does not inherit the  non-negativity of the smoothing kernel $g_{\sigma}$. Consequently, the matrix $g_{\sigma}^{(m)}(t \mtx{I}_n - \mtx{A})$ may become 
(slightly) indefinite. This causes issues when attempting to apply the Nyström++ trace estimator \cref{equ:nystrompp-trace-estimator}. While we have never observed this in our numerical experiments, the Nyström approximation may fail when the SPSD condition is not met; see also the discussion in~\cite{nakatsukasa-2023-randomized-lowrank}. On the theoretical side,~\cref{thm:nystrom-pp}
critically relies on the SPSD condition.

In \cite[Algorithm 6]{braverman-2022-sublinear-timea}, the loss of non-negativity is addressed by introducing Jackson damping, though this comes at the cost of losing the very favorable, exponential convergence. To ensure a non-negative Chebyshev approximation and maintain exponential convergence, we suggest the following simple method: Given an even degree $m$, approximate the square root of the smoothing kernel up to degree $m/2$:
\begin{equation}
    \sqrt{g_{\sigma}}(t - s) \approx \sqrt{g_{\sigma}}^{(\sfrac{m}{2})}(t - s) = \sum_{l=0}^{m/2} \xi_l(t) T_l(s).
    \label{equ:square-root-chebyshev-approximation}
\end{equation}
When squaring this approximation, one obtains a polynomial of degree $m$, which admits an exact approximation
of the form 
$(\sum_{l=0}^{m/2} \xi_l(t) T_l(s))^2 = \sum_{l=0}^{m} \underline{\mu}_l(t) T_l(s)$.
Thus, we obtain the non-negative Chebyshev approximation
\begin{equation}
g_{\sigma}(t - s) \approx    \underline{g}_{\sigma}^{(m)}(t - s) = \sum_{l=0}^{m} \underline{\mu}_l(t) T_l(s).
    \label{equ:non-negative-chebyshev-approximation}
\end{equation}
The corresponding non-negative approximation of the smoothed spectral density,
\begin{equation}
    \phi_{\sigma}(t) \approx \underline{\phi}_{\sigma}^{(m)}(t) := \frac{1}{n} \Trace(\underline{g}_{\sigma}^{(m)}(t\mtx{I}_n - \mtx{A})),
    \label{equ:non-negative-expanded-density}
\end{equation}
satisfies the following bound.

\begin{lemma}[$L^1$-error of non-negative approximate smooth spectral density]\label{lem:non-negative-chebyshev-error}
Let $\phi_{\sigma}$ be the smoothed spectral density defined in~\cref{equ:smooth-spectral-density}
for a symmetric matrix $\mtx{A} \in \mathbb{R}^{n \times n}$ with spectrum contained in $[-1, 1]$. Then the approximation $\underline{\phi}_{\sigma}^{(m)}$ defined in~\cref{equ:non-negative-expanded-density} satisfies
    \begin{equation*}
        \int_{-1}^{1} \big| \phi_{\sigma}(t) - \underline{\phi}_{\sigma}^{(m)}(t) \big|~\mathrm{d}t \leq 4 \sqrt{2} \left(1 + \sigma \sqrt{\pi} \cdot E_{\sqrt{2}\sigma, \sfrac{m}{2}}\right) E_{\sqrt{2}\sigma, \sfrac{m}{2}} =: \underline{E}_{\sigma, m}
        \label{equ:chebyshev-interpolation-sup-error-kernel-nonneg}
    \end{equation*}
    for every \emph{even} degree $m \in \mathbb{N}$, smoothing parameter $\sigma > 0$, and with $E_{\sigma, m}$ defined in \cref{lem:chebyshev-error}.
\end{lemma}
\begin{proof}
    Using that the inequality $| a - b | \leq (2 \sqrt{a} + | \sqrt{a} - \sqrt{b} |)  | \sqrt{a} - \sqrt{b} |$ holds for any non-negative numbers $a, b$, it follows that
    \begin{equation*}
        \big| g_{\sigma} - \underline{g}_{\sigma}^{(m)} \big| \leq \big( 2 \big| \sqrt{g_{\sigma}} \big| + \big| \sqrt{g_{\sigma}} - \sqrt{g_{\sigma}}^{(\sfrac{m}{2})} \big| \big) \big| \sqrt{g_{\sigma}} - \sqrt{g_{\sigma}}^{(\sfrac{m}{2})} \big|.
    \end{equation*}
    Because of $\sqrt{g_{\sigma}} = \sqrt{2 \sigma \sqrt{2 \pi}} \cdot g_{\sqrt{2}\sigma}$,
we can apply \cref{lem:chebyshev-error} with the substitutions
    $\sigma \gets \sqrt{2}\sigma$ and $m \gets m/2$ to bound 
    $\big| \sqrt{g_{\sigma}} - \sqrt{g_{\sigma}}^{(\sfrac{m}{2})} \big|$.
    Together with  $|\sqrt{g_{\sigma}}| \leq 1/\sqrt{\sigma \sqrt{2 \pi}}$, we get
    \begin{align*}
        \big| g_{\sigma}(t - s) - \underline{g}_{\sigma}^{(m)}(t - s) \big| \notag 
        &\leq \Big( 2 \frac{1}{\sqrt{ \sigma \sqrt{2\pi}}} + \sqrt{2 \sigma \sqrt{2 \pi}} \cdot E_{\sqrt{2}\sigma, \sfrac{m}{2}}\Big)\sqrt{2 \sigma \sqrt{2 \pi}} \cdot E_{\sqrt{2}\sigma, \sfrac{m}{2}} \notag \\
        &= 2\sqrt{2} \big(1 + \sigma \sqrt{\pi} \cdot E_{\sqrt{2}\sigma, \sfrac{m}{2}}\big) E_{\sqrt{2}\sigma, \sfrac{m}{2}}
    \end{align*}
    for every $s,t \in [-1, 1]$.
    Bounding the integral as in~\cref{equ:chebyshev-interpolation-spectral-density} then yields the desired result.
\end{proof}

Again, \cref{lem:non-negative-chebyshev-error} implies that choosing a degree $m = \mathcal{O}(\sigma^{-1} \log(\sigma^{-2} \varepsilon^{-1}))$ ensures an error $\mathcal{O}(\varepsilon)$.

As the non-negative approximation~\cref{equ:non-negative-chebyshev-approximation} ``wastes'' half of the degrees of freedom, it can be expected to be less accurate than the standard Chebyshev approximation. On the other hand, $\sqrt{g_{\sigma}} \propto g_{\sqrt{2}\sigma}$ is somewhat easier to approximate because it is less ``pointy" compared to $g_{\sigma}$. Overall, it turns out that one needs to choose the degree of the non-negative approximation roughly a factor $\sqrt{2}$ larger to match the approximation error of the standard Chebyshev approximation (see also \cref{fig:interpolation-issue}). An additional advantage of our approach to preserving non-negativity is that when we apply any of the trace estimators we analyzed in \cref{sec:analysis} to the approximation $\underline{g}_{\sigma}^{(m)}(t\mtx{I} - \mtx{A})$, the resulting spectral density approximation remains non-negative; no post-processing is needed to ensure this property.

\subsubsection{Computation of Chebyshev coefficients}
\label{subsubsec:dct}

It is well known that the coefficients of the Chebyshev approximation \cref{equ:matrix-approximation} can --- in exact arithmetic --- be computed by inverting the (type I) discrete cosine transform (DCT)
\begin{equation*}
    g_{\sigma}(t - \cos(s_i)) = \sum_{l=0}^{m} \mu_l(t) \cos\left(\frac{\pi i l}{m} \right)
    \label{equ:discrete-cosine-transform}
\end{equation*}
between the approximation coefficients $\mu_l(t),~l=0,\dots,m$ and the function evaluations at the Chebyshev points $s_i = \cos(\pi i / m),~i=0,\dots,m$ \cite{baszenski-1997-fast-polynomial, plonka-2018-numerical-fourier, trefethen-2020-approximation-theory}. This can be accomplished in $\mathcal{O}(m \log(m))$ operations.

A Chebyshev approximation can be squared in $\mathcal{O}(m \log(m))$ operations by using the DCT and its inverse \cite{baszenski-1997-fast-polynomial}; see \cref{alg:chebyshev-squaring}. This allows one to efficiently compute the non-negative Chebyshev approximation \cref{equ:non-negative-chebyshev-approximation}, and it will find further use in \cref{subsec:chebyshev-nystrom}. 

\begin{algorithm}
    \caption{Fast squaring of Chebyshev approximations}
    \label{alg:chebyshev-squaring}
    \hspace*{\algorithmicindent} \textbf{Input:}  Coefficients $\{ \mu_l \}_{l=0}^{m}$ of Chebyshev approximation $\sum_{l=0}^{m} \mu_l T_l$\\
    \hspace*{\algorithmicindent} \textbf{Output:} Coefficients $\{ \nu_l \}_{l=0}^{2m}$ such that $\sum_{l=0}^{2m} \nu_l T_l = (\sum_{l=0}^{m} \mu_l T_l)^2${\\[-11pt]}
    \begin{algorithmic}[1]
    \STATE{Compute $\vct{f} = \operatorname{DCT}(\vct{\mu})$}
    \STATE{Compute $[\nu_0, \dots, \nu_{2m}] = \operatorname{DCT}^{-1}(\vct{f} \odot \vct{f})$}
    \label{lin:inverse-DCT}
    \end{algorithmic}
\end{algorithm}

On line \ref{lin:inverse-DCT} of \cref{alg:chebyshev-squaring}, the symbol $\odot$ denotes the entrywise product. 

In \cite{lin-2017-randomized-estimation}, a slightly different approach to approximating $g_{\sigma}$ is taken: Instead of Chebyshev interpolation,
the kernel $g_{\sigma}$ is approximated by a Chebyshev series~\cite[Chapter 4]{trefethen-2020-approximation-theory} truncated after term $m$. In contrast to interpolation, the expansion coefficients involve integrals, which need to be approximated by quadrature. By choosing $2m+1$ Chebyshev quadrature nodes, this approximation can be carried using the fast Fourier transform (FFT). Although interpolation and truncated expansion satisfy nearly identical convergence bounds~\cite[Theorem 8.2]{trefethen-2020-approximation-theory}, we observed the procedure from~\cite{lin-2017-randomized-estimation} to be significantly slower due to need of applying the FFT (instead of DCT) to longer vectors. See~\cref{tab:chebyshev-timing-interpolation}, where we compare a benevolent implementation of \cite[Algorithm 1]{lin-2017-randomized-estimation} with our approach, both for the standard and the non-negative Chebyshev approximation.

\begin{table}[ht]
    \caption{Runtime in milliseconds for computing the coefficients of a degree-$m$ Chebyshev approximation to the smoothing kernel $g_{\sigma}$ in three different ways. \emph{FFT}: Using~\cite[Algorithm 1]{lin-2017-randomized-estimation}.  \emph{DCT}: Applying the discrete cosine transform directly to $g_\sigma$.  \emph{Non-negative DCT}: Applying the discrete cosine transform to $\sqrt{g_\sigma}$ together with \cref{alg:chebyshev-squaring}. We use $\sigma=0.005$ and various degrees $m$. We average over 7 runs of the algorithms and repeat this $1000$ times to form the mean and standard deviation reported below.}
    \label{tab:chebyshev-timing-interpolation}
    \centering
    \input{tables/interpolation.tex}
\end{table}

The suspiciously longer runtime of the FFT-based procedure \cite[Algorithm 1]{lin-2017-randomized-estimation} for $m=1600$ is an artifact of the poor factorizability of the transform size this algorithm requires.

\subsection{Chebyshev-Nyström++ method for spectral density approximation}
\label{subsec:chebyshev-nystrom}

Combining the results of \cref{sec:analysis} and \cref{subsec:chebyshev-approximation}, we apply the Nyström++ estimator $\Nystrpp{\mtx{\Omega}}{\mtx{\Psi}}$ from~\cref{equ:nystrompp-trace-estimator}
to the non-negative 
Chebyshev approximation, which yields the approximation
\begin{equation}
    \phi_{\sigma}(t) \approx \widetilde{\underline{\phi}}_{\sigma}^{(m)}(t) := \frac{1}{n} \Nystrpp{\mtx{\Omega}}{\mtx{\Psi}}(\mtx{B}(t)), \quad \mtx{B}(t) = \underline{g}_{\sigma}^{(m)}(t \mtx{I}_n - \mtx{A}).
    \label{equ:chebyshev-nystrom-formula}
\end{equation}

\subsubsection{Implementation}
\label{subsubsec:chebyshev-nystrom-implementation}

One usually aims at evaluating~\cref{equ:chebyshev-nystrom-formula} for many values of $t$. In the following, we explain how this can be done efficiently by exploiting 
the well-known three-term recurrence relation of Chebyshev polynomials:
\begin{equation}
    T_0(s) = 1,\quad  T_1(s) = s,\quad 
    T_l(s) = 2 s T_{l-1}(s) - T_{l-2}(s) \text{ for $l \geq 2$.}
    \label{equ:chebyshev-recurrence}
\end{equation}

We start with the first term in the definition~\cref{equ:nystrompp-trace-estimator} of $\Nystrpp{\mtx{\Omega}}{\mtx{\Psi}}$, which requires evaluating the trace of the 
Nyström approximation \cref{equ:nystrom-approximation}:
\begin{equation}
    \frac{1}{n} \Trace \big((\mtx{B}(t) \mtx{\Omega}) (\mtx{\Omega}^{\top} \mtx{B}(t) \mtx{\Omega})^{\dagger} (\mtx{B}(t) \mtx{\Omega})^{\top} \big) =
    \frac{1}{n} \Trace \big((\underbrace{\mtx{\Omega}^{\top} \mtx{B}(t) \mtx{\Omega}}_{=: \mtx{K}_1(t)})^{\dagger} (\underbrace{\mtx{\Omega}^{\top} \mtx{B}(t)^2 \mtx{\Omega}}_{=: \mtx{K}_2(t)}) \big).
    \label{equ:cyclic-property}
\end{equation}
where we used the invariance of the trace under cyclic permutations.
Inserting the Chebyshev approximation  $\mtx{B}(t) = \underline{g}_{\sigma}^{(m)}(t \mtx{I}_n - \mtx{A})$ gives
\begin{equation*}
    \mtx{K}_1(t) = \mtx{\Omega}^{\top} \underline{g}_{\sigma}^{(m)}(t \mtx{I}_n - \mtx{A}) \mtx{\Omega} = \sum_{l=0}^{m} \underline{\mu}_l(t) (\mtx{\Omega}^{\top} T_l(\mtx{A}) \mtx{\Omega}).
\end{equation*}
Observing that $\mtx{\Omega}^{\top} T_l(\mtx{A}) \mtx{\Omega}$ is independent of $t$ allows us to precompute these (small) $n_{\mtx{\Omega}} \times n_{\mtx{\Omega}}$ matrices for $l=0, \dots, m$ using the Chebyshev recurrence \cref{equ:chebyshev-recurrence}. For each value of $t$, we then only need to sum these matrices (multiplied by the corresponding coefficient $\underline{\mu}_l(t)$).

The squaring of Chebyshev approximations discussed in \cref{sec:preservnonneg} allows us to express $\mtx{B}(t)^2 = \underline{g}_{\sigma}^{(m)}(t \mtx{I}_n - \mtx{A})^2 = \sum_{l=0}^{2m} \underline{\nu}_l(t) T_l(\mtx{A})$, where the coefficients $\underline{\nu}_l(t)$ can be cheaply computed from $\underline{\mu}_l(t)$ using~\cref{alg:chebyshev-squaring}. Hence, if we precompute
$\mtx{\Omega}^{\top} T_l(\mtx{A}) \mtx{\Omega}$
also for $l=m+1, \dots, 2m$, this allows us to evaluate $\mtx{K}_2(t)$ cheaply, in an analogous way to $\mtx{K}_1(t)$. Alternatively, \cite{lin-2017-randomized-estimation} proposes to instead approximate the squared smoothing kernel $g_{\sigma}^2$ in a separate ``auxiliary" Chebyshev series truncated after term $2m$. While this direct approximation of $g_{\sigma}^2$ might be more accurate than our squaring approach, it is \emph{inconsistent} in the sense that numerically computing the truncated Chebyshev series of $g_{\sigma}^2$ gives --- in general --- not the same result as squaring the truncated Chebyshev series of $g_{\sigma}$. We have observed that meticulously ensuring consistency between the approximations of $g_{\sigma}$ and $g_{\sigma}^2$ makes a visible difference in terms of accuracy (cf. \cref{fig:interpolation-issue}). This consistency will also allow us to derive bounds for the estimator \cref{equ:chebyshev-nystrom-formula}.

\begin{figure}[ht]
    \centering
    \scalebox{0.8}{\input{plots/interpolation.pgf}}
    \caption{Accuracy when computing the estimate \cref{equ:cyclic-property} in three different ways. \emph{Inconsistent}: Applying an auxiliary truncated Chebyshev series to $g_{\sigma}^2$, as proposed in~\cite{lin-2017-randomized-estimation}. \emph{Consistent}: Squaring the Chebyshev approximation $g_{\sigma}^{(m)}$ using \cref{alg:chebyshev-squaring}. \emph{Non-negative consistent}: Squaring the non-negative Chebyshev approximation $\underline{g}_{\sigma}^{(m)}$ with~\cref{alg:chebyshev-squaring}. As input matrix we use the $1000 \times 1000$ matrix $\mtx{A}$ coming from a three-dimensional finite difference discretization of the Hamiltonian defined in~\cref{equ:electronic-hamiltonian} on a $10 \times 10$ grid. We fix $n_{\mtx{\Omega}} = 120$ to make sure the low-rank approximation is accurate, set $\sigma = 0.005$, and compute the $L^1$ error from the approximation of the spectral density at $n_t = 100$ uniformly spaced values of $t$ in the interval $[-1, 1]$.}
    \label{fig:interpolation-issue}
\end{figure}

To treat the second term in the definition~\cref{equ:nystrompp-trace-estimator} of $\Nystrpp{\mtx{\Omega}}{\mtx{\Psi}}$, we observe that 
\begin{align}
    &\frac{1}{n n_{\mtx{\Psi}}} \Trace \big(\mtx{\Psi}^{\top} (\mtx{B}(t) - (\mtx{B}(t) \mtx{\Omega}) (\mtx{\Omega}^{\top} \mtx{B}(t) \mtx{\Omega})^{\dagger} (\mtx{B}(t) \mtx{\Omega})^{\top}) \mtx{\Psi}^{\top} \big) \notag \\
    &= \frac{1}{n n_{\mtx{\Psi}}} \big( \underbrace{\Trace(\mtx{\Psi}^{\top} \mtx{B}(t)\mtx{\Psi})}_{=: \ell(t)} - \Trace((\underbrace{\mtx{\Psi}^{\top} \mtx{B}(t) \mtx{\Omega}}_{=: \mtx{L}_1(t)^{\top}}) (\underbrace{\mtx{\Omega}^{\top} \mtx{B}(t) \mtx{\Omega}}_{= \mtx{K}_1(t)^{\top}})^{\dagger} (\underbrace{\mtx{\Omega}^{\top} \mtx{B}(t) \mtx{\Psi}^{\top}}_{=: \mtx{L}_1(t)})) \big).
    \label{equ:cyclic-property-hutchinson}
\end{align}
Each of the annotated terms can be evaluated efficiently for many values of $t$, by precomputing quantities using the Chebyshev recurrence 
$\mtx{B}(t) = \underline{g}_{\sigma}^{(m)}(t \mtx{I}_n - \mtx{A})$, as explained above.

\cref{alg:nystrom-chebyshev-pp} summarizes the considerations of this section.

\begin{algorithm}
    \caption{Chebyshev-Nyström++ method}
    \label{alg:nystrom-chebyshev-pp}
    \hspace*{\algorithmicindent} \textbf{Input:}  Symmetric $\mtx{A} \in \mathbb{R}^{n \times n}$ with spectrum in $[-1, 1]$, points $\{t_i\}_{i=1}^{n_t} \subset \mathbb{R}$\\
    \hspace*{\algorithmicindent} \textbf{Parameters:} Even approximation degree $m \in \mathbb{N}$, number of Girard-Hutchinson queries $n_{\mtx{\Psi}} \in \mathbb{N}_0$, Nyström sketch size $n_{\mtx{\Omega}} \in \mathbb{N}_0$,  smoothing parameter $\sigma > 0$\\
    \hspace*{\algorithmicindent} \textbf{Output:} Spectral density approximation $\{\widetilde{\underline{\phi}}_{\sigma}^{(m)}(t_i)\}_{i=1}^{n_t}$
    \begin{algorithmic}[1]
    \FOR{$i = 1, \dots, n_t$}
        \STATE{Compute coefficients $\{\xi_l(t_i)\}_{l=0}^{\sfrac{m}{2}}$ of Chebyshev approx.~\cref{equ:square-root-chebyshev-approximation} to $\sqrt{g_{\sigma}}$ with DCT}
        \STATE{Compute coefficients $\{\underline{\mu}_l(t_i)\}_{l=0}^{m}$ of non-negative Chebyshev approx. to $g_{\sigma}$ by applying \cref{alg:chebyshev-squaring} to $\{\xi_l(t_i)\}_{l=0}^{\sfrac{m}{2}}$}
        \STATE{Compute coefficients $\{\underline{\nu}_l(t_i)\}_{l=0}^{2m}$ of squared non-negative Chebyshev approx. by applying \cref{alg:chebyshev-squaring} to $\{\underline{\mu}_l(t_i)\}_{l=0}^{m}$}
    \ENDFOR
    \STATE{Generate Gaussian random matrices $\mtx{\Omega} \in \mathbb{R}^{n \times n_{\mtx{\Omega}}}$, $\mtx{\Psi} \in \mathbb{R}^{n \times n_{\mtx{\Psi}}}$}
    \STATE{Initialize $[\mtx{V}_1, \mtx{V}_2, \mtx{V}_3] \gets [\mtx{0}_{n \times n_{\mtx{\Omega}}}, \mtx{\mtx{\Omega}}, \mtx{0}_{n \times n_{\mtx{\Omega}}}]$}
    \STATE{Initialize $[\mtx{W}_1, \mtx{W}_2, \mtx{W}_3] \gets [\mtx{0}_{n \times n_{\mtx{\Psi}}}, \mtx{\Psi}, \mtx{0}_{n \times n_{\mtx{\Psi}}}]$}
    \STATE{Initialize $[\mtx{K}_1(t_i), \mtx{K}_2(t_i)] \gets [\mtx{0}_{n_{\mtx{\Omega}} \times n_{\mtx{\Omega}}}, \mtx{0}_{n_{\mtx{\Omega}} \times n_{\mtx{\Omega}}}]$ for all $t_i$}
    \STATE{Initialize $[\mtx{L}_1(t_i), \ell(t_i)] \gets [\mtx{0}_{n_{\mtx{\Omega}} \times n_{\mtx{\Psi}}}, 0]$ for all $t_i$}
    \FOR{$l = 0, \dots, 2m$}
        \STATE{$[\mtx{X}, \mtx{Y}] \gets \mtx{\mtx{\Omega}}^{\top} [\mtx{V}_2, \mtx{W}_2]$  }
        \STATE{ $z \gets \Trace(\mtx{\Psi}^{\top} \mtx{W}_2)$}
        \FOR{$i = 1, \dots, n_t$}
            \IF {$l \leq m$}
                \STATE{$\mtx{K}_1(t_i) \gets \mtx{K}_1(t_i) + \underline{\mu}_l(t_i) \mtx{X}$ \hfill\COMMENT{assembles $\mtx{\Omega}^{\top} \underline{g}_{\sigma}^{(m)}(t\mtx{I}_n - \mtx{A}) \mtx{\Omega}$}}
                \STATE{$\mtx{L}_1(t_i) \gets \mtx{L}_1(t_i) + \underline{\mu}_l(t_i) \mtx{Y}$ \hfill\COMMENT{assembles $\mtx{\Omega}^{\top} \underline{g}_{\sigma}^{(m)}(t\mtx{I}_n - \mtx{A}) \mtx{\Psi}$}}
                \STATE{$\ell(t_i) \gets \ell(t_i) + \underline{\mu}_l(t_i) z$ \hfill\COMMENT{assembles $\Trace(\mtx{\Psi}^{\top} \underline{g}_{\sigma}^{(m)}(t\mtx{I}_n - \mtx{A}) \mtx{\Psi})$}}
            \ENDIF
            \STATE{$\mtx{K}_2(t_i) \gets \mtx{K}_2(t_i) + \underline{\nu}_l(t_i) \mtx{X}$ \hfill\COMMENT{assembles $\mtx{\Omega}^{\top} \underline{g}_{\sigma}^{(m)}(t\mtx{I}_n - \mtx{A})^2 \mtx{\Omega}$}}
        \ENDFOR
        \STATE{$[\mtx{V}_3, \mtx{W}_3] \gets (2 - \delta_{l0}) \mtx{A} [\mtx{V}_2, \mtx{W}_2] - [\mtx{V}_1, \mtx{W}_1]$ \hfill\COMMENT{$\delta_{l0}=1$ if $l=0$ else $\delta_{l0}=0$}}
        \STATE{$[\mtx{V}_1, \mtx{W}_1] \gets [\mtx{V}_2, \mtx{W}_2]$}
        \STATE{$[\mtx{V}_2, \mtx{W}_2] \gets [\mtx{V}_3, \mtx{W}_3]$}
    \ENDFOR
    \FOR{$i = 1, \dots, n_t$}
        \STATE{$\widetilde{\underline{\phi}}_{\sigma}(t_i) \gets \frac{1}{n} \Trace\left( \mtx{K}_1(t_i)^{\dagger}\mtx{K}_2(t_i) \right) + \frac{1}{n n_{\mtx{\Psi}}} \left( \ell(t_i) + \Trace\left( \mtx{L}_1(t_i)^{\top} \mtx{K}_1(t_i)^{\dagger} \mtx{L}_1(t_i) \right)  \right) $} \label{lin:4-nystromchebyshev-nystrom-pp}
    \ENDFOR
    \end{algorithmic}
\end{algorithm}

\cref{alg:nystrom-chebyshev-pp} requires $(2m + 1)(n_{\mtx{\Psi}} + n_{\mtx{\Omega}})$ matrix-vector multiplications with $\mtx{A}$, which can be expected to dominate its cost. For each of the evaluation points $t_1, \dots, t_{n_t}$ another $\mathcal{O}(m \log(m) + m (n_{\mtx{\Omega}}^2 + n_{\mtx{\Omega}} n_{\mtx{\Psi}} + 1) + (n_{\mtx{\Omega}}^3 + n_{\mtx{\Omega}} n_{\mtx{\Psi}}^2 + 1))$ operations are needed, which is independent of $n$ --- thanks to the precomputations discussed above.

\cref{alg:nystrom-chebyshev-pp} also accommodates the other two estimators discussed in \cref{sec:analysis}: For $n_{\mtx{\Omega}} = 0$ we recover the Girard-Hutchinson estimator~\cref{equ:hutchinson-trace-estimator}. For $n_{\mtx{\Psi}} = 0$ we recover the Nyström estimator~\cref{equ:nystrom-trace-estimator}, provided that $0/0$ is evaluated to $0$ in the second term on line \ref{lin:4-nystromchebyshev-nystrom-pp}. In this case, \cref{alg:nystrom-chebyshev-pp} essentially amounts to the ``spectrum sweeping" method from \cite[Algorithm 5]{lin-2017-randomized-estimation}, featuring a significantly lower computational complexity.

The two pseudoinverses on line \ref{lin:4-nystromchebyshev-nystrom-pp} need to be computed with care. One way to ensure numerical stability is to 
consider an associated eigenvalue problem and threshold the smallest eigenvalues~\cite{lin-2017-randomized-estimation, epperly-2022-theory-quantum}. However, this makes the algorithm significantly more complicated, and we found that it improves accuracy only slightly. Instead, we use the standard NumPy least-squares solver when applying $\mtx{K}_1(t_i)^{\dagger}$, which is based on the LAPACK DGELSD routine, and likewise truncates the smallest eigenvalues by treating them as zero if they are smaller than a certain threshold. We have found a threshold of $10^{-5}$ times the largest eigenvalue to give the best results in our numerical experiments. 

To ensure numerical accuracy, we found it crucial to separately treat regions
where the spectral density nearly vanishes, because there are no eigenvalues in or close to the region.
Especially for small smoothing parameters $\sigma$, the rapid decay of the Gaussian $g_{\sigma}$ will cause $\mtx{K}_1(t) \approx \mtx{\Omega}^{\top} g_{\sigma}(t\mtx{I}_n - \mtx{A}) \mtx{\Omega}$ to be close to zero for parameter values $t$ located in such regions. Hence, on line \ref{lin:4-nystromchebyshev-nystrom-pp} we would compute the pseudoinverse of an almost zero matrix, which can be extremely unstable; see \cref{fig:zerocheck}. To detect this situation, we can cheaply estimate the nuclear norm with the Girard-Hutchinson estimate $\Trace(\mtx{K}_1(t)) / n_{\mtx{\Omega}}$. If this estimate is small (in our experiments, smaller than $10^{-5}$), we use it as an indicator that $t$ is within a region with a nearly vanishing spectral density, and therefore set $\widetilde{\underline{\phi}}_{\sigma}^{(m)}(t) = 0$ instead of computing the full expression on line \ref{lin:4-nystromchebyshev-nystrom-pp}.

\begin{figure}[ht]
    \centering
    \scalebox{0.8}{\input{plots/zerocheck.pgf}}
    \caption{The difference the non-zero check can make when approximating a spectral density using \cref{alg:nystrom-chebyshev-pp}. As input matrix we use the $1000 \times 1000$ matrix $\mtx{A}$ coming from a three-dimensional finite difference discretization of the Hamiltonian defined in~\cref{equ:electronic-hamiltonian} on a $10 \times 10$ grid. We ran the Chebyshev-Nyström++ method with and without estimating if the matrix function is zero before taking its pseudo-inverse. We let $m=2000$, $n_{\mtx{\Omega}}=80$, $n_{\mtx{\Psi}}=0$, and $\sigma = 0.004$. We compute the $L^1$ error from the approximation of the spectral density at $n_t = 100$ uniformly spaced values of $t$ in the interval $[-1, 1]$.}
    \label{fig:zerocheck}
\end{figure}

\color{black}

\subsubsection{Analysis}
\label{subsubsec:chebyshev-nystrom-analysis}

In this section, we analyze the accuracy of the Chebyshev-Nyström++ estimate returned by~\cref{alg:nystrom-chebyshev-pp}.
\begin{theorem}[Error of Chebyshev-Nyström++ estimate for smoothed spectral density]\label{thm:chebyshev-nystrom}
    For a smoothing parameter $\sigma > 0$, let $\phi_{\sigma}$ be the smoothed spectral density~\cref{equ:smooth-spectral-density} of a symmetric matrix $\mtx{A} \in \mathbb{R}^{n \times n}$ with spectrum contained in $[-1, 1]$. Consider $\varepsilon \in (0, 3)$ and $\delta \in (0, 1]$.
    Then the Chebyshev-Nyström++ estimate $\widetilde{\underline{\phi}}_{\sigma}^{(m)}$ 
    satisfies the error bound
    \begin{equation*}
        \int_{-1}^{1} \left| \phi_{\sigma}(t) - \widetilde{\underline{\phi}}_{\sigma}^{(m)}(t) \right|~\mathrm{d}t \leq \varepsilon
        \label{equ:chebyshev-nystrom-error}
    \end{equation*}
    with probability $\geq 1 - \delta$ for $m = \mathcal{O}(\sigma^{-1} \log(\sigma^{-2} \varepsilon^{-1}))$ and $n_{\mtx{\Omega}} = n_{\mtx{\Psi}} = \mathcal{O}(\varepsilon^{-1}\log(\delta^{-1}))$.
\end{theorem}

\begin{proof}
Combining results derived above, we obtain that 
    \begin{align*}
        &\int_{-1}^{1} \left| \phi_{\sigma}(t) - \widetilde{\underline{\phi}}_{\sigma}^{(m)}(t) \right|~\mathrm{d}t \notag \\
        &\leq \int_{-1}^{1} \left| \phi_{\sigma}(t) - \underline{\phi}_{\sigma}^{(m)}(t) \right|~\mathrm{d}t + \int_{-1}^{1} \left| \underline{\phi}_{\sigma}^{(m)}(t) - \widetilde{\underline{\phi}}_{\sigma}^{(m)}(t) \right|~\mathrm{d}t &&  \notag \\
        &\leq \int_{-1}^{1} \left| \phi_{\sigma}(t) - \underline{\phi}_{\sigma}^{(m)}(t) \right|~\mathrm{d}t + \widetilde{\varepsilon} \int_{-1}^{1} \left| \underline{\phi}_{\sigma}^{(m)}(t) \right|~\mathrm{d}t && \text{(\cref{thm:nystrom-pp}, $\mathbb{P} \geq 1 - \delta$)} \notag \\
        &\leq (1 + \widetilde{\varepsilon}) \int_{-1}^{1} \left| \phi_{\sigma}(t) - \underline{\phi}_{\sigma}^{(m)}(t) \right|~\mathrm{d}t + \widetilde{\varepsilon} \int_{-1}^{1} \left| \phi_{\sigma}(t) \right|~\mathrm{d}t &&  \notag \\
        &\leq (1 + \widetilde{\varepsilon}) \underline{E}_{\sigma, m} + \widetilde{\varepsilon}. && \text{(\cref{lem:non-negative-chebyshev-error})} \notag 
    \end{align*}
    From \cref{lem:non-negative-chebyshev-error} and the discussion after \cref{lem:chebyshev-error} we see that the choice $m = \mathcal{O}(\sigma^{-1} \log(\sigma^{-2} \varepsilon^{-1}))$ ensures $\underline{E}_{\sigma, m} \leq \widetilde{\varepsilon}$, and thus
    the result follows by setting $\varepsilon = 3 \widetilde{\varepsilon}$.
\end{proof}

This result is in close analogy to the corresponding result for cumulative spectral distribution approximation with the stochastic Lanczos quadrature \cite[Theorem 1]{chen-2021-analysis-stochastic}, except it is stated in terms of the $L^1$-norm instead of the Wasserstein-1 metric. Moreover, the overall structure of the two proofs is almost identical.

\cref{thm:chebyshev-nystrom} is extended to matrices with the spectrum contained in an abritrary real interval $[a, b]$ by employing the usual affine transformation $\tau : [a, b] \to [-1, 1]$.

In practice, the $\mathcal{O}(\varepsilon^{-1})$ complexity of $n_{\mtx{\Omega}}$ and $n_{\mtx{\Psi}}$ in \cref{thm:chebyshev-nystrom} is often pessimistic, because  the underlying result of \cref{thm:nystrom-pp} is agnostic to the low-rank structure of the matrix $\mtx{B}(t) = g_{\sigma}(t\mtx{I}_n - \mtx{A})$. The singular values of this matrix are given by
\begin{equation*}
    \sigma_i(t) = g_{\sigma}(t - \lambda_{(i)}) = \frac{1}{\sigma\sqrt{2 \pi}} \exp\Big( -\frac{(t - \lambda_{(i)})^2}{2 \sigma^2} \Big),  
    \label{equ:gaussian-kernel-eigenvalues}
\end{equation*}
where $\lambda_{(1)}, \dots, \lambda_{(n)}$ denote the eigenvalues of $\mtx{A}$ sorted by increasing distance from $t$, such that $\sigma_1(t) \geq \dots \geq \sigma_n(t) \ge 0$. An eigenvalue $\lambda_{(i)}$ distant from $t$ results in an exponentially small singular value $\sigma_i(t)$. In particular, for a small value of $\sigma$, all singular values corresponding to eigenvalues outside a certain window $[t - d, t + d]$ are negligibly small; see \cref{fig:numerical-rank}.
\begin{figure}[ht]
    \centering
    \input{figures/numerical-rank.tex}
    \caption{Singular values $\sigma_i$ of $g_{\sigma}(t\mtx{I}_n - \mtx{A})$ corresponding to eigenvalues $\lambda_{(i)}$ of $\mtx{A}$ which lie further away from $t$ than a certain distance $d$ are negligibly small.}
    \label{fig:numerical-rank}
\end{figure}
By choosing $n_{\mtx{\Omega}}$ such that the sum $\sum_{i=n_{\mtx{\Omega}}+5}^{n} \sigma_i(t)$ becomes negligibly small,  \cref{thm:nystrom} implies an excellent approximation with high probability. In particular, if we assume the eigenvalues of $\mtx{A}$ to be roughly uniformly distributed, a small calculation shows that for a choice $\sigma \approx \varepsilon$ (as is motivated in \cref{subsec:spectral-density}), $n_{\mtx{\Omega}} = \mathcal{O}(n \varepsilon \sqrt{\log(\varepsilon^{-1})}) \approx \mathcal{O}(n \varepsilon)$ will ensure this.

It is worth noting that for $n_{\mtx{\Omega}} = 0$ (that is, no low-rank approximation), the Chebyshev-Nyström++ estimator is closely related to the stochastic Lanczos quadrature (SLQ) \cite{chen-2021-analysis-stochastic}. In fact, it fits into the framework introduced in \cite{chen-2025-randomized-matrixfree}, and only distinguishes itself from SLQ by using Clenshaw-Curtis quadrature instead of Gaussian quadrature. As a consequence, it is expected to converge slightly slower than SLQ when approximating the smoothed spectral density \cite{trefethen-2008-gauss-quadrature}. However, its connection to the DCT makes it a more favorable choice when the approximation is evaluated in multiple values of the parameter $t$, which may easily offset the slightly slower convergence.

%% file: tables/interpolation.tex
\centering
\renewcommand{\arraystretch}{1.2}
\begin{tabular}{@{}lcccc@{}}
\toprule
 & $m=800$ & $m=1600$ & $m=2400$ & $m=3200$\\
\midrule
FFT & $101.4$ $\pm$ $0.2$ & $273.9$ $\pm$ $0.6$ & $214.6$ $\pm$ $0.9$ & $301.5$ $\pm$ $1.1$ \\
DCT & $57.4$ $\pm$ $0.1$ & $86.3$ $\pm$ $0.4$ & $115.4$ $\pm$ $0.2$ & $143.5$ $\pm$ $0.3$ \\
non-negative DCT & $123.7$ $\pm$ $0.3$ & $156.6$ $\pm$ $0.2$ & $188.4$ $\pm$ $0.4$ & $218.0$ $\pm$ $0.8$ \\
\bottomrule
\end{tabular}

%% file: plots/interpolation.pgf
\begingroup%
\makeatletter%
\begin{pgfpicture}%
\pgfpathrectangle{\pgfpointorigin}{\pgfqpoint{5.340163in}{2.959073in}}%
\pgfusepath{use as bounding box, clip}%
\begin{pgfscope}%
\pgfsetbuttcap%
\pgfsetmiterjoin%
\definecolor{currentfill}{rgb}{1.000000,1.000000,1.000000}%
\pgfsetfillcolor{currentfill}%
\pgfsetlinewidth{0.000000pt}%
\definecolor{currentstroke}{rgb}{1.000000,1.000000,1.000000}%
\pgfsetstrokecolor{currentstroke}%
\pgfsetdash{}{0pt}%
\pgfpathmoveto{\pgfqpoint{0.000000in}{-0.000000in}}%
\pgfpathlineto{\pgfqpoint{5.340163in}{-0.000000in}}%
\pgfpathlineto{\pgfqpoint{5.340163in}{2.959073in}}%
\pgfpathlineto{\pgfqpoint{0.000000in}{2.959073in}}%
\pgfpathlineto{\pgfqpoint{0.000000in}{-0.000000in}}%
\pgfpathclose%
\pgfusepath{fill}%
\end{pgfscope}%
\begin{pgfscope}%
\pgfsetbuttcap%
\pgfsetmiterjoin%
\definecolor{currentfill}{rgb}{1.000000,1.000000,1.000000}%
\pgfsetfillcolor{currentfill}%
\pgfsetlinewidth{0.000000pt}%
\definecolor{currentstroke}{rgb}{0.000000,0.000000,0.000000}%
\pgfsetstrokecolor{currentstroke}%
\pgfsetstrokeopacity{0.000000}%
\pgfsetdash{}{0pt}%
\pgfpathmoveto{\pgfqpoint{0.721913in}{0.549073in}}%
\pgfpathlineto{\pgfqpoint{5.240163in}{0.549073in}}%
\pgfpathlineto{\pgfqpoint{5.240163in}{2.859073in}}%
\pgfpathlineto{\pgfqpoint{0.721913in}{2.859073in}}%
\pgfpathlineto{\pgfqpoint{0.721913in}{0.549073in}}%
\pgfpathclose%
\pgfusepath{fill}%
\end{pgfscope}%
\begin{pgfscope}%
\pgfpathrectangle{\pgfqpoint{0.721913in}{0.549073in}}{\pgfqpoint{4.518250in}{2.310000in}}%
\pgfusepath{clip}%
\pgfsetrectcap%
\pgfsetroundjoin%
\pgfsetlinewidth{0.250937pt}%
\definecolor{currentstroke}{rgb}{0.000000,0.000000,0.000000}%
\pgfsetstrokecolor{currentstroke}%
\pgfsetstrokeopacity{0.200000}%
\pgfsetdash{}{0pt}%
\pgfpathmoveto{\pgfqpoint{1.422166in}{0.549073in}}%
\pgfpathlineto{\pgfqpoint{1.422166in}{2.859073in}}%
\pgfusepath{stroke}%
\end{pgfscope}%
\begin{pgfscope}%
\pgfsetbuttcap%
\pgfsetroundjoin%
\definecolor{currentfill}{rgb}{0.000000,0.000000,0.000000}%
\pgfsetfillcolor{currentfill}%
\pgfsetlinewidth{0.803000pt}%
\definecolor{currentstroke}{rgb}{0.000000,0.000000,0.000000}%
\pgfsetstrokecolor{currentstroke}%
\pgfsetdash{}{0pt}%
\pgfsys@defobject{currentmarker}{\pgfqpoint{0.000000in}{-0.048611in}}{\pgfqpoint{0.000000in}{0.000000in}}{%
\pgfpathmoveto{\pgfqpoint{0.000000in}{0.000000in}}%
\pgfpathlineto{\pgfqpoint{0.000000in}{-0.048611in}}%
\pgfusepath{stroke,fill}%
}%
\begin{pgfscope}%
\pgfsys@transformshift{1.422166in}{0.549073in}%
\pgfsys@useobject{currentmarker}{}%
\end{pgfscope}%
\end{pgfscope}%
\begin{pgfscope}%
\definecolor{textcolor}{rgb}{0.000000,0.000000,0.000000}%
\pgfsetstrokecolor{textcolor}%
\pgfsetfillcolor{textcolor}%
\pgftext[x=1.422166in,y=0.451851in,,top]{\color{textcolor}{\rmfamily\fontsize{12.000000}{14.400000}\selectfont\catcode`\^=\active\def^{\ifmmode\sp\else\^{}\fi}\catcode`\%=\active\def
\end{pgfscope}%
\begin{pgfscope}%
\pgfpathrectangle{\pgfqpoint{0.721913in}{0.549073in}}{\pgfqpoint{4.518250in}{2.310000in}}%
\pgfusepath{clip}%
\pgfsetrectcap%
\pgfsetroundjoin%
\pgfsetlinewidth{0.250937pt}%
\definecolor{currentstroke}{rgb}{0.000000,0.000000,0.000000}%
\pgfsetstrokecolor{currentstroke}%
\pgfsetstrokeopacity{0.200000}%
\pgfsetdash{}{0pt}%
\pgfpathmoveto{\pgfqpoint{4.276939in}{0.549073in}}%
\pgfpathlineto{\pgfqpoint{4.276939in}{2.859073in}}%
\pgfusepath{stroke}%
\end{pgfscope}%
\begin{pgfscope}%
\pgfsetbuttcap%
\pgfsetroundjoin%
\definecolor{currentfill}{rgb}{0.000000,0.000000,0.000000}%
\pgfsetfillcolor{currentfill}%
\pgfsetlinewidth{0.803000pt}%
\definecolor{currentstroke}{rgb}{0.000000,0.000000,0.000000}%
\pgfsetstrokecolor{currentstroke}%
\pgfsetdash{}{0pt}%
\pgfsys@defobject{currentmarker}{\pgfqpoint{0.000000in}{-0.048611in}}{\pgfqpoint{0.000000in}{0.000000in}}{%
\pgfpathmoveto{\pgfqpoint{0.000000in}{0.000000in}}%
\pgfpathlineto{\pgfqpoint{0.000000in}{-0.048611in}}%
\pgfusepath{stroke,fill}%
}%
\begin{pgfscope}%
\pgfsys@transformshift{4.276939in}{0.549073in}%
\pgfsys@useobject{currentmarker}{}%
\end{pgfscope}%
\end{pgfscope}%
\begin{pgfscope}%
\definecolor{textcolor}{rgb}{0.000000,0.000000,0.000000}%
\pgfsetstrokecolor{textcolor}%
\pgfsetfillcolor{textcolor}%
\pgftext[x=4.276939in,y=0.451851in,,top]{\color{textcolor}{\rmfamily\fontsize{12.000000}{14.400000}\selectfont\catcode`\^=\active\def^{\ifmmode\sp\else\^{}\fi}\catcode`\%=\active\def
\end{pgfscope}%
\begin{pgfscope}%
\pgfpathrectangle{\pgfqpoint{0.721913in}{0.549073in}}{\pgfqpoint{4.518250in}{2.310000in}}%
\pgfusepath{clip}%
\pgfsetrectcap%
\pgfsetroundjoin%
\pgfsetlinewidth{0.250937pt}%
\definecolor{currentstroke}{rgb}{0.000000,0.000000,0.000000}%
\pgfsetstrokecolor{currentstroke}%
\pgfsetstrokeopacity{0.200000}%
\pgfsetdash{}{0pt}%
\pgfpathmoveto{\pgfqpoint{0.788838in}{0.549073in}}%
\pgfpathlineto{\pgfqpoint{0.788838in}{2.859073in}}%
\pgfusepath{stroke}%
\end{pgfscope}%
\begin{pgfscope}%
\pgfsetbuttcap%
\pgfsetroundjoin%
\definecolor{currentfill}{rgb}{0.000000,0.000000,0.000000}%
\pgfsetfillcolor{currentfill}%
\pgfsetlinewidth{0.602250pt}%
\definecolor{currentstroke}{rgb}{0.000000,0.000000,0.000000}%
\pgfsetstrokecolor{currentstroke}%
\pgfsetdash{}{0pt}%
\pgfsys@defobject{currentmarker}{\pgfqpoint{0.000000in}{-0.027778in}}{\pgfqpoint{0.000000in}{0.000000in}}{%
\pgfpathmoveto{\pgfqpoint{0.000000in}{0.000000in}}%
\pgfpathlineto{\pgfqpoint{0.000000in}{-0.027778in}}%
\pgfusepath{stroke,fill}%
}%
\begin{pgfscope}%
\pgfsys@transformshift{0.788838in}{0.549073in}%
\pgfsys@useobject{currentmarker}{}%
\end{pgfscope}%
\end{pgfscope}%
\begin{pgfscope}%
\pgfpathrectangle{\pgfqpoint{0.721913in}{0.549073in}}{\pgfqpoint{4.518250in}{2.310000in}}%
\pgfusepath{clip}%
\pgfsetrectcap%
\pgfsetroundjoin%
\pgfsetlinewidth{0.250937pt}%
\definecolor{currentstroke}{rgb}{0.000000,0.000000,0.000000}%
\pgfsetstrokecolor{currentstroke}%
\pgfsetstrokeopacity{0.200000}%
\pgfsetdash{}{0pt}%
\pgfpathmoveto{\pgfqpoint{0.979956in}{0.549073in}}%
\pgfpathlineto{\pgfqpoint{0.979956in}{2.859073in}}%
\pgfusepath{stroke}%
\end{pgfscope}%
\begin{pgfscope}%
\pgfsetbuttcap%
\pgfsetroundjoin%
\definecolor{currentfill}{rgb}{0.000000,0.000000,0.000000}%
\pgfsetfillcolor{currentfill}%
\pgfsetlinewidth{0.602250pt}%
\definecolor{currentstroke}{rgb}{0.000000,0.000000,0.000000}%
\pgfsetstrokecolor{currentstroke}%
\pgfsetdash{}{0pt}%
\pgfsys@defobject{currentmarker}{\pgfqpoint{0.000000in}{-0.027778in}}{\pgfqpoint{0.000000in}{0.000000in}}{%
\pgfpathmoveto{\pgfqpoint{0.000000in}{0.000000in}}%
\pgfpathlineto{\pgfqpoint{0.000000in}{-0.027778in}}%
\pgfusepath{stroke,fill}%
}%
\begin{pgfscope}%
\pgfsys@transformshift{0.979956in}{0.549073in}%
\pgfsys@useobject{currentmarker}{}%
\end{pgfscope}%
\end{pgfscope}%
\begin{pgfscope}%
\pgfpathrectangle{\pgfqpoint{0.721913in}{0.549073in}}{\pgfqpoint{4.518250in}{2.310000in}}%
\pgfusepath{clip}%
\pgfsetrectcap%
\pgfsetroundjoin%
\pgfsetlinewidth{0.250937pt}%
\definecolor{currentstroke}{rgb}{0.000000,0.000000,0.000000}%
\pgfsetstrokecolor{currentstroke}%
\pgfsetstrokeopacity{0.200000}%
\pgfsetdash{}{0pt}%
\pgfpathmoveto{\pgfqpoint{1.145509in}{0.549073in}}%
\pgfpathlineto{\pgfqpoint{1.145509in}{2.859073in}}%
\pgfusepath{stroke}%
\end{pgfscope}%
\begin{pgfscope}%
\pgfsetbuttcap%
\pgfsetroundjoin%
\definecolor{currentfill}{rgb}{0.000000,0.000000,0.000000}%
\pgfsetfillcolor{currentfill}%
\pgfsetlinewidth{0.602250pt}%
\definecolor{currentstroke}{rgb}{0.000000,0.000000,0.000000}%
\pgfsetstrokecolor{currentstroke}%
\pgfsetdash{}{0pt}%
\pgfsys@defobject{currentmarker}{\pgfqpoint{0.000000in}{-0.027778in}}{\pgfqpoint{0.000000in}{0.000000in}}{%
\pgfpathmoveto{\pgfqpoint{0.000000in}{0.000000in}}%
\pgfpathlineto{\pgfqpoint{0.000000in}{-0.027778in}}%
\pgfusepath{stroke,fill}%
}%
\begin{pgfscope}%
\pgfsys@transformshift{1.145509in}{0.549073in}%
\pgfsys@useobject{currentmarker}{}%
\end{pgfscope}%
\end{pgfscope}%
\begin{pgfscope}%
\pgfpathrectangle{\pgfqpoint{0.721913in}{0.549073in}}{\pgfqpoint{4.518250in}{2.310000in}}%
\pgfusepath{clip}%
\pgfsetrectcap%
\pgfsetroundjoin%
\pgfsetlinewidth{0.250937pt}%
\definecolor{currentstroke}{rgb}{0.000000,0.000000,0.000000}%
\pgfsetstrokecolor{currentstroke}%
\pgfsetstrokeopacity{0.200000}%
\pgfsetdash{}{0pt}%
\pgfpathmoveto{\pgfqpoint{1.291538in}{0.549073in}}%
\pgfpathlineto{\pgfqpoint{1.291538in}{2.859073in}}%
\pgfusepath{stroke}%
\end{pgfscope}%
\begin{pgfscope}%
\pgfsetbuttcap%
\pgfsetroundjoin%
\definecolor{currentfill}{rgb}{0.000000,0.000000,0.000000}%
\pgfsetfillcolor{currentfill}%
\pgfsetlinewidth{0.602250pt}%
\definecolor{currentstroke}{rgb}{0.000000,0.000000,0.000000}%
\pgfsetstrokecolor{currentstroke}%
\pgfsetdash{}{0pt}%
\pgfsys@defobject{currentmarker}{\pgfqpoint{0.000000in}{-0.027778in}}{\pgfqpoint{0.000000in}{0.000000in}}{%
\pgfpathmoveto{\pgfqpoint{0.000000in}{0.000000in}}%
\pgfpathlineto{\pgfqpoint{0.000000in}{-0.027778in}}%
\pgfusepath{stroke,fill}%
}%
\begin{pgfscope}%
\pgfsys@transformshift{1.291538in}{0.549073in}%
\pgfsys@useobject{currentmarker}{}%
\end{pgfscope}%
\end{pgfscope}%
\begin{pgfscope}%
\pgfpathrectangle{\pgfqpoint{0.721913in}{0.549073in}}{\pgfqpoint{4.518250in}{2.310000in}}%
\pgfusepath{clip}%
\pgfsetrectcap%
\pgfsetroundjoin%
\pgfsetlinewidth{0.250937pt}%
\definecolor{currentstroke}{rgb}{0.000000,0.000000,0.000000}%
\pgfsetstrokecolor{currentstroke}%
\pgfsetstrokeopacity{0.200000}%
\pgfsetdash{}{0pt}%
\pgfpathmoveto{\pgfqpoint{2.281538in}{0.549073in}}%
\pgfpathlineto{\pgfqpoint{2.281538in}{2.859073in}}%
\pgfusepath{stroke}%
\end{pgfscope}%
\begin{pgfscope}%
\pgfsetbuttcap%
\pgfsetroundjoin%
\definecolor{currentfill}{rgb}{0.000000,0.000000,0.000000}%
\pgfsetfillcolor{currentfill}%
\pgfsetlinewidth{0.602250pt}%
\definecolor{currentstroke}{rgb}{0.000000,0.000000,0.000000}%
\pgfsetstrokecolor{currentstroke}%
\pgfsetdash{}{0pt}%
\pgfsys@defobject{currentmarker}{\pgfqpoint{0.000000in}{-0.027778in}}{\pgfqpoint{0.000000in}{0.000000in}}{%
\pgfpathmoveto{\pgfqpoint{0.000000in}{0.000000in}}%
\pgfpathlineto{\pgfqpoint{0.000000in}{-0.027778in}}%
\pgfusepath{stroke,fill}%
}%
\begin{pgfscope}%
\pgfsys@transformshift{2.281538in}{0.549073in}%
\pgfsys@useobject{currentmarker}{}%
\end{pgfscope}%
\end{pgfscope}%
\begin{pgfscope}%
\pgfpathrectangle{\pgfqpoint{0.721913in}{0.549073in}}{\pgfqpoint{4.518250in}{2.310000in}}%
\pgfusepath{clip}%
\pgfsetrectcap%
\pgfsetroundjoin%
\pgfsetlinewidth{0.250937pt}%
\definecolor{currentstroke}{rgb}{0.000000,0.000000,0.000000}%
\pgfsetstrokecolor{currentstroke}%
\pgfsetstrokeopacity{0.200000}%
\pgfsetdash{}{0pt}%
\pgfpathmoveto{\pgfqpoint{2.784238in}{0.549073in}}%
\pgfpathlineto{\pgfqpoint{2.784238in}{2.859073in}}%
\pgfusepath{stroke}%
\end{pgfscope}%
\begin{pgfscope}%
\pgfsetbuttcap%
\pgfsetroundjoin%
\definecolor{currentfill}{rgb}{0.000000,0.000000,0.000000}%
\pgfsetfillcolor{currentfill}%
\pgfsetlinewidth{0.602250pt}%
\definecolor{currentstroke}{rgb}{0.000000,0.000000,0.000000}%
\pgfsetstrokecolor{currentstroke}%
\pgfsetdash{}{0pt}%
\pgfsys@defobject{currentmarker}{\pgfqpoint{0.000000in}{-0.027778in}}{\pgfqpoint{0.000000in}{0.000000in}}{%
\pgfpathmoveto{\pgfqpoint{0.000000in}{0.000000in}}%
\pgfpathlineto{\pgfqpoint{0.000000in}{-0.027778in}}%
\pgfusepath{stroke,fill}%
}%
\begin{pgfscope}%
\pgfsys@transformshift{2.784238in}{0.549073in}%
\pgfsys@useobject{currentmarker}{}%
\end{pgfscope}%
\end{pgfscope}%
\begin{pgfscope}%
\pgfpathrectangle{\pgfqpoint{0.721913in}{0.549073in}}{\pgfqpoint{4.518250in}{2.310000in}}%
\pgfusepath{clip}%
\pgfsetrectcap%
\pgfsetroundjoin%
\pgfsetlinewidth{0.250937pt}%
\definecolor{currentstroke}{rgb}{0.000000,0.000000,0.000000}%
\pgfsetstrokecolor{currentstroke}%
\pgfsetstrokeopacity{0.200000}%
\pgfsetdash{}{0pt}%
\pgfpathmoveto{\pgfqpoint{3.140910in}{0.549073in}}%
\pgfpathlineto{\pgfqpoint{3.140910in}{2.859073in}}%
\pgfusepath{stroke}%
\end{pgfscope}%
\begin{pgfscope}%
\pgfsetbuttcap%
\pgfsetroundjoin%
\definecolor{currentfill}{rgb}{0.000000,0.000000,0.000000}%
\pgfsetfillcolor{currentfill}%
\pgfsetlinewidth{0.602250pt}%
\definecolor{currentstroke}{rgb}{0.000000,0.000000,0.000000}%
\pgfsetstrokecolor{currentstroke}%
\pgfsetdash{}{0pt}%
\pgfsys@defobject{currentmarker}{\pgfqpoint{0.000000in}{-0.027778in}}{\pgfqpoint{0.000000in}{0.000000in}}{%
\pgfpathmoveto{\pgfqpoint{0.000000in}{0.000000in}}%
\pgfpathlineto{\pgfqpoint{0.000000in}{-0.027778in}}%
\pgfusepath{stroke,fill}%
}%
\begin{pgfscope}%
\pgfsys@transformshift{3.140910in}{0.549073in}%
\pgfsys@useobject{currentmarker}{}%
\end{pgfscope}%
\end{pgfscope}%
\begin{pgfscope}%
\pgfpathrectangle{\pgfqpoint{0.721913in}{0.549073in}}{\pgfqpoint{4.518250in}{2.310000in}}%
\pgfusepath{clip}%
\pgfsetrectcap%
\pgfsetroundjoin%
\pgfsetlinewidth{0.250937pt}%
\definecolor{currentstroke}{rgb}{0.000000,0.000000,0.000000}%
\pgfsetstrokecolor{currentstroke}%
\pgfsetstrokeopacity{0.200000}%
\pgfsetdash{}{0pt}%
\pgfpathmoveto{\pgfqpoint{3.417566in}{0.549073in}}%
\pgfpathlineto{\pgfqpoint{3.417566in}{2.859073in}}%
\pgfusepath{stroke}%
\end{pgfscope}%
\begin{pgfscope}%
\pgfsetbuttcap%
\pgfsetroundjoin%
\definecolor{currentfill}{rgb}{0.000000,0.000000,0.000000}%
\pgfsetfillcolor{currentfill}%
\pgfsetlinewidth{0.602250pt}%
\definecolor{currentstroke}{rgb}{0.000000,0.000000,0.000000}%
\pgfsetstrokecolor{currentstroke}%
\pgfsetdash{}{0pt}%
\pgfsys@defobject{currentmarker}{\pgfqpoint{0.000000in}{-0.027778in}}{\pgfqpoint{0.000000in}{0.000000in}}{%
\pgfpathmoveto{\pgfqpoint{0.000000in}{0.000000in}}%
\pgfpathlineto{\pgfqpoint{0.000000in}{-0.027778in}}%
\pgfusepath{stroke,fill}%
}%
\begin{pgfscope}%
\pgfsys@transformshift{3.417566in}{0.549073in}%
\pgfsys@useobject{currentmarker}{}%
\end{pgfscope}%
\end{pgfscope}%
\begin{pgfscope}%
\pgfpathrectangle{\pgfqpoint{0.721913in}{0.549073in}}{\pgfqpoint{4.518250in}{2.310000in}}%
\pgfusepath{clip}%
\pgfsetrectcap%
\pgfsetroundjoin%
\pgfsetlinewidth{0.250937pt}%
\definecolor{currentstroke}{rgb}{0.000000,0.000000,0.000000}%
\pgfsetstrokecolor{currentstroke}%
\pgfsetstrokeopacity{0.200000}%
\pgfsetdash{}{0pt}%
\pgfpathmoveto{\pgfqpoint{3.643611in}{0.549073in}}%
\pgfpathlineto{\pgfqpoint{3.643611in}{2.859073in}}%
\pgfusepath{stroke}%
\end{pgfscope}%
\begin{pgfscope}%
\pgfsetbuttcap%
\pgfsetroundjoin%
\definecolor{currentfill}{rgb}{0.000000,0.000000,0.000000}%
\pgfsetfillcolor{currentfill}%
\pgfsetlinewidth{0.602250pt}%
\definecolor{currentstroke}{rgb}{0.000000,0.000000,0.000000}%
\pgfsetstrokecolor{currentstroke}%
\pgfsetdash{}{0pt}%
\pgfsys@defobject{currentmarker}{\pgfqpoint{0.000000in}{-0.027778in}}{\pgfqpoint{0.000000in}{0.000000in}}{%
\pgfpathmoveto{\pgfqpoint{0.000000in}{0.000000in}}%
\pgfpathlineto{\pgfqpoint{0.000000in}{-0.027778in}}%
\pgfusepath{stroke,fill}%
}%
\begin{pgfscope}%
\pgfsys@transformshift{3.643611in}{0.549073in}%
\pgfsys@useobject{currentmarker}{}%
\end{pgfscope}%
\end{pgfscope}%
\begin{pgfscope}%
\pgfpathrectangle{\pgfqpoint{0.721913in}{0.549073in}}{\pgfqpoint{4.518250in}{2.310000in}}%
\pgfusepath{clip}%
\pgfsetrectcap%
\pgfsetroundjoin%
\pgfsetlinewidth{0.250937pt}%
\definecolor{currentstroke}{rgb}{0.000000,0.000000,0.000000}%
\pgfsetstrokecolor{currentstroke}%
\pgfsetstrokeopacity{0.200000}%
\pgfsetdash{}{0pt}%
\pgfpathmoveto{\pgfqpoint{3.834729in}{0.549073in}}%
\pgfpathlineto{\pgfqpoint{3.834729in}{2.859073in}}%
\pgfusepath{stroke}%
\end{pgfscope}%
\begin{pgfscope}%
\pgfsetbuttcap%
\pgfsetroundjoin%
\definecolor{currentfill}{rgb}{0.000000,0.000000,0.000000}%
\pgfsetfillcolor{currentfill}%
\pgfsetlinewidth{0.602250pt}%
\definecolor{currentstroke}{rgb}{0.000000,0.000000,0.000000}%
\pgfsetstrokecolor{currentstroke}%
\pgfsetdash{}{0pt}%
\pgfsys@defobject{currentmarker}{\pgfqpoint{0.000000in}{-0.027778in}}{\pgfqpoint{0.000000in}{0.000000in}}{%
\pgfpathmoveto{\pgfqpoint{0.000000in}{0.000000in}}%
\pgfpathlineto{\pgfqpoint{0.000000in}{-0.027778in}}%
\pgfusepath{stroke,fill}%
}%
\begin{pgfscope}%
\pgfsys@transformshift{3.834729in}{0.549073in}%
\pgfsys@useobject{currentmarker}{}%
\end{pgfscope}%
\end{pgfscope}%
\begin{pgfscope}%
\pgfpathrectangle{\pgfqpoint{0.721913in}{0.549073in}}{\pgfqpoint{4.518250in}{2.310000in}}%
\pgfusepath{clip}%
\pgfsetrectcap%
\pgfsetroundjoin%
\pgfsetlinewidth{0.250937pt}%
\definecolor{currentstroke}{rgb}{0.000000,0.000000,0.000000}%
\pgfsetstrokecolor{currentstroke}%
\pgfsetstrokeopacity{0.200000}%
\pgfsetdash{}{0pt}%
\pgfpathmoveto{\pgfqpoint{4.000283in}{0.549073in}}%
\pgfpathlineto{\pgfqpoint{4.000283in}{2.859073in}}%
\pgfusepath{stroke}%
\end{pgfscope}%
\begin{pgfscope}%
\pgfsetbuttcap%
\pgfsetroundjoin%
\definecolor{currentfill}{rgb}{0.000000,0.000000,0.000000}%
\pgfsetfillcolor{currentfill}%
\pgfsetlinewidth{0.602250pt}%
\definecolor{currentstroke}{rgb}{0.000000,0.000000,0.000000}%
\pgfsetstrokecolor{currentstroke}%
\pgfsetdash{}{0pt}%
\pgfsys@defobject{currentmarker}{\pgfqpoint{0.000000in}{-0.027778in}}{\pgfqpoint{0.000000in}{0.000000in}}{%
\pgfpathmoveto{\pgfqpoint{0.000000in}{0.000000in}}%
\pgfpathlineto{\pgfqpoint{0.000000in}{-0.027778in}}%
\pgfusepath{stroke,fill}%
}%
\begin{pgfscope}%
\pgfsys@transformshift{4.000283in}{0.549073in}%
\pgfsys@useobject{currentmarker}{}%
\end{pgfscope}%
\end{pgfscope}%
\begin{pgfscope}%
\pgfpathrectangle{\pgfqpoint{0.721913in}{0.549073in}}{\pgfqpoint{4.518250in}{2.310000in}}%
\pgfusepath{clip}%
\pgfsetrectcap%
\pgfsetroundjoin%
\pgfsetlinewidth{0.250937pt}%
\definecolor{currentstroke}{rgb}{0.000000,0.000000,0.000000}%
\pgfsetstrokecolor{currentstroke}%
\pgfsetstrokeopacity{0.200000}%
\pgfsetdash{}{0pt}%
\pgfpathmoveto{\pgfqpoint{4.146311in}{0.549073in}}%
\pgfpathlineto{\pgfqpoint{4.146311in}{2.859073in}}%
\pgfusepath{stroke}%
\end{pgfscope}%
\begin{pgfscope}%
\pgfsetbuttcap%
\pgfsetroundjoin%
\definecolor{currentfill}{rgb}{0.000000,0.000000,0.000000}%
\pgfsetfillcolor{currentfill}%
\pgfsetlinewidth{0.602250pt}%
\definecolor{currentstroke}{rgb}{0.000000,0.000000,0.000000}%
\pgfsetstrokecolor{currentstroke}%
\pgfsetdash{}{0pt}%
\pgfsys@defobject{currentmarker}{\pgfqpoint{0.000000in}{-0.027778in}}{\pgfqpoint{0.000000in}{0.000000in}}{%
\pgfpathmoveto{\pgfqpoint{0.000000in}{0.000000in}}%
\pgfpathlineto{\pgfqpoint{0.000000in}{-0.027778in}}%
\pgfusepath{stroke,fill}%
}%
\begin{pgfscope}%
\pgfsys@transformshift{4.146311in}{0.549073in}%
\pgfsys@useobject{currentmarker}{}%
\end{pgfscope}%
\end{pgfscope}%
\begin{pgfscope}%
\pgfpathrectangle{\pgfqpoint{0.721913in}{0.549073in}}{\pgfqpoint{4.518250in}{2.310000in}}%
\pgfusepath{clip}%
\pgfsetrectcap%
\pgfsetroundjoin%
\pgfsetlinewidth{0.250937pt}%
\definecolor{currentstroke}{rgb}{0.000000,0.000000,0.000000}%
\pgfsetstrokecolor{currentstroke}%
\pgfsetstrokeopacity{0.200000}%
\pgfsetdash{}{0pt}%
\pgfpathmoveto{\pgfqpoint{5.136311in}{0.549073in}}%
\pgfpathlineto{\pgfqpoint{5.136311in}{2.859073in}}%
\pgfusepath{stroke}%
\end{pgfscope}%
\begin{pgfscope}%
\pgfsetbuttcap%
\pgfsetroundjoin%
\definecolor{currentfill}{rgb}{0.000000,0.000000,0.000000}%
\pgfsetfillcolor{currentfill}%
\pgfsetlinewidth{0.602250pt}%
\definecolor{currentstroke}{rgb}{0.000000,0.000000,0.000000}%
\pgfsetstrokecolor{currentstroke}%
\pgfsetdash{}{0pt}%
\pgfsys@defobject{currentmarker}{\pgfqpoint{0.000000in}{-0.027778in}}{\pgfqpoint{0.000000in}{0.000000in}}{%
\pgfpathmoveto{\pgfqpoint{0.000000in}{0.000000in}}%
\pgfpathlineto{\pgfqpoint{0.000000in}{-0.027778in}}%
\pgfusepath{stroke,fill}%
}%
\begin{pgfscope}%
\pgfsys@transformshift{5.136311in}{0.549073in}%
\pgfsys@useobject{currentmarker}{}%
\end{pgfscope}%
\end{pgfscope}%
\begin{pgfscope}%
\definecolor{textcolor}{rgb}{0.000000,0.000000,0.000000}%
\pgfsetstrokecolor{textcolor}%
\pgfsetfillcolor{textcolor}%
\pgftext[x=2.981038in,y=0.248148in,,top]{\color{textcolor}{\rmfamily\fontsize{12.000000}{14.400000}\selectfont\catcode`\^=\active\def^{\ifmmode\sp\else\^{}\fi}\catcode`\%=\active\def
\end{pgfscope}%
\begin{pgfscope}%
\pgfpathrectangle{\pgfqpoint{0.721913in}{0.549073in}}{\pgfqpoint{4.518250in}{2.310000in}}%
\pgfusepath{clip}%
\pgfsetrectcap%
\pgfsetroundjoin%
\pgfsetlinewidth{0.250937pt}%
\definecolor{currentstroke}{rgb}{0.000000,0.000000,0.000000}%
\pgfsetstrokecolor{currentstroke}%
\pgfsetstrokeopacity{0.200000}%
\pgfsetdash{}{0pt}%
\pgfpathmoveto{\pgfqpoint{0.721913in}{0.822963in}}%
\pgfpathlineto{\pgfqpoint{5.240163in}{0.822963in}}%
\pgfusepath{stroke}%
\end{pgfscope}%
\begin{pgfscope}%
\pgfsetbuttcap%
\pgfsetroundjoin%
\definecolor{currentfill}{rgb}{0.000000,0.000000,0.000000}%
\pgfsetfillcolor{currentfill}%
\pgfsetlinewidth{0.803000pt}%
\definecolor{currentstroke}{rgb}{0.000000,0.000000,0.000000}%
\pgfsetstrokecolor{currentstroke}%
\pgfsetdash{}{0pt}%
\pgfsys@defobject{currentmarker}{\pgfqpoint{-0.048611in}{0.000000in}}{\pgfqpoint{-0.000000in}{0.000000in}}{%
\pgfpathmoveto{\pgfqpoint{-0.000000in}{0.000000in}}%
\pgfpathlineto{\pgfqpoint{-0.048611in}{0.000000in}}%
\pgfusepath{stroke,fill}%
}%
\begin{pgfscope}%
\pgfsys@transformshift{0.721913in}{0.822963in}%
\pgfsys@useobject{currentmarker}{}%
\end{pgfscope}%
\end{pgfscope}%
\begin{pgfscope}%
\definecolor{textcolor}{rgb}{0.000000,0.000000,0.000000}%
\pgfsetstrokecolor{textcolor}%
\pgfsetfillcolor{textcolor}%
\pgftext[x=0.303703in, y=0.765093in, left, base]{\color{textcolor}{\rmfamily\fontsize{12.000000}{14.400000}\selectfont\catcode`\^=\active\def^{\ifmmode\sp\else\^{}\fi}\catcode`\%=\active\def
\end{pgfscope}%
\begin{pgfscope}%
\pgfpathrectangle{\pgfqpoint{0.721913in}{0.549073in}}{\pgfqpoint{4.518250in}{2.310000in}}%
\pgfusepath{clip}%
\pgfsetrectcap%
\pgfsetroundjoin%
\pgfsetlinewidth{0.250937pt}%
\definecolor{currentstroke}{rgb}{0.000000,0.000000,0.000000}%
\pgfsetstrokecolor{currentstroke}%
\pgfsetstrokeopacity{0.200000}%
\pgfsetdash{}{0pt}%
\pgfpathmoveto{\pgfqpoint{0.721913in}{1.249888in}}%
\pgfpathlineto{\pgfqpoint{5.240163in}{1.249888in}}%
\pgfusepath{stroke}%
\end{pgfscope}%
\begin{pgfscope}%
\pgfsetbuttcap%
\pgfsetroundjoin%
\definecolor{currentfill}{rgb}{0.000000,0.000000,0.000000}%
\pgfsetfillcolor{currentfill}%
\pgfsetlinewidth{0.803000pt}%
\definecolor{currentstroke}{rgb}{0.000000,0.000000,0.000000}%
\pgfsetstrokecolor{currentstroke}%
\pgfsetdash{}{0pt}%
\pgfsys@defobject{currentmarker}{\pgfqpoint{-0.048611in}{0.000000in}}{\pgfqpoint{-0.000000in}{0.000000in}}{%
\pgfpathmoveto{\pgfqpoint{-0.000000in}{0.000000in}}%
\pgfpathlineto{\pgfqpoint{-0.048611in}{0.000000in}}%
\pgfusepath{stroke,fill}%
}%
\begin{pgfscope}%
\pgfsys@transformshift{0.721913in}{1.249888in}%
\pgfsys@useobject{currentmarker}{}%
\end{pgfscope}%
\end{pgfscope}%
\begin{pgfscope}%
\definecolor{textcolor}{rgb}{0.000000,0.000000,0.000000}%
\pgfsetstrokecolor{textcolor}%
\pgfsetfillcolor{textcolor}%
\pgftext[x=0.303703in, y=1.192018in, left, base]{\color{textcolor}{\rmfamily\fontsize{12.000000}{14.400000}\selectfont\catcode`\^=\active\def^{\ifmmode\sp\else\^{}\fi}\catcode`\%=\active\def
\end{pgfscope}%
\begin{pgfscope}%
\pgfpathrectangle{\pgfqpoint{0.721913in}{0.549073in}}{\pgfqpoint{4.518250in}{2.310000in}}%
\pgfusepath{clip}%
\pgfsetrectcap%
\pgfsetroundjoin%
\pgfsetlinewidth{0.250937pt}%
\definecolor{currentstroke}{rgb}{0.000000,0.000000,0.000000}%
\pgfsetstrokecolor{currentstroke}%
\pgfsetstrokeopacity{0.200000}%
\pgfsetdash{}{0pt}%
\pgfpathmoveto{\pgfqpoint{0.721913in}{1.676813in}}%
\pgfpathlineto{\pgfqpoint{5.240163in}{1.676813in}}%
\pgfusepath{stroke}%
\end{pgfscope}%
\begin{pgfscope}%
\pgfsetbuttcap%
\pgfsetroundjoin%
\definecolor{currentfill}{rgb}{0.000000,0.000000,0.000000}%
\pgfsetfillcolor{currentfill}%
\pgfsetlinewidth{0.803000pt}%
\definecolor{currentstroke}{rgb}{0.000000,0.000000,0.000000}%
\pgfsetstrokecolor{currentstroke}%
\pgfsetdash{}{0pt}%
\pgfsys@defobject{currentmarker}{\pgfqpoint{-0.048611in}{0.000000in}}{\pgfqpoint{-0.000000in}{0.000000in}}{%
\pgfpathmoveto{\pgfqpoint{-0.000000in}{0.000000in}}%
\pgfpathlineto{\pgfqpoint{-0.048611in}{0.000000in}}%
\pgfusepath{stroke,fill}%
}%
\begin{pgfscope}%
\pgfsys@transformshift{0.721913in}{1.676813in}%
\pgfsys@useobject{currentmarker}{}%
\end{pgfscope}%
\end{pgfscope}%
\begin{pgfscope}%
\definecolor{textcolor}{rgb}{0.000000,0.000000,0.000000}%
\pgfsetstrokecolor{textcolor}%
\pgfsetfillcolor{textcolor}%
\pgftext[x=0.303703in, y=1.618943in, left, base]{\color{textcolor}{\rmfamily\fontsize{12.000000}{14.400000}\selectfont\catcode`\^=\active\def^{\ifmmode\sp\else\^{}\fi}\catcode`\%=\active\def
\end{pgfscope}%
\begin{pgfscope}%
\pgfpathrectangle{\pgfqpoint{0.721913in}{0.549073in}}{\pgfqpoint{4.518250in}{2.310000in}}%
\pgfusepath{clip}%
\pgfsetrectcap%
\pgfsetroundjoin%
\pgfsetlinewidth{0.250937pt}%
\definecolor{currentstroke}{rgb}{0.000000,0.000000,0.000000}%
\pgfsetstrokecolor{currentstroke}%
\pgfsetstrokeopacity{0.200000}%
\pgfsetdash{}{0pt}%
\pgfpathmoveto{\pgfqpoint{0.721913in}{2.103738in}}%
\pgfpathlineto{\pgfqpoint{5.240163in}{2.103738in}}%
\pgfusepath{stroke}%
\end{pgfscope}%
\begin{pgfscope}%
\pgfsetbuttcap%
\pgfsetroundjoin%
\definecolor{currentfill}{rgb}{0.000000,0.000000,0.000000}%
\pgfsetfillcolor{currentfill}%
\pgfsetlinewidth{0.803000pt}%
\definecolor{currentstroke}{rgb}{0.000000,0.000000,0.000000}%
\pgfsetstrokecolor{currentstroke}%
\pgfsetdash{}{0pt}%
\pgfsys@defobject{currentmarker}{\pgfqpoint{-0.048611in}{0.000000in}}{\pgfqpoint{-0.000000in}{0.000000in}}{%
\pgfpathmoveto{\pgfqpoint{-0.000000in}{0.000000in}}%
\pgfpathlineto{\pgfqpoint{-0.048611in}{0.000000in}}%
\pgfusepath{stroke,fill}%
}%
\begin{pgfscope}%
\pgfsys@transformshift{0.721913in}{2.103738in}%
\pgfsys@useobject{currentmarker}{}%
\end{pgfscope}%
\end{pgfscope}%
\begin{pgfscope}%
\definecolor{textcolor}{rgb}{0.000000,0.000000,0.000000}%
\pgfsetstrokecolor{textcolor}%
\pgfsetfillcolor{textcolor}%
\pgftext[x=0.395525in, y=2.045868in, left, base]{\color{textcolor}{\rmfamily\fontsize{12.000000}{14.400000}\selectfont\catcode`\^=\active\def^{\ifmmode\sp\else\^{}\fi}\catcode`\%=\active\def
\end{pgfscope}%
\begin{pgfscope}%
\pgfpathrectangle{\pgfqpoint{0.721913in}{0.549073in}}{\pgfqpoint{4.518250in}{2.310000in}}%
\pgfusepath{clip}%
\pgfsetrectcap%
\pgfsetroundjoin%
\pgfsetlinewidth{0.250937pt}%
\definecolor{currentstroke}{rgb}{0.000000,0.000000,0.000000}%
\pgfsetstrokecolor{currentstroke}%
\pgfsetstrokeopacity{0.200000}%
\pgfsetdash{}{0pt}%
\pgfpathmoveto{\pgfqpoint{0.721913in}{2.530663in}}%
\pgfpathlineto{\pgfqpoint{5.240163in}{2.530663in}}%
\pgfusepath{stroke}%
\end{pgfscope}%
\begin{pgfscope}%
\pgfsetbuttcap%
\pgfsetroundjoin%
\definecolor{currentfill}{rgb}{0.000000,0.000000,0.000000}%
\pgfsetfillcolor{currentfill}%
\pgfsetlinewidth{0.803000pt}%
\definecolor{currentstroke}{rgb}{0.000000,0.000000,0.000000}%
\pgfsetstrokecolor{currentstroke}%
\pgfsetdash{}{0pt}%
\pgfsys@defobject{currentmarker}{\pgfqpoint{-0.048611in}{0.000000in}}{\pgfqpoint{-0.000000in}{0.000000in}}{%
\pgfpathmoveto{\pgfqpoint{-0.000000in}{0.000000in}}%
\pgfpathlineto{\pgfqpoint{-0.048611in}{0.000000in}}%
\pgfusepath{stroke,fill}%
}%
\begin{pgfscope}%
\pgfsys@transformshift{0.721913in}{2.530663in}%
\pgfsys@useobject{currentmarker}{}%
\end{pgfscope}%
\end{pgfscope}%
\begin{pgfscope}%
\definecolor{textcolor}{rgb}{0.000000,0.000000,0.000000}%
\pgfsetstrokecolor{textcolor}%
\pgfsetfillcolor{textcolor}%
\pgftext[x=0.395525in, y=2.472793in, left, base]{\color{textcolor}{\rmfamily\fontsize{12.000000}{14.400000}\selectfont\catcode`\^=\active\def^{\ifmmode\sp\else\^{}\fi}\catcode`\%=\active\def
\end{pgfscope}%
\begin{pgfscope}%
\definecolor{textcolor}{rgb}{0.000000,0.000000,0.000000}%
\pgfsetstrokecolor{textcolor}%
\pgfsetfillcolor{textcolor}%
\pgftext[x=0.248148in,y=1.704073in,,bottom,rotate=90.000000]{\color{textcolor}{\rmfamily\fontsize{12.000000}{14.400000}\selectfont\catcode`\^=\active\def^{\ifmmode\sp\else\^{}\fi}\catcode`\%=\active\def
\end{pgfscope}%
\begin{pgfscope}%
\pgfpathrectangle{\pgfqpoint{0.721913in}{0.549073in}}{\pgfqpoint{4.518250in}{2.310000in}}%
\pgfusepath{clip}%
\pgfsetrectcap%
\pgfsetroundjoin%
\pgfsetlinewidth{1.505625pt}%
\definecolor{currentstroke}{rgb}{0.392157,0.560784,1.000000}%
\pgfsetstrokecolor{currentstroke}%
\pgfsetdash{}{0pt}%
\pgfpathmoveto{\pgfqpoint{0.829491in}{2.666573in}}%
\pgfpathlineto{\pgfqpoint{1.562672in}{2.156909in}}%
\pgfpathlineto{\pgfqpoint{2.269077in}{2.009321in}}%
\pgfpathlineto{\pgfqpoint{2.989445in}{1.828200in}}%
\pgfpathlineto{\pgfqpoint{3.704102in}{1.411721in}}%
\pgfpathlineto{\pgfqpoint{4.419657in}{0.741665in}}%
\pgfpathlineto{\pgfqpoint{5.132586in}{0.741575in}}%
\pgfusepath{stroke}%
\end{pgfscope}%
\begin{pgfscope}%
\pgfpathrectangle{\pgfqpoint{0.721913in}{0.549073in}}{\pgfqpoint{4.518250in}{2.310000in}}%
\pgfusepath{clip}%
\pgfsetbuttcap%
\pgfsetroundjoin%
\definecolor{currentfill}{rgb}{0.392157,0.560784,1.000000}%
\pgfsetfillcolor{currentfill}%
\pgfsetlinewidth{1.003750pt}%
\definecolor{currentstroke}{rgb}{0.392157,0.560784,1.000000}%
\pgfsetstrokecolor{currentstroke}%
\pgfsetdash{}{0pt}%
\pgfsys@defobject{currentmarker}{\pgfqpoint{-0.041667in}{-0.041667in}}{\pgfqpoint{0.041667in}{0.041667in}}{%
\pgfpathmoveto{\pgfqpoint{0.000000in}{-0.041667in}}%
\pgfpathcurveto{\pgfqpoint{0.011050in}{-0.041667in}}{\pgfqpoint{0.021649in}{-0.037276in}}{\pgfqpoint{0.029463in}{-0.029463in}}%
\pgfpathcurveto{\pgfqpoint{0.037276in}{-0.021649in}}{\pgfqpoint{0.041667in}{-0.011050in}}{\pgfqpoint{0.041667in}{0.000000in}}%
\pgfpathcurveto{\pgfqpoint{0.041667in}{0.011050in}}{\pgfqpoint{0.037276in}{0.021649in}}{\pgfqpoint{0.029463in}{0.029463in}}%
\pgfpathcurveto{\pgfqpoint{0.021649in}{0.037276in}}{\pgfqpoint{0.011050in}{0.041667in}}{\pgfqpoint{0.000000in}{0.041667in}}%
\pgfpathcurveto{\pgfqpoint{-0.011050in}{0.041667in}}{\pgfqpoint{-0.021649in}{0.037276in}}{\pgfqpoint{-0.029463in}{0.029463in}}%
\pgfpathcurveto{\pgfqpoint{-0.037276in}{0.021649in}}{\pgfqpoint{-0.041667in}{0.011050in}}{\pgfqpoint{-0.041667in}{0.000000in}}%
\pgfpathcurveto{\pgfqpoint{-0.041667in}{-0.011050in}}{\pgfqpoint{-0.037276in}{-0.021649in}}{\pgfqpoint{-0.029463in}{-0.029463in}}%
\pgfpathcurveto{\pgfqpoint{-0.021649in}{-0.037276in}}{\pgfqpoint{-0.011050in}{-0.041667in}}{\pgfqpoint{0.000000in}{-0.041667in}}%
\pgfpathlineto{\pgfqpoint{0.000000in}{-0.041667in}}%
\pgfpathclose%
\pgfusepath{stroke,fill}%
}%
\begin{pgfscope}%
\pgfsys@transformshift{0.829491in}{2.666573in}%
\pgfsys@useobject{currentmarker}{}%
\end{pgfscope}%
\begin{pgfscope}%
\pgfsys@transformshift{1.562672in}{2.156909in}%
\pgfsys@useobject{currentmarker}{}%
\end{pgfscope}%
\begin{pgfscope}%
\pgfsys@transformshift{2.269077in}{2.009321in}%
\pgfsys@useobject{currentmarker}{}%
\end{pgfscope}%
\begin{pgfscope}%
\pgfsys@transformshift{2.989445in}{1.828200in}%
\pgfsys@useobject{currentmarker}{}%
\end{pgfscope}%
\begin{pgfscope}%
\pgfsys@transformshift{3.704102in}{1.411721in}%
\pgfsys@useobject{currentmarker}{}%
\end{pgfscope}%
\begin{pgfscope}%
\pgfsys@transformshift{4.419657in}{0.741665in}%
\pgfsys@useobject{currentmarker}{}%
\end{pgfscope}%
\begin{pgfscope}%
\pgfsys@transformshift{5.132586in}{0.741575in}%
\pgfsys@useobject{currentmarker}{}%
\end{pgfscope}%
\end{pgfscope}%
\begin{pgfscope}%
\pgfpathrectangle{\pgfqpoint{0.721913in}{0.549073in}}{\pgfqpoint{4.518250in}{2.310000in}}%
\pgfusepath{clip}%
\pgfsetrectcap%
\pgfsetroundjoin%
\pgfsetlinewidth{1.505625pt}%
\definecolor{currentstroke}{rgb}{0.862745,0.149020,0.498039}%
\pgfsetstrokecolor{currentstroke}%
\pgfsetdash{}{0pt}%
\pgfpathmoveto{\pgfqpoint{0.829491in}{1.951822in}}%
\pgfpathlineto{\pgfqpoint{1.562672in}{1.960635in}}%
\pgfpathlineto{\pgfqpoint{2.269077in}{1.909660in}}%
\pgfpathlineto{\pgfqpoint{2.989445in}{1.715951in}}%
\pgfpathlineto{\pgfqpoint{3.704102in}{1.319000in}}%
\pgfpathlineto{\pgfqpoint{4.419657in}{0.741721in}}%
\pgfpathlineto{\pgfqpoint{5.132586in}{0.741573in}}%
\pgfusepath{stroke}%
\end{pgfscope}%
\begin{pgfscope}%
\pgfpathrectangle{\pgfqpoint{0.721913in}{0.549073in}}{\pgfqpoint{4.518250in}{2.310000in}}%
\pgfusepath{clip}%
\pgfsetbuttcap%
\pgfsetmiterjoin%
\definecolor{currentfill}{rgb}{0.862745,0.149020,0.498039}%
\pgfsetfillcolor{currentfill}%
\pgfsetlinewidth{1.003750pt}%
\definecolor{currentstroke}{rgb}{0.862745,0.149020,0.498039}%
\pgfsetstrokecolor{currentstroke}%
\pgfsetdash{}{0pt}%
\pgfsys@defobject{currentmarker}{\pgfqpoint{-0.041667in}{-0.041667in}}{\pgfqpoint{0.041667in}{0.041667in}}{%
\pgfpathmoveto{\pgfqpoint{-0.041667in}{-0.041667in}}%
\pgfpathlineto{\pgfqpoint{0.041667in}{-0.041667in}}%
\pgfpathlineto{\pgfqpoint{0.041667in}{0.041667in}}%
\pgfpathlineto{\pgfqpoint{-0.041667in}{0.041667in}}%
\pgfpathlineto{\pgfqpoint{-0.041667in}{-0.041667in}}%
\pgfpathclose%
\pgfusepath{stroke,fill}%
}%
\begin{pgfscope}%
\pgfsys@transformshift{0.829491in}{1.951822in}%
\pgfsys@useobject{currentmarker}{}%
\end{pgfscope}%
\begin{pgfscope}%
\pgfsys@transformshift{1.562672in}{1.960635in}%
\pgfsys@useobject{currentmarker}{}%
\end{pgfscope}%
\begin{pgfscope}%
\pgfsys@transformshift{2.269077in}{1.909660in}%
\pgfsys@useobject{currentmarker}{}%
\end{pgfscope}%
\begin{pgfscope}%
\pgfsys@transformshift{2.989445in}{1.715951in}%
\pgfsys@useobject{currentmarker}{}%
\end{pgfscope}%
\begin{pgfscope}%
\pgfsys@transformshift{3.704102in}{1.319000in}%
\pgfsys@useobject{currentmarker}{}%
\end{pgfscope}%
\begin{pgfscope}%
\pgfsys@transformshift{4.419657in}{0.741721in}%
\pgfsys@useobject{currentmarker}{}%
\end{pgfscope}%
\begin{pgfscope}%
\pgfsys@transformshift{5.132586in}{0.741573in}%
\pgfsys@useobject{currentmarker}{}%
\end{pgfscope}%
\end{pgfscope}%
\begin{pgfscope}%
\pgfpathrectangle{\pgfqpoint{0.721913in}{0.549073in}}{\pgfqpoint{4.518250in}{2.310000in}}%
\pgfusepath{clip}%
\pgfsetrectcap%
\pgfsetroundjoin%
\pgfsetlinewidth{1.505625pt}%
\definecolor{currentstroke}{rgb}{1.000000,0.690196,0.000000}%
\pgfsetstrokecolor{currentstroke}%
\pgfsetdash{}{0pt}%
\pgfpathmoveto{\pgfqpoint{0.829491in}{1.925108in}}%
\pgfpathlineto{\pgfqpoint{1.562672in}{1.890578in}}%
\pgfpathlineto{\pgfqpoint{2.269077in}{1.876566in}}%
\pgfpathlineto{\pgfqpoint{2.989445in}{1.772451in}}%
\pgfpathlineto{\pgfqpoint{3.704102in}{1.527958in}}%
\pgfpathlineto{\pgfqpoint{4.419657in}{0.901857in}}%
\pgfpathlineto{\pgfqpoint{5.132586in}{0.741573in}}%
\pgfusepath{stroke}%
\end{pgfscope}%
\begin{pgfscope}%
\pgfpathrectangle{\pgfqpoint{0.721913in}{0.549073in}}{\pgfqpoint{4.518250in}{2.310000in}}%
\pgfusepath{clip}%
\pgfsetbuttcap%
\pgfsetmiterjoin%
\definecolor{currentfill}{rgb}{1.000000,0.690196,0.000000}%
\pgfsetfillcolor{currentfill}%
\pgfsetlinewidth{1.003750pt}%
\definecolor{currentstroke}{rgb}{1.000000,0.690196,0.000000}%
\pgfsetstrokecolor{currentstroke}%
\pgfsetdash{}{0pt}%
\pgfsys@defobject{currentmarker}{\pgfqpoint{-0.035355in}{-0.058926in}}{\pgfqpoint{0.035355in}{0.058926in}}{%
\pgfpathmoveto{\pgfqpoint{-0.000000in}{-0.058926in}}%
\pgfpathlineto{\pgfqpoint{0.035355in}{0.000000in}}%
\pgfpathlineto{\pgfqpoint{0.000000in}{0.058926in}}%
\pgfpathlineto{\pgfqpoint{-0.035355in}{0.000000in}}%
\pgfpathlineto{\pgfqpoint{-0.000000in}{-0.058926in}}%
\pgfpathclose%
\pgfusepath{stroke,fill}%
}%
\begin{pgfscope}%
\pgfsys@transformshift{0.829491in}{1.925108in}%
\pgfsys@useobject{currentmarker}{}%
\end{pgfscope}%
\begin{pgfscope}%
\pgfsys@transformshift{1.562672in}{1.890578in}%
\pgfsys@useobject{currentmarker}{}%
\end{pgfscope}%
\begin{pgfscope}%
\pgfsys@transformshift{2.269077in}{1.876566in}%
\pgfsys@useobject{currentmarker}{}%
\end{pgfscope}%
\begin{pgfscope}%
\pgfsys@transformshift{2.989445in}{1.772451in}%
\pgfsys@useobject{currentmarker}{}%
\end{pgfscope}%
\begin{pgfscope}%
\pgfsys@transformshift{3.704102in}{1.527958in}%
\pgfsys@useobject{currentmarker}{}%
\end{pgfscope}%
\begin{pgfscope}%
\pgfsys@transformshift{4.419657in}{0.901857in}%
\pgfsys@useobject{currentmarker}{}%
\end{pgfscope}%
\begin{pgfscope}%
\pgfsys@transformshift{5.132586in}{0.741573in}%
\pgfsys@useobject{currentmarker}{}%
\end{pgfscope}%
\end{pgfscope}%
\begin{pgfscope}%
\pgfsetrectcap%
\pgfsetmiterjoin%
\pgfsetlinewidth{0.803000pt}%
\definecolor{currentstroke}{rgb}{0.000000,0.000000,0.000000}%
\pgfsetstrokecolor{currentstroke}%
\pgfsetdash{}{0pt}%
\pgfpathmoveto{\pgfqpoint{0.721913in}{0.549073in}}%
\pgfpathlineto{\pgfqpoint{0.721913in}{2.859073in}}%
\pgfusepath{stroke}%
\end{pgfscope}%
\begin{pgfscope}%
\pgfsetrectcap%
\pgfsetmiterjoin%
\pgfsetlinewidth{0.803000pt}%
\definecolor{currentstroke}{rgb}{0.000000,0.000000,0.000000}%
\pgfsetstrokecolor{currentstroke}%
\pgfsetdash{}{0pt}%
\pgfpathmoveto{\pgfqpoint{5.240163in}{0.549073in}}%
\pgfpathlineto{\pgfqpoint{5.240163in}{2.859073in}}%
\pgfusepath{stroke}%
\end{pgfscope}%
\begin{pgfscope}%
\pgfsetrectcap%
\pgfsetmiterjoin%
\pgfsetlinewidth{0.803000pt}%
\definecolor{currentstroke}{rgb}{0.000000,0.000000,0.000000}%
\pgfsetstrokecolor{currentstroke}%
\pgfsetdash{}{0pt}%
\pgfpathmoveto{\pgfqpoint{0.721913in}{0.549073in}}%
\pgfpathlineto{\pgfqpoint{5.240163in}{0.549073in}}%
\pgfusepath{stroke}%
\end{pgfscope}%
\begin{pgfscope}%
\pgfsetrectcap%
\pgfsetmiterjoin%
\pgfsetlinewidth{0.803000pt}%
\definecolor{currentstroke}{rgb}{0.000000,0.000000,0.000000}%
\pgfsetstrokecolor{currentstroke}%
\pgfsetdash{}{0pt}%
\pgfpathmoveto{\pgfqpoint{0.721913in}{2.859073in}}%
\pgfpathlineto{\pgfqpoint{5.240163in}{2.859073in}}%
\pgfusepath{stroke}%
\end{pgfscope}%
\begin{pgfscope}%
\pgfsetbuttcap%
\pgfsetmiterjoin%
\definecolor{currentfill}{rgb}{1.000000,1.000000,1.000000}%
\pgfsetfillcolor{currentfill}%
\pgfsetfillopacity{0.800000}%
\pgfsetlinewidth{1.003750pt}%
\definecolor{currentstroke}{rgb}{0.800000,0.800000,0.800000}%
\pgfsetstrokecolor{currentstroke}%
\pgfsetstrokeopacity{0.800000}%
\pgfsetdash{}{0pt}%
\pgfpathmoveto{\pgfqpoint{2.710032in}{2.028518in}}%
\pgfpathlineto{\pgfqpoint{5.156830in}{2.028518in}}%
\pgfpathlineto{\pgfqpoint{5.156830in}{2.775739in}}%
\pgfpathlineto{\pgfqpoint{2.710032in}{2.775739in}}%
\pgfpathlineto{\pgfqpoint{2.710032in}{2.028518in}}%
\pgfpathclose%
\pgfusepath{stroke,fill}%
\end{pgfscope}%
\begin{pgfscope}%
\pgfsetrectcap%
\pgfsetroundjoin%
\pgfsetlinewidth{1.505625pt}%
\definecolor{currentstroke}{rgb}{0.392157,0.560784,1.000000}%
\pgfsetstrokecolor{currentstroke}%
\pgfsetdash{}{0pt}%
\pgfpathmoveto{\pgfqpoint{2.776699in}{2.650739in}}%
\pgfpathlineto{\pgfqpoint{2.943366in}{2.650739in}}%
\pgfpathlineto{\pgfqpoint{3.110032in}{2.650739in}}%
\pgfusepath{stroke}%
\end{pgfscope}%
\begin{pgfscope}%
\pgfsetbuttcap%
\pgfsetroundjoin%
\definecolor{currentfill}{rgb}{0.392157,0.560784,1.000000}%
\pgfsetfillcolor{currentfill}%
\pgfsetlinewidth{1.003750pt}%
\definecolor{currentstroke}{rgb}{0.392157,0.560784,1.000000}%
\pgfsetstrokecolor{currentstroke}%
\pgfsetdash{}{0pt}%
\pgfsys@defobject{currentmarker}{\pgfqpoint{-0.031250in}{-0.031250in}}{\pgfqpoint{0.031250in}{0.031250in}}{%
\pgfpathmoveto{\pgfqpoint{0.000000in}{-0.031250in}}%
\pgfpathcurveto{\pgfqpoint{0.008288in}{-0.031250in}}{\pgfqpoint{0.016237in}{-0.027957in}}{\pgfqpoint{0.022097in}{-0.022097in}}%
\pgfpathcurveto{\pgfqpoint{0.027957in}{-0.016237in}}{\pgfqpoint{0.031250in}{-0.008288in}}{\pgfqpoint{0.031250in}{0.000000in}}%
\pgfpathcurveto{\pgfqpoint{0.031250in}{0.008288in}}{\pgfqpoint{0.027957in}{0.016237in}}{\pgfqpoint{0.022097in}{0.022097in}}%
\pgfpathcurveto{\pgfqpoint{0.016237in}{0.027957in}}{\pgfqpoint{0.008288in}{0.031250in}}{\pgfqpoint{0.000000in}{0.031250in}}%
\pgfpathcurveto{\pgfqpoint{-0.008288in}{0.031250in}}{\pgfqpoint{-0.016237in}{0.027957in}}{\pgfqpoint{-0.022097in}{0.022097in}}%
\pgfpathcurveto{\pgfqpoint{-0.027957in}{0.016237in}}{\pgfqpoint{-0.031250in}{0.008288in}}{\pgfqpoint{-0.031250in}{0.000000in}}%
\pgfpathcurveto{\pgfqpoint{-0.031250in}{-0.008288in}}{\pgfqpoint{-0.027957in}{-0.016237in}}{\pgfqpoint{-0.022097in}{-0.022097in}}%
\pgfpathcurveto{\pgfqpoint{-0.016237in}{-0.027957in}}{\pgfqpoint{-0.008288in}{-0.031250in}}{\pgfqpoint{0.000000in}{-0.031250in}}%
\pgfpathlineto{\pgfqpoint{0.000000in}{-0.031250in}}%
\pgfpathclose%
\pgfusepath{stroke,fill}%
}%
\begin{pgfscope}%
\pgfsys@transformshift{2.943366in}{2.650739in}%
\pgfsys@useobject{currentmarker}{}%
\end{pgfscope}%
\end{pgfscope}%
\begin{pgfscope}%
\definecolor{textcolor}{rgb}{0.000000,0.000000,0.000000}%
\pgfsetstrokecolor{textcolor}%
\pgfsetfillcolor{textcolor}%
\pgftext[x=3.243366in,y=2.592406in,left,base]{\color{textcolor}{\rmfamily\fontsize{12.000000}{14.400000}\selectfont\catcode`\^=\active\def^{\ifmmode\sp\else\^{}\fi}\catcode`\%=\active\def
\end{pgfscope}%
\begin{pgfscope}%
\pgfsetrectcap%
\pgfsetroundjoin%
\pgfsetlinewidth{1.505625pt}%
\definecolor{currentstroke}{rgb}{0.862745,0.149020,0.498039}%
\pgfsetstrokecolor{currentstroke}%
\pgfsetdash{}{0pt}%
\pgfpathmoveto{\pgfqpoint{2.776699in}{2.418332in}}%
\pgfpathlineto{\pgfqpoint{2.943366in}{2.418332in}}%
\pgfpathlineto{\pgfqpoint{3.110032in}{2.418332in}}%
\pgfusepath{stroke}%
\end{pgfscope}%
\begin{pgfscope}%
\pgfsetbuttcap%
\pgfsetmiterjoin%
\definecolor{currentfill}{rgb}{0.862745,0.149020,0.498039}%
\pgfsetfillcolor{currentfill}%
\pgfsetlinewidth{1.003750pt}%
\definecolor{currentstroke}{rgb}{0.862745,0.149020,0.498039}%
\pgfsetstrokecolor{currentstroke}%
\pgfsetdash{}{0pt}%
\pgfsys@defobject{currentmarker}{\pgfqpoint{-0.031250in}{-0.031250in}}{\pgfqpoint{0.031250in}{0.031250in}}{%
\pgfpathmoveto{\pgfqpoint{-0.031250in}{-0.031250in}}%
\pgfpathlineto{\pgfqpoint{0.031250in}{-0.031250in}}%
\pgfpathlineto{\pgfqpoint{0.031250in}{0.031250in}}%
\pgfpathlineto{\pgfqpoint{-0.031250in}{0.031250in}}%
\pgfpathlineto{\pgfqpoint{-0.031250in}{-0.031250in}}%
\pgfpathclose%
\pgfusepath{stroke,fill}%
}%
\begin{pgfscope}%
\pgfsys@transformshift{2.943366in}{2.418332in}%
\pgfsys@useobject{currentmarker}{}%
\end{pgfscope}%
\end{pgfscope}%
\begin{pgfscope}%
\definecolor{textcolor}{rgb}{0.000000,0.000000,0.000000}%
\pgfsetstrokecolor{textcolor}%
\pgfsetfillcolor{textcolor}%
\pgftext[x=3.243366in,y=2.359999in,left,base]{\color{textcolor}{\rmfamily\fontsize{12.000000}{14.400000}\selectfont\catcode`\^=\active\def^{\ifmmode\sp\else\^{}\fi}\catcode`\%=\active\def
\end{pgfscope}%
\begin{pgfscope}%
\pgfsetrectcap%
\pgfsetroundjoin%
\pgfsetlinewidth{1.505625pt}%
\definecolor{currentstroke}{rgb}{1.000000,0.690196,0.000000}%
\pgfsetstrokecolor{currentstroke}%
\pgfsetdash{}{0pt}%
\pgfpathmoveto{\pgfqpoint{2.776699in}{2.185925in}}%
\pgfpathlineto{\pgfqpoint{2.943366in}{2.185925in}}%
\pgfpathlineto{\pgfqpoint{3.110032in}{2.185925in}}%
\pgfusepath{stroke}%
\end{pgfscope}%
\begin{pgfscope}%
\pgfsetbuttcap%
\pgfsetmiterjoin%
\definecolor{currentfill}{rgb}{1.000000,0.690196,0.000000}%
\pgfsetfillcolor{currentfill}%
\pgfsetlinewidth{1.003750pt}%
\definecolor{currentstroke}{rgb}{1.000000,0.690196,0.000000}%
\pgfsetstrokecolor{currentstroke}%
\pgfsetdash{}{0pt}%
\pgfsys@defobject{currentmarker}{\pgfqpoint{-0.026517in}{-0.044194in}}{\pgfqpoint{0.026517in}{0.044194in}}{%
\pgfpathmoveto{\pgfqpoint{-0.000000in}{-0.044194in}}%
\pgfpathlineto{\pgfqpoint{0.026517in}{0.000000in}}%
\pgfpathlineto{\pgfqpoint{0.000000in}{0.044194in}}%
\pgfpathlineto{\pgfqpoint{-0.026517in}{0.000000in}}%
\pgfpathlineto{\pgfqpoint{-0.000000in}{-0.044194in}}%
\pgfpathclose%
\pgfusepath{stroke,fill}%
}%
\begin{pgfscope}%
\pgfsys@transformshift{2.943366in}{2.185925in}%
\pgfsys@useobject{currentmarker}{}%
\end{pgfscope}%
\end{pgfscope}%
\begin{pgfscope}%
\definecolor{textcolor}{rgb}{0.000000,0.000000,0.000000}%
\pgfsetstrokecolor{textcolor}%
\pgfsetfillcolor{textcolor}%
\pgftext[x=3.243366in,y=2.127592in,left,base]{\color{textcolor}{\rmfamily\fontsize{12.000000}{14.400000}\selectfont\catcode`\^=\active\def^{\ifmmode\sp\else\^{}\fi}\catcode`\%=\active\def
\end{pgfscope}%
\end{pgfpicture}%
\makeatother%
\endgroup%

%% file: figures/numerical-rank.tex
\begin{tikzpicture}
    \fill[black!15!white] (-2.0, -0.15) rectangle (2.0, 0.15);
    \draw[thick, ->] (-5, 0) to (5, 0);
    \fill[black] (-4, 0) circle (0.075) node[below] {$\lambda_{(7)}$};
    \fill[black] (-3, 0) circle (0.075) node[below] {$\lambda_{(5)}$};
    \fill[black] (-2.35, 0) circle (0.075) node[below] {$\lambda_{(4)}$};
    \fill[black] (-0.8, 0) circle (0.075) node[below] {$\lambda_{(2)}$};
    \fill[black] (0.3, 0) circle (0.075) node[below] {$\lambda_{(1)}$};
    \fill[black] (1, 0) circle (0.075) node[below] {$\lambda_{(3)}$};
    \fill[black] (3, 0) circle (0.075) node[below] {$\lambda_{(6)}$};
    \fill[black] (4.2, 0) circle (0.075) node[below] {$\lambda_{(8)}$};
    \fill[black] (4.8, 0) circle (0.075) node[below] {$\lambda_{(9)}$};
    \fill[black] (-4, 0) circle (0.075) node[above] {$\sigma_{7}$};
    \fill[black] (-3, 0) circle (0.075) node[above] {$\sigma_{5}$};
    \fill[black] (-2.35, 0) circle (0.075) node[above] {$\sigma_{4}$};
    \fill[black] (-0.8, 1.3) circle (0.075) node[above left] {$\sigma_{2}$};
    \fill[black] (0.3, 2.28) circle (0.075) node[above right] {$\sigma_{1}$};
    \fill[black] (1, 0.91) circle (0.075) node[above right] {$\sigma_{3}$};
    \fill[black] (3, 0) circle (0.075) node[above] {$\sigma_{6}$};
    \fill[black] (4.2, 0) circle (0.075) node[above] {$\sigma_{8}$};
    \fill[black] (4.8, 0) circle (0.075) node[above] {$\sigma_{9}$};
    \draw[black, ultra thick] (-0.05, 0.15) to (-0.05, -0.55) node[below] {$t$};
    \draw[black, thick] (-2.0, 0.15) to (-2.0, -0.55) node[below] {$t - d$};
    \draw[black, thick] (2.0, 0.15) to (2.0, -0.55) node[below] {$t + d$};
    \draw[black, thick, dashed] (-0.8, 0) to (-0.8, 1.3);
    \draw[black, thick, dashed] (0.3, 0) to (0.3, 2.28);
    \draw[black, thick, dashed] (1, 0) to (1, 0.91);
    \draw[black!60!white, thick] (-5.000, 0.000) to (-4.899, 0.000) to (-4.798, 0.000) to (-4.697, 0.000) to (-4.596, 0.000) to (-4.495, 0.000) to (-4.394, 0.000) to (-4.293, 0.000) to (-4.192, 0.000) to (-4.091, 0.000) to (-3.990, 0.000) to (-3.889, 0.000) to (-3.788, 0.000) to (-3.687, 0.000) to (-3.586, 0.000) to (-3.485, 0.000) to (-3.384, 0.000) to (-3.283, 0.000) to (-3.182, 0.000) to (-3.081, 0.000) to (-2.980, 0.000) to (-2.879, 0.001) to (-2.778, 0.001) to (-2.677, 0.002) to (-2.576, 0.003) to (-2.475, 0.005) to (-2.374, 0.009) to (-2.273, 0.014) to (-2.172, 0.022) to (-2.071, 0.034) to (-1.970, 0.052) to (-1.869, 0.076) to (-1.768, 0.110) to (-1.667, 0.155) to (-1.566, 0.215) to (-1.465, 0.293) to (-1.364, 0.389) to (-1.263, 0.508) to (-1.162, 0.649) to (-1.061, 0.812) to (-0.960, 0.995) to (-0.859, 1.196) to (-0.758, 1.408) to (-0.657, 1.625) to (-0.556, 1.836) to (-0.455, 2.033) to (-0.354, 2.206) to (-0.253, 2.346) to (-0.152, 2.443) to (-0.051, 2.494) to (0.051, 2.494) to (0.152, 2.443) to (0.253, 2.346) to (0.354, 2.206) to (0.455, 2.033) to (0.556, 1.836) to (0.657, 1.625) to (0.758, 1.408) to (0.859, 1.196) to (0.960, 0.995) to (1.061, 0.812) to (1.162, 0.649) to (1.263, 0.508) to (1.364, 0.389) to (1.465, 0.293) to (1.566, 0.215) to (1.667, 0.155) to (1.768, 0.110) to (1.869, 0.076) to (1.970, 0.052) to (2.071, 0.034) to (2.172, 0.022) to (2.273, 0.014) to (2.374, 0.009) to (2.475, 0.005) to (2.576, 0.003) to (2.677, 0.002) to (2.778, 0.001) to (2.879, 0.001) to (2.980, 0.000) to (3.081, 0.000) to (3.182, 0.000) to (3.283, 0.000) to (3.384, 0.000) to (3.485, 0.000) to (3.586, 0.000) to (3.687, 0.000) to (3.788, 0.000) to (3.889, 0.000) to (3.990, 0.000) to (4.091, 0.000) to (4.192, 0.000) to (4.293, 0.000) to (4.394, 0.000) to (4.495, 0.000) to (4.596, 0.000) to (4.697, 0.000) to (4.798, 0.000) to (4.899, 0.000) to (5.000, 0.000);
    \node[black!60!white] at (-1.3, 2.3) {$g_{\sigma}(t - \cdot)$};
\end{tikzpicture}

%% file: numresults.tex
\section{Numerical experiments}
\label{sec:results}

In this section, we verify the performance of the Chebyshev-Nyström++ algorithm proposed in this paper in application scenarios from electronic structure interaction and neural network optimization.

While theoretically an important tool, we have observed that the preservation of non-negativity discussed in~\cref{sec:preservnonneg} is not needed in practice. The (slight) indefiniteness of the standard Chebyshev approximation $g_{\sigma}^{(m)}(t \mtx{I}_n - \mtx{A})$ from~\cref{equ:matrix-approximation} does not seem to impede the reliability of the Nyström++ estimator. Unless otherwise stated, we transform the input matrices such that their eigenvalues are contained in $[-1, 1]$, by approximating their spectrum with NumPy's Hermitian eigenvalue solver before we apply the methods to them. We compute spectral densities at $n_t = 100$ uniformly distributed values of the parameter $t \in [-1, 1]$. The integral involved in computing the $L^1$-errors is approximated with the composite midpoint quadrature rule using the $n_t$ values of $t$ as nodes. Both parameters introduced at the end of \cref{subsubsec:chebyshev-nystrom-implementation}, the eigenvalue truncation threshold and the parameter for detecting a vanishing spectral density, are fixed to $10^{-5}$.

Our implementations are written in Python 3.12.3 using the packages NumPy 2.0.0 and SciPy 1.14.1. They are executed on a single thread of a GitHub-hosted Ubuntu runner with a 64-bit processor and 16 GB of RAM. 

\subsection{Spectral density for Hamiltonian of electronic structure}
\label{subsec:hamiltonian}

We use the electronic structure interaction example from~\cite[Section 6]{lin-2017-randomized-estimation}, which involves the second order finite difference discretization of the Hamiltonian
\begin{equation}
    \mathcal{H} = - \Delta + V
    \label{equ:electronic-hamiltonian}
\end{equation}
in three dimensions. The potential $V$ interacting with the electrons is generated by Gaussian wells
    $v(r) = v_0 e^{-\lambda r^2}$, 
with $v_0 = -4$ and $\lambda = 8$, centered in cells of side length $L=6$ which are stacked $n_c \in \mathbb{N}$ times in each spatial dimension; see \cref{fig:gaussian-well}. The grid width of the finite difference discretization is fixed to $h=0.6$. For $n_c = 1$, this leads to a sparse matrix of size $1\,000 \times 1\,000$ and for $n_c = 3$ of size $27\,000 \times 27\,000$. This example represents an idealized model for the interaction of nuclei on a regular grid with electrons for a $k$-vector in the center of the first Brillouin zone. The distribution of the eigenvalues of the Hamiltonian --- its spectral density --- allows one to interpret the system's energy levels.

\begin{figure}[ht]
    \begin{subfigure}[b]{0.32\columnwidth}
        \scalebox{0.7}{\input{plots/gaussian-well-1.pgf}}
        \caption{$n_c=1.$}
        \label{fig:gaussian-well-1}
    \end{subfigure}
    \begin{subfigure}[b]{0.32\columnwidth}
        \scalebox{0.7}{\input{plots/gaussian-well-2.pgf}}
        \caption{$n_c=2.$}
        \label{fig:gaussian-well-2}
    \end{subfigure}
    \begin{subfigure}[b]{0.32\columnwidth}
        \scalebox{0.7}{\input{plots/gaussian-well-3.pgf}}
        \caption{$n_c=3.$}
        \label{fig:gaussian-well-3}
    \end{subfigure}
    \caption{Cross-sections of the periodic Gaussian wells potential $V$ for different $n_c$.}
    \label{fig:gaussian-well}
\end{figure}

\cref{fig:convergence} displays the convergence behavior of Chebyshev-Nyström++ applied to the matrix $\mtx{A}$ resulting from the discretization of the Hamiltonian for $n_c = 1$ and $n_c = 3$. Unlike \cite{lin-2017-randomized-estimation}, where just a small portion of the spectral density is considered, we approximate the whole spectrum of $\mtx{A}$ (transformed to $[-1, 1]$) in accordance with our theoretical results (\cref{subsubsec:chebyshev-nystrom-analysis}) and to demonstrate the effectiveness of our implementational particularities (\cref{subsubsec:chebyshev-nystrom-implementation}). As a reference, we use the eigenvalues computed by NumPy's Hermitian eigenvalue solver. As expected from \cref{thm:hutchinson}, Girard-Hutchinson alone ($n_{\mtx{\Omega}} = 0$) requires $n_{\mtx{\Psi}} = \mathcal{O}(\varepsilon^{-2})$ samples to achieve an error of order $\varepsilon$. Because the eigenvalues of the matrix are spread out, we observe the better convergence of the Nyström approximation discussed at the end of \cref{sec:application} once $n_{\mtx{\Omega}}$ is sufficiently large. At this point, we can observe that the errors stagnate around $10^{-6}$ as a consequence of the error made in the Chebyshev approximation. Due to the higher eigenvalue density of $\mtx{A}$ for $n_c = 3$, the stronger convergence can be observed to kick in later in this case, because a larger $n_{\mtx{\Omega}}$ is needed for a good Nyström approximation. In fact, for $n_c = 5$ the eigenvalue density is so high that to profit from this effect, $n_{\mtx{\Omega}}$ would need to be chosen so large that the computations become unfeasible.
\begin{figure}[ht]
    \centering
    \begin{subfigure}[b]{0.495\textwidth}
        \scalebox{0.8}{\input{plots/convergence-1.pgf}}
        \caption{$n_c = 1.$}
        \label{fig:convergence-1}
    \end{subfigure}
    \begin{subfigure}[b]{0.495\textwidth}
        \scalebox{0.8}{\input{plots/convergence-3.pgf}}
        \caption{$n_c = 3.$}
        \label{fig:convergence-3}
    \end{subfigure}
    \caption{Error vs. number of random vectors $n_{\mtx{\Omega}} + n_{\mtx{\Psi}}$ when applying  Chebyshev-Nyström++ to the Hamiltonian matrix from~\cref{subsec:hamiltonian}, with fixed $\sigma=0.005$ and $m=1\,000$, such that the error of the Chebyshev approximation remains sufficiently small according to \cref{fig:chebyshev-heatmap}.}
    \label{fig:convergence}
\end{figure}

For a fixed budget of $n_{\mtx{\Omega}} + n_{\mtx{\Psi}} = 80$ of random vectors, \cref{fig:distribution} shows how the value of the smoothing parameter $\sigma$ impacts the choice between Nyström ($n_{\mtx{\Omega}}$) and 
Girard-Hutchinson ($n_{\mtx{\Psi}}$). As expected, for small $\sigma$ the Nyström approximation alone suffices, but 
this is not an effective choice for large $\sigma$. Choosing $n_{\mtx{\Omega}} = n_{\mtx{\Psi}}$ constitutes a good compromise.

\begin{figure}[ht]
    \centering
    \scalebox{0.8}{ \input{plots/distribution.pgf}}
    \caption{Error vs. $\sigma$ for a fixed budget of $80$ random vectors 
    when applying  Chebyshev-Nyström++ to the Hamiltonian matrix from~\cref{subsec:hamiltonian} with $n_c = 1$.
The Chebyshev approximation error is made negligible by setting $m=16 / \sigma$ (see \cref{lem:non-negative-chebyshev-error}).}
    \label{fig:distribution}
\end{figure}

\subsection{Comparison to Krylov-aware trace estimator}
\label{subsec:krylov-aware}

We notice that the Krylov-aware stochastic trace estimator from \cite[Algorithm 3.1]{chen-2023-krylovaware-stochastic} can also effectively be used in the setting of spectral density approximation. It also samples two Gaussian random matrices $\mtx{\Omega} \in \mathbb{R}^{n \times n_{\mtx{\Omega}}}$ and $\mtx{\Psi} \in \mathbb{R}^{n \times n_{\mtx{\Psi}}}$. It first runs block Lanczos algorithm on $\mtx{A}$ with starting block $\mtx{\Omega}$ for $r$ iterations with reorthogonalization and subsequently $q$ without. The columns of $\mtx{\Psi}$ are then projected onto the complement of the block Krylov subspace generated by the Lanczos algorithm, and are afterwards used as starting vectors for $r$ iterations of the Lanczos algorithm on $\mtx{A}$. What is remarkable is that all these heavy computations are completely independent of the smoothing kernel $g_{\sigma}$ and the parameter $t$. Further, some algebraic manipulations show that similarly to the stochastic Lanczos quadrature \cite[Section 3]{ubaru-2017-fast-estimation}, we can compress the estimator into a quadrature. Since this quadrature is independent of $g_{\sigma}$ and $t$, it can be applied to different smoothing kernels $g_{\sigma}$ and to as many values of $t$ as needed --- at little additional cost.

We run our own, faithfully optimized implementation of the Krylov-aware stochastic trace estimator \cite[Algorithm 3.1]{chen-2023-krylovaware-stochastic} for different parameter settings on the example from \cref{subsec:hamiltonian} and plot the approximation errors for logarithmically spaced values of the smoothing parameter $\sigma$ in \cref{fig:krylov-aware-density}. For reference, we also report the error of Chebyshev-Nyström++ for parameters that lead to a comparable run-time. The non-monotonic behavior of the Chebyshev-Nyström++ method in this case arises because, for small values of $\sigma$, the Chebyshev interpolation requires a very high polynomial degree to be accurate, whereas for large $\sigma$, the low-rank approximation becomes less effective. Due to the fixed parameters, the sweet spot lies around $\sigma \approx 0.002$: the low-rank approximation benefits from some underlying low-rank structure in $g_{\sigma}(t \mtx{I} - \mtx{A})$, while the Chebyshev interpolation remains good because $g_{\sigma}$ is not too sharply peaked. While the Krylov-aware estimator is clearly an attractive alternative, we argue that the parameters of the Chebyshev-Nyström++ method can be chosen in a more interpretable way, due to the explicit requirement for $m$ based on $\sigma$ (\cref{fig:chebyshev-heatmap}) and for $n_{\mtx{\Omega}}$ based on the distribution of the eigenvalues (discussed towards the end of \cref{subsubsec:chebyshev-nystrom-analysis}).

\begin{figure}[ht]
    \begin{minipage}[c]{.515\linewidth}
        \centering
        \scalebox{0.8}{\input{plots/krylov_aware_density.pgf}}
    \end{minipage}\hfill%
    \begin{minipage}[c]{.475\linewidth}
        \vspace{-35pt}
        \scalebox{0.8}{\input{tables/krylov_aware_density_KA.tex}}
        \newline
        \vspace{15pt}
        \newline
        \scalebox{0.8}{\input{tables/krylov_aware_density_CN.tex}}
    \end{minipage}
    \caption{For the example from \cref{subsec:hamiltonian}, error vs. $\sigma$ the Krylov-Aware (KA) stochastic trace estimator with different parameter settings, compared to the Chebyshev-Nyström++ (CN++) method.}
    \label{fig:krylov-aware-density}
\end{figure}

\subsection{Spectral density of Hessian matrix of neural network}
\label{subsec:hessian}

When fitting the parameters of a neural network, it can be of interest to determine whether a (local) minimum has been found and whether this minimum is robust, i.e., whether a small perturbation of the parameters leads to a significant decrease in the fit of the neural network. Both properties are reflected in the eigenvalues of the Hessian matrix $\mtx{A}$ with respect to some loss function: if all eigenvalues are non-negative, then the loss attains a local minimum, and if, additionally, all of them are small, we talk about a flat minimum; one where small perturbations of the parameters do not cause large increases of the loss and are therefore associated with better generalization \cite{yao-2020-pyhessian-neural}. On the other hand, large eigenvalues in the spectrum indicate high sensitivity of the loss to small changes in the parameters. Those we want to avoid, which is the reason why we monitor them.

Since neural networks are usually parametrized by a large number of parameters, assembling the Hessian matrix explicitly and then computing its eigenvalues is clearly infeasible. On the other hand, one can cheaply compute matrix-vector products $\mtx{A} \vct{x}$ with the Hessian, at a cost that scales proportionally with the number of parameters \cite{pearlmutter-1994-fast-exact}.

To demonstrate that our algorithm can effectively be applied in the setting of neural network optimization, we approximate the spectral density of a small fully connected convolutional neural network with $6\,782$ parameters. We train this network on the handwritten digit classification task given by the MNIST dataset\footnote{Handwritten digit classification; taken from \url{http://yann.lecun.com/exdb/mnist}.} in PyTorch 2.4.1 in a standard way. We use the bound from \cite[Theorem 1]{zhou-2011-bounding-spectrum} to estimate bounds on the spectrum.

We plot the spectral density of the Hessian matrix of the untrained neural network, and at different stages of training in \cref{fig:hessian-density}, as well as the corresponding mean squared error loss. It can be observed that the eigenvalues creep towards zero as the training proceeds. Furthermore, despite the loss steadily decreasing, the presence of relatively large eigenvalues in some epochs indicates a sharp minimum of the network, hence, unfavorable generalization properties.

\begin{figure}
    \begin{minipage}[c]{.49\linewidth}
        \centering
        \scalebox{0.8}{\input{plots/hessian_density.pgf}}
    \end{minipage}\hfill%
    \begin{minipage}[c]{.49\linewidth}
        \centering
        \scalebox{0.8}{\input{plots/hessian_density_loss.pgf}}
    \end{minipage}
    \caption{Mean squared error training loss (right) and corresponding approximate spectral density of the Hessian matrix of a fully connected convolutional neural network in different epochs of training on the MNIST dataset (left). The spectral density is approximated in $n_t=150$ uniformly spaced points using Chebyshev-Nyström++ with parameters $n_{\mtx{\Omega}} = 10$, $n_{\mtx{\Psi}} = 10$, $m = 1000$, and $\sigma = 0.005$.}
    \label{fig:hessian-density}
\end{figure}

%% file: plots/convergence-1.pgf
\begingroup%
\makeatletter%
\begin{pgfpicture}%
\pgfpathrectangle{\pgfpointorigin}{\pgfqpoint{3.146913in}{2.959073in}}%
\pgfusepath{use as bounding box, clip}%
\begin{pgfscope}%
\pgfsetbuttcap%
\pgfsetmiterjoin%
\definecolor{currentfill}{rgb}{1.000000,1.000000,1.000000}%
\pgfsetfillcolor{currentfill}%
\pgfsetlinewidth{0.000000pt}%
\definecolor{currentstroke}{rgb}{1.000000,1.000000,1.000000}%
\pgfsetstrokecolor{currentstroke}%
\pgfsetdash{}{0pt}%
\pgfpathmoveto{\pgfqpoint{0.000000in}{-0.000000in}}%
\pgfpathlineto{\pgfqpoint{3.146913in}{-0.000000in}}%
\pgfpathlineto{\pgfqpoint{3.146913in}{2.959073in}}%
\pgfpathlineto{\pgfqpoint{0.000000in}{2.959073in}}%
\pgfpathlineto{\pgfqpoint{0.000000in}{-0.000000in}}%
\pgfpathclose%
\pgfusepath{fill}%
\end{pgfscope}%
\begin{pgfscope}%
\pgfsetbuttcap%
\pgfsetmiterjoin%
\definecolor{currentfill}{rgb}{1.000000,1.000000,1.000000}%
\pgfsetfillcolor{currentfill}%
\pgfsetlinewidth{0.000000pt}%
\definecolor{currentstroke}{rgb}{0.000000,0.000000,0.000000}%
\pgfsetstrokecolor{currentstroke}%
\pgfsetstrokeopacity{0.000000}%
\pgfsetdash{}{0pt}%
\pgfpathmoveto{\pgfqpoint{0.721913in}{0.549073in}}%
\pgfpathlineto{\pgfqpoint{3.046913in}{0.549073in}}%
\pgfpathlineto{\pgfqpoint{3.046913in}{2.859073in}}%
\pgfpathlineto{\pgfqpoint{0.721913in}{2.859073in}}%
\pgfpathlineto{\pgfqpoint{0.721913in}{0.549073in}}%
\pgfpathclose%
\pgfusepath{fill}%
\end{pgfscope}%
\begin{pgfscope}%
\pgfpathrectangle{\pgfqpoint{0.721913in}{0.549073in}}{\pgfqpoint{2.325000in}{2.310000in}}%
\pgfusepath{clip}%
\pgfsetrectcap%
\pgfsetroundjoin%
\pgfsetlinewidth{0.250937pt}%
\definecolor{currentstroke}{rgb}{0.000000,0.000000,0.000000}%
\pgfsetstrokecolor{currentstroke}%
\pgfsetstrokeopacity{0.200000}%
\pgfsetdash{}{0pt}%
\pgfpathmoveto{\pgfqpoint{1.020314in}{0.549073in}}%
\pgfpathlineto{\pgfqpoint{1.020314in}{2.859073in}}%
\pgfusepath{stroke}%
\end{pgfscope}%
\begin{pgfscope}%
\pgfsetbuttcap%
\pgfsetroundjoin%
\definecolor{currentfill}{rgb}{0.000000,0.000000,0.000000}%
\pgfsetfillcolor{currentfill}%
\pgfsetlinewidth{0.803000pt}%
\definecolor{currentstroke}{rgb}{0.000000,0.000000,0.000000}%
\pgfsetstrokecolor{currentstroke}%
\pgfsetdash{}{0pt}%
\pgfsys@defobject{currentmarker}{\pgfqpoint{0.000000in}{-0.048611in}}{\pgfqpoint{0.000000in}{0.000000in}}{%
\pgfpathmoveto{\pgfqpoint{0.000000in}{0.000000in}}%
\pgfpathlineto{\pgfqpoint{0.000000in}{-0.048611in}}%
\pgfusepath{stroke,fill}%
}%
\begin{pgfscope}%
\pgfsys@transformshift{1.020314in}{0.549073in}%
\pgfsys@useobject{currentmarker}{}%
\end{pgfscope}%
\end{pgfscope}%
\begin{pgfscope}%
\definecolor{textcolor}{rgb}{0.000000,0.000000,0.000000}%
\pgfsetstrokecolor{textcolor}%
\pgfsetfillcolor{textcolor}%
\pgftext[x=1.020314in,y=0.451851in,,top]{\color{textcolor}{\rmfamily\fontsize{12.000000}{14.400000}\selectfont\catcode`\^=\active\def^{\ifmmode\sp\else\^{}\fi}\catcode`\%=\active\def
\end{pgfscope}%
\begin{pgfscope}%
\pgfpathrectangle{\pgfqpoint{0.721913in}{0.549073in}}{\pgfqpoint{2.325000in}{2.310000in}}%
\pgfusepath{clip}%
\pgfsetrectcap%
\pgfsetroundjoin%
\pgfsetlinewidth{0.250937pt}%
\definecolor{currentstroke}{rgb}{0.000000,0.000000,0.000000}%
\pgfsetstrokecolor{currentstroke}%
\pgfsetstrokeopacity{0.200000}%
\pgfsetdash{}{0pt}%
\pgfpathmoveto{\pgfqpoint{2.115850in}{0.549073in}}%
\pgfpathlineto{\pgfqpoint{2.115850in}{2.859073in}}%
\pgfusepath{stroke}%
\end{pgfscope}%
\begin{pgfscope}%
\pgfsetbuttcap%
\pgfsetroundjoin%
\definecolor{currentfill}{rgb}{0.000000,0.000000,0.000000}%
\pgfsetfillcolor{currentfill}%
\pgfsetlinewidth{0.803000pt}%
\definecolor{currentstroke}{rgb}{0.000000,0.000000,0.000000}%
\pgfsetstrokecolor{currentstroke}%
\pgfsetdash{}{0pt}%
\pgfsys@defobject{currentmarker}{\pgfqpoint{0.000000in}{-0.048611in}}{\pgfqpoint{0.000000in}{0.000000in}}{%
\pgfpathmoveto{\pgfqpoint{0.000000in}{0.000000in}}%
\pgfpathlineto{\pgfqpoint{0.000000in}{-0.048611in}}%
\pgfusepath{stroke,fill}%
}%
\begin{pgfscope}%
\pgfsys@transformshift{2.115850in}{0.549073in}%
\pgfsys@useobject{currentmarker}{}%
\end{pgfscope}%
\end{pgfscope}%
\begin{pgfscope}%
\definecolor{textcolor}{rgb}{0.000000,0.000000,0.000000}%
\pgfsetstrokecolor{textcolor}%
\pgfsetfillcolor{textcolor}%
\pgftext[x=2.115850in,y=0.451851in,,top]{\color{textcolor}{\rmfamily\fontsize{12.000000}{14.400000}\selectfont\catcode`\^=\active\def^{\ifmmode\sp\else\^{}\fi}\catcode`\%=\active\def
\end{pgfscope}%
\begin{pgfscope}%
\pgfsetbuttcap%
\pgfsetroundjoin%
\definecolor{currentfill}{rgb}{0.000000,0.000000,0.000000}%
\pgfsetfillcolor{currentfill}%
\pgfsetlinewidth{0.602250pt}%
\definecolor{currentstroke}{rgb}{0.000000,0.000000,0.000000}%
\pgfsetstrokecolor{currentstroke}%
\pgfsetdash{}{0pt}%
\pgfsys@defobject{currentmarker}{\pgfqpoint{0.000000in}{-0.027778in}}{\pgfqpoint{0.000000in}{0.000000in}}{%
\pgfpathmoveto{\pgfqpoint{0.000000in}{0.000000in}}%
\pgfpathlineto{\pgfqpoint{0.000000in}{-0.027778in}}%
\pgfusepath{stroke,fill}%
}%
\begin{pgfscope}%
\pgfsys@transformshift{0.777271in}{0.549073in}%
\pgfsys@useobject{currentmarker}{}%
\end{pgfscope}%
\end{pgfscope}%
\begin{pgfscope}%
\pgfsetbuttcap%
\pgfsetroundjoin%
\definecolor{currentfill}{rgb}{0.000000,0.000000,0.000000}%
\pgfsetfillcolor{currentfill}%
\pgfsetlinewidth{0.602250pt}%
\definecolor{currentstroke}{rgb}{0.000000,0.000000,0.000000}%
\pgfsetstrokecolor{currentstroke}%
\pgfsetdash{}{0pt}%
\pgfsys@defobject{currentmarker}{\pgfqpoint{0.000000in}{-0.027778in}}{\pgfqpoint{0.000000in}{0.000000in}}{%
\pgfpathmoveto{\pgfqpoint{0.000000in}{0.000000in}}%
\pgfpathlineto{\pgfqpoint{0.000000in}{-0.027778in}}%
\pgfusepath{stroke,fill}%
}%
\begin{pgfscope}%
\pgfsys@transformshift{0.850613in}{0.549073in}%
\pgfsys@useobject{currentmarker}{}%
\end{pgfscope}%
\end{pgfscope}%
\begin{pgfscope}%
\pgfsetbuttcap%
\pgfsetroundjoin%
\definecolor{currentfill}{rgb}{0.000000,0.000000,0.000000}%
\pgfsetfillcolor{currentfill}%
\pgfsetlinewidth{0.602250pt}%
\definecolor{currentstroke}{rgb}{0.000000,0.000000,0.000000}%
\pgfsetstrokecolor{currentstroke}%
\pgfsetdash{}{0pt}%
\pgfsys@defobject{currentmarker}{\pgfqpoint{0.000000in}{-0.027778in}}{\pgfqpoint{0.000000in}{0.000000in}}{%
\pgfpathmoveto{\pgfqpoint{0.000000in}{0.000000in}}%
\pgfpathlineto{\pgfqpoint{0.000000in}{-0.027778in}}%
\pgfusepath{stroke,fill}%
}%
\begin{pgfscope}%
\pgfsys@transformshift{0.914146in}{0.549073in}%
\pgfsys@useobject{currentmarker}{}%
\end{pgfscope}%
\end{pgfscope}%
\begin{pgfscope}%
\pgfsetbuttcap%
\pgfsetroundjoin%
\definecolor{currentfill}{rgb}{0.000000,0.000000,0.000000}%
\pgfsetfillcolor{currentfill}%
\pgfsetlinewidth{0.602250pt}%
\definecolor{currentstroke}{rgb}{0.000000,0.000000,0.000000}%
\pgfsetstrokecolor{currentstroke}%
\pgfsetdash{}{0pt}%
\pgfsys@defobject{currentmarker}{\pgfqpoint{0.000000in}{-0.027778in}}{\pgfqpoint{0.000000in}{0.000000in}}{%
\pgfpathmoveto{\pgfqpoint{0.000000in}{0.000000in}}%
\pgfpathlineto{\pgfqpoint{0.000000in}{-0.027778in}}%
\pgfusepath{stroke,fill}%
}%
\begin{pgfscope}%
\pgfsys@transformshift{0.970185in}{0.549073in}%
\pgfsys@useobject{currentmarker}{}%
\end{pgfscope}%
\end{pgfscope}%
\begin{pgfscope}%
\pgfsetbuttcap%
\pgfsetroundjoin%
\definecolor{currentfill}{rgb}{0.000000,0.000000,0.000000}%
\pgfsetfillcolor{currentfill}%
\pgfsetlinewidth{0.602250pt}%
\definecolor{currentstroke}{rgb}{0.000000,0.000000,0.000000}%
\pgfsetstrokecolor{currentstroke}%
\pgfsetdash{}{0pt}%
\pgfsys@defobject{currentmarker}{\pgfqpoint{0.000000in}{-0.027778in}}{\pgfqpoint{0.000000in}{0.000000in}}{%
\pgfpathmoveto{\pgfqpoint{0.000000in}{0.000000in}}%
\pgfpathlineto{\pgfqpoint{0.000000in}{-0.027778in}}%
\pgfusepath{stroke,fill}%
}%
\begin{pgfscope}%
\pgfsys@transformshift{1.350103in}{0.549073in}%
\pgfsys@useobject{currentmarker}{}%
\end{pgfscope}%
\end{pgfscope}%
\begin{pgfscope}%
\pgfsetbuttcap%
\pgfsetroundjoin%
\definecolor{currentfill}{rgb}{0.000000,0.000000,0.000000}%
\pgfsetfillcolor{currentfill}%
\pgfsetlinewidth{0.602250pt}%
\definecolor{currentstroke}{rgb}{0.000000,0.000000,0.000000}%
\pgfsetstrokecolor{currentstroke}%
\pgfsetdash{}{0pt}%
\pgfsys@defobject{currentmarker}{\pgfqpoint{0.000000in}{-0.027778in}}{\pgfqpoint{0.000000in}{0.000000in}}{%
\pgfpathmoveto{\pgfqpoint{0.000000in}{0.000000in}}%
\pgfpathlineto{\pgfqpoint{0.000000in}{-0.027778in}}%
\pgfusepath{stroke,fill}%
}%
\begin{pgfscope}%
\pgfsys@transformshift{1.543017in}{0.549073in}%
\pgfsys@useobject{currentmarker}{}%
\end{pgfscope}%
\end{pgfscope}%
\begin{pgfscope}%
\pgfsetbuttcap%
\pgfsetroundjoin%
\definecolor{currentfill}{rgb}{0.000000,0.000000,0.000000}%
\pgfsetfillcolor{currentfill}%
\pgfsetlinewidth{0.602250pt}%
\definecolor{currentstroke}{rgb}{0.000000,0.000000,0.000000}%
\pgfsetstrokecolor{currentstroke}%
\pgfsetdash{}{0pt}%
\pgfsys@defobject{currentmarker}{\pgfqpoint{0.000000in}{-0.027778in}}{\pgfqpoint{0.000000in}{0.000000in}}{%
\pgfpathmoveto{\pgfqpoint{0.000000in}{0.000000in}}%
\pgfpathlineto{\pgfqpoint{0.000000in}{-0.027778in}}%
\pgfusepath{stroke,fill}%
}%
\begin{pgfscope}%
\pgfsys@transformshift{1.679892in}{0.549073in}%
\pgfsys@useobject{currentmarker}{}%
\end{pgfscope}%
\end{pgfscope}%
\begin{pgfscope}%
\pgfsetbuttcap%
\pgfsetroundjoin%
\definecolor{currentfill}{rgb}{0.000000,0.000000,0.000000}%
\pgfsetfillcolor{currentfill}%
\pgfsetlinewidth{0.602250pt}%
\definecolor{currentstroke}{rgb}{0.000000,0.000000,0.000000}%
\pgfsetstrokecolor{currentstroke}%
\pgfsetdash{}{0pt}%
\pgfsys@defobject{currentmarker}{\pgfqpoint{0.000000in}{-0.027778in}}{\pgfqpoint{0.000000in}{0.000000in}}{%
\pgfpathmoveto{\pgfqpoint{0.000000in}{0.000000in}}%
\pgfpathlineto{\pgfqpoint{0.000000in}{-0.027778in}}%
\pgfusepath{stroke,fill}%
}%
\begin{pgfscope}%
\pgfsys@transformshift{1.786061in}{0.549073in}%
\pgfsys@useobject{currentmarker}{}%
\end{pgfscope}%
\end{pgfscope}%
\begin{pgfscope}%
\pgfsetbuttcap%
\pgfsetroundjoin%
\definecolor{currentfill}{rgb}{0.000000,0.000000,0.000000}%
\pgfsetfillcolor{currentfill}%
\pgfsetlinewidth{0.602250pt}%
\definecolor{currentstroke}{rgb}{0.000000,0.000000,0.000000}%
\pgfsetstrokecolor{currentstroke}%
\pgfsetdash{}{0pt}%
\pgfsys@defobject{currentmarker}{\pgfqpoint{0.000000in}{-0.027778in}}{\pgfqpoint{0.000000in}{0.000000in}}{%
\pgfpathmoveto{\pgfqpoint{0.000000in}{0.000000in}}%
\pgfpathlineto{\pgfqpoint{0.000000in}{-0.027778in}}%
\pgfusepath{stroke,fill}%
}%
\begin{pgfscope}%
\pgfsys@transformshift{1.872807in}{0.549073in}%
\pgfsys@useobject{currentmarker}{}%
\end{pgfscope}%
\end{pgfscope}%
\begin{pgfscope}%
\pgfsetbuttcap%
\pgfsetroundjoin%
\definecolor{currentfill}{rgb}{0.000000,0.000000,0.000000}%
\pgfsetfillcolor{currentfill}%
\pgfsetlinewidth{0.602250pt}%
\definecolor{currentstroke}{rgb}{0.000000,0.000000,0.000000}%
\pgfsetstrokecolor{currentstroke}%
\pgfsetdash{}{0pt}%
\pgfsys@defobject{currentmarker}{\pgfqpoint{0.000000in}{-0.027778in}}{\pgfqpoint{0.000000in}{0.000000in}}{%
\pgfpathmoveto{\pgfqpoint{0.000000in}{0.000000in}}%
\pgfpathlineto{\pgfqpoint{0.000000in}{-0.027778in}}%
\pgfusepath{stroke,fill}%
}%
\begin{pgfscope}%
\pgfsys@transformshift{1.946149in}{0.549073in}%
\pgfsys@useobject{currentmarker}{}%
\end{pgfscope}%
\end{pgfscope}%
\begin{pgfscope}%
\pgfsetbuttcap%
\pgfsetroundjoin%
\definecolor{currentfill}{rgb}{0.000000,0.000000,0.000000}%
\pgfsetfillcolor{currentfill}%
\pgfsetlinewidth{0.602250pt}%
\definecolor{currentstroke}{rgb}{0.000000,0.000000,0.000000}%
\pgfsetstrokecolor{currentstroke}%
\pgfsetdash{}{0pt}%
\pgfsys@defobject{currentmarker}{\pgfqpoint{0.000000in}{-0.027778in}}{\pgfqpoint{0.000000in}{0.000000in}}{%
\pgfpathmoveto{\pgfqpoint{0.000000in}{0.000000in}}%
\pgfpathlineto{\pgfqpoint{0.000000in}{-0.027778in}}%
\pgfusepath{stroke,fill}%
}%
\begin{pgfscope}%
\pgfsys@transformshift{2.009682in}{0.549073in}%
\pgfsys@useobject{currentmarker}{}%
\end{pgfscope}%
\end{pgfscope}%
\begin{pgfscope}%
\pgfsetbuttcap%
\pgfsetroundjoin%
\definecolor{currentfill}{rgb}{0.000000,0.000000,0.000000}%
\pgfsetfillcolor{currentfill}%
\pgfsetlinewidth{0.602250pt}%
\definecolor{currentstroke}{rgb}{0.000000,0.000000,0.000000}%
\pgfsetstrokecolor{currentstroke}%
\pgfsetdash{}{0pt}%
\pgfsys@defobject{currentmarker}{\pgfqpoint{0.000000in}{-0.027778in}}{\pgfqpoint{0.000000in}{0.000000in}}{%
\pgfpathmoveto{\pgfqpoint{0.000000in}{0.000000in}}%
\pgfpathlineto{\pgfqpoint{0.000000in}{-0.027778in}}%
\pgfusepath{stroke,fill}%
}%
\begin{pgfscope}%
\pgfsys@transformshift{2.065721in}{0.549073in}%
\pgfsys@useobject{currentmarker}{}%
\end{pgfscope}%
\end{pgfscope}%
\begin{pgfscope}%
\pgfsetbuttcap%
\pgfsetroundjoin%
\definecolor{currentfill}{rgb}{0.000000,0.000000,0.000000}%
\pgfsetfillcolor{currentfill}%
\pgfsetlinewidth{0.602250pt}%
\definecolor{currentstroke}{rgb}{0.000000,0.000000,0.000000}%
\pgfsetstrokecolor{currentstroke}%
\pgfsetdash{}{0pt}%
\pgfsys@defobject{currentmarker}{\pgfqpoint{0.000000in}{-0.027778in}}{\pgfqpoint{0.000000in}{0.000000in}}{%
\pgfpathmoveto{\pgfqpoint{0.000000in}{0.000000in}}%
\pgfpathlineto{\pgfqpoint{0.000000in}{-0.027778in}}%
\pgfusepath{stroke,fill}%
}%
\begin{pgfscope}%
\pgfsys@transformshift{2.445639in}{0.549073in}%
\pgfsys@useobject{currentmarker}{}%
\end{pgfscope}%
\end{pgfscope}%
\begin{pgfscope}%
\pgfsetbuttcap%
\pgfsetroundjoin%
\definecolor{currentfill}{rgb}{0.000000,0.000000,0.000000}%
\pgfsetfillcolor{currentfill}%
\pgfsetlinewidth{0.602250pt}%
\definecolor{currentstroke}{rgb}{0.000000,0.000000,0.000000}%
\pgfsetstrokecolor{currentstroke}%
\pgfsetdash{}{0pt}%
\pgfsys@defobject{currentmarker}{\pgfqpoint{0.000000in}{-0.027778in}}{\pgfqpoint{0.000000in}{0.000000in}}{%
\pgfpathmoveto{\pgfqpoint{0.000000in}{0.000000in}}%
\pgfpathlineto{\pgfqpoint{0.000000in}{-0.027778in}}%
\pgfusepath{stroke,fill}%
}%
\begin{pgfscope}%
\pgfsys@transformshift{2.638553in}{0.549073in}%
\pgfsys@useobject{currentmarker}{}%
\end{pgfscope}%
\end{pgfscope}%
\begin{pgfscope}%
\pgfsetbuttcap%
\pgfsetroundjoin%
\definecolor{currentfill}{rgb}{0.000000,0.000000,0.000000}%
\pgfsetfillcolor{currentfill}%
\pgfsetlinewidth{0.602250pt}%
\definecolor{currentstroke}{rgb}{0.000000,0.000000,0.000000}%
\pgfsetstrokecolor{currentstroke}%
\pgfsetdash{}{0pt}%
\pgfsys@defobject{currentmarker}{\pgfqpoint{0.000000in}{-0.027778in}}{\pgfqpoint{0.000000in}{0.000000in}}{%
\pgfpathmoveto{\pgfqpoint{0.000000in}{0.000000in}}%
\pgfpathlineto{\pgfqpoint{0.000000in}{-0.027778in}}%
\pgfusepath{stroke,fill}%
}%
\begin{pgfscope}%
\pgfsys@transformshift{2.775428in}{0.549073in}%
\pgfsys@useobject{currentmarker}{}%
\end{pgfscope}%
\end{pgfscope}%
\begin{pgfscope}%
\pgfsetbuttcap%
\pgfsetroundjoin%
\definecolor{currentfill}{rgb}{0.000000,0.000000,0.000000}%
\pgfsetfillcolor{currentfill}%
\pgfsetlinewidth{0.602250pt}%
\definecolor{currentstroke}{rgb}{0.000000,0.000000,0.000000}%
\pgfsetstrokecolor{currentstroke}%
\pgfsetdash{}{0pt}%
\pgfsys@defobject{currentmarker}{\pgfqpoint{0.000000in}{-0.027778in}}{\pgfqpoint{0.000000in}{0.000000in}}{%
\pgfpathmoveto{\pgfqpoint{0.000000in}{0.000000in}}%
\pgfpathlineto{\pgfqpoint{0.000000in}{-0.027778in}}%
\pgfusepath{stroke,fill}%
}%
\begin{pgfscope}%
\pgfsys@transformshift{2.881597in}{0.549073in}%
\pgfsys@useobject{currentmarker}{}%
\end{pgfscope}%
\end{pgfscope}%
\begin{pgfscope}%
\pgfsetbuttcap%
\pgfsetroundjoin%
\definecolor{currentfill}{rgb}{0.000000,0.000000,0.000000}%
\pgfsetfillcolor{currentfill}%
\pgfsetlinewidth{0.602250pt}%
\definecolor{currentstroke}{rgb}{0.000000,0.000000,0.000000}%
\pgfsetstrokecolor{currentstroke}%
\pgfsetdash{}{0pt}%
\pgfsys@defobject{currentmarker}{\pgfqpoint{0.000000in}{-0.027778in}}{\pgfqpoint{0.000000in}{0.000000in}}{%
\pgfpathmoveto{\pgfqpoint{0.000000in}{0.000000in}}%
\pgfpathlineto{\pgfqpoint{0.000000in}{-0.027778in}}%
\pgfusepath{stroke,fill}%
}%
\begin{pgfscope}%
\pgfsys@transformshift{2.968343in}{0.549073in}%
\pgfsys@useobject{currentmarker}{}%
\end{pgfscope}%
\end{pgfscope}%
\begin{pgfscope}%
\pgfsetbuttcap%
\pgfsetroundjoin%
\definecolor{currentfill}{rgb}{0.000000,0.000000,0.000000}%
\pgfsetfillcolor{currentfill}%
\pgfsetlinewidth{0.602250pt}%
\definecolor{currentstroke}{rgb}{0.000000,0.000000,0.000000}%
\pgfsetstrokecolor{currentstroke}%
\pgfsetdash{}{0pt}%
\pgfsys@defobject{currentmarker}{\pgfqpoint{0.000000in}{-0.027778in}}{\pgfqpoint{0.000000in}{0.000000in}}{%
\pgfpathmoveto{\pgfqpoint{0.000000in}{0.000000in}}%
\pgfpathlineto{\pgfqpoint{0.000000in}{-0.027778in}}%
\pgfusepath{stroke,fill}%
}%
\begin{pgfscope}%
\pgfsys@transformshift{3.041685in}{0.549073in}%
\pgfsys@useobject{currentmarker}{}%
\end{pgfscope}%
\end{pgfscope}%
\begin{pgfscope}%
\definecolor{textcolor}{rgb}{0.000000,0.000000,0.000000}%
\pgfsetstrokecolor{textcolor}%
\pgfsetfillcolor{textcolor}%
\pgftext[x=1.884413in,y=0.248148in,,top]{\color{textcolor}{\rmfamily\fontsize{12.000000}{14.400000}\selectfont\catcode`\^=\active\def^{\ifmmode\sp\else\^{}\fi}\catcode`\%=\active\def
\end{pgfscope}%
\begin{pgfscope}%
\pgfpathrectangle{\pgfqpoint{0.721913in}{0.549073in}}{\pgfqpoint{2.325000in}{2.310000in}}%
\pgfusepath{clip}%
\pgfsetrectcap%
\pgfsetroundjoin%
\pgfsetlinewidth{0.250937pt}%
\definecolor{currentstroke}{rgb}{0.000000,0.000000,0.000000}%
\pgfsetstrokecolor{currentstroke}%
\pgfsetstrokeopacity{0.200000}%
\pgfsetdash{}{0pt}%
\pgfpathmoveto{\pgfqpoint{0.721913in}{0.867122in}}%
\pgfpathlineto{\pgfqpoint{3.046913in}{0.867122in}}%
\pgfusepath{stroke}%
\end{pgfscope}%
\begin{pgfscope}%
\pgfsetbuttcap%
\pgfsetroundjoin%
\definecolor{currentfill}{rgb}{0.000000,0.000000,0.000000}%
\pgfsetfillcolor{currentfill}%
\pgfsetlinewidth{0.803000pt}%
\definecolor{currentstroke}{rgb}{0.000000,0.000000,0.000000}%
\pgfsetstrokecolor{currentstroke}%
\pgfsetdash{}{0pt}%
\pgfsys@defobject{currentmarker}{\pgfqpoint{-0.048611in}{0.000000in}}{\pgfqpoint{-0.000000in}{0.000000in}}{%
\pgfpathmoveto{\pgfqpoint{-0.000000in}{0.000000in}}%
\pgfpathlineto{\pgfqpoint{-0.048611in}{0.000000in}}%
\pgfusepath{stroke,fill}%
}%
\begin{pgfscope}%
\pgfsys@transformshift{0.721913in}{0.867122in}%
\pgfsys@useobject{currentmarker}{}%
\end{pgfscope}%
\end{pgfscope}%
\begin{pgfscope}%
\definecolor{textcolor}{rgb}{0.000000,0.000000,0.000000}%
\pgfsetstrokecolor{textcolor}%
\pgfsetfillcolor{textcolor}%
\pgftext[x=0.303703in, y=0.809252in, left, base]{\color{textcolor}{\rmfamily\fontsize{12.000000}{14.400000}\selectfont\catcode`\^=\active\def^{\ifmmode\sp\else\^{}\fi}\catcode`\%=\active\def
\end{pgfscope}%
\begin{pgfscope}%
\pgfpathrectangle{\pgfqpoint{0.721913in}{0.549073in}}{\pgfqpoint{2.325000in}{2.310000in}}%
\pgfusepath{clip}%
\pgfsetrectcap%
\pgfsetroundjoin%
\pgfsetlinewidth{0.250937pt}%
\definecolor{currentstroke}{rgb}{0.000000,0.000000,0.000000}%
\pgfsetstrokecolor{currentstroke}%
\pgfsetstrokeopacity{0.200000}%
\pgfsetdash{}{0pt}%
\pgfpathmoveto{\pgfqpoint{0.721913in}{1.483027in}}%
\pgfpathlineto{\pgfqpoint{3.046913in}{1.483027in}}%
\pgfusepath{stroke}%
\end{pgfscope}%
\begin{pgfscope}%
\pgfsetbuttcap%
\pgfsetroundjoin%
\definecolor{currentfill}{rgb}{0.000000,0.000000,0.000000}%
\pgfsetfillcolor{currentfill}%
\pgfsetlinewidth{0.803000pt}%
\definecolor{currentstroke}{rgb}{0.000000,0.000000,0.000000}%
\pgfsetstrokecolor{currentstroke}%
\pgfsetdash{}{0pt}%
\pgfsys@defobject{currentmarker}{\pgfqpoint{-0.048611in}{0.000000in}}{\pgfqpoint{-0.000000in}{0.000000in}}{%
\pgfpathmoveto{\pgfqpoint{-0.000000in}{0.000000in}}%
\pgfpathlineto{\pgfqpoint{-0.048611in}{0.000000in}}%
\pgfusepath{stroke,fill}%
}%
\begin{pgfscope}%
\pgfsys@transformshift{0.721913in}{1.483027in}%
\pgfsys@useobject{currentmarker}{}%
\end{pgfscope}%
\end{pgfscope}%
\begin{pgfscope}%
\definecolor{textcolor}{rgb}{0.000000,0.000000,0.000000}%
\pgfsetstrokecolor{textcolor}%
\pgfsetfillcolor{textcolor}%
\pgftext[x=0.303703in, y=1.425157in, left, base]{\color{textcolor}{\rmfamily\fontsize{12.000000}{14.400000}\selectfont\catcode`\^=\active\def^{\ifmmode\sp\else\^{}\fi}\catcode`\%=\active\def
\end{pgfscope}%
\begin{pgfscope}%
\pgfpathrectangle{\pgfqpoint{0.721913in}{0.549073in}}{\pgfqpoint{2.325000in}{2.310000in}}%
\pgfusepath{clip}%
\pgfsetrectcap%
\pgfsetroundjoin%
\pgfsetlinewidth{0.250937pt}%
\definecolor{currentstroke}{rgb}{0.000000,0.000000,0.000000}%
\pgfsetstrokecolor{currentstroke}%
\pgfsetstrokeopacity{0.200000}%
\pgfsetdash{}{0pt}%
\pgfpathmoveto{\pgfqpoint{0.721913in}{2.098932in}}%
\pgfpathlineto{\pgfqpoint{3.046913in}{2.098932in}}%
\pgfusepath{stroke}%
\end{pgfscope}%
\begin{pgfscope}%
\pgfsetbuttcap%
\pgfsetroundjoin%
\definecolor{currentfill}{rgb}{0.000000,0.000000,0.000000}%
\pgfsetfillcolor{currentfill}%
\pgfsetlinewidth{0.803000pt}%
\definecolor{currentstroke}{rgb}{0.000000,0.000000,0.000000}%
\pgfsetstrokecolor{currentstroke}%
\pgfsetdash{}{0pt}%
\pgfsys@defobject{currentmarker}{\pgfqpoint{-0.048611in}{0.000000in}}{\pgfqpoint{-0.000000in}{0.000000in}}{%
\pgfpathmoveto{\pgfqpoint{-0.000000in}{0.000000in}}%
\pgfpathlineto{\pgfqpoint{-0.048611in}{0.000000in}}%
\pgfusepath{stroke,fill}%
}%
\begin{pgfscope}%
\pgfsys@transformshift{0.721913in}{2.098932in}%
\pgfsys@useobject{currentmarker}{}%
\end{pgfscope}%
\end{pgfscope}%
\begin{pgfscope}%
\definecolor{textcolor}{rgb}{0.000000,0.000000,0.000000}%
\pgfsetstrokecolor{textcolor}%
\pgfsetfillcolor{textcolor}%
\pgftext[x=0.303703in, y=2.041062in, left, base]{\color{textcolor}{\rmfamily\fontsize{12.000000}{14.400000}\selectfont\catcode`\^=\active\def^{\ifmmode\sp\else\^{}\fi}\catcode`\%=\active\def
\end{pgfscope}%
\begin{pgfscope}%
\pgfpathrectangle{\pgfqpoint{0.721913in}{0.549073in}}{\pgfqpoint{2.325000in}{2.310000in}}%
\pgfusepath{clip}%
\pgfsetrectcap%
\pgfsetroundjoin%
\pgfsetlinewidth{0.250937pt}%
\definecolor{currentstroke}{rgb}{0.000000,0.000000,0.000000}%
\pgfsetstrokecolor{currentstroke}%
\pgfsetstrokeopacity{0.200000}%
\pgfsetdash{}{0pt}%
\pgfpathmoveto{\pgfqpoint{0.721913in}{2.714837in}}%
\pgfpathlineto{\pgfqpoint{3.046913in}{2.714837in}}%
\pgfusepath{stroke}%
\end{pgfscope}%
\begin{pgfscope}%
\pgfsetbuttcap%
\pgfsetroundjoin%
\definecolor{currentfill}{rgb}{0.000000,0.000000,0.000000}%
\pgfsetfillcolor{currentfill}%
\pgfsetlinewidth{0.803000pt}%
\definecolor{currentstroke}{rgb}{0.000000,0.000000,0.000000}%
\pgfsetstrokecolor{currentstroke}%
\pgfsetdash{}{0pt}%
\pgfsys@defobject{currentmarker}{\pgfqpoint{-0.048611in}{0.000000in}}{\pgfqpoint{-0.000000in}{0.000000in}}{%
\pgfpathmoveto{\pgfqpoint{-0.000000in}{0.000000in}}%
\pgfpathlineto{\pgfqpoint{-0.048611in}{0.000000in}}%
\pgfusepath{stroke,fill}%
}%
\begin{pgfscope}%
\pgfsys@transformshift{0.721913in}{2.714837in}%
\pgfsys@useobject{currentmarker}{}%
\end{pgfscope}%
\end{pgfscope}%
\begin{pgfscope}%
\definecolor{textcolor}{rgb}{0.000000,0.000000,0.000000}%
\pgfsetstrokecolor{textcolor}%
\pgfsetfillcolor{textcolor}%
\pgftext[x=0.395525in, y=2.656967in, left, base]{\color{textcolor}{\rmfamily\fontsize{12.000000}{14.400000}\selectfont\catcode`\^=\active\def^{\ifmmode\sp\else\^{}\fi}\catcode`\%=\active\def
\end{pgfscope}%
\begin{pgfscope}%
\definecolor{textcolor}{rgb}{0.000000,0.000000,0.000000}%
\pgfsetstrokecolor{textcolor}%
\pgfsetfillcolor{textcolor}%
\pgftext[x=0.248147in,y=1.704073in,,bottom,rotate=90.000000]{\color{textcolor}{\rmfamily\fontsize{12.000000}{14.400000}\selectfont\catcode`\^=\active\def^{\ifmmode\sp\else\^{}\fi}\catcode`\%=\active\def
\end{pgfscope}%
\begin{pgfscope}%
\pgfpathrectangle{\pgfqpoint{0.721913in}{0.549073in}}{\pgfqpoint{2.325000in}{2.310000in}}%
\pgfusepath{clip}%
\pgfsetrectcap%
\pgfsetroundjoin%
\pgfsetlinewidth{1.505625pt}%
\definecolor{currentstroke}{rgb}{0.392157,0.560784,1.000000}%
\pgfsetstrokecolor{currentstroke}%
\pgfsetdash{}{0pt}%
\pgfpathmoveto{\pgfqpoint{0.777271in}{2.426349in}}%
\pgfpathlineto{\pgfqpoint{1.145143in}{2.356873in}}%
\pgfpathlineto{\pgfqpoint{1.526888in}{2.309200in}}%
\pgfpathlineto{\pgfqpoint{1.896020in}{2.264836in}}%
\pgfpathlineto{\pgfqpoint{2.258635in}{2.224582in}}%
\pgfpathlineto{\pgfqpoint{2.625694in}{2.135301in}}%
\pgfpathlineto{\pgfqpoint{2.991556in}{2.122197in}}%
\pgfusepath{stroke}%
\end{pgfscope}%
\begin{pgfscope}%
\pgfpathrectangle{\pgfqpoint{0.721913in}{0.549073in}}{\pgfqpoint{2.325000in}{2.310000in}}%
\pgfusepath{clip}%
\pgfsetbuttcap%
\pgfsetroundjoin%
\definecolor{currentfill}{rgb}{0.392157,0.560784,1.000000}%
\pgfsetfillcolor{currentfill}%
\pgfsetlinewidth{1.003750pt}%
\definecolor{currentstroke}{rgb}{0.392157,0.560784,1.000000}%
\pgfsetstrokecolor{currentstroke}%
\pgfsetdash{}{0pt}%
\pgfsys@defobject{currentmarker}{\pgfqpoint{-0.041667in}{-0.041667in}}{\pgfqpoint{0.041667in}{0.041667in}}{%
\pgfpathmoveto{\pgfqpoint{0.000000in}{-0.041667in}}%
\pgfpathcurveto{\pgfqpoint{0.011050in}{-0.041667in}}{\pgfqpoint{0.021649in}{-0.037276in}}{\pgfqpoint{0.029463in}{-0.029463in}}%
\pgfpathcurveto{\pgfqpoint{0.037276in}{-0.021649in}}{\pgfqpoint{0.041667in}{-0.011050in}}{\pgfqpoint{0.041667in}{0.000000in}}%
\pgfpathcurveto{\pgfqpoint{0.041667in}{0.011050in}}{\pgfqpoint{0.037276in}{0.021649in}}{\pgfqpoint{0.029463in}{0.029463in}}%
\pgfpathcurveto{\pgfqpoint{0.021649in}{0.037276in}}{\pgfqpoint{0.011050in}{0.041667in}}{\pgfqpoint{0.000000in}{0.041667in}}%
\pgfpathcurveto{\pgfqpoint{-0.011050in}{0.041667in}}{\pgfqpoint{-0.021649in}{0.037276in}}{\pgfqpoint{-0.029463in}{0.029463in}}%
\pgfpathcurveto{\pgfqpoint{-0.037276in}{0.021649in}}{\pgfqpoint{-0.041667in}{0.011050in}}{\pgfqpoint{-0.041667in}{0.000000in}}%
\pgfpathcurveto{\pgfqpoint{-0.041667in}{-0.011050in}}{\pgfqpoint{-0.037276in}{-0.021649in}}{\pgfqpoint{-0.029463in}{-0.029463in}}%
\pgfpathcurveto{\pgfqpoint{-0.021649in}{-0.037276in}}{\pgfqpoint{-0.011050in}{-0.041667in}}{\pgfqpoint{0.000000in}{-0.041667in}}%
\pgfpathlineto{\pgfqpoint{0.000000in}{-0.041667in}}%
\pgfpathclose%
\pgfusepath{stroke,fill}%
}%
\begin{pgfscope}%
\pgfsys@transformshift{0.777271in}{2.426349in}%
\pgfsys@useobject{currentmarker}{}%
\end{pgfscope}%
\begin{pgfscope}%
\pgfsys@transformshift{1.145143in}{2.356873in}%
\pgfsys@useobject{currentmarker}{}%
\end{pgfscope}%
\begin{pgfscope}%
\pgfsys@transformshift{1.526888in}{2.309200in}%
\pgfsys@useobject{currentmarker}{}%
\end{pgfscope}%
\begin{pgfscope}%
\pgfsys@transformshift{1.896020in}{2.264836in}%
\pgfsys@useobject{currentmarker}{}%
\end{pgfscope}%
\begin{pgfscope}%
\pgfsys@transformshift{2.258635in}{2.224582in}%
\pgfsys@useobject{currentmarker}{}%
\end{pgfscope}%
\begin{pgfscope}%
\pgfsys@transformshift{2.625694in}{2.135301in}%
\pgfsys@useobject{currentmarker}{}%
\end{pgfscope}%
\begin{pgfscope}%
\pgfsys@transformshift{2.991556in}{2.122197in}%
\pgfsys@useobject{currentmarker}{}%
\end{pgfscope}%
\end{pgfscope}%
\begin{pgfscope}%
\pgfpathrectangle{\pgfqpoint{0.721913in}{0.549073in}}{\pgfqpoint{2.325000in}{2.310000in}}%
\pgfusepath{clip}%
\pgfsetrectcap%
\pgfsetroundjoin%
\pgfsetlinewidth{1.505625pt}%
\definecolor{currentstroke}{rgb}{0.862745,0.149020,0.498039}%
\pgfsetstrokecolor{currentstroke}%
\pgfsetdash{}{0pt}%
\pgfpathmoveto{\pgfqpoint{0.777271in}{2.410049in}}%
\pgfpathlineto{\pgfqpoint{1.145143in}{2.327028in}}%
\pgfpathlineto{\pgfqpoint{1.526888in}{2.247549in}}%
\pgfpathlineto{\pgfqpoint{1.896020in}{2.068361in}}%
\pgfpathlineto{\pgfqpoint{2.258635in}{1.755726in}}%
\pgfpathlineto{\pgfqpoint{2.625694in}{0.772535in}}%
\pgfpathlineto{\pgfqpoint{2.991556in}{0.741573in}}%
\pgfusepath{stroke}%
\end{pgfscope}%
\begin{pgfscope}%
\pgfpathrectangle{\pgfqpoint{0.721913in}{0.549073in}}{\pgfqpoint{2.325000in}{2.310000in}}%
\pgfusepath{clip}%
\pgfsetbuttcap%
\pgfsetmiterjoin%
\definecolor{currentfill}{rgb}{0.862745,0.149020,0.498039}%
\pgfsetfillcolor{currentfill}%
\pgfsetlinewidth{1.003750pt}%
\definecolor{currentstroke}{rgb}{0.862745,0.149020,0.498039}%
\pgfsetstrokecolor{currentstroke}%
\pgfsetdash{}{0pt}%
\pgfsys@defobject{currentmarker}{\pgfqpoint{-0.041667in}{-0.041667in}}{\pgfqpoint{0.041667in}{0.041667in}}{%
\pgfpathmoveto{\pgfqpoint{-0.041667in}{-0.041667in}}%
\pgfpathlineto{\pgfqpoint{0.041667in}{-0.041667in}}%
\pgfpathlineto{\pgfqpoint{0.041667in}{0.041667in}}%
\pgfpathlineto{\pgfqpoint{-0.041667in}{0.041667in}}%
\pgfpathlineto{\pgfqpoint{-0.041667in}{-0.041667in}}%
\pgfpathclose%
\pgfusepath{stroke,fill}%
}%
\begin{pgfscope}%
\pgfsys@transformshift{0.777271in}{2.410049in}%
\pgfsys@useobject{currentmarker}{}%
\end{pgfscope}%
\begin{pgfscope}%
\pgfsys@transformshift{1.145143in}{2.327028in}%
\pgfsys@useobject{currentmarker}{}%
\end{pgfscope}%
\begin{pgfscope}%
\pgfsys@transformshift{1.526888in}{2.247549in}%
\pgfsys@useobject{currentmarker}{}%
\end{pgfscope}%
\begin{pgfscope}%
\pgfsys@transformshift{1.896020in}{2.068361in}%
\pgfsys@useobject{currentmarker}{}%
\end{pgfscope}%
\begin{pgfscope}%
\pgfsys@transformshift{2.258635in}{1.755726in}%
\pgfsys@useobject{currentmarker}{}%
\end{pgfscope}%
\begin{pgfscope}%
\pgfsys@transformshift{2.625694in}{0.772535in}%
\pgfsys@useobject{currentmarker}{}%
\end{pgfscope}%
\begin{pgfscope}%
\pgfsys@transformshift{2.991556in}{0.741573in}%
\pgfsys@useobject{currentmarker}{}%
\end{pgfscope}%
\end{pgfscope}%
\begin{pgfscope}%
\pgfpathrectangle{\pgfqpoint{0.721913in}{0.549073in}}{\pgfqpoint{2.325000in}{2.310000in}}%
\pgfusepath{clip}%
\pgfsetrectcap%
\pgfsetroundjoin%
\pgfsetlinewidth{1.505625pt}%
\definecolor{currentstroke}{rgb}{1.000000,0.690196,0.000000}%
\pgfsetstrokecolor{currentstroke}%
\pgfsetdash{}{0pt}%
\pgfpathmoveto{\pgfqpoint{0.777271in}{2.666573in}}%
\pgfpathlineto{\pgfqpoint{1.145143in}{2.607370in}}%
\pgfpathlineto{\pgfqpoint{1.526888in}{2.493234in}}%
\pgfpathlineto{\pgfqpoint{1.896020in}{2.240948in}}%
\pgfpathlineto{\pgfqpoint{2.258635in}{1.023687in}}%
\pgfpathlineto{\pgfqpoint{2.625694in}{0.900995in}}%
\pgfpathlineto{\pgfqpoint{2.991556in}{0.860471in}}%
\pgfusepath{stroke}%
\end{pgfscope}%
\begin{pgfscope}%
\pgfpathrectangle{\pgfqpoint{0.721913in}{0.549073in}}{\pgfqpoint{2.325000in}{2.310000in}}%
\pgfusepath{clip}%
\pgfsetbuttcap%
\pgfsetmiterjoin%
\definecolor{currentfill}{rgb}{1.000000,0.690196,0.000000}%
\pgfsetfillcolor{currentfill}%
\pgfsetlinewidth{1.003750pt}%
\definecolor{currentstroke}{rgb}{1.000000,0.690196,0.000000}%
\pgfsetstrokecolor{currentstroke}%
\pgfsetdash{}{0pt}%
\pgfsys@defobject{currentmarker}{\pgfqpoint{-0.035355in}{-0.058926in}}{\pgfqpoint{0.035355in}{0.058926in}}{%
\pgfpathmoveto{\pgfqpoint{-0.000000in}{-0.058926in}}%
\pgfpathlineto{\pgfqpoint{0.035355in}{0.000000in}}%
\pgfpathlineto{\pgfqpoint{0.000000in}{0.058926in}}%
\pgfpathlineto{\pgfqpoint{-0.035355in}{0.000000in}}%
\pgfpathlineto{\pgfqpoint{-0.000000in}{-0.058926in}}%
\pgfpathclose%
\pgfusepath{stroke,fill}%
}%
\begin{pgfscope}%
\pgfsys@transformshift{0.777271in}{2.666573in}%
\pgfsys@useobject{currentmarker}{}%
\end{pgfscope}%
\begin{pgfscope}%
\pgfsys@transformshift{1.145143in}{2.607370in}%
\pgfsys@useobject{currentmarker}{}%
\end{pgfscope}%
\begin{pgfscope}%
\pgfsys@transformshift{1.526888in}{2.493234in}%
\pgfsys@useobject{currentmarker}{}%
\end{pgfscope}%
\begin{pgfscope}%
\pgfsys@transformshift{1.896020in}{2.240948in}%
\pgfsys@useobject{currentmarker}{}%
\end{pgfscope}%
\begin{pgfscope}%
\pgfsys@transformshift{2.258635in}{1.023687in}%
\pgfsys@useobject{currentmarker}{}%
\end{pgfscope}%
\begin{pgfscope}%
\pgfsys@transformshift{2.625694in}{0.900995in}%
\pgfsys@useobject{currentmarker}{}%
\end{pgfscope}%
\begin{pgfscope}%
\pgfsys@transformshift{2.991556in}{0.860471in}%
\pgfsys@useobject{currentmarker}{}%
\end{pgfscope}%
\end{pgfscope}%
\begin{pgfscope}%
\pgfpathrectangle{\pgfqpoint{0.721913in}{0.549073in}}{\pgfqpoint{2.325000in}{2.310000in}}%
\pgfusepath{clip}%
\pgfsetbuttcap%
\pgfsetroundjoin%
\pgfsetlinewidth{1.505625pt}%
\definecolor{currentstroke}{rgb}{0.478431,0.478431,0.478431}%
\pgfsetstrokecolor{currentstroke}%
\pgfsetstrokeopacity{0.500000}%
\pgfsetdash{{5.550000pt}{2.400000pt}}{0.000000pt}%
\pgfpathmoveto{\pgfqpoint{0.777271in}{2.334798in}}%
\pgfpathlineto{\pgfqpoint{1.145143in}{2.231390in}}%
\pgfpathlineto{\pgfqpoint{1.526888in}{2.124082in}}%
\pgfpathlineto{\pgfqpoint{1.896020in}{2.020320in}}%
\pgfpathlineto{\pgfqpoint{2.258635in}{1.918390in}}%
\pgfpathlineto{\pgfqpoint{2.625694in}{1.815211in}}%
\pgfpathlineto{\pgfqpoint{2.991556in}{1.712368in}}%
\pgfusepath{stroke}%
\end{pgfscope}%
\begin{pgfscope}%
\pgfpathrectangle{\pgfqpoint{0.721913in}{0.549073in}}{\pgfqpoint{2.325000in}{2.310000in}}%
\pgfusepath{clip}%
\pgfsetbuttcap%
\pgfsetroundjoin%
\pgfsetlinewidth{1.505625pt}%
\definecolor{currentstroke}{rgb}{0.478431,0.478431,0.478431}%
\pgfsetstrokecolor{currentstroke}%
\pgfsetstrokeopacity{0.500000}%
\pgfsetdash{{5.550000pt}{2.400000pt}}{0.000000pt}%
\pgfpathmoveto{\pgfqpoint{0.777271in}{2.488226in}}%
\pgfpathlineto{\pgfqpoint{1.145143in}{2.436522in}}%
\pgfpathlineto{\pgfqpoint{1.526888in}{2.382868in}}%
\pgfpathlineto{\pgfqpoint{1.896020in}{2.330987in}}%
\pgfpathlineto{\pgfqpoint{2.258635in}{2.280022in}}%
\pgfpathlineto{\pgfqpoint{2.625694in}{2.228432in}}%
\pgfpathlineto{\pgfqpoint{2.991556in}{2.177011in}}%
\pgfusepath{stroke}%
\end{pgfscope}%
\begin{pgfscope}%
\pgfsetrectcap%
\pgfsetmiterjoin%
\pgfsetlinewidth{0.803000pt}%
\definecolor{currentstroke}{rgb}{0.000000,0.000000,0.000000}%
\pgfsetstrokecolor{currentstroke}%
\pgfsetdash{}{0pt}%
\pgfpathmoveto{\pgfqpoint{0.721913in}{0.549073in}}%
\pgfpathlineto{\pgfqpoint{0.721913in}{2.859073in}}%
\pgfusepath{stroke}%
\end{pgfscope}%
\begin{pgfscope}%
\pgfsetrectcap%
\pgfsetmiterjoin%
\pgfsetlinewidth{0.803000pt}%
\definecolor{currentstroke}{rgb}{0.000000,0.000000,0.000000}%
\pgfsetstrokecolor{currentstroke}%
\pgfsetdash{}{0pt}%
\pgfpathmoveto{\pgfqpoint{3.046913in}{0.549073in}}%
\pgfpathlineto{\pgfqpoint{3.046913in}{2.859073in}}%
\pgfusepath{stroke}%
\end{pgfscope}%
\begin{pgfscope}%
\pgfsetrectcap%
\pgfsetmiterjoin%
\pgfsetlinewidth{0.803000pt}%
\definecolor{currentstroke}{rgb}{0.000000,0.000000,0.000000}%
\pgfsetstrokecolor{currentstroke}%
\pgfsetdash{}{0pt}%
\pgfpathmoveto{\pgfqpoint{0.721913in}{0.549073in}}%
\pgfpathlineto{\pgfqpoint{3.046913in}{0.549073in}}%
\pgfusepath{stroke}%
\end{pgfscope}%
\begin{pgfscope}%
\pgfsetrectcap%
\pgfsetmiterjoin%
\pgfsetlinewidth{0.803000pt}%
\definecolor{currentstroke}{rgb}{0.000000,0.000000,0.000000}%
\pgfsetstrokecolor{currentstroke}%
\pgfsetdash{}{0pt}%
\pgfpathmoveto{\pgfqpoint{0.721913in}{2.859073in}}%
\pgfpathlineto{\pgfqpoint{3.046913in}{2.859073in}}%
\pgfusepath{stroke}%
\end{pgfscope}%
\begin{pgfscope}%
\definecolor{textcolor}{rgb}{0.478431,0.478431,0.478431}%
\pgfsetstrokecolor{textcolor}%
\pgfsetfillcolor{textcolor}%
\pgftext[x=2.009682in,y=2.006229in,left,base]{\color{textcolor}{\rmfamily\fontsize{12.000000}{14.400000}\selectfont\catcode`\^=\active\def^{\ifmmode\sp\else\^{}\fi}\catcode`\%=\active\def
\end{pgfscope}%
\begin{pgfscope}%
\definecolor{textcolor}{rgb}{0.478431,0.478431,0.478431}%
\pgfsetstrokecolor{textcolor}%
\pgfsetfillcolor{textcolor}%
\pgftext[x=1.872807in,y=2.400024in,left,base]{\color{textcolor}{\rmfamily\fontsize{12.000000}{14.400000}\selectfont\catcode`\^=\active\def^{\ifmmode\sp\else\^{}\fi}\catcode`\%=\active\def
\end{pgfscope}%
\begin{pgfscope}%
\pgfsetbuttcap%
\pgfsetmiterjoin%
\definecolor{currentfill}{rgb}{1.000000,1.000000,1.000000}%
\pgfsetfillcolor{currentfill}%
\pgfsetfillopacity{0.800000}%
\pgfsetlinewidth{1.003750pt}%
\definecolor{currentstroke}{rgb}{0.800000,0.800000,0.800000}%
\pgfsetstrokecolor{currentstroke}%
\pgfsetstrokeopacity{0.800000}%
\pgfsetdash{}{0pt}%
\pgfpathmoveto{\pgfqpoint{0.805247in}{0.632406in}}%
\pgfpathlineto{\pgfqpoint{2.037651in}{0.632406in}}%
\pgfpathlineto{\pgfqpoint{2.037651in}{1.379627in}}%
\pgfpathlineto{\pgfqpoint{0.805247in}{1.379627in}}%
\pgfpathlineto{\pgfqpoint{0.805247in}{0.632406in}}%
\pgfpathclose%
\pgfusepath{stroke,fill}%
\end{pgfscope}%
\begin{pgfscope}%
\pgfsetrectcap%
\pgfsetroundjoin%
\pgfsetlinewidth{1.505625pt}%
\definecolor{currentstroke}{rgb}{0.392157,0.560784,1.000000}%
\pgfsetstrokecolor{currentstroke}%
\pgfsetdash{}{0pt}%
\pgfpathmoveto{\pgfqpoint{0.871913in}{1.254627in}}%
\pgfpathlineto{\pgfqpoint{1.038580in}{1.254627in}}%
\pgfpathlineto{\pgfqpoint{1.205247in}{1.254627in}}%
\pgfusepath{stroke}%
\end{pgfscope}%
\begin{pgfscope}%
\pgfsetbuttcap%
\pgfsetroundjoin%
\definecolor{currentfill}{rgb}{0.392157,0.560784,1.000000}%
\pgfsetfillcolor{currentfill}%
\pgfsetlinewidth{1.003750pt}%
\definecolor{currentstroke}{rgb}{0.392157,0.560784,1.000000}%
\pgfsetstrokecolor{currentstroke}%
\pgfsetdash{}{0pt}%
\pgfsys@defobject{currentmarker}{\pgfqpoint{-0.031250in}{-0.031250in}}{\pgfqpoint{0.031250in}{0.031250in}}{%
\pgfpathmoveto{\pgfqpoint{0.000000in}{-0.031250in}}%
\pgfpathcurveto{\pgfqpoint{0.008288in}{-0.031250in}}{\pgfqpoint{0.016237in}{-0.027957in}}{\pgfqpoint{0.022097in}{-0.022097in}}%
\pgfpathcurveto{\pgfqpoint{0.027957in}{-0.016237in}}{\pgfqpoint{0.031250in}{-0.008288in}}{\pgfqpoint{0.031250in}{0.000000in}}%
\pgfpathcurveto{\pgfqpoint{0.031250in}{0.008288in}}{\pgfqpoint{0.027957in}{0.016237in}}{\pgfqpoint{0.022097in}{0.022097in}}%
\pgfpathcurveto{\pgfqpoint{0.016237in}{0.027957in}}{\pgfqpoint{0.008288in}{0.031250in}}{\pgfqpoint{0.000000in}{0.031250in}}%
\pgfpathcurveto{\pgfqpoint{-0.008288in}{0.031250in}}{\pgfqpoint{-0.016237in}{0.027957in}}{\pgfqpoint{-0.022097in}{0.022097in}}%
\pgfpathcurveto{\pgfqpoint{-0.027957in}{0.016237in}}{\pgfqpoint{-0.031250in}{0.008288in}}{\pgfqpoint{-0.031250in}{0.000000in}}%
\pgfpathcurveto{\pgfqpoint{-0.031250in}{-0.008288in}}{\pgfqpoint{-0.027957in}{-0.016237in}}{\pgfqpoint{-0.022097in}{-0.022097in}}%
\pgfpathcurveto{\pgfqpoint{-0.016237in}{-0.027957in}}{\pgfqpoint{-0.008288in}{-0.031250in}}{\pgfqpoint{0.000000in}{-0.031250in}}%
\pgfpathlineto{\pgfqpoint{0.000000in}{-0.031250in}}%
\pgfpathclose%
\pgfusepath{stroke,fill}%
}%
\begin{pgfscope}%
\pgfsys@transformshift{1.038580in}{1.254627in}%
\pgfsys@useobject{currentmarker}{}%
\end{pgfscope}%
\end{pgfscope}%
\begin{pgfscope}%
\definecolor{textcolor}{rgb}{0.000000,0.000000,0.000000}%
\pgfsetstrokecolor{textcolor}%
\pgfsetfillcolor{textcolor}%
\pgftext[x=1.338580in,y=1.196294in,left,base]{\color{textcolor}{\rmfamily\fontsize{12.000000}{14.400000}\selectfont\catcode`\^=\active\def^{\ifmmode\sp\else\^{}\fi}\catcode`\%=\active\def
\end{pgfscope}%
\begin{pgfscope}%
\pgfsetrectcap%
\pgfsetroundjoin%
\pgfsetlinewidth{1.505625pt}%
\definecolor{currentstroke}{rgb}{0.862745,0.149020,0.498039}%
\pgfsetstrokecolor{currentstroke}%
\pgfsetdash{}{0pt}%
\pgfpathmoveto{\pgfqpoint{0.871913in}{1.022220in}}%
\pgfpathlineto{\pgfqpoint{1.038580in}{1.022220in}}%
\pgfpathlineto{\pgfqpoint{1.205247in}{1.022220in}}%
\pgfusepath{stroke}%
\end{pgfscope}%
\begin{pgfscope}%
\pgfsetbuttcap%
\pgfsetmiterjoin%
\definecolor{currentfill}{rgb}{0.862745,0.149020,0.498039}%
\pgfsetfillcolor{currentfill}%
\pgfsetlinewidth{1.003750pt}%
\definecolor{currentstroke}{rgb}{0.862745,0.149020,0.498039}%
\pgfsetstrokecolor{currentstroke}%
\pgfsetdash{}{0pt}%
\pgfsys@defobject{currentmarker}{\pgfqpoint{-0.031250in}{-0.031250in}}{\pgfqpoint{0.031250in}{0.031250in}}{%
\pgfpathmoveto{\pgfqpoint{-0.031250in}{-0.031250in}}%
\pgfpathlineto{\pgfqpoint{0.031250in}{-0.031250in}}%
\pgfpathlineto{\pgfqpoint{0.031250in}{0.031250in}}%
\pgfpathlineto{\pgfqpoint{-0.031250in}{0.031250in}}%
\pgfpathlineto{\pgfqpoint{-0.031250in}{-0.031250in}}%
\pgfpathclose%
\pgfusepath{stroke,fill}%
}%
\begin{pgfscope}%
\pgfsys@transformshift{1.038580in}{1.022220in}%
\pgfsys@useobject{currentmarker}{}%
\end{pgfscope}%
\end{pgfscope}%
\begin{pgfscope}%
\definecolor{textcolor}{rgb}{0.000000,0.000000,0.000000}%
\pgfsetstrokecolor{textcolor}%
\pgfsetfillcolor{textcolor}%
\pgftext[x=1.338580in,y=0.963887in,left,base]{\color{textcolor}{\rmfamily\fontsize{12.000000}{14.400000}\selectfont\catcode`\^=\active\def^{\ifmmode\sp\else\^{}\fi}\catcode`\%=\active\def
\end{pgfscope}%
\begin{pgfscope}%
\pgfsetrectcap%
\pgfsetroundjoin%
\pgfsetlinewidth{1.505625pt}%
\definecolor{currentstroke}{rgb}{1.000000,0.690196,0.000000}%
\pgfsetstrokecolor{currentstroke}%
\pgfsetdash{}{0pt}%
\pgfpathmoveto{\pgfqpoint{0.871913in}{0.789813in}}%
\pgfpathlineto{\pgfqpoint{1.038580in}{0.789813in}}%
\pgfpathlineto{\pgfqpoint{1.205247in}{0.789813in}}%
\pgfusepath{stroke}%
\end{pgfscope}%
\begin{pgfscope}%
\pgfsetbuttcap%
\pgfsetmiterjoin%
\definecolor{currentfill}{rgb}{1.000000,0.690196,0.000000}%
\pgfsetfillcolor{currentfill}%
\pgfsetlinewidth{1.003750pt}%
\definecolor{currentstroke}{rgb}{1.000000,0.690196,0.000000}%
\pgfsetstrokecolor{currentstroke}%
\pgfsetdash{}{0pt}%
\pgfsys@defobject{currentmarker}{\pgfqpoint{-0.026517in}{-0.044194in}}{\pgfqpoint{0.026517in}{0.044194in}}{%
\pgfpathmoveto{\pgfqpoint{-0.000000in}{-0.044194in}}%
\pgfpathlineto{\pgfqpoint{0.026517in}{0.000000in}}%
\pgfpathlineto{\pgfqpoint{0.000000in}{0.044194in}}%
\pgfpathlineto{\pgfqpoint{-0.026517in}{0.000000in}}%
\pgfpathlineto{\pgfqpoint{-0.000000in}{-0.044194in}}%
\pgfpathclose%
\pgfusepath{stroke,fill}%
}%
\begin{pgfscope}%
\pgfsys@transformshift{1.038580in}{0.789813in}%
\pgfsys@useobject{currentmarker}{}%
\end{pgfscope}%
\end{pgfscope}%
\begin{pgfscope}%
\definecolor{textcolor}{rgb}{0.000000,0.000000,0.000000}%
\pgfsetstrokecolor{textcolor}%
\pgfsetfillcolor{textcolor}%
\pgftext[x=1.338580in,y=0.731480in,left,base]{\color{textcolor}{\rmfamily\fontsize{12.000000}{14.400000}\selectfont\catcode`\^=\active\def^{\ifmmode\sp\else\^{}\fi}\catcode`\%=\active\def
\end{pgfscope}%
\end{pgfpicture}%
\makeatother%
\endgroup%

%% file: plots/convergence-3.pgf
\begingroup%
\makeatletter%
\begin{pgfpicture}%
\pgfpathrectangle{\pgfpointorigin}{\pgfqpoint{3.146913in}{2.959073in}}%
\pgfusepath{use as bounding box, clip}%
\begin{pgfscope}%
\pgfsetbuttcap%
\pgfsetmiterjoin%
\definecolor{currentfill}{rgb}{1.000000,1.000000,1.000000}%
\pgfsetfillcolor{currentfill}%
\pgfsetlinewidth{0.000000pt}%
\definecolor{currentstroke}{rgb}{1.000000,1.000000,1.000000}%
\pgfsetstrokecolor{currentstroke}%
\pgfsetdash{}{0pt}%
\pgfpathmoveto{\pgfqpoint{0.000000in}{-0.000000in}}%
\pgfpathlineto{\pgfqpoint{3.146913in}{-0.000000in}}%
\pgfpathlineto{\pgfqpoint{3.146913in}{2.959073in}}%
\pgfpathlineto{\pgfqpoint{0.000000in}{2.959073in}}%
\pgfpathlineto{\pgfqpoint{0.000000in}{-0.000000in}}%
\pgfpathclose%
\pgfusepath{fill}%
\end{pgfscope}%
\begin{pgfscope}%
\pgfsetbuttcap%
\pgfsetmiterjoin%
\definecolor{currentfill}{rgb}{1.000000,1.000000,1.000000}%
\pgfsetfillcolor{currentfill}%
\pgfsetlinewidth{0.000000pt}%
\definecolor{currentstroke}{rgb}{0.000000,0.000000,0.000000}%
\pgfsetstrokecolor{currentstroke}%
\pgfsetstrokeopacity{0.000000}%
\pgfsetdash{}{0pt}%
\pgfpathmoveto{\pgfqpoint{0.721913in}{0.549073in}}%
\pgfpathlineto{\pgfqpoint{3.046913in}{0.549073in}}%
\pgfpathlineto{\pgfqpoint{3.046913in}{2.859073in}}%
\pgfpathlineto{\pgfqpoint{0.721913in}{2.859073in}}%
\pgfpathlineto{\pgfqpoint{0.721913in}{0.549073in}}%
\pgfpathclose%
\pgfusepath{fill}%
\end{pgfscope}%
\begin{pgfscope}%
\pgfpathrectangle{\pgfqpoint{0.721913in}{0.549073in}}{\pgfqpoint{2.325000in}{2.310000in}}%
\pgfusepath{clip}%
\pgfsetrectcap%
\pgfsetroundjoin%
\pgfsetlinewidth{0.250937pt}%
\definecolor{currentstroke}{rgb}{0.000000,0.000000,0.000000}%
\pgfsetstrokecolor{currentstroke}%
\pgfsetstrokeopacity{0.200000}%
\pgfsetdash{}{0pt}%
\pgfpathmoveto{\pgfqpoint{0.964644in}{0.549073in}}%
\pgfpathlineto{\pgfqpoint{0.964644in}{2.859073in}}%
\pgfusepath{stroke}%
\end{pgfscope}%
\begin{pgfscope}%
\pgfsetbuttcap%
\pgfsetroundjoin%
\definecolor{currentfill}{rgb}{0.000000,0.000000,0.000000}%
\pgfsetfillcolor{currentfill}%
\pgfsetlinewidth{0.803000pt}%
\definecolor{currentstroke}{rgb}{0.000000,0.000000,0.000000}%
\pgfsetstrokecolor{currentstroke}%
\pgfsetdash{}{0pt}%
\pgfsys@defobject{currentmarker}{\pgfqpoint{0.000000in}{-0.048611in}}{\pgfqpoint{0.000000in}{0.000000in}}{%
\pgfpathmoveto{\pgfqpoint{0.000000in}{0.000000in}}%
\pgfpathlineto{\pgfqpoint{0.000000in}{-0.048611in}}%
\pgfusepath{stroke,fill}%
}%
\begin{pgfscope}%
\pgfsys@transformshift{0.964644in}{0.549073in}%
\pgfsys@useobject{currentmarker}{}%
\end{pgfscope}%
\end{pgfscope}%
\begin{pgfscope}%
\definecolor{textcolor}{rgb}{0.000000,0.000000,0.000000}%
\pgfsetstrokecolor{textcolor}%
\pgfsetfillcolor{textcolor}%
\pgftext[x=0.964644in,y=0.451851in,,top]{\color{textcolor}{\rmfamily\fontsize{12.000000}{14.400000}\selectfont\catcode`\^=\active\def^{\ifmmode\sp\else\^{}\fi}\catcode`\%=\active\def
\end{pgfscope}%
\begin{pgfscope}%
\pgfpathrectangle{\pgfqpoint{0.721913in}{0.549073in}}{\pgfqpoint{2.325000in}{2.310000in}}%
\pgfusepath{clip}%
\pgfsetrectcap%
\pgfsetroundjoin%
\pgfsetlinewidth{0.250937pt}%
\definecolor{currentstroke}{rgb}{0.000000,0.000000,0.000000}%
\pgfsetstrokecolor{currentstroke}%
\pgfsetstrokeopacity{0.200000}%
\pgfsetdash{}{0pt}%
\pgfpathmoveto{\pgfqpoint{1.809245in}{0.549073in}}%
\pgfpathlineto{\pgfqpoint{1.809245in}{2.859073in}}%
\pgfusepath{stroke}%
\end{pgfscope}%
\begin{pgfscope}%
\pgfsetbuttcap%
\pgfsetroundjoin%
\definecolor{currentfill}{rgb}{0.000000,0.000000,0.000000}%
\pgfsetfillcolor{currentfill}%
\pgfsetlinewidth{0.803000pt}%
\definecolor{currentstroke}{rgb}{0.000000,0.000000,0.000000}%
\pgfsetstrokecolor{currentstroke}%
\pgfsetdash{}{0pt}%
\pgfsys@defobject{currentmarker}{\pgfqpoint{0.000000in}{-0.048611in}}{\pgfqpoint{0.000000in}{0.000000in}}{%
\pgfpathmoveto{\pgfqpoint{0.000000in}{0.000000in}}%
\pgfpathlineto{\pgfqpoint{0.000000in}{-0.048611in}}%
\pgfusepath{stroke,fill}%
}%
\begin{pgfscope}%
\pgfsys@transformshift{1.809245in}{0.549073in}%
\pgfsys@useobject{currentmarker}{}%
\end{pgfscope}%
\end{pgfscope}%
\begin{pgfscope}%
\definecolor{textcolor}{rgb}{0.000000,0.000000,0.000000}%
\pgfsetstrokecolor{textcolor}%
\pgfsetfillcolor{textcolor}%
\pgftext[x=1.809245in,y=0.451851in,,top]{\color{textcolor}{\rmfamily\fontsize{12.000000}{14.400000}\selectfont\catcode`\^=\active\def^{\ifmmode\sp\else\^{}\fi}\catcode`\%=\active\def
\end{pgfscope}%
\begin{pgfscope}%
\pgfpathrectangle{\pgfqpoint{0.721913in}{0.549073in}}{\pgfqpoint{2.325000in}{2.310000in}}%
\pgfusepath{clip}%
\pgfsetrectcap%
\pgfsetroundjoin%
\pgfsetlinewidth{0.250937pt}%
\definecolor{currentstroke}{rgb}{0.000000,0.000000,0.000000}%
\pgfsetstrokecolor{currentstroke}%
\pgfsetstrokeopacity{0.200000}%
\pgfsetdash{}{0pt}%
\pgfpathmoveto{\pgfqpoint{2.653846in}{0.549073in}}%
\pgfpathlineto{\pgfqpoint{2.653846in}{2.859073in}}%
\pgfusepath{stroke}%
\end{pgfscope}%
\begin{pgfscope}%
\pgfsetbuttcap%
\pgfsetroundjoin%
\definecolor{currentfill}{rgb}{0.000000,0.000000,0.000000}%
\pgfsetfillcolor{currentfill}%
\pgfsetlinewidth{0.803000pt}%
\definecolor{currentstroke}{rgb}{0.000000,0.000000,0.000000}%
\pgfsetstrokecolor{currentstroke}%
\pgfsetdash{}{0pt}%
\pgfsys@defobject{currentmarker}{\pgfqpoint{0.000000in}{-0.048611in}}{\pgfqpoint{0.000000in}{0.000000in}}{%
\pgfpathmoveto{\pgfqpoint{0.000000in}{0.000000in}}%
\pgfpathlineto{\pgfqpoint{0.000000in}{-0.048611in}}%
\pgfusepath{stroke,fill}%
}%
\begin{pgfscope}%
\pgfsys@transformshift{2.653846in}{0.549073in}%
\pgfsys@useobject{currentmarker}{}%
\end{pgfscope}%
\end{pgfscope}%
\begin{pgfscope}%
\definecolor{textcolor}{rgb}{0.000000,0.000000,0.000000}%
\pgfsetstrokecolor{textcolor}%
\pgfsetfillcolor{textcolor}%
\pgftext[x=2.653846in,y=0.451851in,,top]{\color{textcolor}{\rmfamily\fontsize{12.000000}{14.400000}\selectfont\catcode`\^=\active\def^{\ifmmode\sp\else\^{}\fi}\catcode`\%=\active\def
\end{pgfscope}%
\begin{pgfscope}%
\pgfsetbuttcap%
\pgfsetroundjoin%
\definecolor{currentfill}{rgb}{0.000000,0.000000,0.000000}%
\pgfsetfillcolor{currentfill}%
\pgfsetlinewidth{0.602250pt}%
\definecolor{currentstroke}{rgb}{0.000000,0.000000,0.000000}%
\pgfsetstrokecolor{currentstroke}%
\pgfsetdash{}{0pt}%
\pgfsys@defobject{currentmarker}{\pgfqpoint{0.000000in}{-0.027778in}}{\pgfqpoint{0.000000in}{0.000000in}}{%
\pgfpathmoveto{\pgfqpoint{0.000000in}{0.000000in}}%
\pgfpathlineto{\pgfqpoint{0.000000in}{-0.027778in}}%
\pgfusepath{stroke,fill}%
}%
\begin{pgfscope}%
\pgfsys@transformshift{0.777271in}{0.549073in}%
\pgfsys@useobject{currentmarker}{}%
\end{pgfscope}%
\end{pgfscope}%
\begin{pgfscope}%
\pgfsetbuttcap%
\pgfsetroundjoin%
\definecolor{currentfill}{rgb}{0.000000,0.000000,0.000000}%
\pgfsetfillcolor{currentfill}%
\pgfsetlinewidth{0.602250pt}%
\definecolor{currentstroke}{rgb}{0.000000,0.000000,0.000000}%
\pgfsetstrokecolor{currentstroke}%
\pgfsetdash{}{0pt}%
\pgfsys@defobject{currentmarker}{\pgfqpoint{0.000000in}{-0.027778in}}{\pgfqpoint{0.000000in}{0.000000in}}{%
\pgfpathmoveto{\pgfqpoint{0.000000in}{0.000000in}}%
\pgfpathlineto{\pgfqpoint{0.000000in}{-0.027778in}}%
\pgfusepath{stroke,fill}%
}%
\begin{pgfscope}%
\pgfsys@transformshift{0.833814in}{0.549073in}%
\pgfsys@useobject{currentmarker}{}%
\end{pgfscope}%
\end{pgfscope}%
\begin{pgfscope}%
\pgfsetbuttcap%
\pgfsetroundjoin%
\definecolor{currentfill}{rgb}{0.000000,0.000000,0.000000}%
\pgfsetfillcolor{currentfill}%
\pgfsetlinewidth{0.602250pt}%
\definecolor{currentstroke}{rgb}{0.000000,0.000000,0.000000}%
\pgfsetstrokecolor{currentstroke}%
\pgfsetdash{}{0pt}%
\pgfsys@defobject{currentmarker}{\pgfqpoint{0.000000in}{-0.027778in}}{\pgfqpoint{0.000000in}{0.000000in}}{%
\pgfpathmoveto{\pgfqpoint{0.000000in}{0.000000in}}%
\pgfpathlineto{\pgfqpoint{0.000000in}{-0.027778in}}%
\pgfusepath{stroke,fill}%
}%
\begin{pgfscope}%
\pgfsys@transformshift{0.882794in}{0.549073in}%
\pgfsys@useobject{currentmarker}{}%
\end{pgfscope}%
\end{pgfscope}%
\begin{pgfscope}%
\pgfsetbuttcap%
\pgfsetroundjoin%
\definecolor{currentfill}{rgb}{0.000000,0.000000,0.000000}%
\pgfsetfillcolor{currentfill}%
\pgfsetlinewidth{0.602250pt}%
\definecolor{currentstroke}{rgb}{0.000000,0.000000,0.000000}%
\pgfsetstrokecolor{currentstroke}%
\pgfsetdash{}{0pt}%
\pgfsys@defobject{currentmarker}{\pgfqpoint{0.000000in}{-0.027778in}}{\pgfqpoint{0.000000in}{0.000000in}}{%
\pgfpathmoveto{\pgfqpoint{0.000000in}{0.000000in}}%
\pgfpathlineto{\pgfqpoint{0.000000in}{-0.027778in}}%
\pgfusepath{stroke,fill}%
}%
\begin{pgfscope}%
\pgfsys@transformshift{0.925997in}{0.549073in}%
\pgfsys@useobject{currentmarker}{}%
\end{pgfscope}%
\end{pgfscope}%
\begin{pgfscope}%
\pgfsetbuttcap%
\pgfsetroundjoin%
\definecolor{currentfill}{rgb}{0.000000,0.000000,0.000000}%
\pgfsetfillcolor{currentfill}%
\pgfsetlinewidth{0.602250pt}%
\definecolor{currentstroke}{rgb}{0.000000,0.000000,0.000000}%
\pgfsetstrokecolor{currentstroke}%
\pgfsetdash{}{0pt}%
\pgfsys@defobject{currentmarker}{\pgfqpoint{0.000000in}{-0.027778in}}{\pgfqpoint{0.000000in}{0.000000in}}{%
\pgfpathmoveto{\pgfqpoint{0.000000in}{0.000000in}}%
\pgfpathlineto{\pgfqpoint{0.000000in}{-0.027778in}}%
\pgfusepath{stroke,fill}%
}%
\begin{pgfscope}%
\pgfsys@transformshift{1.218894in}{0.549073in}%
\pgfsys@useobject{currentmarker}{}%
\end{pgfscope}%
\end{pgfscope}%
\begin{pgfscope}%
\pgfsetbuttcap%
\pgfsetroundjoin%
\definecolor{currentfill}{rgb}{0.000000,0.000000,0.000000}%
\pgfsetfillcolor{currentfill}%
\pgfsetlinewidth{0.602250pt}%
\definecolor{currentstroke}{rgb}{0.000000,0.000000,0.000000}%
\pgfsetstrokecolor{currentstroke}%
\pgfsetdash{}{0pt}%
\pgfsys@defobject{currentmarker}{\pgfqpoint{0.000000in}{-0.027778in}}{\pgfqpoint{0.000000in}{0.000000in}}{%
\pgfpathmoveto{\pgfqpoint{0.000000in}{0.000000in}}%
\pgfpathlineto{\pgfqpoint{0.000000in}{-0.027778in}}%
\pgfusepath{stroke,fill}%
}%
\begin{pgfscope}%
\pgfsys@transformshift{1.367621in}{0.549073in}%
\pgfsys@useobject{currentmarker}{}%
\end{pgfscope}%
\end{pgfscope}%
\begin{pgfscope}%
\pgfsetbuttcap%
\pgfsetroundjoin%
\definecolor{currentfill}{rgb}{0.000000,0.000000,0.000000}%
\pgfsetfillcolor{currentfill}%
\pgfsetlinewidth{0.602250pt}%
\definecolor{currentstroke}{rgb}{0.000000,0.000000,0.000000}%
\pgfsetstrokecolor{currentstroke}%
\pgfsetdash{}{0pt}%
\pgfsys@defobject{currentmarker}{\pgfqpoint{0.000000in}{-0.027778in}}{\pgfqpoint{0.000000in}{0.000000in}}{%
\pgfpathmoveto{\pgfqpoint{0.000000in}{0.000000in}}%
\pgfpathlineto{\pgfqpoint{0.000000in}{-0.027778in}}%
\pgfusepath{stroke,fill}%
}%
\begin{pgfscope}%
\pgfsys@transformshift{1.473144in}{0.549073in}%
\pgfsys@useobject{currentmarker}{}%
\end{pgfscope}%
\end{pgfscope}%
\begin{pgfscope}%
\pgfsetbuttcap%
\pgfsetroundjoin%
\definecolor{currentfill}{rgb}{0.000000,0.000000,0.000000}%
\pgfsetfillcolor{currentfill}%
\pgfsetlinewidth{0.602250pt}%
\definecolor{currentstroke}{rgb}{0.000000,0.000000,0.000000}%
\pgfsetstrokecolor{currentstroke}%
\pgfsetdash{}{0pt}%
\pgfsys@defobject{currentmarker}{\pgfqpoint{0.000000in}{-0.027778in}}{\pgfqpoint{0.000000in}{0.000000in}}{%
\pgfpathmoveto{\pgfqpoint{0.000000in}{0.000000in}}%
\pgfpathlineto{\pgfqpoint{0.000000in}{-0.027778in}}%
\pgfusepath{stroke,fill}%
}%
\begin{pgfscope}%
\pgfsys@transformshift{1.554995in}{0.549073in}%
\pgfsys@useobject{currentmarker}{}%
\end{pgfscope}%
\end{pgfscope}%
\begin{pgfscope}%
\pgfsetbuttcap%
\pgfsetroundjoin%
\definecolor{currentfill}{rgb}{0.000000,0.000000,0.000000}%
\pgfsetfillcolor{currentfill}%
\pgfsetlinewidth{0.602250pt}%
\definecolor{currentstroke}{rgb}{0.000000,0.000000,0.000000}%
\pgfsetstrokecolor{currentstroke}%
\pgfsetdash{}{0pt}%
\pgfsys@defobject{currentmarker}{\pgfqpoint{0.000000in}{-0.027778in}}{\pgfqpoint{0.000000in}{0.000000in}}{%
\pgfpathmoveto{\pgfqpoint{0.000000in}{0.000000in}}%
\pgfpathlineto{\pgfqpoint{0.000000in}{-0.027778in}}%
\pgfusepath{stroke,fill}%
}%
\begin{pgfscope}%
\pgfsys@transformshift{1.621871in}{0.549073in}%
\pgfsys@useobject{currentmarker}{}%
\end{pgfscope}%
\end{pgfscope}%
\begin{pgfscope}%
\pgfsetbuttcap%
\pgfsetroundjoin%
\definecolor{currentfill}{rgb}{0.000000,0.000000,0.000000}%
\pgfsetfillcolor{currentfill}%
\pgfsetlinewidth{0.602250pt}%
\definecolor{currentstroke}{rgb}{0.000000,0.000000,0.000000}%
\pgfsetstrokecolor{currentstroke}%
\pgfsetdash{}{0pt}%
\pgfsys@defobject{currentmarker}{\pgfqpoint{0.000000in}{-0.027778in}}{\pgfqpoint{0.000000in}{0.000000in}}{%
\pgfpathmoveto{\pgfqpoint{0.000000in}{0.000000in}}%
\pgfpathlineto{\pgfqpoint{0.000000in}{-0.027778in}}%
\pgfusepath{stroke,fill}%
}%
\begin{pgfscope}%
\pgfsys@transformshift{1.678415in}{0.549073in}%
\pgfsys@useobject{currentmarker}{}%
\end{pgfscope}%
\end{pgfscope}%
\begin{pgfscope}%
\pgfsetbuttcap%
\pgfsetroundjoin%
\definecolor{currentfill}{rgb}{0.000000,0.000000,0.000000}%
\pgfsetfillcolor{currentfill}%
\pgfsetlinewidth{0.602250pt}%
\definecolor{currentstroke}{rgb}{0.000000,0.000000,0.000000}%
\pgfsetstrokecolor{currentstroke}%
\pgfsetdash{}{0pt}%
\pgfsys@defobject{currentmarker}{\pgfqpoint{0.000000in}{-0.027778in}}{\pgfqpoint{0.000000in}{0.000000in}}{%
\pgfpathmoveto{\pgfqpoint{0.000000in}{0.000000in}}%
\pgfpathlineto{\pgfqpoint{0.000000in}{-0.027778in}}%
\pgfusepath{stroke,fill}%
}%
\begin{pgfscope}%
\pgfsys@transformshift{1.727395in}{0.549073in}%
\pgfsys@useobject{currentmarker}{}%
\end{pgfscope}%
\end{pgfscope}%
\begin{pgfscope}%
\pgfsetbuttcap%
\pgfsetroundjoin%
\definecolor{currentfill}{rgb}{0.000000,0.000000,0.000000}%
\pgfsetfillcolor{currentfill}%
\pgfsetlinewidth{0.602250pt}%
\definecolor{currentstroke}{rgb}{0.000000,0.000000,0.000000}%
\pgfsetstrokecolor{currentstroke}%
\pgfsetdash{}{0pt}%
\pgfsys@defobject{currentmarker}{\pgfqpoint{0.000000in}{-0.027778in}}{\pgfqpoint{0.000000in}{0.000000in}}{%
\pgfpathmoveto{\pgfqpoint{0.000000in}{0.000000in}}%
\pgfpathlineto{\pgfqpoint{0.000000in}{-0.027778in}}%
\pgfusepath{stroke,fill}%
}%
\begin{pgfscope}%
\pgfsys@transformshift{1.770598in}{0.549073in}%
\pgfsys@useobject{currentmarker}{}%
\end{pgfscope}%
\end{pgfscope}%
\begin{pgfscope}%
\pgfsetbuttcap%
\pgfsetroundjoin%
\definecolor{currentfill}{rgb}{0.000000,0.000000,0.000000}%
\pgfsetfillcolor{currentfill}%
\pgfsetlinewidth{0.602250pt}%
\definecolor{currentstroke}{rgb}{0.000000,0.000000,0.000000}%
\pgfsetstrokecolor{currentstroke}%
\pgfsetdash{}{0pt}%
\pgfsys@defobject{currentmarker}{\pgfqpoint{0.000000in}{-0.027778in}}{\pgfqpoint{0.000000in}{0.000000in}}{%
\pgfpathmoveto{\pgfqpoint{0.000000in}{0.000000in}}%
\pgfpathlineto{\pgfqpoint{0.000000in}{-0.027778in}}%
\pgfusepath{stroke,fill}%
}%
\begin{pgfscope}%
\pgfsys@transformshift{2.063495in}{0.549073in}%
\pgfsys@useobject{currentmarker}{}%
\end{pgfscope}%
\end{pgfscope}%
\begin{pgfscope}%
\pgfsetbuttcap%
\pgfsetroundjoin%
\definecolor{currentfill}{rgb}{0.000000,0.000000,0.000000}%
\pgfsetfillcolor{currentfill}%
\pgfsetlinewidth{0.602250pt}%
\definecolor{currentstroke}{rgb}{0.000000,0.000000,0.000000}%
\pgfsetstrokecolor{currentstroke}%
\pgfsetdash{}{0pt}%
\pgfsys@defobject{currentmarker}{\pgfqpoint{0.000000in}{-0.027778in}}{\pgfqpoint{0.000000in}{0.000000in}}{%
\pgfpathmoveto{\pgfqpoint{0.000000in}{0.000000in}}%
\pgfpathlineto{\pgfqpoint{0.000000in}{-0.027778in}}%
\pgfusepath{stroke,fill}%
}%
\begin{pgfscope}%
\pgfsys@transformshift{2.212222in}{0.549073in}%
\pgfsys@useobject{currentmarker}{}%
\end{pgfscope}%
\end{pgfscope}%
\begin{pgfscope}%
\pgfsetbuttcap%
\pgfsetroundjoin%
\definecolor{currentfill}{rgb}{0.000000,0.000000,0.000000}%
\pgfsetfillcolor{currentfill}%
\pgfsetlinewidth{0.602250pt}%
\definecolor{currentstroke}{rgb}{0.000000,0.000000,0.000000}%
\pgfsetstrokecolor{currentstroke}%
\pgfsetdash{}{0pt}%
\pgfsys@defobject{currentmarker}{\pgfqpoint{0.000000in}{-0.027778in}}{\pgfqpoint{0.000000in}{0.000000in}}{%
\pgfpathmoveto{\pgfqpoint{0.000000in}{0.000000in}}%
\pgfpathlineto{\pgfqpoint{0.000000in}{-0.027778in}}%
\pgfusepath{stroke,fill}%
}%
\begin{pgfscope}%
\pgfsys@transformshift{2.317745in}{0.549073in}%
\pgfsys@useobject{currentmarker}{}%
\end{pgfscope}%
\end{pgfscope}%
\begin{pgfscope}%
\pgfsetbuttcap%
\pgfsetroundjoin%
\definecolor{currentfill}{rgb}{0.000000,0.000000,0.000000}%
\pgfsetfillcolor{currentfill}%
\pgfsetlinewidth{0.602250pt}%
\definecolor{currentstroke}{rgb}{0.000000,0.000000,0.000000}%
\pgfsetstrokecolor{currentstroke}%
\pgfsetdash{}{0pt}%
\pgfsys@defobject{currentmarker}{\pgfqpoint{0.000000in}{-0.027778in}}{\pgfqpoint{0.000000in}{0.000000in}}{%
\pgfpathmoveto{\pgfqpoint{0.000000in}{0.000000in}}%
\pgfpathlineto{\pgfqpoint{0.000000in}{-0.027778in}}%
\pgfusepath{stroke,fill}%
}%
\begin{pgfscope}%
\pgfsys@transformshift{2.399595in}{0.549073in}%
\pgfsys@useobject{currentmarker}{}%
\end{pgfscope}%
\end{pgfscope}%
\begin{pgfscope}%
\pgfsetbuttcap%
\pgfsetroundjoin%
\definecolor{currentfill}{rgb}{0.000000,0.000000,0.000000}%
\pgfsetfillcolor{currentfill}%
\pgfsetlinewidth{0.602250pt}%
\definecolor{currentstroke}{rgb}{0.000000,0.000000,0.000000}%
\pgfsetstrokecolor{currentstroke}%
\pgfsetdash{}{0pt}%
\pgfsys@defobject{currentmarker}{\pgfqpoint{0.000000in}{-0.027778in}}{\pgfqpoint{0.000000in}{0.000000in}}{%
\pgfpathmoveto{\pgfqpoint{0.000000in}{0.000000in}}%
\pgfpathlineto{\pgfqpoint{0.000000in}{-0.027778in}}%
\pgfusepath{stroke,fill}%
}%
\begin{pgfscope}%
\pgfsys@transformshift{2.466472in}{0.549073in}%
\pgfsys@useobject{currentmarker}{}%
\end{pgfscope}%
\end{pgfscope}%
\begin{pgfscope}%
\pgfsetbuttcap%
\pgfsetroundjoin%
\definecolor{currentfill}{rgb}{0.000000,0.000000,0.000000}%
\pgfsetfillcolor{currentfill}%
\pgfsetlinewidth{0.602250pt}%
\definecolor{currentstroke}{rgb}{0.000000,0.000000,0.000000}%
\pgfsetstrokecolor{currentstroke}%
\pgfsetdash{}{0pt}%
\pgfsys@defobject{currentmarker}{\pgfqpoint{0.000000in}{-0.027778in}}{\pgfqpoint{0.000000in}{0.000000in}}{%
\pgfpathmoveto{\pgfqpoint{0.000000in}{0.000000in}}%
\pgfpathlineto{\pgfqpoint{0.000000in}{-0.027778in}}%
\pgfusepath{stroke,fill}%
}%
\begin{pgfscope}%
\pgfsys@transformshift{2.523015in}{0.549073in}%
\pgfsys@useobject{currentmarker}{}%
\end{pgfscope}%
\end{pgfscope}%
\begin{pgfscope}%
\pgfsetbuttcap%
\pgfsetroundjoin%
\definecolor{currentfill}{rgb}{0.000000,0.000000,0.000000}%
\pgfsetfillcolor{currentfill}%
\pgfsetlinewidth{0.602250pt}%
\definecolor{currentstroke}{rgb}{0.000000,0.000000,0.000000}%
\pgfsetstrokecolor{currentstroke}%
\pgfsetdash{}{0pt}%
\pgfsys@defobject{currentmarker}{\pgfqpoint{0.000000in}{-0.027778in}}{\pgfqpoint{0.000000in}{0.000000in}}{%
\pgfpathmoveto{\pgfqpoint{0.000000in}{0.000000in}}%
\pgfpathlineto{\pgfqpoint{0.000000in}{-0.027778in}}%
\pgfusepath{stroke,fill}%
}%
\begin{pgfscope}%
\pgfsys@transformshift{2.571995in}{0.549073in}%
\pgfsys@useobject{currentmarker}{}%
\end{pgfscope}%
\end{pgfscope}%
\begin{pgfscope}%
\pgfsetbuttcap%
\pgfsetroundjoin%
\definecolor{currentfill}{rgb}{0.000000,0.000000,0.000000}%
\pgfsetfillcolor{currentfill}%
\pgfsetlinewidth{0.602250pt}%
\definecolor{currentstroke}{rgb}{0.000000,0.000000,0.000000}%
\pgfsetstrokecolor{currentstroke}%
\pgfsetdash{}{0pt}%
\pgfsys@defobject{currentmarker}{\pgfqpoint{0.000000in}{-0.027778in}}{\pgfqpoint{0.000000in}{0.000000in}}{%
\pgfpathmoveto{\pgfqpoint{0.000000in}{0.000000in}}%
\pgfpathlineto{\pgfqpoint{0.000000in}{-0.027778in}}%
\pgfusepath{stroke,fill}%
}%
\begin{pgfscope}%
\pgfsys@transformshift{2.615199in}{0.549073in}%
\pgfsys@useobject{currentmarker}{}%
\end{pgfscope}%
\end{pgfscope}%
\begin{pgfscope}%
\pgfsetbuttcap%
\pgfsetroundjoin%
\definecolor{currentfill}{rgb}{0.000000,0.000000,0.000000}%
\pgfsetfillcolor{currentfill}%
\pgfsetlinewidth{0.602250pt}%
\definecolor{currentstroke}{rgb}{0.000000,0.000000,0.000000}%
\pgfsetstrokecolor{currentstroke}%
\pgfsetdash{}{0pt}%
\pgfsys@defobject{currentmarker}{\pgfqpoint{0.000000in}{-0.027778in}}{\pgfqpoint{0.000000in}{0.000000in}}{%
\pgfpathmoveto{\pgfqpoint{0.000000in}{0.000000in}}%
\pgfpathlineto{\pgfqpoint{0.000000in}{-0.027778in}}%
\pgfusepath{stroke,fill}%
}%
\begin{pgfscope}%
\pgfsys@transformshift{2.908096in}{0.549073in}%
\pgfsys@useobject{currentmarker}{}%
\end{pgfscope}%
\end{pgfscope}%
\begin{pgfscope}%
\definecolor{textcolor}{rgb}{0.000000,0.000000,0.000000}%
\pgfsetstrokecolor{textcolor}%
\pgfsetfillcolor{textcolor}%
\pgftext[x=1.884413in,y=0.248148in,,top]{\color{textcolor}{\rmfamily\fontsize{12.000000}{14.400000}\selectfont\catcode`\^=\active\def^{\ifmmode\sp\else\^{}\fi}\catcode`\%=\active\def
\end{pgfscope}%
\begin{pgfscope}%
\pgfpathrectangle{\pgfqpoint{0.721913in}{0.549073in}}{\pgfqpoint{2.325000in}{2.310000in}}%
\pgfusepath{clip}%
\pgfsetrectcap%
\pgfsetroundjoin%
\pgfsetlinewidth{0.250937pt}%
\definecolor{currentstroke}{rgb}{0.000000,0.000000,0.000000}%
\pgfsetstrokecolor{currentstroke}%
\pgfsetstrokeopacity{0.200000}%
\pgfsetdash{}{0pt}%
\pgfpathmoveto{\pgfqpoint{0.721913in}{0.856046in}}%
\pgfpathlineto{\pgfqpoint{3.046913in}{0.856046in}}%
\pgfusepath{stroke}%
\end{pgfscope}%
\begin{pgfscope}%
\pgfsetbuttcap%
\pgfsetroundjoin%
\definecolor{currentfill}{rgb}{0.000000,0.000000,0.000000}%
\pgfsetfillcolor{currentfill}%
\pgfsetlinewidth{0.803000pt}%
\definecolor{currentstroke}{rgb}{0.000000,0.000000,0.000000}%
\pgfsetstrokecolor{currentstroke}%
\pgfsetdash{}{0pt}%
\pgfsys@defobject{currentmarker}{\pgfqpoint{-0.048611in}{0.000000in}}{\pgfqpoint{-0.000000in}{0.000000in}}{%
\pgfpathmoveto{\pgfqpoint{-0.000000in}{0.000000in}}%
\pgfpathlineto{\pgfqpoint{-0.048611in}{0.000000in}}%
\pgfusepath{stroke,fill}%
}%
\begin{pgfscope}%
\pgfsys@transformshift{0.721913in}{0.856046in}%
\pgfsys@useobject{currentmarker}{}%
\end{pgfscope}%
\end{pgfscope}%
\begin{pgfscope}%
\definecolor{textcolor}{rgb}{0.000000,0.000000,0.000000}%
\pgfsetstrokecolor{textcolor}%
\pgfsetfillcolor{textcolor}%
\pgftext[x=0.303703in, y=0.798176in, left, base]{\color{textcolor}{\rmfamily\fontsize{12.000000}{14.400000}\selectfont\catcode`\^=\active\def^{\ifmmode\sp\else\^{}\fi}\catcode`\%=\active\def
\end{pgfscope}%
\begin{pgfscope}%
\pgfpathrectangle{\pgfqpoint{0.721913in}{0.549073in}}{\pgfqpoint{2.325000in}{2.310000in}}%
\pgfusepath{clip}%
\pgfsetrectcap%
\pgfsetroundjoin%
\pgfsetlinewidth{0.250937pt}%
\definecolor{currentstroke}{rgb}{0.000000,0.000000,0.000000}%
\pgfsetstrokecolor{currentstroke}%
\pgfsetstrokeopacity{0.200000}%
\pgfsetdash{}{0pt}%
\pgfpathmoveto{\pgfqpoint{0.721913in}{1.460880in}}%
\pgfpathlineto{\pgfqpoint{3.046913in}{1.460880in}}%
\pgfusepath{stroke}%
\end{pgfscope}%
\begin{pgfscope}%
\pgfsetbuttcap%
\pgfsetroundjoin%
\definecolor{currentfill}{rgb}{0.000000,0.000000,0.000000}%
\pgfsetfillcolor{currentfill}%
\pgfsetlinewidth{0.803000pt}%
\definecolor{currentstroke}{rgb}{0.000000,0.000000,0.000000}%
\pgfsetstrokecolor{currentstroke}%
\pgfsetdash{}{0pt}%
\pgfsys@defobject{currentmarker}{\pgfqpoint{-0.048611in}{0.000000in}}{\pgfqpoint{-0.000000in}{0.000000in}}{%
\pgfpathmoveto{\pgfqpoint{-0.000000in}{0.000000in}}%
\pgfpathlineto{\pgfqpoint{-0.048611in}{0.000000in}}%
\pgfusepath{stroke,fill}%
}%
\begin{pgfscope}%
\pgfsys@transformshift{0.721913in}{1.460880in}%
\pgfsys@useobject{currentmarker}{}%
\end{pgfscope}%
\end{pgfscope}%
\begin{pgfscope}%
\definecolor{textcolor}{rgb}{0.000000,0.000000,0.000000}%
\pgfsetstrokecolor{textcolor}%
\pgfsetfillcolor{textcolor}%
\pgftext[x=0.303703in, y=1.403010in, left, base]{\color{textcolor}{\rmfamily\fontsize{12.000000}{14.400000}\selectfont\catcode`\^=\active\def^{\ifmmode\sp\else\^{}\fi}\catcode`\%=\active\def
\end{pgfscope}%
\begin{pgfscope}%
\pgfpathrectangle{\pgfqpoint{0.721913in}{0.549073in}}{\pgfqpoint{2.325000in}{2.310000in}}%
\pgfusepath{clip}%
\pgfsetrectcap%
\pgfsetroundjoin%
\pgfsetlinewidth{0.250937pt}%
\definecolor{currentstroke}{rgb}{0.000000,0.000000,0.000000}%
\pgfsetstrokecolor{currentstroke}%
\pgfsetstrokeopacity{0.200000}%
\pgfsetdash{}{0pt}%
\pgfpathmoveto{\pgfqpoint{0.721913in}{2.065714in}}%
\pgfpathlineto{\pgfqpoint{3.046913in}{2.065714in}}%
\pgfusepath{stroke}%
\end{pgfscope}%
\begin{pgfscope}%
\pgfsetbuttcap%
\pgfsetroundjoin%
\definecolor{currentfill}{rgb}{0.000000,0.000000,0.000000}%
\pgfsetfillcolor{currentfill}%
\pgfsetlinewidth{0.803000pt}%
\definecolor{currentstroke}{rgb}{0.000000,0.000000,0.000000}%
\pgfsetstrokecolor{currentstroke}%
\pgfsetdash{}{0pt}%
\pgfsys@defobject{currentmarker}{\pgfqpoint{-0.048611in}{0.000000in}}{\pgfqpoint{-0.000000in}{0.000000in}}{%
\pgfpathmoveto{\pgfqpoint{-0.000000in}{0.000000in}}%
\pgfpathlineto{\pgfqpoint{-0.048611in}{0.000000in}}%
\pgfusepath{stroke,fill}%
}%
\begin{pgfscope}%
\pgfsys@transformshift{0.721913in}{2.065714in}%
\pgfsys@useobject{currentmarker}{}%
\end{pgfscope}%
\end{pgfscope}%
\begin{pgfscope}%
\definecolor{textcolor}{rgb}{0.000000,0.000000,0.000000}%
\pgfsetstrokecolor{textcolor}%
\pgfsetfillcolor{textcolor}%
\pgftext[x=0.303703in, y=2.007844in, left, base]{\color{textcolor}{\rmfamily\fontsize{12.000000}{14.400000}\selectfont\catcode`\^=\active\def^{\ifmmode\sp\else\^{}\fi}\catcode`\%=\active\def
\end{pgfscope}%
\begin{pgfscope}%
\pgfpathrectangle{\pgfqpoint{0.721913in}{0.549073in}}{\pgfqpoint{2.325000in}{2.310000in}}%
\pgfusepath{clip}%
\pgfsetrectcap%
\pgfsetroundjoin%
\pgfsetlinewidth{0.250937pt}%
\definecolor{currentstroke}{rgb}{0.000000,0.000000,0.000000}%
\pgfsetstrokecolor{currentstroke}%
\pgfsetstrokeopacity{0.200000}%
\pgfsetdash{}{0pt}%
\pgfpathmoveto{\pgfqpoint{0.721913in}{2.670548in}}%
\pgfpathlineto{\pgfqpoint{3.046913in}{2.670548in}}%
\pgfusepath{stroke}%
\end{pgfscope}%
\begin{pgfscope}%
\pgfsetbuttcap%
\pgfsetroundjoin%
\definecolor{currentfill}{rgb}{0.000000,0.000000,0.000000}%
\pgfsetfillcolor{currentfill}%
\pgfsetlinewidth{0.803000pt}%
\definecolor{currentstroke}{rgb}{0.000000,0.000000,0.000000}%
\pgfsetstrokecolor{currentstroke}%
\pgfsetdash{}{0pt}%
\pgfsys@defobject{currentmarker}{\pgfqpoint{-0.048611in}{0.000000in}}{\pgfqpoint{-0.000000in}{0.000000in}}{%
\pgfpathmoveto{\pgfqpoint{-0.000000in}{0.000000in}}%
\pgfpathlineto{\pgfqpoint{-0.048611in}{0.000000in}}%
\pgfusepath{stroke,fill}%
}%
\begin{pgfscope}%
\pgfsys@transformshift{0.721913in}{2.670548in}%
\pgfsys@useobject{currentmarker}{}%
\end{pgfscope}%
\end{pgfscope}%
\begin{pgfscope}%
\definecolor{textcolor}{rgb}{0.000000,0.000000,0.000000}%
\pgfsetstrokecolor{textcolor}%
\pgfsetfillcolor{textcolor}%
\pgftext[x=0.395525in, y=2.612678in, left, base]{\color{textcolor}{\rmfamily\fontsize{12.000000}{14.400000}\selectfont\catcode`\^=\active\def^{\ifmmode\sp\else\^{}\fi}\catcode`\%=\active\def
\end{pgfscope}%
\begin{pgfscope}%
\definecolor{textcolor}{rgb}{0.000000,0.000000,0.000000}%
\pgfsetstrokecolor{textcolor}%
\pgfsetfillcolor{textcolor}%
\pgftext[x=0.248147in,y=1.704073in,,bottom,rotate=90.000000]{\color{textcolor}{\rmfamily\fontsize{12.000000}{14.400000}\selectfont\catcode`\^=\active\def^{\ifmmode\sp\else\^{}\fi}\catcode`\%=\active\def
\end{pgfscope}%
\begin{pgfscope}%
\pgfpathrectangle{\pgfqpoint{0.721913in}{0.549073in}}{\pgfqpoint{2.325000in}{2.310000in}}%
\pgfusepath{clip}%
\pgfsetrectcap%
\pgfsetroundjoin%
\pgfsetlinewidth{1.505625pt}%
\definecolor{currentstroke}{rgb}{0.392157,0.560784,1.000000}%
\pgfsetstrokecolor{currentstroke}%
\pgfsetdash{}{0pt}%
\pgfpathmoveto{\pgfqpoint{0.777271in}{2.204699in}}%
\pgfpathlineto{\pgfqpoint{1.159282in}{2.114072in}}%
\pgfpathlineto{\pgfqpoint{1.524410in}{2.059564in}}%
\pgfpathlineto{\pgfqpoint{1.891095in}{2.010133in}}%
\pgfpathlineto{\pgfqpoint{2.259210in}{1.934435in}}%
\pgfpathlineto{\pgfqpoint{2.625645in}{1.861712in}}%
\pgfpathlineto{\pgfqpoint{2.991556in}{1.804047in}}%
\pgfusepath{stroke}%
\end{pgfscope}%
\begin{pgfscope}%
\pgfpathrectangle{\pgfqpoint{0.721913in}{0.549073in}}{\pgfqpoint{2.325000in}{2.310000in}}%
\pgfusepath{clip}%
\pgfsetbuttcap%
\pgfsetroundjoin%
\definecolor{currentfill}{rgb}{0.392157,0.560784,1.000000}%
\pgfsetfillcolor{currentfill}%
\pgfsetlinewidth{1.003750pt}%
\definecolor{currentstroke}{rgb}{0.392157,0.560784,1.000000}%
\pgfsetstrokecolor{currentstroke}%
\pgfsetdash{}{0pt}%
\pgfsys@defobject{currentmarker}{\pgfqpoint{-0.041667in}{-0.041667in}}{\pgfqpoint{0.041667in}{0.041667in}}{%
\pgfpathmoveto{\pgfqpoint{0.000000in}{-0.041667in}}%
\pgfpathcurveto{\pgfqpoint{0.011050in}{-0.041667in}}{\pgfqpoint{0.021649in}{-0.037276in}}{\pgfqpoint{0.029463in}{-0.029463in}}%
\pgfpathcurveto{\pgfqpoint{0.037276in}{-0.021649in}}{\pgfqpoint{0.041667in}{-0.011050in}}{\pgfqpoint{0.041667in}{0.000000in}}%
\pgfpathcurveto{\pgfqpoint{0.041667in}{0.011050in}}{\pgfqpoint{0.037276in}{0.021649in}}{\pgfqpoint{0.029463in}{0.029463in}}%
\pgfpathcurveto{\pgfqpoint{0.021649in}{0.037276in}}{\pgfqpoint{0.011050in}{0.041667in}}{\pgfqpoint{0.000000in}{0.041667in}}%
\pgfpathcurveto{\pgfqpoint{-0.011050in}{0.041667in}}{\pgfqpoint{-0.021649in}{0.037276in}}{\pgfqpoint{-0.029463in}{0.029463in}}%
\pgfpathcurveto{\pgfqpoint{-0.037276in}{0.021649in}}{\pgfqpoint{-0.041667in}{0.011050in}}{\pgfqpoint{-0.041667in}{0.000000in}}%
\pgfpathcurveto{\pgfqpoint{-0.041667in}{-0.011050in}}{\pgfqpoint{-0.037276in}{-0.021649in}}{\pgfqpoint{-0.029463in}{-0.029463in}}%
\pgfpathcurveto{\pgfqpoint{-0.021649in}{-0.037276in}}{\pgfqpoint{-0.011050in}{-0.041667in}}{\pgfqpoint{0.000000in}{-0.041667in}}%
\pgfpathlineto{\pgfqpoint{0.000000in}{-0.041667in}}%
\pgfpathclose%
\pgfusepath{stroke,fill}%
}%
\begin{pgfscope}%
\pgfsys@transformshift{0.777271in}{2.204699in}%
\pgfsys@useobject{currentmarker}{}%
\end{pgfscope}%
\begin{pgfscope}%
\pgfsys@transformshift{1.159282in}{2.114072in}%
\pgfsys@useobject{currentmarker}{}%
\end{pgfscope}%
\begin{pgfscope}%
\pgfsys@transformshift{1.524410in}{2.059564in}%
\pgfsys@useobject{currentmarker}{}%
\end{pgfscope}%
\begin{pgfscope}%
\pgfsys@transformshift{1.891095in}{2.010133in}%
\pgfsys@useobject{currentmarker}{}%
\end{pgfscope}%
\begin{pgfscope}%
\pgfsys@transformshift{2.259210in}{1.934435in}%
\pgfsys@useobject{currentmarker}{}%
\end{pgfscope}%
\begin{pgfscope}%
\pgfsys@transformshift{2.625645in}{1.861712in}%
\pgfsys@useobject{currentmarker}{}%
\end{pgfscope}%
\begin{pgfscope}%
\pgfsys@transformshift{2.991556in}{1.804047in}%
\pgfsys@useobject{currentmarker}{}%
\end{pgfscope}%
\end{pgfscope}%
\begin{pgfscope}%
\pgfpathrectangle{\pgfqpoint{0.721913in}{0.549073in}}{\pgfqpoint{2.325000in}{2.310000in}}%
\pgfusepath{clip}%
\pgfsetrectcap%
\pgfsetroundjoin%
\pgfsetlinewidth{1.505625pt}%
\definecolor{currentstroke}{rgb}{0.862745,0.149020,0.498039}%
\pgfsetstrokecolor{currentstroke}%
\pgfsetdash{}{0pt}%
\pgfpathmoveto{\pgfqpoint{0.777271in}{2.240870in}}%
\pgfpathlineto{\pgfqpoint{1.159282in}{2.156063in}}%
\pgfpathlineto{\pgfqpoint{1.524410in}{2.106084in}}%
\pgfpathlineto{\pgfqpoint{1.891095in}{2.020666in}}%
\pgfpathlineto{\pgfqpoint{2.259210in}{1.894368in}}%
\pgfpathlineto{\pgfqpoint{2.625645in}{1.657607in}}%
\pgfpathlineto{\pgfqpoint{2.991556in}{0.741573in}}%
\pgfusepath{stroke}%
\end{pgfscope}%
\begin{pgfscope}%
\pgfpathrectangle{\pgfqpoint{0.721913in}{0.549073in}}{\pgfqpoint{2.325000in}{2.310000in}}%
\pgfusepath{clip}%
\pgfsetbuttcap%
\pgfsetmiterjoin%
\definecolor{currentfill}{rgb}{0.862745,0.149020,0.498039}%
\pgfsetfillcolor{currentfill}%
\pgfsetlinewidth{1.003750pt}%
\definecolor{currentstroke}{rgb}{0.862745,0.149020,0.498039}%
\pgfsetstrokecolor{currentstroke}%
\pgfsetdash{}{0pt}%
\pgfsys@defobject{currentmarker}{\pgfqpoint{-0.041667in}{-0.041667in}}{\pgfqpoint{0.041667in}{0.041667in}}{%
\pgfpathmoveto{\pgfqpoint{-0.041667in}{-0.041667in}}%
\pgfpathlineto{\pgfqpoint{0.041667in}{-0.041667in}}%
\pgfpathlineto{\pgfqpoint{0.041667in}{0.041667in}}%
\pgfpathlineto{\pgfqpoint{-0.041667in}{0.041667in}}%
\pgfpathlineto{\pgfqpoint{-0.041667in}{-0.041667in}}%
\pgfpathclose%
\pgfusepath{stroke,fill}%
}%
\begin{pgfscope}%
\pgfsys@transformshift{0.777271in}{2.240870in}%
\pgfsys@useobject{currentmarker}{}%
\end{pgfscope}%
\begin{pgfscope}%
\pgfsys@transformshift{1.159282in}{2.156063in}%
\pgfsys@useobject{currentmarker}{}%
\end{pgfscope}%
\begin{pgfscope}%
\pgfsys@transformshift{1.524410in}{2.106084in}%
\pgfsys@useobject{currentmarker}{}%
\end{pgfscope}%
\begin{pgfscope}%
\pgfsys@transformshift{1.891095in}{2.020666in}%
\pgfsys@useobject{currentmarker}{}%
\end{pgfscope}%
\begin{pgfscope}%
\pgfsys@transformshift{2.259210in}{1.894368in}%
\pgfsys@useobject{currentmarker}{}%
\end{pgfscope}%
\begin{pgfscope}%
\pgfsys@transformshift{2.625645in}{1.657607in}%
\pgfsys@useobject{currentmarker}{}%
\end{pgfscope}%
\begin{pgfscope}%
\pgfsys@transformshift{2.991556in}{0.741573in}%
\pgfsys@useobject{currentmarker}{}%
\end{pgfscope}%
\end{pgfscope}%
\begin{pgfscope}%
\pgfpathrectangle{\pgfqpoint{0.721913in}{0.549073in}}{\pgfqpoint{2.325000in}{2.310000in}}%
\pgfusepath{clip}%
\pgfsetrectcap%
\pgfsetroundjoin%
\pgfsetlinewidth{1.505625pt}%
\definecolor{currentstroke}{rgb}{1.000000,0.690196,0.000000}%
\pgfsetstrokecolor{currentstroke}%
\pgfsetdash{}{0pt}%
\pgfpathmoveto{\pgfqpoint{0.777271in}{2.666573in}}%
\pgfpathlineto{\pgfqpoint{1.159282in}{2.660893in}}%
\pgfpathlineto{\pgfqpoint{1.524410in}{2.645765in}}%
\pgfpathlineto{\pgfqpoint{1.891095in}{2.603867in}}%
\pgfpathlineto{\pgfqpoint{2.259210in}{2.484962in}}%
\pgfpathlineto{\pgfqpoint{2.625645in}{1.897006in}}%
\pgfpathlineto{\pgfqpoint{2.991556in}{1.100817in}}%
\pgfusepath{stroke}%
\end{pgfscope}%
\begin{pgfscope}%
\pgfpathrectangle{\pgfqpoint{0.721913in}{0.549073in}}{\pgfqpoint{2.325000in}{2.310000in}}%
\pgfusepath{clip}%
\pgfsetbuttcap%
\pgfsetmiterjoin%
\definecolor{currentfill}{rgb}{1.000000,0.690196,0.000000}%
\pgfsetfillcolor{currentfill}%
\pgfsetlinewidth{1.003750pt}%
\definecolor{currentstroke}{rgb}{1.000000,0.690196,0.000000}%
\pgfsetstrokecolor{currentstroke}%
\pgfsetdash{}{0pt}%
\pgfsys@defobject{currentmarker}{\pgfqpoint{-0.035355in}{-0.058926in}}{\pgfqpoint{0.035355in}{0.058926in}}{%
\pgfpathmoveto{\pgfqpoint{-0.000000in}{-0.058926in}}%
\pgfpathlineto{\pgfqpoint{0.035355in}{0.000000in}}%
\pgfpathlineto{\pgfqpoint{0.000000in}{0.058926in}}%
\pgfpathlineto{\pgfqpoint{-0.035355in}{0.000000in}}%
\pgfpathlineto{\pgfqpoint{-0.000000in}{-0.058926in}}%
\pgfpathclose%
\pgfusepath{stroke,fill}%
}%
\begin{pgfscope}%
\pgfsys@transformshift{0.777271in}{2.666573in}%
\pgfsys@useobject{currentmarker}{}%
\end{pgfscope}%
\begin{pgfscope}%
\pgfsys@transformshift{1.159282in}{2.660893in}%
\pgfsys@useobject{currentmarker}{}%
\end{pgfscope}%
\begin{pgfscope}%
\pgfsys@transformshift{1.524410in}{2.645765in}%
\pgfsys@useobject{currentmarker}{}%
\end{pgfscope}%
\begin{pgfscope}%
\pgfsys@transformshift{1.891095in}{2.603867in}%
\pgfsys@useobject{currentmarker}{}%
\end{pgfscope}%
\begin{pgfscope}%
\pgfsys@transformshift{2.259210in}{2.484962in}%
\pgfsys@useobject{currentmarker}{}%
\end{pgfscope}%
\begin{pgfscope}%
\pgfsys@transformshift{2.625645in}{1.897006in}%
\pgfsys@useobject{currentmarker}{}%
\end{pgfscope}%
\begin{pgfscope}%
\pgfsys@transformshift{2.991556in}{1.100817in}%
\pgfsys@useobject{currentmarker}{}%
\end{pgfscope}%
\end{pgfscope}%
\begin{pgfscope}%
\pgfsetrectcap%
\pgfsetmiterjoin%
\pgfsetlinewidth{0.803000pt}%
\definecolor{currentstroke}{rgb}{0.000000,0.000000,0.000000}%
\pgfsetstrokecolor{currentstroke}%
\pgfsetdash{}{0pt}%
\pgfpathmoveto{\pgfqpoint{0.721913in}{0.549073in}}%
\pgfpathlineto{\pgfqpoint{0.721913in}{2.859073in}}%
\pgfusepath{stroke}%
\end{pgfscope}%
\begin{pgfscope}%
\pgfsetrectcap%
\pgfsetmiterjoin%
\pgfsetlinewidth{0.803000pt}%
\definecolor{currentstroke}{rgb}{0.000000,0.000000,0.000000}%
\pgfsetstrokecolor{currentstroke}%
\pgfsetdash{}{0pt}%
\pgfpathmoveto{\pgfqpoint{3.046913in}{0.549073in}}%
\pgfpathlineto{\pgfqpoint{3.046913in}{2.859073in}}%
\pgfusepath{stroke}%
\end{pgfscope}%
\begin{pgfscope}%
\pgfsetrectcap%
\pgfsetmiterjoin%
\pgfsetlinewidth{0.803000pt}%
\definecolor{currentstroke}{rgb}{0.000000,0.000000,0.000000}%
\pgfsetstrokecolor{currentstroke}%
\pgfsetdash{}{0pt}%
\pgfpathmoveto{\pgfqpoint{0.721913in}{0.549073in}}%
\pgfpathlineto{\pgfqpoint{3.046913in}{0.549073in}}%
\pgfusepath{stroke}%
\end{pgfscope}%
\begin{pgfscope}%
\pgfsetrectcap%
\pgfsetmiterjoin%
\pgfsetlinewidth{0.803000pt}%
\definecolor{currentstroke}{rgb}{0.000000,0.000000,0.000000}%
\pgfsetstrokecolor{currentstroke}%
\pgfsetdash{}{0pt}%
\pgfpathmoveto{\pgfqpoint{0.721913in}{2.859073in}}%
\pgfpathlineto{\pgfqpoint{3.046913in}{2.859073in}}%
\pgfusepath{stroke}%
\end{pgfscope}%
\begin{pgfscope}%
\pgfsetbuttcap%
\pgfsetmiterjoin%
\definecolor{currentfill}{rgb}{1.000000,1.000000,1.000000}%
\pgfsetfillcolor{currentfill}%
\pgfsetfillopacity{0.800000}%
\pgfsetlinewidth{1.003750pt}%
\definecolor{currentstroke}{rgb}{0.800000,0.800000,0.800000}%
\pgfsetstrokecolor{currentstroke}%
\pgfsetstrokeopacity{0.800000}%
\pgfsetdash{}{0pt}%
\pgfpathmoveto{\pgfqpoint{0.805247in}{0.632406in}}%
\pgfpathlineto{\pgfqpoint{2.037651in}{0.632406in}}%
\pgfpathlineto{\pgfqpoint{2.037651in}{1.379627in}}%
\pgfpathlineto{\pgfqpoint{0.805247in}{1.379627in}}%
\pgfpathlineto{\pgfqpoint{0.805247in}{0.632406in}}%
\pgfpathclose%
\pgfusepath{stroke,fill}%
\end{pgfscope}%
\begin{pgfscope}%
\pgfsetrectcap%
\pgfsetroundjoin%
\pgfsetlinewidth{1.505625pt}%
\definecolor{currentstroke}{rgb}{0.392157,0.560784,1.000000}%
\pgfsetstrokecolor{currentstroke}%
\pgfsetdash{}{0pt}%
\pgfpathmoveto{\pgfqpoint{0.871913in}{1.254627in}}%
\pgfpathlineto{\pgfqpoint{1.038580in}{1.254627in}}%
\pgfpathlineto{\pgfqpoint{1.205247in}{1.254627in}}%
\pgfusepath{stroke}%
\end{pgfscope}%
\begin{pgfscope}%
\pgfsetbuttcap%
\pgfsetroundjoin%
\definecolor{currentfill}{rgb}{0.392157,0.560784,1.000000}%
\pgfsetfillcolor{currentfill}%
\pgfsetlinewidth{1.003750pt}%
\definecolor{currentstroke}{rgb}{0.392157,0.560784,1.000000}%
\pgfsetstrokecolor{currentstroke}%
\pgfsetdash{}{0pt}%
\pgfsys@defobject{currentmarker}{\pgfqpoint{-0.031250in}{-0.031250in}}{\pgfqpoint{0.031250in}{0.031250in}}{%
\pgfpathmoveto{\pgfqpoint{0.000000in}{-0.031250in}}%
\pgfpathcurveto{\pgfqpoint{0.008288in}{-0.031250in}}{\pgfqpoint{0.016237in}{-0.027957in}}{\pgfqpoint{0.022097in}{-0.022097in}}%
\pgfpathcurveto{\pgfqpoint{0.027957in}{-0.016237in}}{\pgfqpoint{0.031250in}{-0.008288in}}{\pgfqpoint{0.031250in}{0.000000in}}%
\pgfpathcurveto{\pgfqpoint{0.031250in}{0.008288in}}{\pgfqpoint{0.027957in}{0.016237in}}{\pgfqpoint{0.022097in}{0.022097in}}%
\pgfpathcurveto{\pgfqpoint{0.016237in}{0.027957in}}{\pgfqpoint{0.008288in}{0.031250in}}{\pgfqpoint{0.000000in}{0.031250in}}%
\pgfpathcurveto{\pgfqpoint{-0.008288in}{0.031250in}}{\pgfqpoint{-0.016237in}{0.027957in}}{\pgfqpoint{-0.022097in}{0.022097in}}%
\pgfpathcurveto{\pgfqpoint{-0.027957in}{0.016237in}}{\pgfqpoint{-0.031250in}{0.008288in}}{\pgfqpoint{-0.031250in}{0.000000in}}%
\pgfpathcurveto{\pgfqpoint{-0.031250in}{-0.008288in}}{\pgfqpoint{-0.027957in}{-0.016237in}}{\pgfqpoint{-0.022097in}{-0.022097in}}%
\pgfpathcurveto{\pgfqpoint{-0.016237in}{-0.027957in}}{\pgfqpoint{-0.008288in}{-0.031250in}}{\pgfqpoint{0.000000in}{-0.031250in}}%
\pgfpathlineto{\pgfqpoint{0.000000in}{-0.031250in}}%
\pgfpathclose%
\pgfusepath{stroke,fill}%
}%
\begin{pgfscope}%
\pgfsys@transformshift{1.038580in}{1.254627in}%
\pgfsys@useobject{currentmarker}{}%
\end{pgfscope}%
\end{pgfscope}%
\begin{pgfscope}%
\definecolor{textcolor}{rgb}{0.000000,0.000000,0.000000}%
\pgfsetstrokecolor{textcolor}%
\pgfsetfillcolor{textcolor}%
\pgftext[x=1.338580in,y=1.196294in,left,base]{\color{textcolor}{\rmfamily\fontsize{12.000000}{14.400000}\selectfont\catcode`\^=\active\def^{\ifmmode\sp\else\^{}\fi}\catcode`\%=\active\def
\end{pgfscope}%
\begin{pgfscope}%
\pgfsetrectcap%
\pgfsetroundjoin%
\pgfsetlinewidth{1.505625pt}%
\definecolor{currentstroke}{rgb}{0.862745,0.149020,0.498039}%
\pgfsetstrokecolor{currentstroke}%
\pgfsetdash{}{0pt}%
\pgfpathmoveto{\pgfqpoint{0.871913in}{1.022220in}}%
\pgfpathlineto{\pgfqpoint{1.038580in}{1.022220in}}%
\pgfpathlineto{\pgfqpoint{1.205247in}{1.022220in}}%
\pgfusepath{stroke}%
\end{pgfscope}%
\begin{pgfscope}%
\pgfsetbuttcap%
\pgfsetmiterjoin%
\definecolor{currentfill}{rgb}{0.862745,0.149020,0.498039}%
\pgfsetfillcolor{currentfill}%
\pgfsetlinewidth{1.003750pt}%
\definecolor{currentstroke}{rgb}{0.862745,0.149020,0.498039}%
\pgfsetstrokecolor{currentstroke}%
\pgfsetdash{}{0pt}%
\pgfsys@defobject{currentmarker}{\pgfqpoint{-0.031250in}{-0.031250in}}{\pgfqpoint{0.031250in}{0.031250in}}{%
\pgfpathmoveto{\pgfqpoint{-0.031250in}{-0.031250in}}%
\pgfpathlineto{\pgfqpoint{0.031250in}{-0.031250in}}%
\pgfpathlineto{\pgfqpoint{0.031250in}{0.031250in}}%
\pgfpathlineto{\pgfqpoint{-0.031250in}{0.031250in}}%
\pgfpathlineto{\pgfqpoint{-0.031250in}{-0.031250in}}%
\pgfpathclose%
\pgfusepath{stroke,fill}%
}%
\begin{pgfscope}%
\pgfsys@transformshift{1.038580in}{1.022220in}%
\pgfsys@useobject{currentmarker}{}%
\end{pgfscope}%
\end{pgfscope}%
\begin{pgfscope}%
\definecolor{textcolor}{rgb}{0.000000,0.000000,0.000000}%
\pgfsetstrokecolor{textcolor}%
\pgfsetfillcolor{textcolor}%
\pgftext[x=1.338580in,y=0.963887in,left,base]{\color{textcolor}{\rmfamily\fontsize{12.000000}{14.400000}\selectfont\catcode`\^=\active\def^{\ifmmode\sp\else\^{}\fi}\catcode`\%=\active\def
\end{pgfscope}%
\begin{pgfscope}%
\pgfsetrectcap%
\pgfsetroundjoin%
\pgfsetlinewidth{1.505625pt}%
\definecolor{currentstroke}{rgb}{1.000000,0.690196,0.000000}%
\pgfsetstrokecolor{currentstroke}%
\pgfsetdash{}{0pt}%
\pgfpathmoveto{\pgfqpoint{0.871913in}{0.789813in}}%
\pgfpathlineto{\pgfqpoint{1.038580in}{0.789813in}}%
\pgfpathlineto{\pgfqpoint{1.205247in}{0.789813in}}%
\pgfusepath{stroke}%
\end{pgfscope}%
\begin{pgfscope}%
\pgfsetbuttcap%
\pgfsetmiterjoin%
\definecolor{currentfill}{rgb}{1.000000,0.690196,0.000000}%
\pgfsetfillcolor{currentfill}%
\pgfsetlinewidth{1.003750pt}%
\definecolor{currentstroke}{rgb}{1.000000,0.690196,0.000000}%
\pgfsetstrokecolor{currentstroke}%
\pgfsetdash{}{0pt}%
\pgfsys@defobject{currentmarker}{\pgfqpoint{-0.026517in}{-0.044194in}}{\pgfqpoint{0.026517in}{0.044194in}}{%
\pgfpathmoveto{\pgfqpoint{-0.000000in}{-0.044194in}}%
\pgfpathlineto{\pgfqpoint{0.026517in}{0.000000in}}%
\pgfpathlineto{\pgfqpoint{0.000000in}{0.044194in}}%
\pgfpathlineto{\pgfqpoint{-0.026517in}{0.000000in}}%
\pgfpathlineto{\pgfqpoint{-0.000000in}{-0.044194in}}%
\pgfpathclose%
\pgfusepath{stroke,fill}%
}%
\begin{pgfscope}%
\pgfsys@transformshift{1.038580in}{0.789813in}%
\pgfsys@useobject{currentmarker}{}%
\end{pgfscope}%
\end{pgfscope}%
\begin{pgfscope}%
\definecolor{textcolor}{rgb}{0.000000,0.000000,0.000000}%
\pgfsetstrokecolor{textcolor}%
\pgfsetfillcolor{textcolor}%
\pgftext[x=1.338580in,y=0.731480in,left,base]{\color{textcolor}{\rmfamily\fontsize{12.000000}{14.400000}\selectfont\catcode`\^=\active\def^{\ifmmode\sp\else\^{}\fi}\catcode`\%=\active\def
\end{pgfscope}%
\end{pgfpicture}%
\makeatother%
\endgroup%

%% file: tables/krylov_aware_density_KA.tex
\centering
\renewcommand{\arraystretch}{1.2}
\begin{tabular}{@{}lccccc@{}}
\toprule
 & $n_{\mtx{\Omega}}$ & $n_{\mtx{\Psi}}$ & $q$ & $r$ & time (s)\\
\midrule
KA (I) & $10$ & $60$ & $20$ & $10$ & $3.42$ \\
KA (II) & $20$ & $30$ & $20$ & $10$ & $2.12$ \\
KA (III) & $10$ & $30$ & $40$ & $10$ & $3.92$ \\
KA (IV) & $10$ & $30$ & $20$ & $20$ & $2.20$ \\
\bottomrule
\end{tabular}

%% file: tables/krylov_aware_density_CN.tex
\centering
\renewcommand{\arraystretch}{1.2}
\begin{tabular}{@{}lcccc@{}}
\toprule
 & $n_{\mtx{\Omega}}$ & $n_{\mtx{\Psi}}$ & $m$ & time (s)\\
\midrule
CN++ & $20$ & $20$ & $2000$ & $3.77$ \\
\bottomrule
\end{tabular}

%% file: plots/hessian_density.pgf
\begingroup%
\makeatletter%
\begin{pgfpicture}%
\pgfpathrectangle{\pgfpointorigin}{\pgfqpoint{3.052968in}{2.959073in}}%
\pgfusepath{use as bounding box, clip}%
\begin{pgfscope}%
\pgfsetbuttcap%
\pgfsetmiterjoin%
\definecolor{currentfill}{rgb}{1.000000,1.000000,1.000000}%
\pgfsetfillcolor{currentfill}%
\pgfsetlinewidth{0.000000pt}%
\definecolor{currentstroke}{rgb}{1.000000,1.000000,1.000000}%
\pgfsetstrokecolor{currentstroke}%
\pgfsetdash{}{0pt}%
\pgfpathmoveto{\pgfqpoint{0.000000in}{-0.000000in}}%
\pgfpathlineto{\pgfqpoint{3.052968in}{-0.000000in}}%
\pgfpathlineto{\pgfqpoint{3.052968in}{2.959073in}}%
\pgfpathlineto{\pgfqpoint{0.000000in}{2.959073in}}%
\pgfpathlineto{\pgfqpoint{0.000000in}{-0.000000in}}%
\pgfpathclose%
\pgfusepath{fill}%
\end{pgfscope}%
\begin{pgfscope}%
\pgfsetbuttcap%
\pgfsetmiterjoin%
\definecolor{currentfill}{rgb}{1.000000,1.000000,1.000000}%
\pgfsetfillcolor{currentfill}%
\pgfsetlinewidth{0.000000pt}%
\definecolor{currentstroke}{rgb}{0.000000,0.000000,0.000000}%
\pgfsetstrokecolor{currentstroke}%
\pgfsetstrokeopacity{0.000000}%
\pgfsetdash{}{0pt}%
\pgfpathmoveto{\pgfqpoint{0.627968in}{0.549073in}}%
\pgfpathlineto{\pgfqpoint{2.952968in}{0.549073in}}%
\pgfpathlineto{\pgfqpoint{2.952968in}{2.859073in}}%
\pgfpathlineto{\pgfqpoint{0.627968in}{2.859073in}}%
\pgfpathlineto{\pgfqpoint{0.627968in}{0.549073in}}%
\pgfpathclose%
\pgfusepath{fill}%
\end{pgfscope}%
\begin{pgfscope}%
\pgfpathrectangle{\pgfqpoint{0.627968in}{0.549073in}}{\pgfqpoint{2.325000in}{2.310000in}}%
\pgfusepath{clip}%
\pgfsetrectcap%
\pgfsetroundjoin%
\pgfsetlinewidth{0.250937pt}%
\definecolor{currentstroke}{rgb}{0.000000,0.000000,0.000000}%
\pgfsetstrokecolor{currentstroke}%
\pgfsetstrokeopacity{0.200000}%
\pgfsetdash{}{0pt}%
\pgfpathmoveto{\pgfqpoint{1.555865in}{0.549073in}}%
\pgfpathlineto{\pgfqpoint{1.555865in}{2.859073in}}%
\pgfusepath{stroke}%
\end{pgfscope}%
\begin{pgfscope}%
\pgfsetbuttcap%
\pgfsetroundjoin%
\definecolor{currentfill}{rgb}{0.000000,0.000000,0.000000}%
\pgfsetfillcolor{currentfill}%
\pgfsetlinewidth{0.803000pt}%
\definecolor{currentstroke}{rgb}{0.000000,0.000000,0.000000}%
\pgfsetstrokecolor{currentstroke}%
\pgfsetdash{}{0pt}%
\pgfsys@defobject{currentmarker}{\pgfqpoint{0.000000in}{-0.048611in}}{\pgfqpoint{0.000000in}{0.000000in}}{%
\pgfpathmoveto{\pgfqpoint{0.000000in}{0.000000in}}%
\pgfpathlineto{\pgfqpoint{0.000000in}{-0.048611in}}%
\pgfusepath{stroke,fill}%
}%
\begin{pgfscope}%
\pgfsys@transformshift{1.555865in}{0.549073in}%
\pgfsys@useobject{currentmarker}{}%
\end{pgfscope}%
\end{pgfscope}%
\begin{pgfscope}%
\definecolor{textcolor}{rgb}{0.000000,0.000000,0.000000}%
\pgfsetstrokecolor{textcolor}%
\pgfsetfillcolor{textcolor}%
\pgftext[x=1.555865in,y=0.451851in,,top]{\color{textcolor}{\rmfamily\fontsize{12.000000}{14.400000}\selectfont\catcode`\^=\active\def^{\ifmmode\sp\else\^{}\fi}\catcode`\%=\active\def
\end{pgfscope}%
\begin{pgfscope}%
\pgfpathrectangle{\pgfqpoint{0.627968in}{0.549073in}}{\pgfqpoint{2.325000in}{2.310000in}}%
\pgfusepath{clip}%
\pgfsetrectcap%
\pgfsetroundjoin%
\pgfsetlinewidth{0.250937pt}%
\definecolor{currentstroke}{rgb}{0.000000,0.000000,0.000000}%
\pgfsetstrokecolor{currentstroke}%
\pgfsetstrokeopacity{0.200000}%
\pgfsetdash{}{0pt}%
\pgfpathmoveto{\pgfqpoint{2.544946in}{0.549073in}}%
\pgfpathlineto{\pgfqpoint{2.544946in}{2.859073in}}%
\pgfusepath{stroke}%
\end{pgfscope}%
\begin{pgfscope}%
\pgfsetbuttcap%
\pgfsetroundjoin%
\definecolor{currentfill}{rgb}{0.000000,0.000000,0.000000}%
\pgfsetfillcolor{currentfill}%
\pgfsetlinewidth{0.803000pt}%
\definecolor{currentstroke}{rgb}{0.000000,0.000000,0.000000}%
\pgfsetstrokecolor{currentstroke}%
\pgfsetdash{}{0pt}%
\pgfsys@defobject{currentmarker}{\pgfqpoint{0.000000in}{-0.048611in}}{\pgfqpoint{0.000000in}{0.000000in}}{%
\pgfpathmoveto{\pgfqpoint{0.000000in}{0.000000in}}%
\pgfpathlineto{\pgfqpoint{0.000000in}{-0.048611in}}%
\pgfusepath{stroke,fill}%
}%
\begin{pgfscope}%
\pgfsys@transformshift{2.544946in}{0.549073in}%
\pgfsys@useobject{currentmarker}{}%
\end{pgfscope}%
\end{pgfscope}%
\begin{pgfscope}%
\definecolor{textcolor}{rgb}{0.000000,0.000000,0.000000}%
\pgfsetstrokecolor{textcolor}%
\pgfsetfillcolor{textcolor}%
\pgftext[x=2.544946in,y=0.451851in,,top]{\color{textcolor}{\rmfamily\fontsize{12.000000}{14.400000}\selectfont\catcode`\^=\active\def^{\ifmmode\sp\else\^{}\fi}\catcode`\%=\active\def
\end{pgfscope}%
\begin{pgfscope}%
\definecolor{textcolor}{rgb}{0.000000,0.000000,0.000000}%
\pgfsetstrokecolor{textcolor}%
\pgfsetfillcolor{textcolor}%
\pgftext[x=1.790468in,y=0.248148in,,top]{\color{textcolor}{\rmfamily\fontsize{12.000000}{14.400000}\selectfont\catcode`\^=\active\def^{\ifmmode\sp\else\^{}\fi}\catcode`\%=\active\def
\end{pgfscope}%
\begin{pgfscope}%
\pgfpathrectangle{\pgfqpoint{0.627968in}{0.549073in}}{\pgfqpoint{2.325000in}{2.310000in}}%
\pgfusepath{clip}%
\pgfsetrectcap%
\pgfsetroundjoin%
\pgfsetlinewidth{0.250937pt}%
\definecolor{currentstroke}{rgb}{0.000000,0.000000,0.000000}%
\pgfsetstrokecolor{currentstroke}%
\pgfsetstrokeopacity{0.200000}%
\pgfsetdash{}{0pt}%
\pgfpathmoveto{\pgfqpoint{0.627968in}{0.741573in}}%
\pgfpathlineto{\pgfqpoint{2.952968in}{0.741573in}}%
\pgfusepath{stroke}%
\end{pgfscope}%
\begin{pgfscope}%
\pgfsetbuttcap%
\pgfsetroundjoin%
\definecolor{currentfill}{rgb}{0.000000,0.000000,0.000000}%
\pgfsetfillcolor{currentfill}%
\pgfsetlinewidth{0.803000pt}%
\definecolor{currentstroke}{rgb}{0.000000,0.000000,0.000000}%
\pgfsetstrokecolor{currentstroke}%
\pgfsetdash{}{0pt}%
\pgfsys@defobject{currentmarker}{\pgfqpoint{-0.048611in}{0.000000in}}{\pgfqpoint{-0.000000in}{0.000000in}}{%
\pgfpathmoveto{\pgfqpoint{-0.000000in}{0.000000in}}%
\pgfpathlineto{\pgfqpoint{-0.048611in}{0.000000in}}%
\pgfusepath{stroke,fill}%
}%
\begin{pgfscope}%
\pgfsys@transformshift{0.627968in}{0.741573in}%
\pgfsys@useobject{currentmarker}{}%
\end{pgfscope}%
\end{pgfscope}%
\begin{pgfscope}%
\definecolor{textcolor}{rgb}{0.000000,0.000000,0.000000}%
\pgfsetstrokecolor{textcolor}%
\pgfsetfillcolor{textcolor}%
\pgftext[x=0.322222in, y=0.683703in, left, base]{\color{textcolor}{\rmfamily\fontsize{12.000000}{14.400000}\selectfont\catcode`\^=\active\def^{\ifmmode\sp\else\^{}\fi}\catcode`\%=\active\def
\end{pgfscope}%
\begin{pgfscope}%
\pgfpathrectangle{\pgfqpoint{0.627968in}{0.549073in}}{\pgfqpoint{2.325000in}{2.310000in}}%
\pgfusepath{clip}%
\pgfsetrectcap%
\pgfsetroundjoin%
\pgfsetlinewidth{0.250937pt}%
\definecolor{currentstroke}{rgb}{0.000000,0.000000,0.000000}%
\pgfsetstrokecolor{currentstroke}%
\pgfsetstrokeopacity{0.200000}%
\pgfsetdash{}{0pt}%
\pgfpathmoveto{\pgfqpoint{0.627968in}{1.497651in}}%
\pgfpathlineto{\pgfqpoint{2.952968in}{1.497651in}}%
\pgfusepath{stroke}%
\end{pgfscope}%
\begin{pgfscope}%
\pgfsetbuttcap%
\pgfsetroundjoin%
\definecolor{currentfill}{rgb}{0.000000,0.000000,0.000000}%
\pgfsetfillcolor{currentfill}%
\pgfsetlinewidth{0.803000pt}%
\definecolor{currentstroke}{rgb}{0.000000,0.000000,0.000000}%
\pgfsetstrokecolor{currentstroke}%
\pgfsetdash{}{0pt}%
\pgfsys@defobject{currentmarker}{\pgfqpoint{-0.048611in}{0.000000in}}{\pgfqpoint{-0.000000in}{0.000000in}}{%
\pgfpathmoveto{\pgfqpoint{-0.000000in}{0.000000in}}%
\pgfpathlineto{\pgfqpoint{-0.048611in}{0.000000in}}%
\pgfusepath{stroke,fill}%
}%
\begin{pgfscope}%
\pgfsys@transformshift{0.627968in}{1.497651in}%
\pgfsys@useobject{currentmarker}{}%
\end{pgfscope}%
\end{pgfscope}%
\begin{pgfscope}%
\definecolor{textcolor}{rgb}{0.000000,0.000000,0.000000}%
\pgfsetstrokecolor{textcolor}%
\pgfsetfillcolor{textcolor}%
\pgftext[x=0.322222in, y=1.439780in, left, base]{\color{textcolor}{\rmfamily\fontsize{12.000000}{14.400000}\selectfont\catcode`\^=\active\def^{\ifmmode\sp\else\^{}\fi}\catcode`\%=\active\def
\end{pgfscope}%
\begin{pgfscope}%
\pgfpathrectangle{\pgfqpoint{0.627968in}{0.549073in}}{\pgfqpoint{2.325000in}{2.310000in}}%
\pgfusepath{clip}%
\pgfsetrectcap%
\pgfsetroundjoin%
\pgfsetlinewidth{0.250937pt}%
\definecolor{currentstroke}{rgb}{0.000000,0.000000,0.000000}%
\pgfsetstrokecolor{currentstroke}%
\pgfsetstrokeopacity{0.200000}%
\pgfsetdash{}{0pt}%
\pgfpathmoveto{\pgfqpoint{0.627968in}{2.253729in}}%
\pgfpathlineto{\pgfqpoint{2.952968in}{2.253729in}}%
\pgfusepath{stroke}%
\end{pgfscope}%
\begin{pgfscope}%
\pgfsetbuttcap%
\pgfsetroundjoin%
\definecolor{currentfill}{rgb}{0.000000,0.000000,0.000000}%
\pgfsetfillcolor{currentfill}%
\pgfsetlinewidth{0.803000pt}%
\definecolor{currentstroke}{rgb}{0.000000,0.000000,0.000000}%
\pgfsetstrokecolor{currentstroke}%
\pgfsetdash{}{0pt}%
\pgfsys@defobject{currentmarker}{\pgfqpoint{-0.048611in}{0.000000in}}{\pgfqpoint{-0.000000in}{0.000000in}}{%
\pgfpathmoveto{\pgfqpoint{-0.000000in}{0.000000in}}%
\pgfpathlineto{\pgfqpoint{-0.048611in}{0.000000in}}%
\pgfusepath{stroke,fill}%
}%
\begin{pgfscope}%
\pgfsys@transformshift{0.627968in}{2.253729in}%
\pgfsys@useobject{currentmarker}{}%
\end{pgfscope}%
\end{pgfscope}%
\begin{pgfscope}%
\definecolor{textcolor}{rgb}{0.000000,0.000000,0.000000}%
\pgfsetstrokecolor{textcolor}%
\pgfsetfillcolor{textcolor}%
\pgftext[x=0.322222in, y=2.195858in, left, base]{\color{textcolor}{\rmfamily\fontsize{12.000000}{14.400000}\selectfont\catcode`\^=\active\def^{\ifmmode\sp\else\^{}\fi}\catcode`\%=\active\def
\end{pgfscope}%
\begin{pgfscope}%
\definecolor{textcolor}{rgb}{0.000000,0.000000,0.000000}%
\pgfsetstrokecolor{textcolor}%
\pgfsetfillcolor{textcolor}%
\pgftext[x=0.266667in,y=1.704073in,,bottom,rotate=90.000000]{\color{textcolor}{\rmfamily\fontsize{12.000000}{14.400000}\selectfont\catcode`\^=\active\def^{\ifmmode\sp\else\^{}\fi}\catcode`\%=\active\def
\end{pgfscope}%
\begin{pgfscope}%
\pgfpathrectangle{\pgfqpoint{0.627968in}{0.549073in}}{\pgfqpoint{2.325000in}{2.310000in}}%
\pgfusepath{clip}%
\pgfsetrectcap%
\pgfsetroundjoin%
\pgfsetlinewidth{1.505625pt}%
\definecolor{currentstroke}{rgb}{0.392157,0.560784,1.000000}%
\pgfsetstrokecolor{currentstroke}%
\pgfsetdash{}{0pt}%
\pgfpathmoveto{\pgfqpoint{0.683326in}{2.666573in}}%
\pgfpathlineto{\pgfqpoint{0.698187in}{2.335899in}}%
\pgfpathlineto{\pgfqpoint{0.713048in}{2.103922in}}%
\pgfpathlineto{\pgfqpoint{0.727909in}{2.185464in}}%
\pgfpathlineto{\pgfqpoint{0.742770in}{2.059793in}}%
\pgfpathlineto{\pgfqpoint{0.757631in}{1.819076in}}%
\pgfpathlineto{\pgfqpoint{0.772491in}{1.869803in}}%
\pgfpathlineto{\pgfqpoint{0.787352in}{1.599727in}}%
\pgfpathlineto{\pgfqpoint{0.802213in}{1.270806in}}%
\pgfpathlineto{\pgfqpoint{0.817074in}{1.163032in}}%
\pgfpathlineto{\pgfqpoint{0.831935in}{1.113466in}}%
\pgfpathlineto{\pgfqpoint{0.846796in}{0.975844in}}%
\pgfpathlineto{\pgfqpoint{0.861657in}{0.932377in}}%
\pgfpathlineto{\pgfqpoint{0.876518in}{0.972453in}}%
\pgfpathlineto{\pgfqpoint{0.891379in}{0.993572in}}%
\pgfpathlineto{\pgfqpoint{0.906240in}{0.870943in}}%
\pgfpathlineto{\pgfqpoint{0.921101in}{0.797100in}}%
\pgfpathlineto{\pgfqpoint{0.935962in}{0.830977in}}%
\pgfpathlineto{\pgfqpoint{0.950823in}{0.782515in}}%
\pgfpathlineto{\pgfqpoint{0.965684in}{0.771482in}}%
\pgfpathlineto{\pgfqpoint{0.980545in}{0.827135in}}%
\pgfpathlineto{\pgfqpoint{0.995406in}{0.795853in}}%
\pgfpathlineto{\pgfqpoint{1.010267in}{0.772870in}}%
\pgfpathlineto{\pgfqpoint{1.025128in}{0.791679in}}%
\pgfpathlineto{\pgfqpoint{1.039989in}{0.900218in}}%
\pgfpathlineto{\pgfqpoint{1.054850in}{0.864255in}}%
\pgfpathlineto{\pgfqpoint{1.069711in}{0.811674in}}%
\pgfpathlineto{\pgfqpoint{1.084572in}{0.829462in}}%
\pgfpathlineto{\pgfqpoint{1.099433in}{0.775016in}}%
\pgfpathlineto{\pgfqpoint{1.114294in}{0.784884in}}%
\pgfpathlineto{\pgfqpoint{1.129155in}{0.834565in}}%
\pgfpathlineto{\pgfqpoint{1.144016in}{0.822284in}}%
\pgfpathlineto{\pgfqpoint{1.158877in}{0.842125in}}%
\pgfpathlineto{\pgfqpoint{1.173738in}{0.838038in}}%
\pgfpathlineto{\pgfqpoint{1.188599in}{0.826957in}}%
\pgfpathlineto{\pgfqpoint{1.203460in}{0.767723in}}%
\pgfpathlineto{\pgfqpoint{1.218321in}{0.784932in}}%
\pgfpathlineto{\pgfqpoint{1.233182in}{0.830818in}}%
\pgfpathlineto{\pgfqpoint{1.248043in}{0.790423in}}%
\pgfpathlineto{\pgfqpoint{1.262904in}{0.813273in}}%
\pgfpathlineto{\pgfqpoint{1.277765in}{0.824743in}}%
\pgfpathlineto{\pgfqpoint{1.292626in}{0.812564in}}%
\pgfpathlineto{\pgfqpoint{1.307487in}{0.827805in}}%
\pgfpathlineto{\pgfqpoint{1.322348in}{0.768098in}}%
\pgfpathlineto{\pgfqpoint{1.337209in}{0.785723in}}%
\pgfpathlineto{\pgfqpoint{1.352070in}{0.743635in}}%
\pgfpathlineto{\pgfqpoint{1.366931in}{0.771681in}}%
\pgfpathlineto{\pgfqpoint{1.381792in}{0.828818in}}%
\pgfpathlineto{\pgfqpoint{1.396653in}{0.791302in}}%
\pgfpathlineto{\pgfqpoint{1.411514in}{0.747149in}}%
\pgfpathlineto{\pgfqpoint{1.426375in}{1.131801in}}%
\pgfpathlineto{\pgfqpoint{1.441235in}{0.741573in}}%
\pgfpathlineto{\pgfqpoint{2.585531in}{0.741573in}}%
\pgfpathlineto{\pgfqpoint{2.600392in}{0.743545in}}%
\pgfpathlineto{\pgfqpoint{2.615253in}{0.771116in}}%
\pgfpathlineto{\pgfqpoint{2.630114in}{0.828591in}}%
\pgfpathlineto{\pgfqpoint{2.644975in}{0.791990in}}%
\pgfpathlineto{\pgfqpoint{2.659836in}{0.747318in}}%
\pgfpathlineto{\pgfqpoint{2.674697in}{1.383691in}}%
\pgfpathlineto{\pgfqpoint{2.689558in}{0.741573in}}%
\pgfpathlineto{\pgfqpoint{2.897611in}{0.741573in}}%
\pgfpathlineto{\pgfqpoint{2.897611in}{0.741573in}}%
\pgfusepath{stroke}%
\end{pgfscope}%
\begin{pgfscope}%
\pgfpathrectangle{\pgfqpoint{0.627968in}{0.549073in}}{\pgfqpoint{2.325000in}{2.310000in}}%
\pgfusepath{clip}%
\pgfsetrectcap%
\pgfsetroundjoin%
\pgfsetlinewidth{1.505625pt}%
\definecolor{currentstroke}{rgb}{0.470588,0.368627,0.941176}%
\pgfsetstrokecolor{currentstroke}%
\pgfsetdash{}{0pt}%
\pgfpathmoveto{\pgfqpoint{0.683326in}{1.519977in}}%
\pgfpathlineto{\pgfqpoint{0.698187in}{1.328463in}}%
\pgfpathlineto{\pgfqpoint{0.713048in}{1.367577in}}%
\pgfpathlineto{\pgfqpoint{0.727909in}{1.280477in}}%
\pgfpathlineto{\pgfqpoint{0.742770in}{0.963614in}}%
\pgfpathlineto{\pgfqpoint{0.757631in}{0.893562in}}%
\pgfpathlineto{\pgfqpoint{0.772491in}{0.895132in}}%
\pgfpathlineto{\pgfqpoint{0.787352in}{0.864968in}}%
\pgfpathlineto{\pgfqpoint{0.802213in}{0.806503in}}%
\pgfpathlineto{\pgfqpoint{0.817074in}{0.762134in}}%
\pgfpathlineto{\pgfqpoint{0.831935in}{0.811090in}}%
\pgfpathlineto{\pgfqpoint{0.846796in}{0.817339in}}%
\pgfpathlineto{\pgfqpoint{0.861657in}{0.757853in}}%
\pgfpathlineto{\pgfqpoint{0.876518in}{0.771689in}}%
\pgfpathlineto{\pgfqpoint{0.891379in}{0.741573in}}%
\pgfpathlineto{\pgfqpoint{0.950823in}{0.741573in}}%
\pgfpathlineto{\pgfqpoint{0.965684in}{0.753887in}}%
\pgfpathlineto{\pgfqpoint{0.980545in}{0.774787in}}%
\pgfpathlineto{\pgfqpoint{0.995406in}{0.837906in}}%
\pgfpathlineto{\pgfqpoint{1.010267in}{0.846729in}}%
\pgfpathlineto{\pgfqpoint{1.025128in}{0.831392in}}%
\pgfpathlineto{\pgfqpoint{1.039989in}{0.795417in}}%
\pgfpathlineto{\pgfqpoint{1.054850in}{0.864527in}}%
\pgfpathlineto{\pgfqpoint{1.069711in}{0.882862in}}%
\pgfpathlineto{\pgfqpoint{1.084572in}{0.791499in}}%
\pgfpathlineto{\pgfqpoint{1.099433in}{0.749102in}}%
\pgfpathlineto{\pgfqpoint{1.114294in}{0.768202in}}%
\pgfpathlineto{\pgfqpoint{1.129155in}{0.832325in}}%
\pgfpathlineto{\pgfqpoint{1.144016in}{0.844286in}}%
\pgfpathlineto{\pgfqpoint{1.158877in}{0.837864in}}%
\pgfpathlineto{\pgfqpoint{1.173738in}{0.801607in}}%
\pgfpathlineto{\pgfqpoint{1.188599in}{0.830382in}}%
\pgfpathlineto{\pgfqpoint{1.203460in}{0.793143in}}%
\pgfpathlineto{\pgfqpoint{1.218321in}{0.747821in}}%
\pgfpathlineto{\pgfqpoint{1.233182in}{0.746213in}}%
\pgfpathlineto{\pgfqpoint{1.262904in}{0.830011in}}%
\pgfpathlineto{\pgfqpoint{1.277765in}{0.775773in}}%
\pgfpathlineto{\pgfqpoint{1.292626in}{0.744174in}}%
\pgfpathlineto{\pgfqpoint{1.307487in}{0.741573in}}%
\pgfpathlineto{\pgfqpoint{1.589845in}{0.742276in}}%
\pgfpathlineto{\pgfqpoint{1.604706in}{0.758068in}}%
\pgfpathlineto{\pgfqpoint{1.619567in}{0.817589in}}%
\pgfpathlineto{\pgfqpoint{1.634428in}{0.810481in}}%
\pgfpathlineto{\pgfqpoint{1.649289in}{0.753860in}}%
\pgfpathlineto{\pgfqpoint{1.664150in}{0.804141in}}%
\pgfpathlineto{\pgfqpoint{1.679011in}{0.741573in}}%
\pgfpathlineto{\pgfqpoint{2.897611in}{0.741573in}}%
\pgfpathlineto{\pgfqpoint{2.897611in}{0.741573in}}%
\pgfusepath{stroke}%
\end{pgfscope}%
\begin{pgfscope}%
\pgfpathrectangle{\pgfqpoint{0.627968in}{0.549073in}}{\pgfqpoint{2.325000in}{2.310000in}}%
\pgfusepath{clip}%
\pgfsetrectcap%
\pgfsetroundjoin%
\pgfsetlinewidth{1.505625pt}%
\definecolor{currentstroke}{rgb}{0.862745,0.149020,0.498039}%
\pgfsetstrokecolor{currentstroke}%
\pgfsetdash{}{0pt}%
\pgfpathmoveto{\pgfqpoint{0.683326in}{1.581543in}}%
\pgfpathlineto{\pgfqpoint{0.698187in}{1.346106in}}%
\pgfpathlineto{\pgfqpoint{0.727909in}{1.121802in}}%
\pgfpathlineto{\pgfqpoint{0.742770in}{0.922405in}}%
\pgfpathlineto{\pgfqpoint{0.757631in}{0.831439in}}%
\pgfpathlineto{\pgfqpoint{0.772491in}{0.826947in}}%
\pgfpathlineto{\pgfqpoint{0.787352in}{0.812840in}}%
\pgfpathlineto{\pgfqpoint{0.802213in}{0.762332in}}%
\pgfpathlineto{\pgfqpoint{0.817074in}{0.795348in}}%
\pgfpathlineto{\pgfqpoint{0.831935in}{0.832952in}}%
\pgfpathlineto{\pgfqpoint{0.846796in}{0.817841in}}%
\pgfpathlineto{\pgfqpoint{0.861657in}{0.830776in}}%
\pgfpathlineto{\pgfqpoint{0.876518in}{0.772366in}}%
\pgfpathlineto{\pgfqpoint{0.891379in}{0.754603in}}%
\pgfpathlineto{\pgfqpoint{0.906240in}{0.741573in}}%
\pgfpathlineto{\pgfqpoint{0.950823in}{0.741573in}}%
\pgfpathlineto{\pgfqpoint{0.965684in}{0.748368in}}%
\pgfpathlineto{\pgfqpoint{0.980545in}{0.761104in}}%
\pgfpathlineto{\pgfqpoint{0.995406in}{0.861579in}}%
\pgfpathlineto{\pgfqpoint{1.010267in}{0.901627in}}%
\pgfpathlineto{\pgfqpoint{1.025128in}{0.789454in}}%
\pgfpathlineto{\pgfqpoint{1.039989in}{0.768733in}}%
\pgfpathlineto{\pgfqpoint{1.054850in}{0.825722in}}%
\pgfpathlineto{\pgfqpoint{1.084572in}{0.776577in}}%
\pgfpathlineto{\pgfqpoint{1.099433in}{0.842146in}}%
\pgfpathlineto{\pgfqpoint{1.114294in}{0.868824in}}%
\pgfpathlineto{\pgfqpoint{1.129155in}{0.826307in}}%
\pgfpathlineto{\pgfqpoint{1.144016in}{0.807291in}}%
\pgfpathlineto{\pgfqpoint{1.158877in}{0.830454in}}%
\pgfpathlineto{\pgfqpoint{1.173738in}{0.797837in}}%
\pgfpathlineto{\pgfqpoint{1.188599in}{0.829584in}}%
\pgfpathlineto{\pgfqpoint{1.203460in}{0.794342in}}%
\pgfpathlineto{\pgfqpoint{1.218321in}{0.747367in}}%
\pgfpathlineto{\pgfqpoint{1.233182in}{0.822711in}}%
\pgfpathlineto{\pgfqpoint{1.248043in}{0.741573in}}%
\pgfpathlineto{\pgfqpoint{1.634428in}{0.741573in}}%
\pgfpathlineto{\pgfqpoint{1.649289in}{0.816639in}}%
\pgfpathlineto{\pgfqpoint{1.664150in}{0.752369in}}%
\pgfpathlineto{\pgfqpoint{1.679011in}{0.807266in}}%
\pgfpathlineto{\pgfqpoint{1.693872in}{0.820204in}}%
\pgfpathlineto{\pgfqpoint{1.708733in}{0.760086in}}%
\pgfpathlineto{\pgfqpoint{1.723594in}{0.742430in}}%
\pgfpathlineto{\pgfqpoint{1.753316in}{0.741573in}}%
\pgfpathlineto{\pgfqpoint{2.897611in}{0.741573in}}%
\pgfpathlineto{\pgfqpoint{2.897611in}{0.741573in}}%
\pgfusepath{stroke}%
\end{pgfscope}%
\begin{pgfscope}%
\pgfpathrectangle{\pgfqpoint{0.627968in}{0.549073in}}{\pgfqpoint{2.325000in}{2.310000in}}%
\pgfusepath{clip}%
\pgfsetrectcap%
\pgfsetroundjoin%
\pgfsetlinewidth{1.505625pt}%
\definecolor{currentstroke}{rgb}{0.996078,0.380392,0.000000}%
\pgfsetstrokecolor{currentstroke}%
\pgfsetdash{}{0pt}%
\pgfpathmoveto{\pgfqpoint{0.683326in}{1.432125in}}%
\pgfpathlineto{\pgfqpoint{0.698187in}{1.350751in}}%
\pgfpathlineto{\pgfqpoint{0.713048in}{1.247086in}}%
\pgfpathlineto{\pgfqpoint{0.727909in}{1.042743in}}%
\pgfpathlineto{\pgfqpoint{0.742770in}{0.882458in}}%
\pgfpathlineto{\pgfqpoint{0.757631in}{0.843069in}}%
\pgfpathlineto{\pgfqpoint{0.772491in}{0.845930in}}%
\pgfpathlineto{\pgfqpoint{0.787352in}{0.950384in}}%
\pgfpathlineto{\pgfqpoint{0.802213in}{0.972232in}}%
\pgfpathlineto{\pgfqpoint{0.817074in}{0.850237in}}%
\pgfpathlineto{\pgfqpoint{0.831935in}{0.831480in}}%
\pgfpathlineto{\pgfqpoint{0.846796in}{0.770269in}}%
\pgfpathlineto{\pgfqpoint{0.861657in}{0.766646in}}%
\pgfpathlineto{\pgfqpoint{0.876518in}{0.741573in}}%
\pgfpathlineto{\pgfqpoint{0.906240in}{0.741573in}}%
\pgfpathlineto{\pgfqpoint{0.921101in}{0.854486in}}%
\pgfpathlineto{\pgfqpoint{0.935962in}{0.759687in}}%
\pgfpathlineto{\pgfqpoint{0.950823in}{0.820542in}}%
\pgfpathlineto{\pgfqpoint{0.980545in}{0.830787in}}%
\pgfpathlineto{\pgfqpoint{1.010267in}{0.795749in}}%
\pgfpathlineto{\pgfqpoint{1.025128in}{0.830700in}}%
\pgfpathlineto{\pgfqpoint{1.039989in}{0.778492in}}%
\pgfpathlineto{\pgfqpoint{1.054850in}{0.764254in}}%
\pgfpathlineto{\pgfqpoint{1.069711in}{0.821723in}}%
\pgfpathlineto{\pgfqpoint{1.084572in}{0.807099in}}%
\pgfpathlineto{\pgfqpoint{1.099433in}{0.779875in}}%
\pgfpathlineto{\pgfqpoint{1.114294in}{0.828321in}}%
\pgfpathlineto{\pgfqpoint{1.129155in}{0.794108in}}%
\pgfpathlineto{\pgfqpoint{1.144016in}{0.763903in}}%
\pgfpathlineto{\pgfqpoint{1.158877in}{0.817857in}}%
\pgfpathlineto{\pgfqpoint{1.173738in}{0.825048in}}%
\pgfpathlineto{\pgfqpoint{1.188599in}{0.826367in}}%
\pgfpathlineto{\pgfqpoint{1.203460in}{0.814928in}}%
\pgfpathlineto{\pgfqpoint{1.218321in}{0.756045in}}%
\pgfpathlineto{\pgfqpoint{1.233182in}{0.848461in}}%
\pgfpathlineto{\pgfqpoint{1.248043in}{0.741573in}}%
\pgfpathlineto{\pgfqpoint{1.649289in}{0.741573in}}%
\pgfpathlineto{\pgfqpoint{1.664150in}{0.745424in}}%
\pgfpathlineto{\pgfqpoint{1.679011in}{0.783286in}}%
\pgfpathlineto{\pgfqpoint{1.693872in}{0.830435in}}%
\pgfpathlineto{\pgfqpoint{1.708733in}{0.778811in}}%
\pgfpathlineto{\pgfqpoint{1.723594in}{0.744642in}}%
\pgfpathlineto{\pgfqpoint{1.738455in}{0.741573in}}%
\pgfpathlineto{\pgfqpoint{2.897611in}{0.741573in}}%
\pgfpathlineto{\pgfqpoint{2.897611in}{0.741573in}}%
\pgfusepath{stroke}%
\end{pgfscope}%
\begin{pgfscope}%
\pgfpathrectangle{\pgfqpoint{0.627968in}{0.549073in}}{\pgfqpoint{2.325000in}{2.310000in}}%
\pgfusepath{clip}%
\pgfsetrectcap%
\pgfsetroundjoin%
\pgfsetlinewidth{1.505625pt}%
\definecolor{currentstroke}{rgb}{1.000000,0.690196,0.000000}%
\pgfsetstrokecolor{currentstroke}%
\pgfsetdash{}{0pt}%
\pgfpathmoveto{\pgfqpoint{0.683326in}{1.498554in}}%
\pgfpathlineto{\pgfqpoint{0.698187in}{1.250357in}}%
\pgfpathlineto{\pgfqpoint{0.713048in}{1.090641in}}%
\pgfpathlineto{\pgfqpoint{0.727909in}{1.032426in}}%
\pgfpathlineto{\pgfqpoint{0.742770in}{0.892449in}}%
\pgfpathlineto{\pgfqpoint{0.757631in}{0.888617in}}%
\pgfpathlineto{\pgfqpoint{0.772491in}{1.004405in}}%
\pgfpathlineto{\pgfqpoint{0.787352in}{0.867701in}}%
\pgfpathlineto{\pgfqpoint{0.802213in}{0.825145in}}%
\pgfpathlineto{\pgfqpoint{0.817074in}{0.869058in}}%
\pgfpathlineto{\pgfqpoint{0.831935in}{0.846060in}}%
\pgfpathlineto{\pgfqpoint{0.846796in}{0.772179in}}%
\pgfpathlineto{\pgfqpoint{0.861657in}{0.751736in}}%
\pgfpathlineto{\pgfqpoint{0.876518in}{0.768301in}}%
\pgfpathlineto{\pgfqpoint{0.891379in}{0.830187in}}%
\pgfpathlineto{\pgfqpoint{0.906240in}{0.832546in}}%
\pgfpathlineto{\pgfqpoint{0.921101in}{0.838935in}}%
\pgfpathlineto{\pgfqpoint{0.935962in}{0.812028in}}%
\pgfpathlineto{\pgfqpoint{0.950823in}{0.832119in}}%
\pgfpathlineto{\pgfqpoint{0.965684in}{0.801355in}}%
\pgfpathlineto{\pgfqpoint{0.980545in}{0.804957in}}%
\pgfpathlineto{\pgfqpoint{0.995406in}{0.825793in}}%
\pgfpathlineto{\pgfqpoint{1.010267in}{0.765985in}}%
\pgfpathlineto{\pgfqpoint{1.025128in}{0.767009in}}%
\pgfpathlineto{\pgfqpoint{1.039989in}{0.750999in}}%
\pgfpathlineto{\pgfqpoint{1.054850in}{0.803862in}}%
\pgfpathlineto{\pgfqpoint{1.069711in}{0.822634in}}%
\pgfpathlineto{\pgfqpoint{1.084572in}{0.762339in}}%
\pgfpathlineto{\pgfqpoint{1.099433in}{0.782066in}}%
\pgfpathlineto{\pgfqpoint{1.114294in}{0.768880in}}%
\pgfpathlineto{\pgfqpoint{1.129155in}{0.827224in}}%
\pgfpathlineto{\pgfqpoint{1.144016in}{0.795033in}}%
\pgfpathlineto{\pgfqpoint{1.158877in}{0.755964in}}%
\pgfpathlineto{\pgfqpoint{1.173738in}{0.799735in}}%
\pgfpathlineto{\pgfqpoint{1.188599in}{0.827310in}}%
\pgfpathlineto{\pgfqpoint{1.203460in}{0.795621in}}%
\pgfpathlineto{\pgfqpoint{1.218321in}{0.830216in}}%
\pgfpathlineto{\pgfqpoint{1.233182in}{0.791102in}}%
\pgfpathlineto{\pgfqpoint{1.248043in}{0.748020in}}%
\pgfpathlineto{\pgfqpoint{1.262904in}{1.360714in}}%
\pgfpathlineto{\pgfqpoint{1.277765in}{0.741573in}}%
\pgfpathlineto{\pgfqpoint{1.634428in}{0.741573in}}%
\pgfpathlineto{\pgfqpoint{1.649289in}{0.820789in}}%
\pgfpathlineto{\pgfqpoint{1.664150in}{0.755950in}}%
\pgfpathlineto{\pgfqpoint{1.679011in}{0.814329in}}%
\pgfpathlineto{\pgfqpoint{1.693872in}{0.813996in}}%
\pgfpathlineto{\pgfqpoint{1.708733in}{0.755753in}}%
\pgfpathlineto{\pgfqpoint{1.723594in}{0.819994in}}%
\pgfpathlineto{\pgfqpoint{1.738455in}{0.741573in}}%
\pgfpathlineto{\pgfqpoint{2.897611in}{0.741573in}}%
\pgfpathlineto{\pgfqpoint{2.897611in}{0.741573in}}%
\pgfusepath{stroke}%
\end{pgfscope}%
\begin{pgfscope}%
\pgfsetrectcap%
\pgfsetmiterjoin%
\pgfsetlinewidth{0.803000pt}%
\definecolor{currentstroke}{rgb}{0.000000,0.000000,0.000000}%
\pgfsetstrokecolor{currentstroke}%
\pgfsetdash{}{0pt}%
\pgfpathmoveto{\pgfqpoint{0.627968in}{0.549073in}}%
\pgfpathlineto{\pgfqpoint{0.627968in}{2.859073in}}%
\pgfusepath{stroke}%
\end{pgfscope}%
\begin{pgfscope}%
\pgfsetrectcap%
\pgfsetmiterjoin%
\pgfsetlinewidth{0.803000pt}%
\definecolor{currentstroke}{rgb}{0.000000,0.000000,0.000000}%
\pgfsetstrokecolor{currentstroke}%
\pgfsetdash{}{0pt}%
\pgfpathmoveto{\pgfqpoint{2.952968in}{0.549073in}}%
\pgfpathlineto{\pgfqpoint{2.952968in}{2.859073in}}%
\pgfusepath{stroke}%
\end{pgfscope}%
\begin{pgfscope}%
\pgfsetrectcap%
\pgfsetmiterjoin%
\pgfsetlinewidth{0.803000pt}%
\definecolor{currentstroke}{rgb}{0.000000,0.000000,0.000000}%
\pgfsetstrokecolor{currentstroke}%
\pgfsetdash{}{0pt}%
\pgfpathmoveto{\pgfqpoint{0.627968in}{0.549073in}}%
\pgfpathlineto{\pgfqpoint{2.952968in}{0.549073in}}%
\pgfusepath{stroke}%
\end{pgfscope}%
\begin{pgfscope}%
\pgfsetrectcap%
\pgfsetmiterjoin%
\pgfsetlinewidth{0.803000pt}%
\definecolor{currentstroke}{rgb}{0.000000,0.000000,0.000000}%
\pgfsetstrokecolor{currentstroke}%
\pgfsetdash{}{0pt}%
\pgfpathmoveto{\pgfqpoint{0.627968in}{2.859073in}}%
\pgfpathlineto{\pgfqpoint{2.952968in}{2.859073in}}%
\pgfusepath{stroke}%
\end{pgfscope}%
\begin{pgfscope}%
\pgfsetbuttcap%
\pgfsetmiterjoin%
\definecolor{currentfill}{rgb}{1.000000,1.000000,1.000000}%
\pgfsetfillcolor{currentfill}%
\pgfsetfillopacity{0.800000}%
\pgfsetlinewidth{1.003750pt}%
\definecolor{currentstroke}{rgb}{0.800000,0.800000,0.800000}%
\pgfsetstrokecolor{currentstroke}%
\pgfsetstrokeopacity{0.800000}%
\pgfsetdash{}{0pt}%
\pgfpathmoveto{\pgfqpoint{1.585132in}{1.563704in}}%
\pgfpathlineto{\pgfqpoint{2.869635in}{1.563704in}}%
\pgfpathlineto{\pgfqpoint{2.869635in}{2.775739in}}%
\pgfpathlineto{\pgfqpoint{1.585132in}{2.775739in}}%
\pgfpathlineto{\pgfqpoint{1.585132in}{1.563704in}}%
\pgfpathclose%
\pgfusepath{stroke,fill}%
\end{pgfscope}%
\begin{pgfscope}%
\pgfsetrectcap%
\pgfsetroundjoin%
\pgfsetlinewidth{1.505625pt}%
\definecolor{currentstroke}{rgb}{0.392157,0.560784,1.000000}%
\pgfsetstrokecolor{currentstroke}%
\pgfsetdash{}{0pt}%
\pgfpathmoveto{\pgfqpoint{1.651799in}{2.650739in}}%
\pgfpathlineto{\pgfqpoint{1.818466in}{2.650739in}}%
\pgfpathlineto{\pgfqpoint{1.985132in}{2.650739in}}%
\pgfusepath{stroke}%
\end{pgfscope}%
\begin{pgfscope}%
\definecolor{textcolor}{rgb}{0.000000,0.000000,0.000000}%
\pgfsetstrokecolor{textcolor}%
\pgfsetfillcolor{textcolor}%
\pgftext[x=2.118466in,y=2.592406in,left,base]{\color{textcolor}{\rmfamily\fontsize{12.000000}{14.400000}\selectfont\catcode`\^=\active\def^{\ifmmode\sp\else\^{}\fi}\catcode`\%=\active\def
\end{pgfscope}%
\begin{pgfscope}%
\pgfsetrectcap%
\pgfsetroundjoin%
\pgfsetlinewidth{1.505625pt}%
\definecolor{currentstroke}{rgb}{0.470588,0.368627,0.941176}%
\pgfsetstrokecolor{currentstroke}%
\pgfsetdash{}{0pt}%
\pgfpathmoveto{\pgfqpoint{1.651799in}{2.418332in}}%
\pgfpathlineto{\pgfqpoint{1.818466in}{2.418332in}}%
\pgfpathlineto{\pgfqpoint{1.985132in}{2.418332in}}%
\pgfusepath{stroke}%
\end{pgfscope}%
\begin{pgfscope}%
\definecolor{textcolor}{rgb}{0.000000,0.000000,0.000000}%
\pgfsetstrokecolor{textcolor}%
\pgfsetfillcolor{textcolor}%
\pgftext[x=2.118466in,y=2.359999in,left,base]{\color{textcolor}{\rmfamily\fontsize{12.000000}{14.400000}\selectfont\catcode`\^=\active\def^{\ifmmode\sp\else\^{}\fi}\catcode`\%=\active\def
\end{pgfscope}%
\begin{pgfscope}%
\pgfsetrectcap%
\pgfsetroundjoin%
\pgfsetlinewidth{1.505625pt}%
\definecolor{currentstroke}{rgb}{0.862745,0.149020,0.498039}%
\pgfsetstrokecolor{currentstroke}%
\pgfsetdash{}{0pt}%
\pgfpathmoveto{\pgfqpoint{1.651799in}{2.185925in}}%
\pgfpathlineto{\pgfqpoint{1.818466in}{2.185925in}}%
\pgfpathlineto{\pgfqpoint{1.985132in}{2.185925in}}%
\pgfusepath{stroke}%
\end{pgfscope}%
\begin{pgfscope}%
\definecolor{textcolor}{rgb}{0.000000,0.000000,0.000000}%
\pgfsetstrokecolor{textcolor}%
\pgfsetfillcolor{textcolor}%
\pgftext[x=2.118466in,y=2.127592in,left,base]{\color{textcolor}{\rmfamily\fontsize{12.000000}{14.400000}\selectfont\catcode`\^=\active\def^{\ifmmode\sp\else\^{}\fi}\catcode`\%=\active\def
\end{pgfscope}%
\begin{pgfscope}%
\pgfsetrectcap%
\pgfsetroundjoin%
\pgfsetlinewidth{1.505625pt}%
\definecolor{currentstroke}{rgb}{0.996078,0.380392,0.000000}%
\pgfsetstrokecolor{currentstroke}%
\pgfsetdash{}{0pt}%
\pgfpathmoveto{\pgfqpoint{1.651799in}{1.953518in}}%
\pgfpathlineto{\pgfqpoint{1.818466in}{1.953518in}}%
\pgfpathlineto{\pgfqpoint{1.985132in}{1.953518in}}%
\pgfusepath{stroke}%
\end{pgfscope}%
\begin{pgfscope}%
\definecolor{textcolor}{rgb}{0.000000,0.000000,0.000000}%
\pgfsetstrokecolor{textcolor}%
\pgfsetfillcolor{textcolor}%
\pgftext[x=2.118466in,y=1.895185in,left,base]{\color{textcolor}{\rmfamily\fontsize{12.000000}{14.400000}\selectfont\catcode`\^=\active\def^{\ifmmode\sp\else\^{}\fi}\catcode`\%=\active\def
\end{pgfscope}%
\begin{pgfscope}%
\pgfsetrectcap%
\pgfsetroundjoin%
\pgfsetlinewidth{1.505625pt}%
\definecolor{currentstroke}{rgb}{1.000000,0.690196,0.000000}%
\pgfsetstrokecolor{currentstroke}%
\pgfsetdash{}{0pt}%
\pgfpathmoveto{\pgfqpoint{1.651799in}{1.721111in}}%
\pgfpathlineto{\pgfqpoint{1.818466in}{1.721111in}}%
\pgfpathlineto{\pgfqpoint{1.985132in}{1.721111in}}%
\pgfusepath{stroke}%
\end{pgfscope}%
\begin{pgfscope}%
\definecolor{textcolor}{rgb}{0.000000,0.000000,0.000000}%
\pgfsetstrokecolor{textcolor}%
\pgfsetfillcolor{textcolor}%
\pgftext[x=2.118466in,y=1.662778in,left,base]{\color{textcolor}{\rmfamily\fontsize{12.000000}{14.400000}\selectfont\catcode`\^=\active\def^{\ifmmode\sp\else\^{}\fi}\catcode`\%=\active\def
\end{pgfscope}%
\end{pgfpicture}%
\makeatother%
\endgroup%

%% file: plots/hessian_density_loss.pgf
\begingroup%
\makeatletter%
\begin{pgfpicture}%
\pgfpathrectangle{\pgfpointorigin}{\pgfqpoint{3.034449in}{2.959073in}}%
\pgfusepath{use as bounding box, clip}%
\begin{pgfscope}%
\pgfsetbuttcap%
\pgfsetmiterjoin%
\definecolor{currentfill}{rgb}{1.000000,1.000000,1.000000}%
\pgfsetfillcolor{currentfill}%
\pgfsetlinewidth{0.000000pt}%
\definecolor{currentstroke}{rgb}{1.000000,1.000000,1.000000}%
\pgfsetstrokecolor{currentstroke}%
\pgfsetdash{}{0pt}%
\pgfpathmoveto{\pgfqpoint{0.000000in}{-0.000000in}}%
\pgfpathlineto{\pgfqpoint{3.034449in}{-0.000000in}}%
\pgfpathlineto{\pgfqpoint{3.034449in}{2.959073in}}%
\pgfpathlineto{\pgfqpoint{0.000000in}{2.959073in}}%
\pgfpathlineto{\pgfqpoint{0.000000in}{-0.000000in}}%
\pgfpathclose%
\pgfusepath{fill}%
\end{pgfscope}%
\begin{pgfscope}%
\pgfsetbuttcap%
\pgfsetmiterjoin%
\definecolor{currentfill}{rgb}{1.000000,1.000000,1.000000}%
\pgfsetfillcolor{currentfill}%
\pgfsetlinewidth{0.000000pt}%
\definecolor{currentstroke}{rgb}{0.000000,0.000000,0.000000}%
\pgfsetstrokecolor{currentstroke}%
\pgfsetstrokeopacity{0.000000}%
\pgfsetdash{}{0pt}%
\pgfpathmoveto{\pgfqpoint{0.609449in}{0.549073in}}%
\pgfpathlineto{\pgfqpoint{2.934449in}{0.549073in}}%
\pgfpathlineto{\pgfqpoint{2.934449in}{2.859073in}}%
\pgfpathlineto{\pgfqpoint{0.609449in}{2.859073in}}%
\pgfpathlineto{\pgfqpoint{0.609449in}{0.549073in}}%
\pgfpathclose%
\pgfusepath{fill}%
\end{pgfscope}%
\begin{pgfscope}%
\pgfpathrectangle{\pgfqpoint{0.609449in}{0.549073in}}{\pgfqpoint{2.325000in}{2.310000in}}%
\pgfusepath{clip}%
\pgfsetrectcap%
\pgfsetroundjoin%
\pgfsetlinewidth{0.250937pt}%
\definecolor{currentstroke}{rgb}{0.000000,0.000000,0.000000}%
\pgfsetstrokecolor{currentstroke}%
\pgfsetstrokeopacity{0.200000}%
\pgfsetdash{}{0pt}%
\pgfpathmoveto{\pgfqpoint{0.664806in}{0.549073in}}%
\pgfpathlineto{\pgfqpoint{0.664806in}{2.859073in}}%
\pgfusepath{stroke}%
\end{pgfscope}%
\begin{pgfscope}%
\pgfsetbuttcap%
\pgfsetroundjoin%
\definecolor{currentfill}{rgb}{0.000000,0.000000,0.000000}%
\pgfsetfillcolor{currentfill}%
\pgfsetlinewidth{0.803000pt}%
\definecolor{currentstroke}{rgb}{0.000000,0.000000,0.000000}%
\pgfsetstrokecolor{currentstroke}%
\pgfsetdash{}{0pt}%
\pgfsys@defobject{currentmarker}{\pgfqpoint{0.000000in}{-0.048611in}}{\pgfqpoint{0.000000in}{0.000000in}}{%
\pgfpathmoveto{\pgfqpoint{0.000000in}{0.000000in}}%
\pgfpathlineto{\pgfqpoint{0.000000in}{-0.048611in}}%
\pgfusepath{stroke,fill}%
}%
\begin{pgfscope}%
\pgfsys@transformshift{0.664806in}{0.549073in}%
\pgfsys@useobject{currentmarker}{}%
\end{pgfscope}%
\end{pgfscope}%
\begin{pgfscope}%
\definecolor{textcolor}{rgb}{0.000000,0.000000,0.000000}%
\pgfsetstrokecolor{textcolor}%
\pgfsetfillcolor{textcolor}%
\pgftext[x=0.664806in,y=0.451851in,,top]{\color{textcolor}{\rmfamily\fontsize{12.000000}{14.400000}\selectfont\catcode`\^=\active\def^{\ifmmode\sp\else\^{}\fi}\catcode`\%=\active\def
\end{pgfscope}%
\begin{pgfscope}%
\pgfpathrectangle{\pgfqpoint{0.609449in}{0.549073in}}{\pgfqpoint{2.325000in}{2.310000in}}%
\pgfusepath{clip}%
\pgfsetrectcap%
\pgfsetroundjoin%
\pgfsetlinewidth{0.250937pt}%
\definecolor{currentstroke}{rgb}{0.000000,0.000000,0.000000}%
\pgfsetstrokecolor{currentstroke}%
\pgfsetstrokeopacity{0.200000}%
\pgfsetdash{}{0pt}%
\pgfpathmoveto{\pgfqpoint{1.218378in}{0.549073in}}%
\pgfpathlineto{\pgfqpoint{1.218378in}{2.859073in}}%
\pgfusepath{stroke}%
\end{pgfscope}%
\begin{pgfscope}%
\pgfsetbuttcap%
\pgfsetroundjoin%
\definecolor{currentfill}{rgb}{0.000000,0.000000,0.000000}%
\pgfsetfillcolor{currentfill}%
\pgfsetlinewidth{0.803000pt}%
\definecolor{currentstroke}{rgb}{0.000000,0.000000,0.000000}%
\pgfsetstrokecolor{currentstroke}%
\pgfsetdash{}{0pt}%
\pgfsys@defobject{currentmarker}{\pgfqpoint{0.000000in}{-0.048611in}}{\pgfqpoint{0.000000in}{0.000000in}}{%
\pgfpathmoveto{\pgfqpoint{0.000000in}{0.000000in}}%
\pgfpathlineto{\pgfqpoint{0.000000in}{-0.048611in}}%
\pgfusepath{stroke,fill}%
}%
\begin{pgfscope}%
\pgfsys@transformshift{1.218378in}{0.549073in}%
\pgfsys@useobject{currentmarker}{}%
\end{pgfscope}%
\end{pgfscope}%
\begin{pgfscope}%
\definecolor{textcolor}{rgb}{0.000000,0.000000,0.000000}%
\pgfsetstrokecolor{textcolor}%
\pgfsetfillcolor{textcolor}%
\pgftext[x=1.218378in,y=0.451851in,,top]{\color{textcolor}{\rmfamily\fontsize{12.000000}{14.400000}\selectfont\catcode`\^=\active\def^{\ifmmode\sp\else\^{}\fi}\catcode`\%=\active\def
\end{pgfscope}%
\begin{pgfscope}%
\pgfpathrectangle{\pgfqpoint{0.609449in}{0.549073in}}{\pgfqpoint{2.325000in}{2.310000in}}%
\pgfusepath{clip}%
\pgfsetrectcap%
\pgfsetroundjoin%
\pgfsetlinewidth{0.250937pt}%
\definecolor{currentstroke}{rgb}{0.000000,0.000000,0.000000}%
\pgfsetstrokecolor{currentstroke}%
\pgfsetstrokeopacity{0.200000}%
\pgfsetdash{}{0pt}%
\pgfpathmoveto{\pgfqpoint{1.771949in}{0.549073in}}%
\pgfpathlineto{\pgfqpoint{1.771949in}{2.859073in}}%
\pgfusepath{stroke}%
\end{pgfscope}%
\begin{pgfscope}%
\pgfsetbuttcap%
\pgfsetroundjoin%
\definecolor{currentfill}{rgb}{0.000000,0.000000,0.000000}%
\pgfsetfillcolor{currentfill}%
\pgfsetlinewidth{0.803000pt}%
\definecolor{currentstroke}{rgb}{0.000000,0.000000,0.000000}%
\pgfsetstrokecolor{currentstroke}%
\pgfsetdash{}{0pt}%
\pgfsys@defobject{currentmarker}{\pgfqpoint{0.000000in}{-0.048611in}}{\pgfqpoint{0.000000in}{0.000000in}}{%
\pgfpathmoveto{\pgfqpoint{0.000000in}{0.000000in}}%
\pgfpathlineto{\pgfqpoint{0.000000in}{-0.048611in}}%
\pgfusepath{stroke,fill}%
}%
\begin{pgfscope}%
\pgfsys@transformshift{1.771949in}{0.549073in}%
\pgfsys@useobject{currentmarker}{}%
\end{pgfscope}%
\end{pgfscope}%
\begin{pgfscope}%
\definecolor{textcolor}{rgb}{0.000000,0.000000,0.000000}%
\pgfsetstrokecolor{textcolor}%
\pgfsetfillcolor{textcolor}%
\pgftext[x=1.771949in,y=0.451851in,,top]{\color{textcolor}{\rmfamily\fontsize{12.000000}{14.400000}\selectfont\catcode`\^=\active\def^{\ifmmode\sp\else\^{}\fi}\catcode`\%=\active\def
\end{pgfscope}%
\begin{pgfscope}%
\pgfpathrectangle{\pgfqpoint{0.609449in}{0.549073in}}{\pgfqpoint{2.325000in}{2.310000in}}%
\pgfusepath{clip}%
\pgfsetrectcap%
\pgfsetroundjoin%
\pgfsetlinewidth{0.250937pt}%
\definecolor{currentstroke}{rgb}{0.000000,0.000000,0.000000}%
\pgfsetstrokecolor{currentstroke}%
\pgfsetstrokeopacity{0.200000}%
\pgfsetdash{}{0pt}%
\pgfpathmoveto{\pgfqpoint{2.325521in}{0.549073in}}%
\pgfpathlineto{\pgfqpoint{2.325521in}{2.859073in}}%
\pgfusepath{stroke}%
\end{pgfscope}%
\begin{pgfscope}%
\pgfsetbuttcap%
\pgfsetroundjoin%
\definecolor{currentfill}{rgb}{0.000000,0.000000,0.000000}%
\pgfsetfillcolor{currentfill}%
\pgfsetlinewidth{0.803000pt}%
\definecolor{currentstroke}{rgb}{0.000000,0.000000,0.000000}%
\pgfsetstrokecolor{currentstroke}%
\pgfsetdash{}{0pt}%
\pgfsys@defobject{currentmarker}{\pgfqpoint{0.000000in}{-0.048611in}}{\pgfqpoint{0.000000in}{0.000000in}}{%
\pgfpathmoveto{\pgfqpoint{0.000000in}{0.000000in}}%
\pgfpathlineto{\pgfqpoint{0.000000in}{-0.048611in}}%
\pgfusepath{stroke,fill}%
}%
\begin{pgfscope}%
\pgfsys@transformshift{2.325521in}{0.549073in}%
\pgfsys@useobject{currentmarker}{}%
\end{pgfscope}%
\end{pgfscope}%
\begin{pgfscope}%
\definecolor{textcolor}{rgb}{0.000000,0.000000,0.000000}%
\pgfsetstrokecolor{textcolor}%
\pgfsetfillcolor{textcolor}%
\pgftext[x=2.325521in,y=0.451851in,,top]{\color{textcolor}{\rmfamily\fontsize{12.000000}{14.400000}\selectfont\catcode`\^=\active\def^{\ifmmode\sp\else\^{}\fi}\catcode`\%=\active\def
\end{pgfscope}%
\begin{pgfscope}%
\pgfpathrectangle{\pgfqpoint{0.609449in}{0.549073in}}{\pgfqpoint{2.325000in}{2.310000in}}%
\pgfusepath{clip}%
\pgfsetrectcap%
\pgfsetroundjoin%
\pgfsetlinewidth{0.250937pt}%
\definecolor{currentstroke}{rgb}{0.000000,0.000000,0.000000}%
\pgfsetstrokecolor{currentstroke}%
\pgfsetstrokeopacity{0.200000}%
\pgfsetdash{}{0pt}%
\pgfpathmoveto{\pgfqpoint{2.879092in}{0.549073in}}%
\pgfpathlineto{\pgfqpoint{2.879092in}{2.859073in}}%
\pgfusepath{stroke}%
\end{pgfscope}%
\begin{pgfscope}%
\pgfsetbuttcap%
\pgfsetroundjoin%
\definecolor{currentfill}{rgb}{0.000000,0.000000,0.000000}%
\pgfsetfillcolor{currentfill}%
\pgfsetlinewidth{0.803000pt}%
\definecolor{currentstroke}{rgb}{0.000000,0.000000,0.000000}%
\pgfsetstrokecolor{currentstroke}%
\pgfsetdash{}{0pt}%
\pgfsys@defobject{currentmarker}{\pgfqpoint{0.000000in}{-0.048611in}}{\pgfqpoint{0.000000in}{0.000000in}}{%
\pgfpathmoveto{\pgfqpoint{0.000000in}{0.000000in}}%
\pgfpathlineto{\pgfqpoint{0.000000in}{-0.048611in}}%
\pgfusepath{stroke,fill}%
}%
\begin{pgfscope}%
\pgfsys@transformshift{2.879092in}{0.549073in}%
\pgfsys@useobject{currentmarker}{}%
\end{pgfscope}%
\end{pgfscope}%
\begin{pgfscope}%
\definecolor{textcolor}{rgb}{0.000000,0.000000,0.000000}%
\pgfsetstrokecolor{textcolor}%
\pgfsetfillcolor{textcolor}%
\pgftext[x=2.879092in,y=0.451851in,,top]{\color{textcolor}{\rmfamily\fontsize{12.000000}{14.400000}\selectfont\catcode`\^=\active\def^{\ifmmode\sp\else\^{}\fi}\catcode`\%=\active\def
\end{pgfscope}%
\begin{pgfscope}%
\definecolor{textcolor}{rgb}{0.000000,0.000000,0.000000}%
\pgfsetstrokecolor{textcolor}%
\pgfsetfillcolor{textcolor}%
\pgftext[x=1.771949in,y=0.248148in,,top]{\color{textcolor}{\rmfamily\fontsize{12.000000}{14.400000}\selectfont\catcode`\^=\active\def^{\ifmmode\sp\else\^{}\fi}\catcode`\%=\active\def
\end{pgfscope}%
\begin{pgfscope}%
\pgfpathrectangle{\pgfqpoint{0.609449in}{0.549073in}}{\pgfqpoint{2.325000in}{2.310000in}}%
\pgfusepath{clip}%
\pgfsetrectcap%
\pgfsetroundjoin%
\pgfsetlinewidth{0.250937pt}%
\definecolor{currentstroke}{rgb}{0.000000,0.000000,0.000000}%
\pgfsetstrokecolor{currentstroke}%
\pgfsetstrokeopacity{0.200000}%
\pgfsetdash{}{0pt}%
\pgfpathmoveto{\pgfqpoint{0.609449in}{0.575455in}}%
\pgfpathlineto{\pgfqpoint{2.934449in}{0.575455in}}%
\pgfusepath{stroke}%
\end{pgfscope}%
\begin{pgfscope}%
\pgfsetbuttcap%
\pgfsetroundjoin%
\definecolor{currentfill}{rgb}{0.000000,0.000000,0.000000}%
\pgfsetfillcolor{currentfill}%
\pgfsetlinewidth{0.803000pt}%
\definecolor{currentstroke}{rgb}{0.000000,0.000000,0.000000}%
\pgfsetstrokecolor{currentstroke}%
\pgfsetdash{}{0pt}%
\pgfsys@defobject{currentmarker}{\pgfqpoint{-0.048611in}{0.000000in}}{\pgfqpoint{-0.000000in}{0.000000in}}{%
\pgfpathmoveto{\pgfqpoint{-0.000000in}{0.000000in}}%
\pgfpathlineto{\pgfqpoint{-0.048611in}{0.000000in}}%
\pgfusepath{stroke,fill}%
}%
\begin{pgfscope}%
\pgfsys@transformshift{0.609449in}{0.575455in}%
\pgfsys@useobject{currentmarker}{}%
\end{pgfscope}%
\end{pgfscope}%
\begin{pgfscope}%
\definecolor{textcolor}{rgb}{0.000000,0.000000,0.000000}%
\pgfsetstrokecolor{textcolor}%
\pgfsetfillcolor{textcolor}%
\pgftext[x=0.303703in, y=0.517585in, left, base]{\color{textcolor}{\rmfamily\fontsize{12.000000}{14.400000}\selectfont\catcode`\^=\active\def^{\ifmmode\sp\else\^{}\fi}\catcode`\%=\active\def
\end{pgfscope}%
\begin{pgfscope}%
\pgfpathrectangle{\pgfqpoint{0.609449in}{0.549073in}}{\pgfqpoint{2.325000in}{2.310000in}}%
\pgfusepath{clip}%
\pgfsetrectcap%
\pgfsetroundjoin%
\pgfsetlinewidth{0.250937pt}%
\definecolor{currentstroke}{rgb}{0.000000,0.000000,0.000000}%
\pgfsetstrokecolor{currentstroke}%
\pgfsetstrokeopacity{0.200000}%
\pgfsetdash{}{0pt}%
\pgfpathmoveto{\pgfqpoint{0.609449in}{1.311946in}}%
\pgfpathlineto{\pgfqpoint{2.934449in}{1.311946in}}%
\pgfusepath{stroke}%
\end{pgfscope}%
\begin{pgfscope}%
\pgfsetbuttcap%
\pgfsetroundjoin%
\definecolor{currentfill}{rgb}{0.000000,0.000000,0.000000}%
\pgfsetfillcolor{currentfill}%
\pgfsetlinewidth{0.803000pt}%
\definecolor{currentstroke}{rgb}{0.000000,0.000000,0.000000}%
\pgfsetstrokecolor{currentstroke}%
\pgfsetdash{}{0pt}%
\pgfsys@defobject{currentmarker}{\pgfqpoint{-0.048611in}{0.000000in}}{\pgfqpoint{-0.000000in}{0.000000in}}{%
\pgfpathmoveto{\pgfqpoint{-0.000000in}{0.000000in}}%
\pgfpathlineto{\pgfqpoint{-0.048611in}{0.000000in}}%
\pgfusepath{stroke,fill}%
}%
\begin{pgfscope}%
\pgfsys@transformshift{0.609449in}{1.311946in}%
\pgfsys@useobject{currentmarker}{}%
\end{pgfscope}%
\end{pgfscope}%
\begin{pgfscope}%
\definecolor{textcolor}{rgb}{0.000000,0.000000,0.000000}%
\pgfsetstrokecolor{textcolor}%
\pgfsetfillcolor{textcolor}%
\pgftext[x=0.303703in, y=1.254076in, left, base]{\color{textcolor}{\rmfamily\fontsize{12.000000}{14.400000}\selectfont\catcode`\^=\active\def^{\ifmmode\sp\else\^{}\fi}\catcode`\%=\active\def
\end{pgfscope}%
\begin{pgfscope}%
\pgfpathrectangle{\pgfqpoint{0.609449in}{0.549073in}}{\pgfqpoint{2.325000in}{2.310000in}}%
\pgfusepath{clip}%
\pgfsetrectcap%
\pgfsetroundjoin%
\pgfsetlinewidth{0.250937pt}%
\definecolor{currentstroke}{rgb}{0.000000,0.000000,0.000000}%
\pgfsetstrokecolor{currentstroke}%
\pgfsetstrokeopacity{0.200000}%
\pgfsetdash{}{0pt}%
\pgfpathmoveto{\pgfqpoint{0.609449in}{2.048438in}}%
\pgfpathlineto{\pgfqpoint{2.934449in}{2.048438in}}%
\pgfusepath{stroke}%
\end{pgfscope}%
\begin{pgfscope}%
\pgfsetbuttcap%
\pgfsetroundjoin%
\definecolor{currentfill}{rgb}{0.000000,0.000000,0.000000}%
\pgfsetfillcolor{currentfill}%
\pgfsetlinewidth{0.803000pt}%
\definecolor{currentstroke}{rgb}{0.000000,0.000000,0.000000}%
\pgfsetstrokecolor{currentstroke}%
\pgfsetdash{}{0pt}%
\pgfsys@defobject{currentmarker}{\pgfqpoint{-0.048611in}{0.000000in}}{\pgfqpoint{-0.000000in}{0.000000in}}{%
\pgfpathmoveto{\pgfqpoint{-0.000000in}{0.000000in}}%
\pgfpathlineto{\pgfqpoint{-0.048611in}{0.000000in}}%
\pgfusepath{stroke,fill}%
}%
\begin{pgfscope}%
\pgfsys@transformshift{0.609449in}{2.048438in}%
\pgfsys@useobject{currentmarker}{}%
\end{pgfscope}%
\end{pgfscope}%
\begin{pgfscope}%
\definecolor{textcolor}{rgb}{0.000000,0.000000,0.000000}%
\pgfsetstrokecolor{textcolor}%
\pgfsetfillcolor{textcolor}%
\pgftext[x=0.303703in, y=1.990568in, left, base]{\color{textcolor}{\rmfamily\fontsize{12.000000}{14.400000}\selectfont\catcode`\^=\active\def^{\ifmmode\sp\else\^{}\fi}\catcode`\%=\active\def
\end{pgfscope}%
\begin{pgfscope}%
\pgfpathrectangle{\pgfqpoint{0.609449in}{0.549073in}}{\pgfqpoint{2.325000in}{2.310000in}}%
\pgfusepath{clip}%
\pgfsetrectcap%
\pgfsetroundjoin%
\pgfsetlinewidth{0.250937pt}%
\definecolor{currentstroke}{rgb}{0.000000,0.000000,0.000000}%
\pgfsetstrokecolor{currentstroke}%
\pgfsetstrokeopacity{0.200000}%
\pgfsetdash{}{0pt}%
\pgfpathmoveto{\pgfqpoint{0.609449in}{2.784929in}}%
\pgfpathlineto{\pgfqpoint{2.934449in}{2.784929in}}%
\pgfusepath{stroke}%
\end{pgfscope}%
\begin{pgfscope}%
\pgfsetbuttcap%
\pgfsetroundjoin%
\definecolor{currentfill}{rgb}{0.000000,0.000000,0.000000}%
\pgfsetfillcolor{currentfill}%
\pgfsetlinewidth{0.803000pt}%
\definecolor{currentstroke}{rgb}{0.000000,0.000000,0.000000}%
\pgfsetstrokecolor{currentstroke}%
\pgfsetdash{}{0pt}%
\pgfsys@defobject{currentmarker}{\pgfqpoint{-0.048611in}{0.000000in}}{\pgfqpoint{-0.000000in}{0.000000in}}{%
\pgfpathmoveto{\pgfqpoint{-0.000000in}{0.000000in}}%
\pgfpathlineto{\pgfqpoint{-0.048611in}{0.000000in}}%
\pgfusepath{stroke,fill}%
}%
\begin{pgfscope}%
\pgfsys@transformshift{0.609449in}{2.784929in}%
\pgfsys@useobject{currentmarker}{}%
\end{pgfscope}%
\end{pgfscope}%
\begin{pgfscope}%
\definecolor{textcolor}{rgb}{0.000000,0.000000,0.000000}%
\pgfsetstrokecolor{textcolor}%
\pgfsetfillcolor{textcolor}%
\pgftext[x=0.303703in, y=2.727059in, left, base]{\color{textcolor}{\rmfamily\fontsize{12.000000}{14.400000}\selectfont\catcode`\^=\active\def^{\ifmmode\sp\else\^{}\fi}\catcode`\%=\active\def
\end{pgfscope}%
\begin{pgfscope}%
\definecolor{textcolor}{rgb}{0.000000,0.000000,0.000000}%
\pgfsetstrokecolor{textcolor}%
\pgfsetfillcolor{textcolor}%
\pgftext[x=0.248148in,y=1.704073in,,bottom,rotate=90.000000]{\color{textcolor}{\rmfamily\fontsize{12.000000}{14.400000}\selectfont\catcode`\^=\active\def^{\ifmmode\sp\else\^{}\fi}\catcode`\%=\active\def
\end{pgfscope}%
\begin{pgfscope}%
\pgfpathrectangle{\pgfqpoint{0.609449in}{0.549073in}}{\pgfqpoint{2.325000in}{2.310000in}}%
\pgfusepath{clip}%
\pgfsetrectcap%
\pgfsetroundjoin%
\pgfsetlinewidth{1.505625pt}%
\definecolor{currentstroke}{rgb}{0.000000,0.000000,0.000000}%
\pgfsetstrokecolor{currentstroke}%
\pgfsetdash{}{0pt}%
\pgfpathmoveto{\pgfqpoint{0.664806in}{2.666573in}}%
\pgfpathlineto{\pgfqpoint{0.941592in}{1.144941in}}%
\pgfpathlineto{\pgfqpoint{1.218378in}{0.863584in}}%
\pgfpathlineto{\pgfqpoint{1.495164in}{0.805609in}}%
\pgfpathlineto{\pgfqpoint{1.771949in}{0.779224in}}%
\pgfpathlineto{\pgfqpoint{2.048735in}{0.764086in}}%
\pgfpathlineto{\pgfqpoint{2.325521in}{0.754136in}}%
\pgfpathlineto{\pgfqpoint{2.602306in}{0.747002in}}%
\pgfpathlineto{\pgfqpoint{2.879092in}{0.741573in}}%
\pgfusepath{stroke}%
\end{pgfscope}%
\begin{pgfscope}%
\pgfpathrectangle{\pgfqpoint{0.609449in}{0.549073in}}{\pgfqpoint{2.325000in}{2.310000in}}%
\pgfusepath{clip}%
\pgfsetbuttcap%
\pgfsetmiterjoin%
\definecolor{currentfill}{rgb}{0.000000,0.000000,0.000000}%
\pgfsetfillcolor{currentfill}%
\pgfsetlinewidth{1.003750pt}%
\definecolor{currentstroke}{rgb}{0.000000,0.000000,0.000000}%
\pgfsetstrokecolor{currentstroke}%
\pgfsetdash{}{0pt}%
\pgfsys@defobject{currentmarker}{\pgfqpoint{-0.035355in}{-0.058926in}}{\pgfqpoint{0.035355in}{0.058926in}}{%
\pgfpathmoveto{\pgfqpoint{-0.000000in}{-0.058926in}}%
\pgfpathlineto{\pgfqpoint{0.035355in}{0.000000in}}%
\pgfpathlineto{\pgfqpoint{0.000000in}{0.058926in}}%
\pgfpathlineto{\pgfqpoint{-0.035355in}{0.000000in}}%
\pgfpathlineto{\pgfqpoint{-0.000000in}{-0.058926in}}%
\pgfpathclose%
\pgfusepath{stroke,fill}%
}%
\begin{pgfscope}%
\pgfsys@transformshift{0.664806in}{2.666573in}%
\pgfsys@useobject{currentmarker}{}%
\end{pgfscope}%
\begin{pgfscope}%
\pgfsys@transformshift{0.941592in}{1.144941in}%
\pgfsys@useobject{currentmarker}{}%
\end{pgfscope}%
\begin{pgfscope}%
\pgfsys@transformshift{1.218378in}{0.863584in}%
\pgfsys@useobject{currentmarker}{}%
\end{pgfscope}%
\begin{pgfscope}%
\pgfsys@transformshift{1.495164in}{0.805609in}%
\pgfsys@useobject{currentmarker}{}%
\end{pgfscope}%
\begin{pgfscope}%
\pgfsys@transformshift{1.771949in}{0.779224in}%
\pgfsys@useobject{currentmarker}{}%
\end{pgfscope}%
\begin{pgfscope}%
\pgfsys@transformshift{2.048735in}{0.764086in}%
\pgfsys@useobject{currentmarker}{}%
\end{pgfscope}%
\begin{pgfscope}%
\pgfsys@transformshift{2.325521in}{0.754136in}%
\pgfsys@useobject{currentmarker}{}%
\end{pgfscope}%
\begin{pgfscope}%
\pgfsys@transformshift{2.602306in}{0.747002in}%
\pgfsys@useobject{currentmarker}{}%
\end{pgfscope}%
\begin{pgfscope}%
\pgfsys@transformshift{2.879092in}{0.741573in}%
\pgfsys@useobject{currentmarker}{}%
\end{pgfscope}%
\end{pgfscope}%
\begin{pgfscope}%
\pgfsetrectcap%
\pgfsetmiterjoin%
\pgfsetlinewidth{0.803000pt}%
\definecolor{currentstroke}{rgb}{0.000000,0.000000,0.000000}%
\pgfsetstrokecolor{currentstroke}%
\pgfsetdash{}{0pt}%
\pgfpathmoveto{\pgfqpoint{0.609449in}{0.549073in}}%
\pgfpathlineto{\pgfqpoint{0.609449in}{2.859073in}}%
\pgfusepath{stroke}%
\end{pgfscope}%
\begin{pgfscope}%
\pgfsetrectcap%
\pgfsetmiterjoin%
\pgfsetlinewidth{0.803000pt}%
\definecolor{currentstroke}{rgb}{0.000000,0.000000,0.000000}%
\pgfsetstrokecolor{currentstroke}%
\pgfsetdash{}{0pt}%
\pgfpathmoveto{\pgfqpoint{2.934449in}{0.549073in}}%
\pgfpathlineto{\pgfqpoint{2.934449in}{2.859073in}}%
\pgfusepath{stroke}%
\end{pgfscope}%
\begin{pgfscope}%
\pgfsetrectcap%
\pgfsetmiterjoin%
\pgfsetlinewidth{0.803000pt}%
\definecolor{currentstroke}{rgb}{0.000000,0.000000,0.000000}%
\pgfsetstrokecolor{currentstroke}%
\pgfsetdash{}{0pt}%
\pgfpathmoveto{\pgfqpoint{0.609449in}{0.549073in}}%
\pgfpathlineto{\pgfqpoint{2.934449in}{0.549073in}}%
\pgfusepath{stroke}%
\end{pgfscope}%
\begin{pgfscope}%
\pgfsetrectcap%
\pgfsetmiterjoin%
\pgfsetlinewidth{0.803000pt}%
\definecolor{currentstroke}{rgb}{0.000000,0.000000,0.000000}%
\pgfsetstrokecolor{currentstroke}%
\pgfsetdash{}{0pt}%
\pgfpathmoveto{\pgfqpoint{0.609449in}{2.859073in}}%
\pgfpathlineto{\pgfqpoint{2.934449in}{2.859073in}}%
\pgfusepath{stroke}%
\end{pgfscope}%
\end{pgfpicture}%
\makeatother%
\endgroup%

%% file: conclusion.tex
\section{Conclusion}
\label{sec:conclusion}

We have analyzed three stochastic trace estimators for constant symmetric matrices $\mtx{B} \in \mathbb{R}^{n \times n}$ when they are applied straightforwardly to a parameter-dependent symmetric matrix $\mtx{B}(t) \in \mathbb{R}^{n \times n}$. We have showed that in the $L^1$-norm, these estimators satisfy bounds which parallel existing results for constant matrices. The results for the estimators allow us to analyze the Nyström-Chebyshev++ method, an efficient technique for approximating the spectral density of a symmetric matrix. Further, we have proposed multiple improvements, both algorithmic and theoretical, to an existing implementation of this method.

%% file: appendix.tex
\appendix

\section{Moment bounds for Gaussian random vectors and matrices}

The goal of this section is to establish moment bounds for $\lVert \mtx{A} \mtx{\Omega} \rVert _2^2$ with a fixed matrix $\mtx{A}$ and a Gaussian random matrix $\mtx{\Omega}$. These bounds are needed in the proof of \cref{lem:nystrom}, but may also be of independent interest.
The corresponding result for the Frobenius norm can be found in~\cite[Lemma 3]{kressner-2024-randomized-lowrank}. For the spectral norm, moment bounds for $\lVert \mtx{\Omega} \rVert _2^2$ (that is, $\mtx{A} = \mtx{I}$) are well established \cite{chen-2005-condition-numbers, edelman-1988-eigenvalues-condition, james-1964-distributions-matrix}. In \cite[Lemma B.1]{tropp-2023-randomized-algorithms} bounds on the first- and second-order moments for general $\mtx{A}$ are derived. In the following, we will generalize this result to moments of arbitrarily high order. To do so, we first establish two preliminary results.

\begin{lemma}[Moment bound of chi-squared random variable]\label{lem:gamma}
    Let $X \sim \chi_k^2$, where $\chi_k^2$ denotes the chi-squared distribution with $k \ge 2$ degrees of freedom. Then 
    $
        \mathbb{E}^{{p}}[X] \leq k + p-1
    $ holds for every $p \ge 1$, $p\in \mathbb R$.
\end{lemma}%
\begin{proof}
It is well known that
    $
        \mathbb{E}^{p}[X] = 2 \Big( \frac{\Gamma(k/2+p)}{\Gamma(k/2)} \Big)^{\sfrac{1}{p}}.
    $ For every $\alpha \ge 0$, $\beta \ge 2$ the bound $\frac{\Gamma(\alpha+\beta)}{\Gamma(\alpha + 1)} \le (\alpha+\beta/2)^{\beta-1}$ holds~\cite[Equation 2.2]{laforgia-1984-further-inequalities}. The result of the lemma follows by using this bound with $\alpha = k/2-1$ and $\beta = p+1$.
\end{proof}

\begin{lemma}[Spectral norm moments of Gaussian random vector]\label{lem:spectral-norm-moment-vector}
    Given $\mtx{A} \in \mathbb{R}^{m \times m}$ and a Gaussian random vector $\vct{\omega} \in \mathbb{R}^{m}$, the bound
    \begin{equation*}
        \mathbb{E}^{p}\big[ \lVert \mtx{A} \vct{\omega} \rVert _2^2 \big]
        \leq  (k + p - 1) \Big( \lVert \mtx{A} \rVert _2^2 + \frac{1}{k} \lVert \mtx{A} \rVert _F^2 \Big).        
    \end{equation*}
    holds for every $k \ge 2$, $k\in \mathbb{N}$, and $p \ge 1$, $p\in \mathbb R$.
\end{lemma}%
\begin{proof}
By the unitary invariance of Gaussian random vectors, we may assume w.l.o.g. that $\mtx{A} = \mtx{\Sigma} = \operatorname{diag}(\sigma_1, \dots, \sigma_m)$ with $\sigma_1 \geq \dots \geq \sigma_m \geq 0$.
Following the proof of \cite[Theorem 1]{cohen-2016-optimal-approximate}, we split the singular values into $\ell = \lceil m/k \rceil$ groups of size $k$:
    \begin{equation*}
        \overbrace{\underbrace{\sigma_1, \dots, \sigma_k}_{\leq \sigma_1}}^{\geq \sigma_{k+1}}, \overbrace{\underbrace{\sigma_{k+1}, \dots, \sigma_{2k}}_{\leq \sigma_{k+1}}}^{\geq \sigma_{2k+1}}, \dots, \overbrace{\underbrace{\sigma_{(\ell - 1)k + 1}, \dots, \sigma_{\ell k}}_{\leq \sigma_{(\ell - 1)k + 1}}}^{\geq 0}.
    \end{equation*}
    If $m$ is not a multiple of $k$, we set $\sigma_i = 0$ for $i > m$. Using that
    $\sigma_{ik + 1}^2 \leq ( \sigma_{(i-1)k + 1} + \cdots + \sigma_{ik - 1} + \sigma_{ik} ) / k$ for $i = 1,\ldots, \ell-1$, we get
    \begin{equation}
        \sum_{i=0}^{\ell-1} \sigma_{ik + 1}^2 = \lVert \mtx{\Sigma} \rVert _2^2 + \sum_{i=1}^{\ell-1} \sigma_{ik + 1}^2  \leq \lVert \mtx{\Sigma} \rVert _2^2 + \frac{1}{k} \sum_{j=1}^{(\ell - 1)k} \sigma_j^2 \leq \lVert \mtx{\Sigma} \rVert _2^2 + \frac{1}{k} \lVert \mtx{\Sigma} \rVert _F^2.
        \label{equ:singular-value-group-bound}
    \end{equation}
    This allows us to bound
    \begin{align*}
    \mathbb{E}^{p}\big[ \lVert \mtx{\Sigma} \vct{\omega} \rVert _2^2 \big] & = 
        \mathbb{E}^{p}\Big[ \sum_{i=1}^{m} \sigma_i^2 \omega_i^2 \Big]
        = \mathbb{E}^{p}\Big[ \sum_{i=0}^{\ell - 1} \sum_{j=1}^{k} \sigma_{ik + j}^2 \omega_{ik + j}^2 \Big] \\
        &\leq \mathbb{E}^{p}\Big[ \sum_{i=0}^{\ell - 1} \sigma_{ik + 1}^2 \sum_{j=1}^{k} \omega_{ik + j}^2 \Big] 
        \leq \sum_{i=0}^{\ell - 1} \sigma_{ik + 1}^2 ~ \mathbb{E}^{p}\Big[ \sum_{j=1}^{k} \omega_{ik + j}^2 \Big] \\
        &\leq (k + p - 1) \sum_{i=0}^{\ell - 1} \sigma_{ik + 1}^2 \leq (k + p - 1) \Big( \lVert \mtx{\Sigma} \rVert _2^2 + \frac{1}{k} \lVert \mtx{\Sigma} \rVert _F^2 \Big), 
    \end{align*}
    where we used the triangle inequality, \cref{lem:gamma}, and \cref{equ:singular-value-group-bound} for the last three inequalities.
\end{proof}

\begin{lemma}[Spectral norm moments of Gaussian random matrix]\label{lem:spectral-norm-moment}
    Given $\mtx{A} \in \mathbb{R}^{m \times m}$ and a Gaussian random matrix $\mtx{\Omega} \in \mathbb{R}^{m \times k}$, the bound
    \begin{equation*}
        \mathbb{E}^{p}\left[ \lVert \mtx{A} \mtx{\Omega} \rVert _2^2 \right]
        \leq  (k + p-1) \Big( 2 \lVert \mtx{A} \rVert _2^2 + \frac{1}{k} \lVert \mtx{A} \rVert _F^2 \Big).
    \end{equation*}
    holds for every $k \ge 2$, $k\in \mathbb{N}$, and $p \ge 1$, $p\in \mathbb R$.
\end{lemma}%
\begin{proof}
    By the proof of \cite[Lemma B.1]{tropp-2023-randomized-algorithms}, we have that
    \begin{equation*}
        \mathbb{E}^{p}\big[ \lVert \mtx{A} \mtx{\Omega} \rVert _2^2 \big]
        \leq \lVert \mtx{A} \rVert _2^2 \mathbb{E}^{p}\big[ \lVert \vct{\omega}_1 \rVert _2^2 \big] + \mathbb{E}^{p}\big[ \lVert \mtx{A} \vct{\omega}_2 \rVert _2^2 \big]
    \end{equation*}
    for Gaussian random vectors $\vct{\omega}_1 \in \mathbb{R}^{k}$ and $\vct{\omega}_2 \in \mathbb{R}^{m}$. The claimed result follows from applying \cref{lem:gamma} and \cref{lem:spectral-norm-moment-vector}:
    \begin{equation*}
        \mathbb{E}^{p}\big[ \lVert \mtx{A} \mtx{\Omega} \rVert _2^2 \big]
        \leq (k + p - 1) \lVert \mtx{A} \rVert _2^2  + (k + p - 1) \Big( \lVert \mtx{A} \rVert _2^2 + \frac{1}{k} \lVert \mtx{A} \rVert _F^2 \Big).
    \end{equation*}
\end{proof}

In passing, we note that \cref{lem:spectral-norm-moment}  yields the bound
    \begin{equation*}
        \mathbb{E}^{p}\big[ \lVert \mtx{A} \mtx{\Omega} \rVert _2 \big] = \sqrt{\mathbb{E}^{\sfrac{p}{2}}\big[ \lVert \mtx{A} \mtx{\Omega} \rVert _2^2 \big]} \leq \sqrt{k + p/2-1} \cdot \Big(\sqrt{2} \lVert \mtx{A} \rVert _2 + \frac{1}{\sqrt{k}}\lVert \mtx{A} \rVert_F\Big)
    \end{equation*}
    with $k\ge 2, p\ge 2$. For $p = 2$, this nearly matches the corresponding bound from~\cite[Lemma B.1]{tropp-2023-randomized-algorithms}, except for the additional factor $\sqrt{2}$.

\section{Moment bounds for sub-gamma random variables}

A real centered random variable $X$ is called sub-gamma with parameters $(\nu,c)$ if its moment generating function satisfies
\begin{equation*}
 \mathbb E[ \exp( tX ) ] \le \exp\left( \frac{ t^2 \nu}{2(1-c|t|)} \right) \ \text{for every} \ 0 < |t| < 1/c.
\end{equation*}
Any such $X$ satisfies the following tail bound~\cite[P. 29]{boucheron-2013-concentration-inequalities}:
\begin{equation}
  \mathbb P\big(|X| > \sqrt{2\nu t} + ct\big) \le 2 e^{-t} \ \text{for every} \ t>0.
  \label{equ:subgammatail}
\end{equation}

\begin{lemma}[Moment bounds for sub-gamma random variables]\label{lem:sub-gamma-moments}
    Let $X$ be a centered, $(\nu,c)$-sub-gamma random variable for $\nu,c>0$. Then
    \begin{equation*}
        \mathbb{E}^p[X] \leq 2 \sqrt{2 \nu p} + 4 c p
    \end{equation*}
    holds for any $p \geq 1$, $p \in \mathbb{R}$.
\end{lemma}
\begin{proof}
The statement follows from a straightforward extension of the proof of~\cite[Theorem 2.3]{boucheron-2013-concentration-inequalities} from  even integers $p$ to general $p$. We include the detailed argument for completeness. Using integration by parts, the reparametrization $x = \sqrt{2\nu t} + ct$, and~\cref{equ:subgammatail}, one gets
\begin{align*}
 \mathbb{E}\big[|X|^p\big] &= p \int_0^\infty  |x|^{p-1} \mathbb P\big(|X| > x\big) \,\mathrm{d}x \\
 &=  p \int_0^\infty  (\sqrt{2\nu t} + ct)^{p-1} \mathbb P\big(|X| > \sqrt{2\nu t} + ct \big) \frac{ \sqrt{2\nu t} + 2ct}{2t} \,\mathrm{d}t \\
 &\le p  \int_0^\infty  (\sqrt{2\nu t} + 2ct)^{p}  \frac{e^{-t}}{t} \,\mathrm{d}t.
\end{align*}
Using the inequality $( (a+b)/2 )^p \le (a^p + b^p)/2$ implied by the convexity of $x^p$ on $[0,\infty)$,
we have that
\begin{align*}
 \mathbb{E}\big[|X|^p\big] & \le 
 p 2^{p-1}   \int_0^\infty  \big( (2\nu t)^{\sfrac{p}{2}} + (2ct)^{p}\big)  \frac{e^{-t}}{t} \,\mathrm{d}t \\
 &= p 2^{p-1} \big( (2\nu)^{\sfrac{p}{2}} \Gamma(p/2) + (2c)^p  \Gamma(p) \big) \\
 &= 2^{p} (2\nu)^{\sfrac{p}{2}} \Gamma(p/2+1) + 2^{p-1} (2c)^p  \Gamma(p+1).
\end{align*}
Using $(a+b)^{\sfrac{1}{p}} \le a^{\sfrac{1}{p}} + b^{\sfrac{1}{p}}$, we obtain

\begin{equation*}
\mathbb{E}^p[X] \leq 2 \sqrt{2\nu} \Gamma(p/2+1)^{\sfrac{1}{p}}  + 4 c \Gamma(p+1)^{\sfrac{1}{p}}.
\end{equation*}
Combined with $\Gamma(p/2+1)^{\sfrac{1}{p}} \le \sqrt{p}$ and $\Gamma(p+1)^{\sfrac{1}{p}} \le p$, this completes the proof.
\end{proof}